\documentclass[reqno,11pt,a4paper]{amsart}
\usepackage[utf8]{inputenc}
\usepackage{pifont}
\usepackage[T1]{fontenc}
\usepackage{enumerate}
\usepackage{textcomp}
\usepackage{verbatim}
\usepackage{pgffor}

\usepackage{esint}
\usepackage{pdfsync}
\usepackage[colorlinks=true,linkcolor=blue]{hyperref}%
\usepackage[T1]{fontenc}
\usepackage{pifont}
\usepackage[english]{babel}

\usepackage{amsmath,esint,color}
\usepackage{amsthm,latexsym,epsfig,graphicx,marvosym, mathrsfs,bigints,stmaryrd,textgreek,upgreek}
\usepackage{amsmath,stmaryrd}
\usepackage{relsize}
\usepackage{amsfonts}
\usepackage{amssymb}
\usepackage[headings]{fullpage}

\allowdisplaybreaks

\usepackage{amsmath}
\usepackage{amsfonts}
\usepackage{amssymb}
\usepackage{amsthm}
\newtheorem{theorem}{Theorem}[section]
\newtheorem{proposition}{Proposition}[section]
\newtheorem{lemma}{Lemma}[section]
\newtheorem{remark}{Remark}[section]

\newtheorem{defn}{Definition}[section]

\begin{document}

\title[Global well-posedness of arbitrarily large Lipschitz solutions]{Global well-posedness of arbitrarily large Lipschitz solutions for the  Muskat problem with surface tension. }
\author{{Omar Lazar}}
\address{New York University Abu Dhabi, Abu Dhabi, UAE}
\email{omar.lazar@nyu.edu}

\begin{abstract} We prove a global well-posedness result  for  the 2D Muskat problem with surface tension. Given any regular enough initial data which is small in some critical space but possibly large in Lipschitz, we prove that the associated Cauchy problem has a unique global solution. Our result allows for the slope of the interface between the two fluids to be arbitrarily large. 
\end{abstract}   
\maketitle
\tableofcontents

\section{Introduction}	
The goal of this paper is to study the Muskat problem \cite{Musk} with surface tension and with or without gravity. The Muskat problem models the dynamic between  two incompressible fluids in a porous media. It is analogous to the two phase Hele-Shaw cell \cite{ADP11,AMS20,TA20,Cp93}. We shall study the case where the two fluids have different densities, $\rho_1$ and $\rho_2$ satisfying $\rho_2>\rho_1$, and same constant viscosity denoted by $\nu$.  We assume that the denser fluid namely the one with density $\rho_2$  is below the moving interface and the lighter one, with density  $\rho_1$, is above.  When surface tension is taking into account, the Muskat problem can be formulated as follows (see Matioc \cite{Mat2}),  
\begin{eqnarray} \label{CP0}
(\mathcal{M}_{\sigma}):\left\{ \begin{array}{ll}
        \displaystyle f_t(x,t) =\displaystyle\frac{k}{2\pi\mu} p.v.  \int_{\mathbb R}\frac{y+ f_x(x,t)\left(f(x,t)-f(x-y,t)\right)}{ y^2+\left(f(x,t)-f(x-y,t)\right)^2 }\partial_x \tau_y \left(\sigma\kappa(f)-g\rho f\right)(t,\cdot)\, dy, \nonumber\\
        f(x,0) = f_0(x)
    \end{array}
\right. 
\end{eqnarray}
Where $\kappa(f)=\kappa(f(x,t))$ is the curvature operator of the  moving interface $\{y=f(x,t)+tV\}$ with constant speed $V$, where~${V\in\mathbb R}$. The curvature operator is defined as follows
$$
\kappa(f):=\frac{\partial_{xx} f(x,t)}{1+(\partial_x f(x,t))^2}
$$
For any $y \in \mathbb{R}$, the operator $\tau_y f(\cdot):=f(\cdot-y)$ is the shifting operator. The constant $k$ is the permeability coefficient of the homogeneous porous media and $\mu$ is the constant viscosity of the fluids. We denote by  $\rho$ the difference between the density of the two fluids, that is $\rho:=\rho_2-\rho_1>0$. For the sake of simplicity,  an easy scaling argument allows us to assume that these constants are equal to 1.
  If the two fluids have the same viscosity and same density, it is by now well known that the Muskat equation without surface tension is locally/globally well-posed in various subcritical and critical spaces \cite{SCH, OS, Ngu22, Mat1, GL22, EMW, Chen, EM, DZL, CL21, CG2,CG,ADP11, CGSV, CCGSP,CCGS, CGS,Cam2, Cam1,ambrose,AN0,AN1, AN2, AN3,AO20,AMat22, NPau20}. Many different authors have studied the case of different physical phenomena like instabilities (\cite{FL, CCFG1, CCFG,CFG23,ACF21, LS}), local regularity/singularity (\cite{JS1, JS2, Pad,APW23}), fingering, mixing or nonuniqueness \cite{Otto,CCF21, CFM22, CFP}, self-similar solutions (\cite{GGNP22}), desingularization (\cite{GGNH23}), Muskat bubbles (\cite{GGPS24,GGPS25}) in different physical scenarios like in the one phase flow Muskat problem with different characteristics (\cite{GGPS19, DGNgyu, DGNgyu2}). We refer to the review articles \cite{GrL, G17,TA} for more details and where other free boundary problems like the water waves equation (see \cite{DL,ABZ11}) are also mentioned.  Different type of singularities are possible in the Muskat problem without surface tension but some remained unclear if one considers surface tension effect. We refer to \cite{GGS20, FN,TA20} for some results when surface tension is considered.\\
  
   Regarding the Cauchy problem for the Muskat equation with surface tension, Huy Nguyen \cite{Ngu20} proved global existence for any data in the largest $L^2$-based subcritical Sobolev spaces in a very general setting, namely, for one fluid or two fluids, with or without viscosity jump and any dimension. Using abstract functional analysis for parabolic equations,  A-V Matioc and B-V Matioc \cite{Mat2} (see also \cite{Mat3}) were able to extend this result to a more general class of Sobolev space. They studied the Cauchy problem for any data belonging to the largest subcritical Sobolev space based on $L^p$. Recall that the Muskat problem with surface tension has the following scaling invariance, if $f$ is a solution to \eqref{CP0} with $g=0$, then, for any $\lambda>0$, the function  
\[
f_\lambda(t,x):=\lambda^{-1}f(\lambda^3  t,\lambda x)  
\]  
 is also a solution of  \eqref{CP0} when $g=0$.  For example, $\dot W^{1,\infty}$ and $\dot H^{\frac{3}{2}}$ are example of critical spaces. \\

     The main goal of this paper is to study the Cauchy problem for the Muskat equation in the case of two fluids and when surface tension is taken into account. We consider a smooth enough initial data which is small in some critical  space which allows the data to be arbitrarily large in the Lipschitz space. The main effort is to prove the {\it{a priori}} estimates in the critical Sobolev space $\dot H^{3/2}$ then starting with regular enough data which is in $H^{5/2}$ for instance, one proves global existence and uniqueness. One may consider a regularization of the Muskat equation  (e.g. \cite{AN3}) and then pass to the limit in the equation. As the solutions have enough compactness then the passage to the limit is classical and therefore omitted. As well, by doing the same estimate for the difference of two solutions emanating from the same initial data (in $\dot H^{3/2}$ for instance) the uniqueness can be obtained in a classical manner and therefore we shall only focus on the proof of the {\it{a priori}} estimates. More precisely, we are going to prove the following theorem.
     
\begin{theorem}
There exists a universal constant $C>0,$ such that for any $f_0 \in  H^{\frac{5}{2}}$ satisfying 
\begin{equation}
\left(C\Vert f_0 \Vert_{\dot H^{\frac{3}{2}}}+ C\Vert f_0 \Vert^4_{\dot H^{\frac{3}{2}}} + \Vert f_0 \Vert_{\dot H^{\frac{3}{2}}} \Vert f_0 \Vert_{\dot B^{1}_{\infty,1}}\right)\left(1+\Vert f_0\Vert^2_{\dot W^{1,\infty}}\right)^{\frac{3}{2}}<1,
\end{equation}  
there exists a unique global solution $f$ to the Muskat problem with surface tension which verifies
$$
f \in L^\infty([0,\infty], H^{\frac{5}{2}}) \cap L^2([0,\infty], \dot H^4 \cap \dot H^{3}).
$$
\end{theorem}
\begin{remark} One  readily notices that the homogeneous Besov semi-norm $\dot B^1_{\infty,1}$ and the Lipschitz semi-norm  may be chosen arbitrarily large by choosing the critical Sobolev semi-norm sufficiently small. It is easy to give explicitly  the value of the universal constant $C>0$ but this is not the goal of our result. The constant $C$ in our main theorem is not greater than $\Gamma(3)\Gamma(4)=12$. The data can be as regular as one wants as $H^{5/2}$ can be replace by any more regular Sobolev space with the associated gain of regularity.
\end{remark}
\begin{remark} Our proof is based on a new formulation of the Muskat problem with surface tension which allows to exhibit a nice structure of the type : elliptic term + commutator + remainders, where some of the remainders are quite singular, especially the term $T_4$.  Even though the surface tension brings more dissipation, the structure of the equation makes the Cauchy problem in the critical setting not easy to deal with. In particular, this is because of the fact that it is a weighted dissipation. The proof requires to notice many cancellations, in particular one crucial cancellations in order to deal with the term $T_{4}$ (see Lemma \ref{t.4} ) but also some cancellations which appeared earlier in the paper (see for instance Lemma \ref{id}).
\end{remark}

\section{Notations, Useful identities and Lemmas}
Throughout the paper we shall use the following notations:

 $$\Delta_\alpha f (x,t):=\frac{f(x,t)-f(x-\alpha,t)}{\alpha},$$
 
and 

$$\bar\Delta_\alpha f (x,t):=\frac{f(x,t)-f(x+\alpha,t)}{\alpha}.$$
 
We need also to introduce the sum and difference of these operators denoted respectively as follows 
  $$S_\alpha f (x,t):=\Delta_\alpha f (x,t)+\bar\Delta_\alpha f (x,t)$$ and $$D_\alpha f (x,t):=\Delta_\alpha f (x,t)-\bar\Delta_\alpha f (x,t).$$
  Therefore, 

$$S_\alpha f (x,t):=\frac{2f(x,t)-f(x+\alpha,t)-f(x-\alpha,t)}{\alpha}$$
 and 
$$D_\alpha f (x,t):=\frac{f(x+\alpha,t)-f(x-\alpha,t)}{\alpha}$$
We shall also use the notations 
\begin{equation} \label{Ops}
s_\alpha f (x,t):={2f(x,t)-f(x-\alpha,t)-f(x+\alpha,t)}
\end{equation}
and
\begin{equation} \label{Opd}
d_\alpha f (x,t):={f(x+\alpha,t)-f(x-\alpha,t)}
\end{equation}

The operators $S_\alpha$ and $D_\alpha$ will play a crucial role throughout the paper. In particular, one of the motivation is the link between these operators and the homogeneous Besov spaces. Let us recall the definition of the homogeneous Besov spaces which are actually the original one introduced by Besov.
\begin{defn}
Let $f$ be a tempered distribution so that its Fourier transform is integrable near the origin, $f$ belongs to $\dot B^{s}_{p,q}$ where $(p,q) \in [1,\infty]^2$  and $s \in (0,1)$ if and only if the following quantity is finite, that is
$$
\Vert f \Vert_{\dot B^{s}_{p,q} }=\left(\int \frac{\Vert \delta_{\alpha} f \Vert^q_{L^p}}{\vert \alpha \vert^{sq}} \ \frac{d\alpha}{\vert\alpha\vert} \right)^{\frac{1}{q}}<\infty.
$$
If $q=\infty$, then
$$
\Vert f \Vert_{\dot B^{s}_{p,\infty} }= \sup_{\alpha \in \mathbb R} \frac{\Vert \delta_{\alpha} f \Vert_{L^p}}{\vert \alpha \vert^s}<\infty.
$$
Importantly, the definition of $\dot B^{s}_{p,q}$ depends strongly on the indice of regularity $"s"$ in the following sense: if $s\in [1,2)$ then the above definition remains valid if and only if the first difference operator $\delta_{\alpha}$ is replaced by the second difference operator $s_{\alpha}$ which the operator introduced in \eqref{Ops}.

\end{defn}

Recall that homogeneous Besov spaces generalize the homogeneous Sobolev spaces, in particular,  $\dot B^{s}_{2,2}=\dot H^s$ if $s<1/2$. The following embedding between homogeneous Besov spaces are very useful as it allows to recover a control in terms of homogeneous Sobolev semi-norms  (see e.g. \cite{BCD}, \cite{PG16}). For all $(p_1,p_2,r_1,r_2) \in [1,\infty]^4$ which are such that $p_1 \leq p_2$ and $r_1 \leq r_2$ and for all $s \in \mathbb R$, we have the following continuous embedding
$$
\dot B^{s}_{p_1,r_1}(\mathbb R) \hookrightarrow \dot B^{s-(\frac{1}{p_1}-\frac{1}{p_2})}_{p_1,r_2}(\mathbb R),
$$
 If $p_1\geq2$ and $r>1$, then the left hand side of the above embedding can possibly be a homogeneous Sobolev space. This would allow one to retrieve estimates depending only on wanted critical semi-norms. Sometimes, it is not possible to get $r$ within the range $(1,\infty]$ and one has to deal with the extreme case  $r=1$. In the latter case, one may use real interpolation and again recover some Sobolev estimates (based on $L^2$). More precisely, when $r=1$, we shall use the following real interpolation inequality. There exists a constant $C>0$, such that for all $(s_1,s_2) \in \mathbb R^2$ such that $s_1<s_2$ and all $\theta \in (0,1)$  
  $$
 \Vert f \Vert_{\dot B^{\theta s_1 + (1-\theta) s_2}_{p,1}} \leq \Vert f \Vert^{\theta}_{\dot B^{s_1}_{p,\infty}}  \Vert f \Vert^{\theta}_{\dot B^{s_2}_{p,\infty}}
 $$  
for all $p\in [1,\infty]$. \\

The notation $Z=[X,Y]_{\theta, 1-\theta}$ is the classical (real or complex) interpolation where $X$, $Y$, $Z$ are any homogeneous Besov spaces. It is equivalent to $\Vert . \Vert_Z \leq C \Vert . \Vert^{\theta}_X \Vert . \Vert^{1-\theta}_Y$ where $C$ is a fixed constant. The notation $A \lesssim B$ (resp. $A\approx B$) means that $A$ is less than $B$ (resp. $A$ equals $B$) up to a fixed positive constant.   \\
 
Throughout the paper, we will need to control terms involving   $\partial_{\alpha}D_{\alpha},\partial^2_{\alpha}D_{\alpha}, \partial_{\alpha}S_{\alpha}, \partial^2_{\alpha}S_{\alpha}$  there are some nice cancellations allowing us to get expressions with finite differences. These show somehow some stability of the symmetrization through the action of differentiation. More precisely, we have the following Lemma.
\begin{lemma} \label{id}
The following identities hold,
\begin{eqnarray}  \label{eq}
 \partial_{\alpha}S_\alpha f &=&\bar\Delta_\alpha f_x-\Delta_\alpha f_x-\frac{s_\alpha f}{\alpha^2} \label{ds}\\
\partial^2_{\alpha}S_\alpha f&=&\frac{s_{\alpha} f_{xx}}{\alpha}-\frac{\int_0^\alpha s_\kappa f_{xx}(x) \ d \kappa}{\alpha^2}+\frac{d_{\alpha} f_x}{\alpha^2}+\frac{s_{\alpha} f}{\alpha^3} \label{dds} \\
 \partial_{\alpha}D_\alpha f&=&- \frac{s_\alpha f_x}{\alpha}-\frac{1}{\alpha^2} \ {\int_0^\alpha s_\kappa f_x \ d\kappa} \label{dd}\\
 \partial^2_{\alpha}D_\alpha f&=&\frac{d_\alpha f_{xx}}{\alpha}+2 \frac{s_\alpha f_x}{\alpha^2} + \frac{2}{\alpha^3}\int_0^\alpha s_{\eta}f_x \ d\eta\label{ddd}
 \end{eqnarray}
\end{lemma} 
\noindent {\bf{Proof of Lemma \ref{id}.}} The identities \eqref{ds} and \eqref{dd} have  been already obtained in \cite{CL21} but we recall them for the reader convenience and for completeness. The novelties of this lemma are the new identities for the second derivatives in $\alpha$ which are obtained through important cancellations of the bad terms. Indeed, as we have
\begin{eqnarray*}
S_\alpha f &=& \frac{f(x)-f(x-\alpha)}{\alpha}+\frac{f(x)-f(x+\alpha)}{\alpha}=\frac{2f(x)-f(x-\alpha)-f(x+\alpha)}{\alpha}
\end{eqnarray*}
We infer that,
\begin{eqnarray*} \label{dS}
\partial_{\alpha}S_\alpha f &=&\frac{f_x(x-\alpha)-f_x(x+\alpha)}{\alpha}+\left(\frac{f(x-\alpha)+f(x+\alpha)-2f(x)}{\alpha^2}\right) \nonumber \\
&=& \bar\Delta_\alpha f_x-\Delta_\alpha f_x-\frac{s_\alpha f}{\alpha^2}
\end{eqnarray*}
Next, we want to get a satisfactory expression {\it{i.e.}} involving symmetric terms  for higher derivatives, we find 
\begin{eqnarray} \label{ddS}
\partial^2_{\alpha}S_\alpha f &=&-\left(\frac{f_{xx}(x-\alpha)+f_{xx}(x+\alpha)}{\alpha}\right)+2\left(\frac{f_{xx}(x+\alpha)-f_{xx}(x-\alpha)}{\alpha^2}\right) \nonumber \\
&&\ + \ 2\left(\frac{2f(x)-f(x-\alpha)-f(x+\alpha)}{\alpha^3}\right)
\end{eqnarray}
Since, on one hand,
$$-\left(\frac{f_{xx}(x-\alpha)+f_{xx}(x+\alpha)}{\alpha}\right)=-\left(\frac{f_{xx}(x-\alpha)+f_{xx}(x+\alpha)-2f_{xx}(x)}{\alpha}\right)-2\frac{f_{xx}(x)}{\alpha}$$
and on the other hand,
$$\left(\frac{f_x(x+\alpha)-f_x(x-\alpha)}{\alpha^2}\right)=\left(\frac{\int_0^\alpha f_{xx}(x+\kappa)+f_{xx}(x-\kappa)-2f_{xx}(x) \ d\kappa}{\alpha^2}\right)+\frac{ 2 f_{xx}(x) }{\alpha}
$$
Importantly, one notices that the bad term $2\frac{f''(x) }{\alpha}$ cancels and we find
\begin{eqnarray*}
\partial^2_{\alpha}S_\alpha f &=&-\left(\frac{f_{xx}(x-\alpha)+f_{xx}(x+\alpha)-2f_{xx}(x)}{\alpha}\right) \\
&+&\left(\frac{\int_0^\alpha f_{xx}(x+\eta)+f_{xx}(x-\eta)-2f_{xx}(x) \ d \eta}{\alpha^2}\right)\\
&+& \left(\frac{f_x(x+\alpha)-f_x(x-\alpha)}{\alpha^2}\right) + \ \left(\frac{2f(x)-f(x-\alpha)-f(x+\alpha)}{\alpha^3}\right) \\
&=&\frac{s_{\alpha} f_{xx}}{\alpha}-\frac{\int_0^\alpha s_\kappa f_{xx}(x) \ d \kappa}{\alpha^2}+\frac{d_{\alpha} f_x}{\alpha^2}+\frac{s_{\alpha} f}{\alpha^3}
\end{eqnarray*}
Let us deal not with the operator $D_\alpha$ and its derivatives. We first recall that    
 $$D_\alpha f =\frac{f(x+\alpha)-f(x-\alpha)}{\alpha}=\frac{d_{\alpha} f}{\alpha}$$
One notices that 
$$D_\alpha f=\frac{1}{\alpha}\int_0^\alpha s_\eta f_x \ d\eta + 2 f_x $$
For the derivatives, we have the following identities
\begin{eqnarray*}
\partial_{\alpha}D_\alpha f&=&\frac{f_x(x+\alpha)+f_x(x-\alpha)}{\alpha}+\left(\frac{f(x-\alpha)-f(x+\alpha)}{\alpha^2}\right) \\
&=&\frac{f_x(x+\alpha)+f_x(x-\alpha)-2f_x(x)}{\alpha}+2\frac{f_x(x)}{\alpha} \\
&& + \ \frac{\int_0^\alpha -f_x(x-\kappa)-f_x(x+\kappa)+2f_x(x) \ d\kappa}{\alpha^2}-2\frac{f_x(x)}{\alpha}\\
&=& - \frac{s_\alpha f_x}{\alpha}-\frac{1}{\alpha^2} \ {\int_0^\alpha s_\kappa f_x \ d\kappa}
\end{eqnarray*}
and 
\begin{eqnarray*}
\partial^2_{\alpha}D_\alpha f &=&\frac{f_{xx}(x+\alpha)-f_{xx}(x-\alpha)}{\alpha}-2\left(\frac{f_x(x+\alpha)+f_x(x-\alpha)}{\alpha^2}\right) \\
&&-2\left(\frac{f(x-\alpha)-f(x+\alpha)}{\alpha^3}\right)
\end{eqnarray*}
Then, we write
\begin{eqnarray*}
-2\left(\frac{f_x(x+\alpha)+f_x(x-\alpha)}{\alpha^2}\right)=-2\left(\frac{-2f_x(x)+f_x(x+\alpha)+f_x(x-\alpha)}{\alpha^2}\right)-4\frac{f_x(x)}{\alpha^2}
\end{eqnarray*}
we may also write,
\begin{eqnarray*}
-2\left(\frac{f(x-\alpha)-f(x+\alpha)}{\alpha^3}\right)=-2\left(\frac{\int_0^\alpha -f_x(x-s)-f_x(x+s)+2f_x(x) \ ds}{\alpha^3}\right)+4\frac{f_x(x)}{\alpha^2}
\end{eqnarray*}
Crucially, the bad term $4\frac{f_x(x)}{\alpha^2}$ cancels again and this  allowing us to get the following very useful identity 
\begin{eqnarray*}
\partial^2_{\alpha}D_\alpha f &=&\frac{f_{xx}(x+\alpha)-f_{xx}(x-\alpha)}{\alpha}+2\left(\frac{2f_x(x)-f_x(x+\alpha)-f_x(x-\alpha)}{\alpha^2}\right) \\
&& + \ 2\left(\frac{\int_0^\alpha f_x(x-\kappa)+f_x(x+ \kappa)-2f_x(x) \ d\kappa}{\alpha^3}\right) \\
&=&\frac{d_\alpha f_{xx}}{\alpha}+2 \frac{s_\alpha f_x}{\alpha^2} + \frac{2}{\alpha^3}\int_0^\alpha s_{\eta}f_x \ d\eta.
\end{eqnarray*}
Hence,  Lemma \ref{id} is proven. 
\qed \\

Throughout the article, we shall use the two following operators : the Hilbert transform, and the 
 $$\displaystyle \mathcal{H}f=\frac{1}{\pi} \ p.v. \int \frac{f(x-y)}{y} \ dy,$$
 
The notation $\Lambda^s$ for any $s\in \mathbb R$ is the fractional Laplacian operator defined through the Fourier transform {\it{via}} the formula $\mathcal{F}\left((-\Delta)^{\alpha} f \right)=\vert \xi \vert^{2\alpha} \mathcal{F}(f).$ We shall use the following proposition whose proof can be found in \cite[Proposition 3.8]{DL19}
\begin{proposition} \label{bmo-est} Let $\mathcal{H}$ denotes the Hilbert transform, then for any $p\in (1,\infty)$ and any integers $k\geq0$ and $m\geq0$,  the following estimate holds
$$
\left\Vert\partial_x^k \left[\mathcal{H}, g \right] \partial_x^m f \right\Vert_{L^p} \leq C(k,m,p) \Vert \partial^{k+m}_x g \Vert_{BMO} \Vert f \Vert_{L^p} 
$$
\end{proposition}

\section{A new formulation  using oscillatory integrals}
It is not difficult to check that the Muskat Equation with surface tension $(\mathcal{M}_{\sigma})$ (see \eqref{CP0}) is equivalent to the following one
\begin{eqnarray} \label{CP1}
\left\{ \begin{array}{ll}
        \displaystyle\partial_t f =\displaystyle\frac{k\sigma}{2\pi \mu} \ {\rm{p.v.}} \displaystyle\int\frac{1+\partial_x f(x) \Delta_\alpha f(x)}{1+\Delta^2_\alpha f(x)}\partial_x \left(\frac{\partial_x^2 f(x-\alpha)}{(1+ f(x-\alpha)^2)^{\frac{3}{2}}}\right)\frac{d\alpha}{\alpha} \\
       \ \ \ \ \ \  + \ \displaystyle\frac{g\rho}{\pi} \ {\rm{p.v.}}\int \partial_\alpha \arctan \Delta_\alpha f \ d\alpha, \\
        f(x,0) = f_0(x).
    \end{array}
\right. 
\end{eqnarray}
When $\sigma=0$ one recovers the classical  2D Muskat Equation. As we shall do estimate in the critical Sobolev space $\dot H^{\frac{3}{2}}$ (and also in $\dot H^{\frac{5}{2}}$) it suffices to obtain {\it{a priori}} estimate of the first part of the PDE, that is, for $g=0$. If $\sigma=0$, we refer to C\'ordoba-Lazar \cite{CL21}. Hence, we shall study the following Cauchy problem for 
\begin{eqnarray} \label{CP1}
\left\{ \begin{array}{ll}
        \displaystyle\partial_t f =\displaystyle\frac{\sigma}{\pi} \ \mbox{p.v.} \displaystyle\int\frac{1+\partial_x f(x) \Delta_\alpha f(x)}{1+\Delta^2_\alpha f(x)}\partial_x \left(\frac{\partial_x^2 f(x-\alpha)}{(1+ f_x(x-\alpha)^2)^{\frac{3}{2}}}\right)\frac{d\alpha}{\alpha} \\
        f(x,0) = f_0(x)
    \end{array}
\right. 
\end{eqnarray}
Then, at the end we will add  the  {\it{a priori}} for the Muskat problem when $\sigma=0$. To prove the  {\it{a priori}} for the Cauchy problem \eqref{CP1} for data in the critical Sobolev space and for more regular data, we are going to use the following reformulation of the problem using oscillatory integrals. More precisely, we shall prove the following Proposition.
\begin{proposition} \label{NF} The Cauchy problem \eqref{CP1}  and the one below are equivalent
\begin{eqnarray*} 
\left\{
    \begin{array}{ll}\vspace{0.3cm}
\partial_t f =& \underbrace{\displaystyle\partial_x \left[\mathcal{H}, \int_{0}^{\infty} e^{-\sigma}\cos(\sigma f_x )\cos(\arctan(f_{x}))   \ d \sigma \right]\partial_{xx} f }_{{\rm{commutator}}} \\
\vspace{0.5cm}
& \  - \ \underbrace{\displaystyle\Lambda^{3} f \ \int_{0}^{\infty} e^{-\sigma}\cos(\sigma f_x )\cos(\tau_{\alpha}\arctan(f_{x}))   \ d \sigma}_{{\rm{elliptic \ component}}}     \\
\vspace{0.5cm}
& \ + \ \displaystyle \mathcal{H}f_{xx} \ \partial_{x}\left(\int_{0}^{\infty} e^{-\sigma}\cos(\sigma f_x )\cos(\arctan(f_{x}))  \ d \sigma \right)\\
\vspace{0.5cm}
& \  - \displaystyle\frac{1}{4\pi}\displaystyle\int\int_{0}^{\infty}\int_{0}^{\infty}\int_{0}^{\alpha} \ e^{-\gamma-\sigma} \ \frac{1}{\alpha} s_{\eta} f_{x}  \left(\sin(\gamma\Delta_\alpha f)-\sin(\gamma\bar\Delta_\alpha f)\right) \\  
\vspace{0.5cm}
& \ \ \ \ \times \ \displaystyle\partial_x \left(\frac{\partial_x^2 \tau_{\alpha} f}{\alpha} \cos(\sigma \tau_{\alpha} f_x )\cos(\arctan(\tau_{\alpha} f_{x}))\right) \ d\eta \ d\alpha \ d\gamma \ d\sigma  \\
\vspace{0.5cm}
& \  - \displaystyle\frac{1}{4\pi}\displaystyle\int\int_{0}^{\infty}\int_{0}^{\infty}\int_{0}^{\alpha} \ e^{-\gamma-\sigma} \ \frac{1}{\alpha} s_{\eta} f_{x}  \left(\sin(\gamma\Delta_\alpha f)+\sin(\gamma\bar\Delta_\alpha f)\right)  \\
\vspace{0.5cm}  
& \ \ \ \ \times \ \displaystyle \underbrace{\partial_x \left(\frac{\partial_x^2 \tau_{\alpha} f}{\alpha} \cos(\sigma \tau_{\alpha} f_x )\cos(\arctan(\tau_{\alpha} f_{x}))\right)  \ d\eta \ d\alpha \ d\gamma \ d\sigma}_{{\rm{singular \ term}}}  \\
\vspace{0.5cm}
& \  - \displaystyle\frac{1}{4\pi}\displaystyle\int\int_{0}^{\infty}\int_{0}^{\infty} \ e^{-\gamma-\sigma} \ S_{\alpha} f \left(\sin(\gamma\Delta_\alpha f)-\sin(\gamma\bar\Delta_\alpha f)\right)  \\
\vspace{0.3cm}
& \ \ \ \ \times \  \underbrace{\displaystyle\partial_x \left(\frac{\partial_x^2 \tau_{\alpha} f}{\alpha} \cos(\sigma \tau_{\alpha} f_x )\cos(\arctan(\tau_{\alpha}f_{x}))\right) \ d\alpha \ d\gamma \ d\sigma}_{{\rm{singular \ term}}}  \\
\vspace{0.3cm}
&  \ - \displaystyle\frac{1}{4\pi}\displaystyle\int\int_{0}^{\infty}\int_{0}^{\infty} \ e^{-\gamma-\sigma} \ S_{\alpha} f \left(\sin(\gamma\Delta_\alpha f)+\sin(\gamma\bar\Delta_\alpha f)\right) \\
& \ \ \ \ \times \  \displaystyle\partial_x \left(\frac{\partial_x^2 \tau_{\alpha} f}{\alpha} \cos(\sigma \tau_{\alpha} f_x )\cos(\arctan(\tau_{\alpha}f_{x}))\right) \ d\alpha \ d\gamma \ d\sigma \\ \\
f(x,0) &= f_0(x)
    \end{array}
\right. 
\end{eqnarray*}
\end{proposition}
\noindent {\bf{Proof of Proposition \ref{NF}.}} As we want to move from the class of rational functions of the slope to the class of oscillatory functions which have somehow better additive properties  we use the elementary identities $\int_0^\infty e^{-\gamma} \cos(r\gamma)  \ d\gamma =\frac{1}{1+r^2}$ and $\cos(\arctan r)=\frac{1}{\sqrt{1+r^2}}$. One may check that the evolution equation may be rewritten as          
\begin{eqnarray*}
\partial_t f&=&\frac{1}{\pi}\int\int_{0}^{\infty} \int_{0}^{\infty}   \cos(\gamma \Delta_\alpha f) \ e^{-\gamma-\sigma}  \\
&& \ \times \  \partial_x \left(\frac{\partial_x^2 f(x-\alpha)}{\alpha} \cos(\sigma \tau_\alpha f_x)\cos(\arctan(\tau_\alpha f_{x}))\right)  \ d\gamma \ d\sigma  \ d\alpha  \\
&+& \frac{1}{\pi}  f_{x} \int \Delta_{\alpha}f\ \cos(\gamma\Delta_\alpha f) \ e^{-\gamma-\sigma}  \\
&& \ \times \  \partial_x \left(\frac{\partial_x^2 \tau_\alpha f}{\alpha} \cos(\sigma\tau_\alpha f_x )\cos(\arctan(\tau_\alpha f_{x}))\right) \ d\gamma \ d\sigma  \ d\alpha \\
\end{eqnarray*}
Then, we linearize the first term as the idea to have a singular integral with nice properties plus a controlled term, more precisely, one uses the decomposition  $\cos(\gamma \Delta_\alpha f)=1-2 \sin(\frac{\gamma}{2}\Delta_\alpha f)$ to get that
\begin{eqnarray*}
\partial_t f&=&\underbrace{\frac{1}{\pi}\int\int_{0}^{\infty} \ e^{-\sigma}  \    \partial_x \left(\frac{\partial_x^2  \tau_\alpha f}{\alpha} \cos(\sigma  \tau_\alpha f_x )\cos(\arctan( \tau_\alpha f_{x}))\right) \ d\sigma \ d\alpha}_{\rm Hibert \ transform}    \\
&-&  \frac{2}{\pi}\int\int_{0}^{\infty} \int_{0}^{\infty}  \sin^2 (\frac{1}{2}\Delta_\alpha (\gamma f)(x))  \ e^{-\gamma-\sigma}  \\
&& \ \times \ \partial_x \left(\frac{\partial_x^2  \tau_\alpha f }{\alpha} \cos(\sigma  \tau_\alpha f_x )\cos(\arctan( \tau_\alpha f_{x}))\right) \ d\gamma \ d\sigma\ d\alpha  \\
&+&\frac{1}{\pi}\int \int_{0}^{\infty} \int_{0}^{\infty} f_{x} \Delta_{\alpha}f(x) \ \cos(\gamma\Delta_\alpha f) \ e^{-\gamma-\sigma}  \\
&& \ \times \  \partial_x \left(\frac{\partial_x^2  \tau_\alpha f}{\alpha} \cos(\sigma  \tau_\alpha f_x )\cos(\arctan( \tau_\alpha f_{x}))\right) \ d\gamma \ d\sigma\ d\alpha \\
\end{eqnarray*}
Crucially, one notices that a Hilbert transform appears in the first term above  and we find
\begin{eqnarray*}
\partial_t f &=&\int_{0}^{\infty} e^{-\sigma}    \partial_x \left(\mathcal{H}\left(\partial_x^2 f \cos(\sigma f_x )\cos(\arctan(f_{x}))\right)\right)   \ d\sigma  \\
&-&  \frac{2}{\pi}\int\int_{0}^{\infty}\int_{0}^{\infty}  \sin^2 (\frac{1}{2}\Delta_\alpha f) \ e^{-\gamma-\sigma}  \   \partial_x \left(\frac{\partial_x^2 \tau_\alpha f}{\alpha} \cos(\gamma \sigma \tau_\alpha f_x )\cos(\arctan(\tau_\alpha f_{x}))\right) \ d\gamma \ d\sigma\ d\alpha  \\
&+& \frac{1}{\pi}\int\int_{0}^{\infty}\int_{0}^{\infty} f_x \ \Delta_{\alpha}f(x) \ \cos(\gamma \Delta_\alpha f)  \ e^{-\gamma-\sigma} \ \partial_x \left(\frac{\partial_x^2 \tau_\alpha f}{\alpha} \cos(\sigma \tau_\alpha f_x )\cos(\arctan(\tau_\alpha f_{x}))\right)  \ d\gamma \ d\sigma \ d\alpha
\end{eqnarray*}
Importantly, one observes that the first term above
$$
I_1:= \int_{0}^{\infty} e^{-\sigma}    \partial_x \left(\mathcal{H}\left(\partial_x^2 f \cos(\sigma f_x )\cos(\arctan(f_{x}))\right)\right)   \ d\sigma,  
$$
 can be written as the sum of a  commutator plus a weighted dissipative term. More precisely, we may write 
\begin{eqnarray*}
I_{1} &=& \partial_x \left(\mathcal{H}\left(\partial_{xx} f \int_{0}^{\infty} e^{-\sigma}\cos(\sigma f_x )\cos(\arctan(f_{x}))\ d\sigma\right)  \right)  \\
&=&  \partial_x \left[\mathcal{H},\int_{0}^{\infty}  e^{-\sigma}\cos(\sigma f_x )\cos(\arctan(f_{x}))  \ d \sigma \right] f_{xx} \\
&& \ + \ \partial_x \left(\mathcal{H}  f_{xx}\  \int_{0}^{\infty} e^{-\sigma}\cos(\sigma f_x )\cos(\arctan(f_{x}))  \ d \sigma \right)
\end{eqnarray*}
Therefore, since $\partial_{x}\mathcal{H}f=\Lambda f$ and $\Lambda^2 f= -f_{xx}$, we find
\begin{eqnarray*}
I_{1} &=&  \partial_x \left[\mathcal{H},\int_{0}^{\infty}  e^{-\sigma}\cos(\sigma f_x )\cos(\arctan(f_{x}))   \ d \sigma \right] f_{xx} \\
&+& \mathcal{H}f_{xxx} \ \int_{0}^{\infty} e^{-\sigma}\cos(\sigma f_x )\cos(\arctan(f_{x}))   \ d \sigma \\
&+& \mathcal{H}f_{xx} \ \partial_{x}\left(\int_{0}^{\infty} e^{-\sigma}\cos(\sigma f_x )\cos(\arctan(f_{x}))  \ d \sigma \right) 
\end{eqnarray*}
Hence, one may reformulate the equation as follows,

\begin{eqnarray*}
\partial_t f &=& \partial_x \left[\mathcal{H}, \int_{0}^{\infty} e^{-\sigma}\cos(\sigma f_x )\cos(\arctan(f_{x})) \ d \sigma \right]\partial_{xx} f \\
&-& \Lambda^{3} f \ \int_{0}^{\infty} e^{-\sigma}\cos(\sigma f_x )\cos(\arctan(f_{x}))   \ d \sigma \\
&+& \mathcal{H}f_{xx} \ \partial_{x}\left(\int_{0}^{\infty} e^{-\sigma}\cos(\sigma f_x )\cos(\arctan(f_{x}))  \ d \sigma \right) \\
&-&  \frac{2}{\pi}\int\int_{0}^{\infty} \int_{0}^{\infty}  \sin^2 (\frac{\gamma}{2}\Delta_\alpha f)  \ e^{-\gamma-\sigma}  \\
&&\ \times \   \partial_x \left(\frac{\partial_x^2  \tau_\alpha f}{\alpha} \cos(\sigma \tau_\alpha f_x)\cos(\arctan(\tau_\alpha f_{x}))\right) \ d\alpha \ d\gamma \ d\sigma  \\
&+& \frac{1}{\pi}\int\int_{0}^{\infty}\int_{0}^{\infty} f_{x} \
\Delta_{\alpha}f \ \cos(\gamma\Delta_\alpha f)  \ e^{-\gamma-\sigma}  \\
&&\ \times \ \partial_x \left(\frac{ \partial_x^2 \tau_\alpha f}{\alpha} \cos(\sigma \tau_\alpha f_x)\cos(\arctan(\tau_\alpha f_{x}))\right) \ d\alpha \ d\gamma \ d\sigma 
\end{eqnarray*}
The last term above can be rewritten by using the  fact that 
  $$ \displaystyle \int_0^\infty e^{-\gamma} {\Delta_\alpha f} \ {\cos(\gamma\Delta_{\alpha}f)}   \ d\gamma=  \int_0^\infty e^{-\gamma}\sin(\gamma\Delta_\alpha f) \ d\gamma.$$ 
Therefore, we get
\begin{eqnarray*}
\partial_t f &=& \partial_x \left[\mathcal{H},\int_{0}^{\infty} \int_{0}^{\infty} e^{-\gamma-\sigma}\cos(\sigma f_x )\cos(\arctan(f_{x})) \ d \gamma  \ d \sigma \right] f_{xx} \\
&-& \Lambda^{3} f \ \int_{0}^{\infty} e^{-\sigma}\cos(\sigma f_x )\cos(\arctan(f_{x}))   \ d \sigma \\
&+& \mathcal{H}f_{xx} \ \partial_{x}\left(\int_{0}^{\infty} e^{-\sigma}\cos(\sigma f_x )\cos(\arctan(f_{x}))  \ d \sigma \right) \\
&-&  \frac{2}{\pi}\int\int_{0}^{\infty} \int_{0}^{\infty} e^{-\gamma-\sigma} \sin^2 (\frac{\gamma}{2}\Delta_\alpha f)  \\
 && \ \times \  \partial_x \left(\frac{\partial_x^2  \tau_\alpha f}{\alpha} \cos(\sigma \tau_\alpha f_x)\cos(\arctan(\tau_\alpha f_{x}))\right) \ d\alpha \ d\gamma \ d\sigma  \\
&+& \frac{1}{\pi}\int\int_{0}^{\infty}\int_{0}^{\infty} \ e^{-\gamma-\sigma} \ f_{x}(x)
 \ \sin(\gamma\Delta_\alpha f) \\
 && \ \times \  \partial_x \left(\frac{\partial_x^2 \tau_\alpha f}{\alpha} \cos(\sigma \tau_\alpha f_x)\cos(\arctan(\tau_\alpha f_{x}))\right) \ d\alpha \ d\gamma \ d\sigma. \\
\end{eqnarray*}
One notices that the third term above may be rewritten by using the following useful identity 
$$
-2\int_{0}^{\infty} e^{-\gamma} \sin^2 (\frac{1}{2}\Delta_\alpha (\gamma f)(x)) \ d\gamma= -\Delta_\alpha f \int_{0}^{\infty} e^{-\gamma} \sin (\gamma\Delta_\alpha f) \ d\gamma.
$$
Consequently, we find
\begin{eqnarray*}
\partial_t f &=& \partial_x \left[\mathcal{H},\int_{0}^{\infty} \int_{0}^{\infty} e^{-\gamma-\sigma}\cos(\sigma f_x )\cos(\arctan(f_{x})) \ d \gamma  \ d \sigma \right] f_{xx} \\
&-& \Lambda^{3} f \ \int_{0}^{\infty} \int_{0}^{\infty} e^{-\gamma-\sigma}\cos(\sigma f_x )\cos(\arctan(f_{x})) \ d \gamma  \ d \sigma    \\
&+& \mathcal{H}f_{xx} \ \partial_{x}\left(\int_{0}^{\infty} e^{-\sigma}\cos(\sigma f_x )\cos(\arctan(f_{x}))  \ d \sigma \right) \\
&+& \mathcal{H} f_{xx} \ \int_{0}^{\infty} \int_{0}^{\infty} e^{-\gamma-\sigma}\partial_{x}\left(\cos(\sigma f_x )\cos(\arctan(f_{x}))\right) \ d \gamma  \ d \sigma    \\
&+& \frac{1}{\pi}\int\int_{0}^{\infty}\int_{0}^{\infty} \ e^{-\gamma-\sigma} \  (f_{x}-\Delta_{\alpha}f) \sin(\gamma\Delta_\alpha f) \\
&& \ \times \ \partial_x \left(\frac{\partial_x^2 \tau_\alpha f }{\alpha} \cos(\sigma \tau_\alpha f_x )\cos(\arctan(\tau_\alpha f_{x})\right) \ d\alpha \ d\gamma \ d\sigma. \\
\end{eqnarray*} 
Then, we use the fact that 
\begin{equation}
f_x-\Delta_\alpha f = f_x-\frac{1}{2}S_{\alpha} f -\frac{1}{2}D_{\alpha} f =f_x-\frac{1}{2}S_{\alpha} f -\frac{1}{2\alpha} \int_{0}^{\alpha} s_{\eta} f_{x} \ d\eta - f_x
\end{equation}
where one notices a crucial cancellation of the bad term $f_{x}$.
Hence, the evolution equation may be rewritten as
\begin{eqnarray*}
\partial_t f &=& \partial_x \left[\mathcal{H}, \int_{0}^{\infty} e^{-\sigma}\cos(\sigma f_x )\cos(\arctan(f_{x}))   \ d \sigma \right] f_{xx} \\ 
&-&  \ \Lambda^{3} f \  \int_{0}^{\infty} e^{-\sigma}\cos(\sigma f_x )\cos(\arctan(f_{x}))   \ d \sigma    \\
&+& \mathcal{H} f_{xx} \ \int_{0}^{\infty}  e^{-\sigma}\partial_{x}\left(\cos(\sigma f_x )\cos(\arctan(f_{x}))\right)  \ d \sigma    \\
&-& \frac{1}{2\pi}\int\int_{0}^{\infty}\int_{0}^{\infty}\ e^{-\gamma-\sigma} \frac{1}{\alpha}\int_{0}^{\alpha} s_{\eta} f_{x}  \sin(\gamma\Delta_\alpha f) \\ && \partial_x \left(\frac{\partial_x^2 \tau_\alpha f}{\alpha} \cos(\sigma \tau_\alpha f_x )\cos(\arctan(\tau_\alpha f_{x}))\right) \ d\eta  \ d\gamma \ d\sigma \ d\alpha \\
&-& \frac{1}{2\pi}\int\int_{0}^{\infty}\int_{0}^{\infty} \ e^{-\gamma-\sigma} \ S_{\alpha} f \sin(\gamma\Delta_\alpha f) \\ 
&&  \partial_x \left(\frac{\partial_x^2 \tau_\alpha f}{\alpha} \cos(\sigma \tau_\alpha f_x )\cos(\arctan(\tau_\alpha f_{x}))\right) \ d\gamma \ d\sigma  \ d\alpha. \\
\end{eqnarray*} 

We may further symmetrize the last two terms by simply writting that 
$$\sin(\gamma\Delta_\alpha f)=\frac{1}{2}\left(\sin(\gamma\Delta_\alpha f)-\sin(\gamma\bar\Delta_\alpha f)\right)+\frac{1}{2}\left(\sin(\gamma\Delta_\alpha f)+\sin(\gamma \bar\Delta_\alpha f)\right)$$ 

Finally, the Muskat problem with surface tension may be reformulated as follows

\begin{eqnarray*}
\partial_t f &=& \partial_x \left[\mathcal{H}, \int_{0}^{\infty} e^{-\sigma}\cos(\sigma f_x )\cos(\arctan(f_{x}))   \ d \sigma \right] f_{xx} \\
&-&  \Lambda^{3} f \  \int_{0}^{\infty} e^{-\sigma}\cos(\sigma f_x )\cos(\arctan(f_{x}))  \ d \sigma    \\
&+& \mathcal{H} f_{xx} \ \int_{0}^{\infty}  e^{-\sigma}\partial_{x}\left(\cos(\sigma f_x )\cos(\arctan(f_{x}))\right)  \ d \sigma    \\
&-& \frac{1}{4\pi}\int\int_{0}^{\infty}\int_{0}^{\infty}\int_{0}^{\alpha} \ e^{-\gamma-\sigma} \ \frac{1}{\alpha} s_{\eta} f_{x}  \left(\sin(\gamma\Delta_\alpha f)-\sin(\gamma\bar\Delta_\alpha f)\right) \\  
&& \ \times \ \partial_x \left(\frac{\partial_x^2 \tau_{\alpha} f}{\alpha} \cos(\sigma \tau_{\alpha} f_x )\cos(\arctan(\tau_{\alpha} f_{x}))\right) \ d\eta  \ d\gamma \ d\sigma \ d\alpha \\
&-& \frac{1}{4\pi}\int\int_{0}^{\infty}\int_{0}^{\infty}\int_{0}^{\alpha} \ e^{-\gamma-\sigma} \ \frac{1}{\alpha} s_{\eta} f_{x}  \left(\sin(\gamma\Delta_\alpha f)+\sin(\gamma\bar\Delta_\alpha f)\right) \\  
&& \ \times \ \partial_x \left(\frac{\partial_x^2 \tau_{\alpha} f}{\alpha} \cos(\sigma \tau_{\alpha} f_x )\cos(\arctan(\tau_{\alpha} f_{x}))\right)  \ d\eta  \ d\gamma \ d\sigma \ d\alpha \\
&-& \frac{1}{4\pi}\int\int_{0}^{\infty}\int_{0}^{\infty} \ e^{-\gamma-\sigma} \ S_{\alpha} f \left(\sin(\gamma\Delta_\alpha f)-\sin(\gamma\bar\Delta_\alpha f)\right) \\
&& \ \times \  \partial_x \left(\frac{\partial_x^2 \tau_{\alpha} f}{\alpha} \cos(\sigma \tau_{\alpha} f_x )\cos(\arctan(\tau_{\alpha}f_{x}))\right)  \ d\gamma \ d\sigma \ d\alpha \\
&-& \frac{1}{4\pi}\int\int_{0}^{\infty}\int_{0}^{\infty} \ e^{-\gamma-\sigma} \ S_{\alpha} f \left(\sin(\gamma\Delta_\alpha f)+\sin(\gamma\bar\Delta_\alpha f)\right) \\
&& \ \times \  \partial_x \left(\frac{\partial_x^2 \tau_{\alpha} f}{\alpha} \cos(\sigma \tau_{\alpha} f_x )\cos(\arctan(\tau_{\alpha}f_{x}))\right)  \ d\gamma \ d\sigma \ d\alpha. \\
\end{eqnarray*} 

This ends the proof of the new formulation.

\qed

\section{{\it{A priori}} estimates in the critical Sobolev space}

We shall prove the following regularity criteria in the fully critical setting.

\begin{theorem} \label{H32}
There exists a constant $C>0$ such that for any $f_0 \in \dot H^{\frac{3}{2}} \cap \dot B^1_{\infty,1}$ 
\begin{equation}
\left(C\Vert f_0 \Vert_{\dot H^{\frac{3}{2}}}+ C\Vert f_0 \Vert^3_{\dot H^{\frac{3}{2}}} + \Vert f_0 \Vert_{\dot H^{\frac{3}{2}}} \Vert f_0 \Vert_{\dot B^{1}_{\infty,1}}\right)\left(1+\Vert f_0\Vert^2_{\dot W^{1,\infty}}\right)^{\frac{3}{2}}<1,
\end{equation} 
then,  if the semi-norms $H^{\frac{3}{2}}$ and $\dot B^1_{\infty,1}$ does not blow-up in finite time, we have for all $T>0$
\begin{eqnarray*} 
&&\Vert f(T) \Vert^2_{\dot H^{3/2}}+  \ \int_0^T \frac{\displaystyle \Vert f \Vert_{\dot H^2}^2}{(1+  \Vert f_x \Vert^2_{L^{\infty}})^2} \ ds +   \ \int_0^T \frac{\displaystyle \Vert f \Vert_{\dot H^3}^2}{(1+  \Vert f_x \Vert^2_{L^{\infty}})^{\frac{3}{2}}} \ ds \\
&& \lesssim \Vert f_0 \Vert^2_{\dot H^{3/2}} +{\displaystyle\int_0^T\Vert f \Vert^2_{\dot H^3} \left( \mathcal{P}(\Vert f \Vert_{\dot H^{3/2}})+\Vert f \Vert_{\dot H^{3/2}}  \Vert  f \Vert_{\dot B^1_{\infty,1}}\right) \ ds}+{\displaystyle\int_0^T \mathcal{Q}(\Vert f \Vert_{\dot H^{3/2}}) \Vert f \Vert_{\dot H^{3/2}} \  \ ds},
\end{eqnarray*}

where $\mathcal{P}(X)=X+X^3$ and $\mathcal{Q}(X)=X+X^2.$
 
\end{theorem}

\noindent{{\bf{Proof of Theorem \ref{H32}}}.} We study the time evolution of the critical semi-norm $\dot H^{3/2}$. Set
$\mathcal{A}(t):=\Vert f \Vert_{\dot H^{3/2}}(t)$, then 
\begin{eqnarray*}
\frac{1}{2}\partial_{t}\mathcal{A}^2(t)&=&\int \Lambda^{3/2} f \ \Lambda^{3/2}f_{t} \ dx =\int  \Lambda^{3}f \ f_{t} \ dx. \\
\end{eqnarray*}
Therefore, using the formulation with the oscillatory integrals obtained in Proposition \ref{NF}, one finds 
\begin{eqnarray*}
\frac{1}{2}\partial_{t}\mathcal{A}^2(t) &=&  -\int \vert \Lambda^{3} f \vert^2  \ \int_{0}^{\infty} e^{-\sigma}\cos(\sigma f_x )\cos(\arctan(f_{x}))   \ d \sigma \ dx    \\ 
&+& \int  \Lambda^{3} f \ \mathcal{H} f_{xx} \ \int_{0}^{\infty}  e^{-\sigma}\partial_{x}\left(\cos(\sigma f_x )\cos(\arctan(f_{x}))\right)  \ d \sigma \ dx    \\
&+& \  \int \Lambda^{3} f \ \partial_x \left[\mathcal{H},\int_{0}^{\infty} e^{-\sigma}\cos(\sigma f_x )\cos(\arctan(f_{x}))   \ d \sigma  \right] f_{xx} \ dx \\
&-& \frac{1}{4\pi} \int \int\int_{0}^{\infty}\int_{0}^{\infty} \ e^{-\gamma-\sigma} \  \Lambda^{3} f \ S_{\alpha} f \left(\sin(\gamma\Delta_\alpha f)+\sin(\gamma\bar\Delta_\alpha f)\right)   \\
&& \  \times\ \partial_x \left(\frac{\partial_x^2 \tau_{\alpha} f}{\alpha} \cos(\sigma \tau_{\alpha} f_x )\cos(\arctan(\tau_{\alpha} f_{x}))\right)   \ d\gamma \ d\sigma \ d\alpha \ dx \\
&-& \frac{1}{4\pi} \int \int\int_{0}^{\infty}\int_{0}^{\infty}\int_{0}^{\alpha} \ e^{-\gamma-\sigma} \ \Lambda^{3} f \ \frac{1}{\alpha} s_{\eta} f_{x}  \left(\sin(\gamma\Delta_\alpha f)-\sin(\gamma\bar\Delta_\alpha f)\right)  \\
&& \  \times\ \partial_x \left(\frac{\partial_x^2 \tau_{\alpha} f}{\alpha} \cos(\sigma \tau_{\alpha} f_x )\cos(\arctan(\tau_{\alpha} f_{x}))\right) \ d\eta  \ d\gamma \ d\sigma \ d\alpha \ dx \\
&-& \frac{1}{4\pi} \int\int\int_{0}^{\infty}\ \int_{0}^{\infty}\int_{0}^{\alpha} \ e^{-\gamma-\sigma} \ \Lambda^{3} f \ \frac{1}{\alpha} s_{\eta} f_{x}  \left(\sin(\gamma\Delta_\alpha f)+\sin(\gamma\bar\Delta_\alpha f)\right)    \\
&& \  \times\ \partial_x \left(\frac{\partial_x^2 \tau_{\alpha} f}{\alpha} \cos(\sigma \tau_{\alpha} f_x )\cos(\arctan(\tau_{\alpha} f_{x}))\right) \ d\eta  \ d\gamma \ d\sigma \ d\alpha \ dx \\
&-& \frac{1}{4\pi} \int \int\int_{0}^{\infty}\int_{0}^{\infty} \ e^{-\gamma-\sigma} \ \Lambda^{3} f \ S_{\alpha} f \left(\sin(\gamma\Delta_\alpha f)-\sin(\gamma\bar\Delta_\alpha f)\right)   \\
&& \  \times\ \partial_x \left(\frac{\partial_x^2 \tau_{\alpha} f}{\alpha} \cos(\sigma \tau_{\alpha} f_x )\cos(\arctan(\tau_{\alpha} f_{x}))\right)   \ d\gamma \ d\sigma \ d\alpha \ dx \\
&:=& T_{1, {\rm weighted \ dissip}} + T_{1,{\rm rem}} + \sum_{j=2}^6 T_j
\end{eqnarray*}
\subsection{The weighted dissipation term $T_{1}$}
It is not difficult to see that $T_{1,dissip}$ is indeed a weighted diffusion term. Indeed, since we have 
$$
\int_{0}^{\infty}  e^{-\sigma}\cos(\sigma f_x )\cos(\tau_{\alpha}\arctan(f_{x}))  \ d \sigma =\frac{1}{(1+ f_x^2)^{\frac{3}{2}}}
$$
Hence,
\begin{equation}
T_{1, {\rm weighted \ dissip}}=-\displaystyle\int   \ \frac{\vert\Lambda^3 f \vert^2}{(1+ \vert f_x \vert^2)^{\frac{3}{2}}} \ dx.
\end{equation}

\subsection{The remainder term $T_{1,rem}$}

We want to estimate 

\begin{eqnarray*}
T_{1,rem}&=& \int  \Lambda^{3} f \ \mathcal{H} f_{xx} \ \int_{0}^{\infty}  e^{-\gamma}\partial_{x}\left(\cos(\sigma f_x )\cos(\arctan(f_{x}))\right)  \ d \sigma \ dx.    
\end{eqnarray*}
We notice that 
\begin{eqnarray*}
T_{1,rem}&\lesssim& \Vert \Lambda^{3} f \Vert_{L^2} \Vert \mathcal{H} f_{xx} \Vert_{L^4} \int \Vert \partial_{x}\left(\cos(\sigma f_x )\cos(\arctan(f_{x}))\right)\Vert_{L^4} \ dx \\
&\lesssim&  \Vert \Lambda^{3} f \Vert_{L^2} \Vert \mathcal{H} f_{xx} \Vert_{L^4} \int_0^\infty (1+\sigma)e^{-\sigma}\Vert f_{xx} \Vert_{L^4} \ d\sigma.
\end{eqnarray*}
Using the fact that $\mathcal{H}$ is continuous on $L^p$ for any $p\in(1,\infty)$ together with the Sobolev embedding $\dot H^{\frac{1}{4}} \hookrightarrow L^4$, we find that
\begin{eqnarray*}
T_{1,rem}&\lesssim& \Vert f \Vert_{\dot H^3} \Vert f \Vert^2_{\dot H^{\frac{9}{4}}}
\end{eqnarray*}
Finally, since $\Vert f \Vert_{\dot H^{\frac{9}{4}}} \leq \Vert f \Vert^{\frac{1}{2}}_{\dot H^3} \Vert f \Vert^{\frac{1}{2}}_{\dot H^{\frac{3}{2}}},$ we conclude that
\begin{eqnarray*}
T_{1,rem}&\lesssim&\Vert f \Vert^2_{\dot H^3} \Vert f \Vert_{\dot H^{\frac{3}{2}}} \\
\end{eqnarray*}

We will now estimates the $T_j$ for $j=2,...,6$.

\subsection{Estimate of the commutator $T_2$}
Recall that,
\begin{eqnarray*}
T_2=- \int \Lambda^{3} f \ \partial_x \left[\mathcal{H}, \int_{0}^{\infty} e^{-\sigma}\cos(\sigma f_x )\cos(\arctan(f_{x})) \ d \sigma   \right] f_{xx} \ dx. \\
\end{eqnarray*}
Then, we see that 
\begin{eqnarray*}
T_2\leq \Vert f \Vert_{\dot H^3}    \left \Vert \partial_x \left[\mathcal{H},\int_{0}^{\infty}  e^{-\sigma}\cos(\sigma f_x )\cos(\arctan(f_{x}))   \ d \sigma  \right] f_{xx} \right \Vert_{L^2} \\
\end{eqnarray*}
Using the commutator estimate of Proposition \ref{bmo-est}, we have
\begin{eqnarray*}
T_2\leq \Vert f \Vert_{\dot H^3} \Vert f \Vert_{\dot H^2}   \left \Vert \partial_x \left(\int_{0}^{\infty}  e^{-\sigma}\cos(\sigma f_x )\cos(\arctan(f_{x}))   \ d \sigma  \right)  \right \Vert_{BMO} \\
\end{eqnarray*}
Using the fact that $\dot H^{1/2} \hookrightarrow BMO,$ we find
\begin{eqnarray*}
T_2&\leq& \Vert f \Vert_{\dot H^3} \Vert f \Vert_{\dot H^2} \int_{0}^{\infty} e^{-\sigma} \left\Vert \partial_x\left(\cos(\sigma f_x )\cos(\arctan(f_{x})) \right)\right\Vert_{\dot H^{\frac{1}{2}}}   \ d \sigma \\
&\leq&\Vert f \Vert_{\dot H^3} \Vert f \Vert_{\dot H^2} \int_{0}^{\infty} \sigma e^{-\sigma} \left\Vert f_{xx}\sin(\sigma f_x )\cos(\arctan(f_{x})) \right\Vert_{\dot H^{\frac{1}{2}}}   \ d \sigma \\
&& + \ \Vert f \Vert_{\dot H^3} \Vert f \Vert_{\dot H^2} \int_{0}^{\infty}  e^{-\sigma} \left\Vert \frac{f_{xx}}{1+(f_x)^2}\cos(\sigma f_x )\sin(\arctan(f_{x})) \right\Vert_{\dot H^{\frac{1}{2}}}   \ d \sigma \\
\end{eqnarray*}
The main effort will be devoted to the $\dot H^{1/2}$-estimate of these products. Note that here we want to avoid the use of the control in $L^\infty$ as it seems difficult to get the wanted critical estimate if $L^\infty$ hits $f_{xx}$. Therefore, the fact that $\dot H^{1/2} \cap L^\infty$ is an algebra is not very helpful. As well, using well-known estimates of the type $\Vert fg \Vert_{\dot H^{1/2}}\leq \Vert f \Vert_{\dot W^{1/2,p}}\Vert f \Vert_{\dot W^{1/2,p}} + \Vert g \Vert_{\dot W^{1/2,q}} \Vert f \Vert_{L^q}$ where $p$ and $q$ are conjugate exponents would give rise to subcritical terms (using the fact that the Lipschitz space acts on homogeneous Sobolev spaces of regularity $s \in (0,1)$). 
Instead, we shall write use the point-wise formula for the root of the Laplacian of a product and then estimate its $L^2$ norm. This would give more flexibility not only in spreading accordingly the derivatives but also in the integrability of each term.  More precisely, we shall prove the following lemma. Before stating the lemma, set
 $$\mathcal{N}_1:=\Lambda^{\frac{1}{2}}\left(f_{xx}\sin(\sigma f_x )\cos(\arctan(f_{x})\right) $$ 
and
 $$\mathcal{N}_2:=\Lambda^{\frac{1}{2}}\left(\frac{f_{xx}}{1+(f_x)^2}\sin(\sigma f_x )\cos(\arctan(f_{x})\right) $$
 The main effort will be devoted to the control in $L^2$ of $\mathcal{N}_1$ and $\mathcal{N}_2$
\begin{lemma} \label{Ni} We have the following the control of $\mathcal{N}_1$ and  $\mathcal{N}_2$ in $L^2$. For $i=1,2$, one has
$$
\Vert \mathcal{N}_i \Vert_{L^2} \lesssim  \Vert f\Vert_{\dot H^{\frac{5}{2}}}+\Vert f \Vert_{\dot H^{\frac{9}{4}}} \left\Vert f \right\Vert_{\dot H^{\frac{7}{4}}}+ \sigma\Vert f \Vert_{\dot H^{\frac{5}{2}}} \Vert f \Vert_{\dot H^{\frac{3}{2}}}
$$ 
\end{lemma} 
\noindent {\bf{Proof of Lemma \ref{Ni}}.} One has the identity
\begin{eqnarray*}
 \mathcal{N}_1 &=& \Lambda^{\frac{1}{2}}(f_{xx}) \sin(\sigma f_x )\cos(\arctan(f_{x})) + f_{xx} \ \Lambda^{\frac{1}{2}}(\sin(\sigma f_x )\cos(\arctan(f_{x}))) \\
 &-&C \int \frac{\delta_y f_{xx} \ \delta_y(\sin(\sigma f_x )\cos(\arctan(f_{x})))}{\vert y \vert^{\frac{3}{2}}} \ dy
\end{eqnarray*}
Hence 
\begin{eqnarray*}
 \Vert \mathcal{N}_1 \Vert_{L^2} &\lesssim& \Vert f \Vert_{\dot H^{\frac{5}{2}}} + \Vert f_{xx}\Vert_{L^4} \ \left\Vert\Lambda^{\frac{1}{2}}\sin(\sigma f_x) \right \Vert_{L^4} \\
 &+& \int \frac{\Vert \delta_y f_{xx} \Vert_{L^4} \ \Vert \delta_y(\sin(\sigma f_x ) \cos(\arctan(f_{x})))\Vert_{L^4}}{\vert y \vert^{\frac{3}{2}}} \ dy
\end{eqnarray*}
Since the function $\psi_{\sigma}(\cdot) :=\sin(\sigma \cdot ) \cos(\arctan(\cdot))$ is $\sigma$-Lipschitz, we have, by Sobolev embedding and the fact that  $\dot H^{3/4}$ is stable by composition with any outer Lipschitz function, one gets
\begin{eqnarray*}
\left\Vert\Lambda^{\frac{1}{2}}\sin(\sigma f_x) \right \Vert_{L^4} \lesssim \Vert \sin(\sigma f_x) \Vert_{\dot H^{3/4}}\lesssim \Vert f \Vert_{\dot H^{\frac{7}{4}}}
\end{eqnarray*}
and
\begin{eqnarray*}
\Vert \delta_y(\sin(\sigma f_x )\cos(\arctan(f_{x}))) \Vert_{L^4}\lesssim\Vert \delta_y \psi_\sigma f_x  \Vert_{L^4} \lesssim \sigma\Vert \delta_y  f_x \Vert_{L^4}  
\end{eqnarray*}
Hence, using the two above inequalities, we infer that
\begin{eqnarray*}
 \Vert \mathcal{N}_1 \Vert_{L^2} &\lesssim& \Vert f \Vert_{\dot H^{\frac{5}{2}}} + \Vert f\Vert_{\dot H^{\frac{9}{4}}}\Vert f \Vert_{\dot H^{\frac{7}{4}}}  +\sigma \int \frac{\Vert \delta_y f_{xx} \Vert_{L^4} \ \Vert \delta_y  f_x \Vert_{L^4} }{\vert y \vert^{\frac{3}{2}}} \ dy
\end{eqnarray*}
Therefore,
\begin{eqnarray*}
 \Vert \mathcal{N}_1 \Vert_{L^2} &\lesssim& \Vert f \Vert_{\dot H^{\frac{5}{2}}} + \Vert f\Vert_{\dot H^{\frac{9}{4}}}\Vert f \Vert_{\dot H^{\frac{7}{4}}}  +\sigma \left(\int \frac{\Vert \delta_y f_{xx} \Vert^2_{L^4}  }{\vert y \vert^{\frac{3}{2}}} \ dy\right)^{\frac{1}{2}} \left(\int \frac{ \Vert \delta_y  f_x \Vert^2_{L^4} }{\vert y \vert^{\frac{3}{2}}} \ dy\right)^{\frac{1}{2}} \\
 &\lesssim& \Vert f \Vert_{\dot H^{\frac{5}{2}}} + \Vert f\Vert_{\dot H^{\frac{9}{4}}}\Vert f \Vert_{\dot H^{\frac{7}{4}}}  +\sigma\Vert f \Vert_{\dot B^{\frac{9}{4}}_{4,2}} \Vert f \Vert_{\dot B^{\frac{5}{4}}_{4,2}}
 \end{eqnarray*} 
 Since $\dot H^{\frac{5}{2}} \hookrightarrow \dot B^{\frac{9}{4}}_{4,2}$ and $\dot H^{\frac{3}{2}} \hookrightarrow \dot B^{\frac{5}{4}}_{4,2}$, we find that, 
  \begin{eqnarray*}
\Vert \mathcal{N}_1 \Vert_{L^2}&\lesssim& \Vert f\Vert_{\dot H^{\frac{5}{2}}}+\Vert f \Vert_{\dot H^{\frac{9}{4}}} \left\Vert f \right\Vert_{\dot H^{\frac{7}{4}}}+ \sigma\Vert f \Vert_{\dot H^{\frac{5}{2}}} \Vert f \Vert_{\dot H^{\frac{3}{2}}}
\end{eqnarray*}

The estimate of $\mathcal{N}_2$ is analogous to the one of $\mathcal{N}_1$ as they are the same up to the bounded factor $(1+(f_x)^2)^{-1}$, to get similar type of terms, let us use
$$
(1+(f_x)^2)^{-1}=\int_0^\infty e^{-k} \cos(k f_x) \ dk
$$
Then, $\mathcal{N}_2$ can be written as
\begin{eqnarray*}
\mathcal{N}_2&=&\int_0^\infty e^{-k} \Lambda^{\frac{1}{2}}\left({f_{xx}}\cos(k f_x)\sin(\sigma f_x )\cos(\arctan(f_{x})\right) \ dk \\
&=&\int_0^\infty e^{-k} \cos(k f_x)\sin(\sigma f_x )\cos(\arctan(f_{x})\Lambda^{\frac{1}{2}}\left({f_{xx}}\right) \ dk \\
&+&\int_0^\infty e^{-k} {f_{xx}}\Lambda^{\frac{1}{2}}\left(\cos(k f_x)\sin(\sigma f_x )\cos(\arctan(f_{x})\right) \ dk \\
&-&C \int_0^\infty e^{-k}\int \frac{\delta_y\left( \sin(\sigma f_x )\cos(k f_x)\cos(\arctan(f_{x}))\right) \ \delta_y f_{xx}}{\vert y \vert^{\frac{3}{2}}} \ dk \ dy\\
 \end{eqnarray*}
 Then, we may estimate the $L^2$-norm of $\mathcal{N}_2$ as follows
 \begin{eqnarray*}
\Vert \mathcal{N}_2 \Vert_{L^2}&\lesssim&\int_0^\infty e^{-k} \Vert \Lambda^{\frac{1}{2}}{f_{xx}}\Vert_{L^2} \ dk \\
&+&\int_0^\infty e^{-k} \Vert {f_{xx}} \Vert_{L^4} \left\Vert \Lambda^{\frac{1}{2}}\left(\cos(k f_x)\sin(\sigma f_x )\cos(\arctan(f_{x})\right) \right\Vert_{L^4} \ dk \\
&+& \int_0^\infty e^{-k}\int \frac{\Vert \delta_y\left( \sin(\sigma f_x )\cos(k f_x)\cos(\arctan(f_{x}))\right)\Vert_{L^4} \ \Vert \delta_y f_{xx} \Vert_{L^4}}{\vert y \vert^{\frac{3}{2}}} \ dk \ dy\\
 \end{eqnarray*}
 In order to estimate the second integral appearing in the right hand side, one notices that the function $\phi_{k,\sigma}(\cdot)=\cos(k \cdot)\sin(\sigma \cdot )\cos(\arctan(\cdot))$ is a product of bounded  Lipschitz functions  and it is easy to check that it is a  $k\sigma$-Lipschitz function. In particular, we use again that the composite function $\phi_{k,\sigma}$ with any $\dot H^s$ function is sublinear with a constant which depends only on the outer function (see e.g. \cite{BMey}). Hence, using the Sobolev embedding $\dot H^{1/4} \hookrightarrow L^4$ one finds
 \begin{eqnarray*}
\Vert \mathcal{N}_2 \Vert_{L^2}&\lesssim& \Vert f\Vert_{\dot H^{\frac{5}{2}}} + \int_0^\infty \sigma ke^{-k}\Vert f \Vert_{\dot H^{\frac{9}{4}}} \left\Vert f_x \right\Vert_{\dot H^{\frac{3}{4}}} \ dk +\int_0^\infty \sigma ke^{-k} \ dk \int  \frac{\Vert \delta_y f_x\Vert_{L^4} \ \Vert \delta_y f_{xx} \Vert_{L^4}}{\vert y \vert^{\frac{3}{2}}} \ dy\\
&\lesssim& \Vert f\Vert_{\dot H^{\frac{5}{2}}}+\Vert f \Vert_{\dot H^{\frac{9}{4}}} \left\Vert f_x \right\Vert_{\dot H^{\frac{3}{4}}}+ \sigma\left(\int  \frac{\Vert \delta_y f_{xx}\Vert^2_{L^4}}{\vert y \vert^{\frac{3}{2}}} \ dy\right)^{\frac{1}{2}} \left(\int  \frac{\Vert \delta_y f_x\Vert^2_{L^4}}{\vert y \vert^{\frac{3}{2}}} \ dy\right)^{\frac{1}{2}} \\
&\lesssim& \Vert f\Vert_{\dot H^{\frac{5}{2}}}+\Vert f \Vert_{\dot H^{\frac{9}{4}}} \left\Vert f \right\Vert_{\dot H^{\frac{7}{4}}}+ \sigma \Vert f \Vert_{\dot B^{\frac{9}{4}}_{4,2}} \Vert f \Vert_{\dot B^{\frac{5}{4}}_{4,2}}
 \end{eqnarray*} 
 Hence, since $\dot H^{\frac{5}{2}}\hookrightarrow \dot B^{\frac{9}{4}}_{4,2}$ and $\dot H^{\frac{3}{2}}\hookrightarrow \dot B^{\frac{5}{4}}_{4,2}$, we obtain  
  \begin{eqnarray*}
\Vert \mathcal{N}_2 \Vert_{L^2}&\lesssim& \Vert f\Vert_{\dot H^{\frac{5}{2}}}+\Vert f \Vert_{\dot H^{\frac{9}{4}}} \left\Vert f \right\Vert_{\dot H^{\frac{7}{4}}}+ \sigma\Vert f \Vert_{\dot H^{\frac{5}{2}}} \Vert f \Vert_{\dot H^{\frac{3}{2}}}
\end{eqnarray*}

Now, the estimate of $T_2$ is done through interpolation, indeed using the Lemma \ref{Ni},  we find

\begin{eqnarray*}
T_2&\lesssim& \Vert f \Vert_{\dot H^3} \Vert f \Vert_{\dot H^2}\left(\Vert f\Vert_{\dot H^{\frac{5}{2}}}+\Vert f \Vert_{\dot H^{\frac{9}{4}}} \left\Vert f \right\Vert_{\dot H^{\frac{7}{4}}}+ \Vert f \Vert_{\dot H^{\frac{5}{2}}} \Vert f \Vert_{\dot H^{\frac{3}{2}}}\right) \\
\end{eqnarray*}
Where the $\sigma$ has disappeared as we have integrated against the positive measure $e^{-\sigma} \ d\sigma$ giving therefore rise to well-defined integrals (actually Gamma functions). Finally, using interpolation, we find
\begin{eqnarray*}
T_2&\lesssim& \Vert f \Vert_{\dot H^3} \Vert f \Vert^\frac{1}{3}_{\dot H^3}\Vert f \Vert^{\frac{2}{3}}_{\dot H^\frac{3}{2}}\left(\Vert f \Vert^\frac{2}{3}_{\dot H^3}\Vert f \Vert^{\frac{1}{3}}_{\dot H^\frac{3}{2}}+\Vert f \Vert^\frac{1}{2}_{\dot H^3}\Vert f \Vert^{\frac{1}{2}}_{\dot H^\frac{3}{2}} \Vert f \Vert^\frac{1}{6}_{\dot H^3}\Vert f \Vert^{\frac{5}{6}}_{\dot H^\frac{3}{2}}+ \Vert f \Vert^\frac{2}{3}_{\dot H^3}\Vert f \Vert^{\frac{1}{3}}_{\dot H^\frac{3}{2}} \Vert f \Vert_{\dot H^{\frac{3}{2}}}\right) \\
\end{eqnarray*}
Consequently, we have obtained 
\begin{eqnarray*}
T_2&\lesssim&\Vert f \Vert^2_{\dot H^3} \left(\Vert f\Vert_{\dot H^{\frac{3}{2}}}+\Vert f\Vert^2_{\dot H^{\frac{3}{2}}}\right)
\end{eqnarray*}

 \subsection{Estimate of $T_3$}
 
Recall that,

\begin{eqnarray*}
T_3&=& - \frac{1}{4\pi} \int \int\int_{0}^{\infty}\int_{0}^{\infty}  \ e^{-\gamma-\sigma} \  \Lambda^{3} f \ S_{\alpha} f \left(\sin(\gamma\Delta_\alpha f)+\sin(\gamma\bar\Delta_\alpha f)\right)   \\
&& \  \times\ \partial_x \left(\frac{\partial_x^2 \tau_{\alpha} f}{\alpha} \cos(\sigma \tau_{\alpha} f_x )\cos(\arctan(\tau_{\alpha} f_{x}))\right)   \ d\gamma \ d\sigma \ d\alpha \ dx \\
\end{eqnarray*}
Even though  the translation $\tau_\alpha$ is not regular enough because of the fact that the whole term is highly nonlinear and nonlocal, this term is still regular enough. This is mainly due to the presence of the second finite differences which allow to spread the regularity in $\alpha$ with more flexibility. In particular, one may use Besov estimates with regularity bigger than 1 (and strictly less than 2). Hence,  one does not need to make appear the more regular operator $\Delta_\alpha$ in stead of the translation $\tau_\alpha$ one. Indeed,
by using the following straightforward inequalities for any $\sigma\geq0$, $\gamma\geq0$: $$\left\vert\sin(\gamma\Delta_\alpha f)+\sin(\gamma\bar\Delta_\alpha f)\right \vert\leq \gamma\vert S_{\alpha}f\vert,$$
and $$\vert \partial_x \cos(\sigma \tau_{\alpha} f_x) \vert \leq \sigma \vert \tau_{\alpha}f_{xx} \vert,$$ as well as $$\vert \partial_x\cos(\arctan(\tau_{\alpha} f_{x}))\vert \leq \vert \tau_{\alpha} f_{xx} \vert,$$
one may write that
\begin{eqnarray*}
 T_3 &\lesssim&  \int \int\int_{0}^{\infty}\int_{0}^{\infty}   \ \gamma(1+\sigma) e^{-\gamma-\sigma} \ \vert \Lambda^{3} f\vert \ \vert S_{\alpha} f \vert^2  \left \vert \frac{\partial_x^3 \tau_{\alpha} f}{\alpha} \right \vert  \ d\gamma \ d\sigma \ d\alpha \ dx \\
 &+& \int \int\int_{0}^{\infty}\int_{0}^{\infty}   \ \gamma(1+\sigma) e^{-\gamma-\sigma} \ \vert \Lambda^{3} f\vert \ \vert S_{\alpha} f \vert^2  \left \vert \frac{\partial_x^2 \tau_{\alpha} f}{\alpha} \right \vert \ \vert f_{xx} \vert
   \ d\gamma \ d\sigma \ d\alpha \  dx.
\end{eqnarray*}
Therefore, by using  Sobolev interpolation  together with the embedding $\dot H^{3/2}\hookrightarrow \dot B^1_{\infty,2}$, we get
\begin{eqnarray*}
 T_3  &\lesssim& \Vert \partial^{3}_x f\Vert^2_{L^2} \ \int \frac{\Vert s_{\alpha} f \Vert^2_{L^{\infty}}}{\vert\alpha\vert^3} \ d\alpha + \Vert \Lambda^{3} f\Vert_{L^2} \Vert f_{xx}\Vert^2_{L^4}\ \int \frac{\Vert s_{\alpha} f \Vert_{L^\infty}^2}{\vert\alpha\vert^3} \ d\alpha \\
&\lesssim& (\Vert f\Vert^2_{\dot H^3}+\Vert f\Vert_{\dot H^3} \Vert f\Vert^2_{\dot H^{9/4}})\ \Vert f \Vert^2_{\dot B^1_{\infty,2}} \\
&\lesssim& (\Vert f\Vert^2_{\dot H^3}+\Vert f\Vert^2_{\dot H^3} \Vert f\Vert_{\dot H^{3/2}})\ \Vert f \Vert^2_{\dot H^{3/2}} \\
&\lesssim& \Vert f\Vert^2_{\dot H^3} \left(\Vert f \Vert^2_{\dot H^{3/2}} + \Vert f\Vert^3_{\dot H^{3/2}}\right) \\
\end{eqnarray*}

 \subsection{Estimate of $T_{4}$}

 This term is  much more delicate than $T_{3}$  in the sense that the term with the integral in $\eta$ does not help anymore. Indeed, one is forced to put an artificial weight which gives rise to a polynomial in $\alpha$ and this latter has a bad effect on the kernel which becomes none integrable at infinity. However, we can find very crucial cancellations of the two very problematic terms. These cancellations which are based on the structure of the equation and especially on the shift operator which appears in the Muskat problem with surface tension. Our symmetrization is crucial to see these hidden cancellations allowing us to get a new expression of $T_{4}$ which is more regular. In particular, this new expression will allow us to get nice {\it{a priori}} estimates in the critical Sobolev semi-norm $\dot H^{\frac{3}{2}}$.  \\

 Recall that,
 
  \begin{eqnarray*}
 T_{4} &=&-\frac{1}{4\pi} \int \int\int_{0}^{\infty}\int_{0}^{\infty}\int_{0}^{\alpha} \ e^{-\gamma-\sigma} \ \Lambda^{3} f \ \frac{1}{\alpha} s_{\eta} f_{x}  \left(\sin(\gamma\Delta_\alpha f)-\sin(\gamma\bar\Delta_\alpha f)\right) \\
&& \  \times\ \partial_x \left(\frac{\partial_x^2 \tau_{\alpha} f}{\alpha} \cos(\sigma \tau_{\alpha} f_x )\cos(\arctan(\tau_{\alpha} f_{x}))\right) \ d\eta  \ d\gamma \ d\sigma  \ d\alpha\ dx \\
  \end{eqnarray*}
  We are going to write $T_4$ in a more convenient way by searching some cancellations. More precisely, we shall prove the following new expression of $T_4$
  \begin{lemma} \label{t.4} $T_4$ may be written as
   \begin{eqnarray*}
  T_{4}&=&-\frac{1}{2\pi} \int \int\int_{0}^{\infty}\int_{0}^{\infty} \ e^{-\gamma-\sigma} \ \Lambda^{3} f \ \partial^2_{\alpha}\left[\frac{1}{\alpha}\sin(\frac{\gamma}{2}D_\alpha f)\cos(\frac{\gamma}{2}S_\alpha f)\frac{1}{\alpha} \int_{0}^{\alpha} s_{\eta} f_{x}    \ d\eta \right]\\
&& \  \times\ {  \delta_\alpha f_x } \cos(\sigma \tau_{\alpha} f_x )\cos(\arctan(\tau_{\alpha} f_{x})) \ d\eta  \ d\gamma \ d\sigma  \ d\alpha \ dx  \\
&-&\frac{1}{2\pi} \int \int\int_{0}^{\infty}\int_{0}^{\infty} \ \sigma e^{-\gamma-\sigma} \ \Lambda^{3} f \ \partial_{\alpha}\left[\frac{1}{\alpha}\sin(\frac{\gamma}{2}D_\alpha f)\cos(\frac{\gamma}{2}S_\alpha f)\frac{1}{\alpha} \int_{0}^{\alpha} s_{\eta} f_{x}    \ d\eta \right]\\
&& \  \times\ {  \delta_\alpha f_x } \left(\partial_{\alpha}\tau_{\alpha}f_x\right) \sin(\sigma \tau_{\alpha} f_x )\cos(\arctan(\tau_{\alpha} f_{x})) \ d\eta  \ d\gamma \ d\sigma  \ d\alpha\ dx \nonumber\\
&-&\frac{1}{2\pi} \int \int\int_{0}^{\infty}\int_{0}^{\infty} \ \sigma e^{-\gamma-\sigma} \ \Lambda^{3} f \ \partial_{\alpha}\left[\frac{1}{\alpha}\sin(\frac{\gamma}{2}D_\alpha f)\cos(\frac{\gamma}{2}S_\alpha f)\frac{1}{\alpha} \int_{0}^{\alpha} s_{\eta} f_{x}    \ d\eta \right]\\
&& \  \times\ {  \delta_\alpha f_x } \frac{\partial_{\alpha}\tau_{\alpha}f_x}{1+(\tau_{\alpha} f_{x})^2} \sin(\sigma \tau_{\alpha} f_x )\cos(\arctan(\tau_{\alpha} f_{x})) \ d\eta  \ d\gamma \ d\sigma  \ d\alpha \ dx \nonumber
   \end{eqnarray*}
  \end{lemma}
  \noindent {\bf{Proof of Lemma \ref{t.4}.}} We have
  \begin{eqnarray*}
 T_{4} &=&-\frac{1}{4\pi} \int \int\int_{0}^{\infty}\int_{0}^{\infty}\int_{0}^{\alpha} \ e^{-\gamma-\sigma} \ \Lambda^{3} f \ \frac{1}{\alpha} s_{\eta} f_{x}  \left(\sin(\gamma\Delta_\alpha f)-\sin(\gamma\bar\Delta_\alpha f)\right) \\
&& \  \times\ \frac{\partial_x^3 \tau_{\alpha} f}{\alpha} \cos(\sigma \tau_{\alpha} f_x )\cos(\arctan(\tau_{\alpha} f_{x})) \ d\eta  \ d\gamma \ d\sigma  \ d\alpha \ dx \nonumber \\
&+&\frac{1}{4\pi} \int \int\int_{0}^{\infty}\int_{0}^{\infty}\int_{0}^{\alpha} \ \sigma e^{-\gamma-\sigma} \ \Lambda^{3} f \ \frac{1}{\alpha} s_{\eta} f_{x}  \left(\sin(\gamma\Delta_\alpha f)-\sin(\gamma\bar\Delta_\alpha f)\right) \\
&& \  \times\ \frac{\partial_x^2 \tau_{\alpha} f}{\alpha}  \tau_{\alpha} f_{xx} \sin(\sigma\tau_{\alpha} f_x )\cos(\arctan(\tau_{\alpha} f_{x})) \ d\eta  \ d\gamma \ d\sigma  \ d\alpha\ dx \nonumber\\
&+&\frac{1}{4\pi} \int \int\int_{0}^{\infty}\int_{0}^{\infty}\int_{0}^{\alpha} \ \sigma e^{-\gamma-\sigma} \ \Lambda^{3} f \ \frac{1}{\alpha} s_{\eta} f_{x}  \left(\sin(\gamma\Delta_\alpha f)-\sin(\gamma\bar\Delta_\alpha f)\right) \\
&& \  \times\ \frac{\partial_x^2 \tau_{\alpha} f}{\alpha}  \frac{\tau_{\alpha} f_{xx}}{1+(\tau_{\alpha} f_{x})^2} \cos(\sigma\tau_{\alpha} f_x )\cos(\arctan(\tau_{\alpha} f_{x})) \ d\eta  \ d\gamma \ d\sigma  \ d\alpha \ dx \\
&:=&A+B+C. 
  \end{eqnarray*}
  On the other hand, we may write $\partial_x^3 \tau_\alpha f=-\partial_\alpha \partial_x^2\delta_\alpha f,$ then integrate by parts in $\alpha$, hence the most singular term $A$, that is the first one in the expression of $T_6$ above, may be written as follows
  \begin{eqnarray*}
  A&=&\frac{1}{4\pi} \int \int\int_{0}^{\infty}\int_{0}^{\infty}\int_{0}^{\alpha} \ e^{-\gamma-\sigma} \ \Lambda^{3} f \ \frac{1}{\alpha} s_{\eta} f_{x}  \left(\sin(\gamma\Delta_\alpha f)-\sin(\gamma\bar\Delta_\alpha f)\right) \\
&& \  \times\ \frac{\partial_\alpha \partial_x^2\delta_\alpha f }{\alpha} \cos(\sigma \tau_{\alpha} f_x )\cos(\arctan(\tau_{\alpha} f_{x}))  \ d\eta  \ d\gamma \ d\sigma \ dx \ d\alpha\\
&+&\frac{1}{4\pi} \int \int\int_{0}^{\infty}\int_{0}^{\infty}\int_{0}^{\alpha} \ \sigma e^{-\gamma-\sigma} \ \Lambda^{3} f \ \frac{1}{\alpha} s_{\eta} f_{x}  \left(\sin(\gamma\Delta_\alpha f)-\sin(\gamma\bar\Delta_\alpha f)\right) \\
&& \  \times\ \frac{\partial_x^2\delta_\alpha f }{\alpha} \partial_{\alpha}(\tau_{\alpha} f_x) \ \sin(\sigma \tau_{\alpha} f_x )\cos(\arctan(\tau_{\alpha} f_{x})) \ d\eta  \ d\gamma \ d\sigma \ dx \ d\alpha \\
  \end{eqnarray*}
  Since $\partial_{\alpha}(\tau_{\alpha} f_x)=-\partial_{x}(\tau_{\alpha} f_x)=-\tau_{\alpha} f_{xx}$, we find
  \begin{eqnarray*}
 A&=&\frac{1}{4\pi} \int \int\int_{0}^{\infty}\int_{0}^{\infty}\int_{0}^{\alpha} \ e^{-\gamma-\sigma} \ \Lambda^{3} f \ \frac{1}{\alpha} s_{\eta} f_{x}  \left(\sin(\gamma\Delta_\alpha f)-\sin(\gamma\bar\Delta_\alpha f)\right) \\
&& \  \times\ \frac{\partial_\alpha \partial_x^2\delta_\alpha f }{\alpha} \cos(\sigma \tau_{\alpha} f_x )\cos(\arctan(\tau_{\alpha} f_{x})) 
\ d\eta  \ d\gamma \ d\sigma \ dx \ d\alpha \\
&=&-\frac{1}{4\pi} \int \int\int_{0}^{\infty}\int_{0}^{\infty}\int_{0}^{\alpha} \ \sigma e^{-\gamma-\sigma} \ \Lambda^{3} f \ \frac{1}{\alpha} s_{\eta} f_{x}  \left(\sin(\gamma\Delta_\alpha f)-\sin(\gamma\bar\Delta_\alpha f)\right) \\
&& \  \times\ \frac{\partial_x^2\tau_\alpha f }{\alpha} \tau_{\alpha} f_{xx} \ \sin(\sigma \tau_{\alpha} f_x )\cos(\arctan(\tau_{\alpha} f_{x})) 
\  d\eta  \ d\gamma \ d\sigma \ dx \ d\alpha \\
&-&\frac{1}{4\pi} \int \int\int_{0}^{\infty}\int_{0}^{\infty}\int_{0}^{\alpha} \ \sigma e^{-\gamma-\sigma} \ \Lambda^{3} f \ \frac{1}{\alpha} s_{\eta} f_{x}  \left(\sin(\gamma\Delta_\alpha f)-\sin(\gamma\bar\Delta_\alpha f)\right) \\
&& \  \times\ \frac{\partial_x^2 \tau_{\alpha} f}{\alpha}  \frac{\tau_{\alpha} f_{xx}}{1+(\tau_{\alpha} f_{x})^2} \cos(\sigma\tau_{\alpha} f_x )\cos(\arctan(\tau_{\alpha} f_{x})) \ d\eta  \ d\gamma \ d\sigma \ dx \ d\alpha \\
&-&\frac{1}{4\pi} \int \int\int_{0}^{\infty}\int_{0}^{\infty}\int_{0}^{\alpha} s_{\eta} f_{x} \ e^{-\gamma-\sigma} \ \Lambda^{3} f \ \partial_{\alpha}\left[\frac{1}{\alpha^2}   \left(\sin(\gamma\Delta_\alpha f)-\sin(\gamma\bar\Delta_\alpha f)\right) \right]\\
&& \  \times\ { \partial_x^2 \tau_\alpha f } \cos(\sigma \tau_{\alpha} f_x )\cos(\arctan(\tau_{\alpha} f_{x})) \ d\eta  \ d\gamma \ d\sigma \ dx \ d\alpha \\
  \end{eqnarray*}
  Therefore, replacing this expression of $A$ in the expression of  $T_{4}$ above leads to {{two crucial cancellations of the bad terms}}, namely $B$ and $C$, we get  the following new expression of $T_{4}$
  \begin{eqnarray*}
  T_{4}&=&-\frac{1}{4\pi} \int \int\int_{0}^{\infty}\int_{0}^{\infty} \ e^{-\gamma-\sigma} \ \Lambda^{3} f \ \partial_{\alpha}\left[\int_{0}^{\alpha} s_{\eta} f_{x}\frac{1}{\alpha^2}  \ d\eta  \left(\sin(\gamma\Delta_\alpha f)-\sin(\gamma\bar\Delta_\alpha f)\right)  \right]\\
&& \  \times\ { \partial_x^2 \tau_\alpha f } \cos(\sigma \tau_{\alpha} f_x )\cos(\arctan(\tau_{\alpha} f_{x})) \ d\eta  \ d\gamma \ d\sigma \ dx \ d\alpha \\
   \end{eqnarray*}
     We can further write this term as follows
   \begin{eqnarray*}
  T_{4}&=&-\frac{1}{2\pi} \int \int\int_{0}^{\infty}\int_{0}^{\infty} \ e^{-\gamma-\sigma} \ \Lambda^{3} f \ \partial_{\alpha}\left[\frac{1}{\alpha}\sin(\frac{\gamma}{2}D_\alpha f)\cos(\frac{\gamma}{2}S_\alpha f)\frac{1}{\alpha} \int_{0}^{\alpha} s_{\eta} f_{x}    \ d\eta \right]\\
&& \  \times\ { \partial_x^2 \tau_\alpha f } \cos(\sigma \tau_{\alpha} f_x )\cos(\arctan(\tau_{\alpha} f_{x})) \ d\eta  \ d\gamma \ d\sigma \ dx \ d\alpha \\
   \end{eqnarray*}
   Using $\partial_x^2 \tau_\alpha f=-\partial_\alpha \delta_\alpha f_x$ and then integrating by parts in $\alpha$ give
   \begin{eqnarray*}
  T_{4}&=&-\frac{1}{2\pi} \int \int\int_{0}^{\infty}\int_{0}^{\infty} \ e^{-\gamma-\sigma} \ \Lambda^{3} f \ \partial^2_{\alpha}\left[\frac{1}{\alpha}\sin(\frac{\gamma}{2}D_\alpha f)\cos(\frac{\gamma}{2}S_\alpha f)\frac{1}{\alpha} \int_{0}^{\alpha} s_{\eta} f_{x}    \ d\eta \right]\\
&& \  \times\ {  \delta_\alpha f_x } \cos(\sigma \tau_{\alpha} f_x )\cos(\arctan(\tau_{\alpha} f_{x})) \ d\eta \ d\gamma \ d\sigma \ dx \ d\alpha \\
&+&\frac{1}{2\pi} \int \int\int_{0}^{\infty}\int_{0}^{\infty} \ \sigma e^{-\gamma-\sigma} \ \Lambda^{3} f \ \partial_{\alpha}\left[\frac{1}{\alpha}\sin(\frac{\gamma}{2}D_\alpha f)\cos(\frac{\gamma}{2}S_\alpha f)\frac{1}{\alpha} \int_{0}^{\alpha} s_{\eta} f_{x}    \ d\eta \right]\\
&& \  \times\ {  \delta_\alpha f_x } \left(\partial_{\alpha}\tau_{\alpha}f_x\right) \sin(\sigma \tau_{\alpha} f_x )\cos(\arctan(\tau_{\alpha} f_{x})) \ d\eta  \ d\gamma \ d\sigma \ dx\ d\alpha \\
&+&\frac{1}{2\pi} \int \int\int_{0}^{\infty}\int_{0}^{\infty} \ \sigma e^{-\gamma-\sigma} \ \Lambda^{3} f \ \partial_{\alpha}\left[\frac{1}{\alpha}\sin(\frac{\gamma}{2}D_\alpha f)\cos(\frac{\gamma}{2}S_\alpha f)\frac{1}{\alpha} \int_{0}^{\alpha} s_{\eta} f_{x}    \ d\eta \right]\\
&& \  \times\ {  \delta_\alpha f_x } \frac{\partial_{\alpha}\tau_{\alpha}f_x}{1+(\tau_{\alpha} f_{x})^2} \sin(\sigma \tau_{\alpha} f_x )\cos(\arctan(\tau_{\alpha} f_{x})) \ d\eta  \ d\gamma \ d\sigma \ dx \ d\alpha \\
   \end{eqnarray*}
   Since $\partial_{\alpha}\tau_{\alpha}f_x=-\tau_{\alpha}f_{xx}$, we find
    \begin{eqnarray*}
  T_{4}&=&-\frac{1}{2\pi} \int \int\int_{0}^{\infty}\int_{0}^{\infty} \ e^{-\gamma-\sigma} \ \Lambda^{3} f \ \partial^2_{\alpha}\left[\frac{1}{\alpha}\sin(\frac{\gamma}{2}D_\alpha f)\cos(\frac{\gamma}{2}S_\alpha f)\frac{1}{\alpha} \int_{0}^{\alpha} s_{\eta} f_{x}    \ d\eta \right]\\
&& \  \times\ {  \delta_\alpha f_x } \cos(\sigma \tau_{\alpha} f_x )\cos(\arctan(\tau_{\alpha} f_{x})) \ d\eta  \ d\gamma \ d\sigma \ dx \ d\alpha \\
&-&\frac{1}{2\pi} \int \int\int_{0}^{\infty}\int_{0}^{\infty} \ \sigma e^{-\gamma-\sigma} \ \Lambda^{3} f \ \partial_{\alpha}\left[\frac{1}{\alpha}\sin(\frac{\gamma}{2}D_\alpha f)\cos(\frac{\gamma}{2}S_\alpha f)\frac{1}{\alpha} \int_{0}^{\alpha} s_{\eta} f_{x}    \ d\eta \right]\\
&& \  \times\ {  \delta_\alpha f_x } \left(\partial_{\alpha}\tau_{\alpha}f_x\right) \sin(\sigma \tau_{\alpha} f_x )\cos(\arctan(\tau_{\alpha} f_{x})) \ d\eta  \ d\gamma \ d\sigma \ dx \ d\alpha\\
&-&\frac{1}{2\pi} \int \int\int_{0}^{\infty}\int_{0}^{\infty} \ \sigma e^{-\gamma-\sigma} \ \Lambda^{3} f \ \partial_{\alpha}\left[\frac{1}{\alpha}\sin(\frac{\gamma}{2}D_\alpha f)\cos(\frac{\gamma}{2}S_\alpha f)\frac{1}{\alpha} \int_{0}^{\alpha} s_{\eta} f_{x}    \ d\eta \right]\\
&& \  \times\ {  \delta_\alpha f_x } \frac{\partial_{\alpha}\tau_{\alpha}f_x}{1+(\tau_{\alpha} f_{x})^2} \sin(\sigma \tau_{\alpha} f_x )\cos(\arctan(\tau_{\alpha} f_{x})) \ d\eta  \ d\gamma \ d\sigma \ dx \ d\alpha. \\
   \end{eqnarray*}

    This end the proof of Lemma \ref{t.4}. 
    
    \qed
    
    We are now going to estimate $T_4$ by controlling each of the 3 above terms that we shall denote $T_4:=T_{4,1}+T_{4,2}+T_{4,3}$. More precisely, we are going to prove the following control that we state as a Lemma. 
    \begin{lemma} \label{t41c}
    The following control holds 
    \begin{eqnarray}
T_{4,1}&\lesssim& \Vert f \Vert^2_{\dot H^3}\left(\Vert f \Vert_{\dot H^{\frac{3}{2}}}+\Vert f \Vert^2_{\dot H^{\frac{3}{2}}}+\Vert f \Vert^3_{\dot H^{\frac{3}{2}}}\right)
\end{eqnarray}
    \end{lemma}
    \noindent {\bf{Proof of Lemma \ref{t41c}.}} We shall estimate each of term in the new expression of $T_{4,1}$ obtained in Lemma \ref{t.4}. \\
   
   \noindent $\bullet$ {Estimate of $T_{4,1}$}
    
    Using the formula $\displaystyle D_\alpha f-2f_x=\frac{1}{\alpha} \int_{0}^{\alpha} s_{\eta} f_{x} \ d\eta,$ we infer that for any $k\geq0$, $$\partial^k_\alpha  D_\alpha f=\partial^k_\alpha\left(\frac{1}{\alpha}\int_{0}^{\alpha} s_{\eta} f_{x}
     \ d\eta\right).$$ Hence, since
    \begin{eqnarray*}
   \partial^2_\alpha\left(\frac{1}{\alpha}\int_{0}^{\alpha} s_{\eta} f_{x}  \ d\eta\right)&=&\partial^2_\alpha  D_\alpha f=\frac{d_\alpha f_{xx}}{\alpha}+2 \frac{s_\alpha f_x}{\alpha^2} + \frac{2}{\alpha^3}\int_0^\alpha s_{\eta}f_x \ d\eta
   \end{eqnarray*}
   We find that

   \begin{eqnarray*}
  T_{4,1}&=&-\frac{1}{2\pi} \int \int\int_{0}^{\infty}\int_{0}^{\infty} \ e^{-\gamma-\sigma} \ \Lambda^{3} f \ \partial^2_{\alpha}\left[\frac{1}{\alpha}\sin(\frac{\gamma}{2}D_\alpha f)\cos(\frac{\gamma}{2}S_\alpha f)\frac{1}{\alpha} \int_{0}^{\alpha} s_{\eta} f_{x}    \ d\eta \right]\\
&& \  \times\ {  \delta_\alpha f_x } \cos(\sigma \tau_{\alpha} f_x )\cos(\arctan(\tau_{\alpha} f_{x}))   \ d\gamma \ d\sigma \ dx \ d\alpha \\
\end{eqnarray*}
Therefore, we get
\begin{eqnarray*}
T_{4,1}&=&-\frac{1}{\pi} \int \int\int_{0}^{\infty}\int_{0}^{\infty} \ e^{-\gamma-\sigma} \ \Lambda^{3} f \ \frac{1}{\alpha^3}\sin(\frac{\gamma}{2}D_\alpha f)\cos(\frac{\gamma}{2}S_\alpha f)\frac{1}{\alpha} \int_{0}^{\alpha} s_{\eta} f_{x}    \ d\eta \\
&& \  \times\ {  \delta_\alpha f_x } \cos(\sigma \tau_{\alpha} f_x )\cos(\arctan(\tau_{\alpha} f_{x}))   \ d\gamma \ d\sigma  \ d\alpha\ dx \\
&-&\frac{1}{4\pi} \int \int\int_{0}^{\infty}\int_{0}^{\infty} \ \gamma e^{-\gamma-\sigma} \ \Lambda^{3} f \ \partial^2_\alpha D_\alpha f \cos(\frac{\gamma}{2}D_\alpha f)\cos(\frac{\gamma}{2}S_\alpha f)\frac{1}{\alpha^2} \int_{0}^{\alpha} s_{\eta} f_{x}    \ d\eta \\
&& \  \times\ {  \delta_\alpha f_x } \cos(\sigma \tau_{\alpha} f_x )\cos(\arctan(\tau_{\alpha} f_{x}))   \ d\gamma \ d\sigma  \ d\alpha \ dx\\
&+&\frac{1}{8\pi} \int \int\int_{0}^{\infty}\int_{0}^{\infty} \ \gamma^2 e^{-\gamma-\sigma} \ \Lambda^{3} f \ (\partial_\alpha D_\alpha f)^2 \sin(\frac{\gamma}{2}D_\alpha f)\cos(\frac{\gamma}{2}S_\alpha f)\frac{1}{\alpha^2} \int_{0}^{\alpha} s_{\eta} f_{x}    \ d\eta \\
&& \  \times\ {  \delta_\alpha f_x } \cos(\sigma \tau_{\alpha} f_x )\cos(\arctan(\tau_{\alpha} f_{x}))   \ d\gamma \ d\sigma  \ d\alpha \ dx \\
&+&\frac{1}{4\pi} \int \int\int_{0}^{\infty}\int_{0}^{\infty} \ \gamma e^{-\gamma-\sigma} \ \Lambda^{3} f \ \partial^2_\alpha S_\alpha f \sin(\frac{\gamma}{2}D_\alpha f)\sin(\frac{\gamma}{2}S_\alpha f)\frac{1}{\alpha^2} \int_{0}^{\alpha} s_{\eta} f_{x}    \ d\eta \\
&& \  \times\ {  \delta_\alpha f_x } \cos(\sigma \tau_{\alpha} f_x )\cos(\arctan(\tau_{\alpha} f_{x})) \ d\eta \ d\gamma \ d\sigma  \ d\alpha  \ dx\\
&+&\frac{1}{8\pi} \int \int\int_{0}^{\infty}\int_{0}^{\infty} \ \gamma^2 e^{-\gamma-\sigma} \ \Lambda^{3} f \ (\partial_\alpha S_\alpha f)^2 \sin(\frac{\gamma}{2}D_\alpha f)\cos(\frac{\gamma}{2}S_\alpha f)\frac{1}{\alpha^2} \int_{0}^{\alpha} s_{\eta} f_{x}    \ d\eta \\
&& \  \times\ {  \delta_\alpha f_x } \cos(\sigma \tau_{\alpha} f_x )\cos(\arctan(\tau_{\alpha} f_{x}))   \ d\gamma \ d\sigma \ d\alpha \ dx \\
&-&\frac{1}{2\pi} \int \int\int_{0}^{\infty}\int_{0}^{\infty} \ e^{-\gamma-\sigma} \ \Lambda^{3} f \ \partial^2_{\alpha}D_\alpha f \ \frac{1}{\alpha}\sin(\frac{\gamma}{2}D_\alpha f)\cos(\frac{\gamma}{2}S_\alpha f)\\
&& \  \times\ {  \delta_\alpha f_x } \cos(\sigma \tau_{\alpha} f_x )\cos(\arctan(\tau_{\alpha} f_{x}))   \ d\gamma \ d\sigma  \ d\alpha \ dx\\
&+&\frac{1}{2\pi} \int \int\int_{0}^{\infty}\int_{0}^{\infty} \ \gamma e^{-\gamma-\sigma} \ \Lambda^{3} f \ \partial_{\alpha}D_\alpha f \cos(\frac{\gamma}{2}D_\alpha f)\cos(\frac{\gamma}{2}S_\alpha f)\frac{1}{\alpha^3} \int_{0}^{\alpha} s_{\eta} f_{x}    \ d\eta \\
&& \  \times\ {  \delta_\alpha f_x } \cos(\sigma \tau_{\alpha} f_x )\cos(\arctan(\tau_{\alpha} f_{x}))  \ d\gamma \ d\sigma  \ d\alpha\ dx\\
&-&\frac{1}{2\pi} \int \int\int_{0}^{\infty}\int_{0}^{\infty} \ \gamma e^{-\gamma-\sigma} \ \Lambda^{3} f \ \partial_{\alpha}S_\alpha f \sin(\frac{\gamma}{2}D_\alpha f)\sin(\frac{\gamma}{2}S_\alpha f)\frac{1}{\alpha^3} \int_{0}^{\alpha} s_{\eta} f_{x}    \ d\eta \\
&& \  \times\ {  \delta_\alpha f_x } \cos(\sigma \tau_{\alpha} f_x )\cos(\arctan(\tau_{\alpha} f_{x}))  \ d\gamma \ d\sigma   \ d\alpha \ dx\\
&+&\frac{1}{\pi} \int \int\int_{0}^{\infty}\int_{0}^{\infty} \ e^{-\gamma-\sigma} \ \Lambda^{3} f \ \partial_\alpha D_\alpha f \sin(\frac{\gamma}{2}D_\alpha f)\cos(\frac{\gamma}{2}S_\alpha f)\frac{1}{\alpha^2}  \\
&& \  \times\ {  \delta_\alpha f_x } \cos(\sigma \tau_{\alpha} f_x )\cos(\arctan(\tau_{\alpha} f_{x}))   \ d\gamma \ d\sigma \ dx  \ d\alpha \\
&+&\frac{1}{4\pi} \int \int\int_{0}^{\infty}\int_{0}^{\infty} \ \gamma^2 e^{-\gamma-\sigma} \ \Lambda^{3} f \ \partial_{\alpha} D_\alpha f \ \partial_{\alpha} S_\alpha f \cos(\frac{\gamma}{2}D_\alpha f)\sin(\frac{\gamma}{2}S_\alpha f)\frac{1}{\alpha^2} \int_{0}^{\alpha} s_{\eta} f_{x}    \ d\eta \\
&& \  \times\ {  \delta_\alpha f_x } \cos(\sigma \tau_{\alpha} f_x )\cos(\arctan(\tau_{\alpha} f_{x}))   \ d\gamma \ d\sigma   \ d\alpha \ dx\\
&-&\frac{1}{2\pi} \int \int\int_{0}^{\infty}\int_{0}^{\infty} \ \gamma e^{-\gamma-\sigma} \ \Lambda^{3} f \ \frac{1}{\alpha} \ (\partial_{\alpha} D_\alpha f)^2 \  \cos(\frac{\gamma}{2}D_\alpha f)\cos(\frac{\gamma}{2}S_\alpha f)\\
&& \  \times\ {  \delta_\alpha f_x } \cos(\sigma \tau_{\alpha} f_x )\cos(\arctan(\tau_{\alpha} f_{x}))  \ d\gamma \ d\sigma   \ d\alpha \ dx\\
&+&\frac{1}{2\pi} \int \int\int_{0}^{\infty}\int_{0}^{\infty} \ \gamma e^{-\gamma-\sigma} \ \Lambda^{3} f \ \frac{1}{\alpha}\ \partial_{\alpha} D_\alpha f \ \partial_{\alpha} S_\alpha f \sin(\frac{\gamma}{2}D_\alpha f)\sin(\frac{\gamma}{2}S_\alpha f) \\
&& \  \times\ {  \delta_\alpha f_x } \cos(\sigma \tau_{\alpha} f_x )\cos(\arctan(\tau_{\alpha} f_{x}))   \ d\gamma \ d\sigma  \ d\alpha\ dx \\
&:=& \sum_{i=1}^{12}  T_{4,1,i}
\end{eqnarray*}

$\bullet$ \noindent{{Estimate of the term $T_{4,1,1}$}} \\

This term corresponds to the one where all the derivatives in $\alpha$ falls on the kernel, hence, one needs to use Besov regularity with indices bigger than 1 in order to milder the decay in $\alpha$ and close the estimate in the critical spaces.
\begin{eqnarray*}
T_{4,1,1}&=&-\frac{1}{\pi} \int \int\int_{0}^{\infty}\int_{0}^{\infty} \ e^{-\gamma-\sigma} \ \Lambda^{3} f \ \frac{1}{\alpha^3}\sin(\frac{\gamma}{2}D_\alpha f)\cos(\frac{\gamma}{2}S_\alpha f)\frac{1}{\alpha} \int_{0}^{\alpha} s_{\eta} f_{x}    \ d\eta \\
&& \  \times\ {  \delta_\alpha f_x } \cos(\sigma \tau_{\alpha} f_x )\cos(\arctan(\tau_{\alpha} f_{x}))   \ d\gamma \ d\sigma \ d\alpha \ dx \\
&\lesssim& \Vert f \Vert_{\dot H^3} \int \frac{1}{\alpha^4} \Vert \delta_\alpha f_x \Vert_{L^\infty} \int_{0}^{\alpha} \frac{\Vert s_{\eta} f_{x} \Vert_{L^2}}{\eta^2} \eta^2    \ d\eta \ d\alpha \\
&\lesssim& \Vert f \Vert_{\dot H^3} \int \frac{1}{\alpha^4} \Vert \delta_\alpha f_x \Vert_{L^\infty}  \left(\int_{0}^{\alpha} \frac{\Vert s_{\eta} f_{x} \Vert^2_{L^2}}{\eta^4}    \ d\eta \right)^{\frac{1}{2}} \left(\int_{0}^{\alpha} \eta^4 \ d\eta \right)^{\frac{1}{2}} \ d\alpha \\
&\lesssim& \Vert f \Vert_{\dot H^3} \int \frac{1}{\alpha^4} \Vert \delta_\alpha f_x \Vert_{L^\infty}  \left(\int \frac{\Vert s_{\eta} f_{x} \Vert^2_{L^2}}{\eta^4}    \ d\eta \right)^{\frac{1}{2}} \left(\int_{0}^{\alpha} \eta^4 \ d\eta \right)^{\frac{1}{2}} \ d\alpha \\
&\lesssim& \Vert f \Vert_{\dot H^3} \Vert f \Vert_{\dot H^\frac{5}{2}} \int \frac{1}{\vert\alpha\vert^{\frac{3}{2}}} \Vert \delta_\alpha f_x \Vert_{L^\infty}  \ d\alpha \\
&\lesssim& \Vert f \Vert_{\dot H^3} \Vert f \Vert_{\dot H^\frac{5}{2}} \Vert f \Vert_{\dot B^\frac{3}{2}_{\infty,1}} \\
&\lesssim& \Vert f \Vert_{\dot H^3} \Vert f \Vert_{\dot H^\frac{5}{2}} 
\Vert f \Vert^{\frac{1}{2}}_{\dot B^1_{\infty,\infty}} \Vert f \Vert^{\frac{1}{2}}_{\dot B^2_{\infty,\infty}} \\
&\lesssim& \Vert f \Vert_{\dot H^3} \Vert f \Vert_{\dot H^\frac{5}{2}} 
\Vert f \Vert^{\frac{1}{2}}_{\dot H^\frac{3}{2}} \Vert f \Vert^{\frac{1}{2}}_{\dot H^\frac{5}{2}} \\
&\lesssim& \Vert f \Vert_{\dot H^3} \Vert f \Vert^{\frac{3}{2}}_{\dot H^\frac{5}{2}} 
\Vert f \Vert^{\frac{1}{2}}_{\dot H^\frac{3}{2}}\\
\end{eqnarray*}

where we used  $\dot  B^\frac{3}{2}_{\infty,1}=\left[\dot B^1_{\infty, \infty}, \dot B^2_{\infty, \infty}\right]_{\frac{1}{2},\frac{1}{2}}$ and that $\left[\dot H^3, \dot H^\frac{3}{2} \right]_{\frac{2}{3},\frac{1}{3}}=\dot H^\frac{5}{2} \hookrightarrow \dot B^2_{\infty, \infty}.$
Then, since
$$
\Vert f \Vert^{{\frac{3}{2}}}_{\dot H^\frac{5}{2}} \leq \Vert f \Vert_{\dot H^3} \Vert f \Vert^{\frac{1}{2}}_{\dot H^\frac{3}{2} } 
$$
we conclude that 

\begin{eqnarray*}
T_{4,1,1}&\lesssim&\Vert f \Vert^2_{\dot H^3} \Vert f \Vert_{\dot H^\frac{3}{2} }
\end{eqnarray*}

$\bullet$ {Estimate of the $T_{4,1,2}$} \\

We need to develop this term as follows, using the formula for the second derivative  of $D_\alpha$ with respect to $\alpha$, we may write

\begin{eqnarray*}
T_{4,1,2}&=&-\frac{1}{4\pi} \int \int\int_{0}^{\infty}\int_{0}^{\infty} \ \gamma e^{-\gamma-\sigma} \ \Lambda^{3} f \ \partial^2_\alpha D_\alpha f \cos(\frac{\gamma}{2}D_\alpha f)\cos(\frac{\gamma}{2}S_\alpha f)\frac{1}{\alpha^2} \int_{0}^{\alpha} s_{\eta} f_{x}    \ d\eta \\
&& \  \times\ {  \delta_\alpha f_x } \cos(\sigma \tau_{\alpha} f_x )\cos(\arctan(\tau_{\alpha} f_{x}))   \ d\gamma \ d\sigma  \ d\alpha \ dx\\
&=&-\frac{1}{4\pi} \int \int\int_{0}^{\infty}\int_{0}^{\infty} \ \gamma e^{-\gamma-\sigma} \ \Lambda^{3} f \ \frac{d_\alpha f_{xx}}{\alpha} \cos(\frac{\gamma}{2}D_\alpha f)\cos(\frac{\gamma}{2}S_\alpha f)\frac{1}{\alpha^2} \int_{0}^{\alpha} s_{\eta} f_{x}    \ d\eta \\
&& \  \times\ {  \delta_\alpha f_x } \cos(\sigma \tau_{\alpha} f_x )\cos(\arctan(\tau_{\alpha} f_{x}))   \ d\gamma \ d\sigma  \ d\alpha \ dx\\
&-&\frac{1}{2\pi} \int \int\int_{0}^{\infty}\int_{0}^{\infty} \ \gamma e^{-\gamma-\sigma} \ \Lambda^{3} f \ \frac{s_\alpha f_x}{\alpha^2} \cos(\frac{\gamma}{2}D_\alpha f)\cos(\frac{\gamma}{2}S_\alpha f)\frac{1}{\alpha^2} \int_{0}^{\alpha} s_{\eta} f_{x}    \ d\eta \\
&& \  \times\ {  \delta_\alpha f_x } \cos(\sigma \tau_{\alpha} f_x )\cos(\arctan(\tau_{\alpha} f_{x}))   \ d\gamma \ d\sigma  \ d\alpha \ dx\\
&-&\frac{1}{2\pi} \int \int\int_{0}^{\infty}\int_{0}^{\infty} \ \gamma e^{-\gamma-\sigma} \ \Lambda^{3} f \ \frac{1}{\alpha^3} \int_0^\alpha s_{\kappa}f_x \ d\kappa \cos(\frac{\gamma}{2}D_\alpha f)\cos(\frac{\gamma}{2}S_\alpha f)\frac{1}{\alpha^2} \int_{0}^{\alpha} s_{\eta} f_{x}    \ d\eta \\
&& \  \times\ {  \delta_\alpha f_x } \cos(\sigma \tau_{\alpha} f_x )\cos(\arctan(\tau_{\alpha} f_{x}))   \ d\gamma \ d\sigma  \ d\alpha\ dx \\
&:=& \sum_{i=1}^3 T_{4,1,2,i}
\end{eqnarray*}

$\bullet$ {{Estimate of $T_{4,1,2,1}$}} \\

We want to estimate 

\begin{eqnarray*}
T_{4,1,2,1}&=&-\frac{1}{4\pi} \int \int\int_{0}^{\infty}\int_{0}^{\infty} \ \gamma e^{-\gamma-\sigma} \ \Lambda^{3} f \ \frac{d_\alpha f_{xx}}{\alpha} \cos(\frac{\gamma}{2}D_\alpha f)\cos(\frac{\gamma}{2}S_\alpha f)\frac{1}{\alpha^2} \int_{0}^{\alpha} s_{\eta} f_{x}    \ d\eta \\
&& \  \times\ {  \delta_\alpha f_x } \cos(\sigma \tau_{\alpha} f_x )\cos(\arctan(\tau_{\alpha} f_{x}))  \ d\gamma \ d\sigma  \ d\alpha \ dx \\
\end{eqnarray*}
We have that
\begin{eqnarray*}
T_{4,1,2,1}&\lesssim& \Vert f \Vert_{\dot H^3} \Vert f \Vert_{\dot H^2} \int \frac{\Vert \delta_\alpha f_x \Vert_{L^\infty}}{\vert \alpha \vert^3} \left(\int_{0}^{\alpha} \frac{\Vert s_{\eta} f_{x} \Vert^2_{L^\infty}}{\eta^2} \ d\eta\right)^{\frac{1}{2}} \left(\int_{0}^{\alpha} \eta^2  \ d\eta\right)^{\frac{1}{2}} \ d\alpha \\
&\lesssim& \Vert f \Vert_{\dot H^3} \Vert f \Vert_{\dot H^2} \int \frac{\Vert \delta_\alpha f_x \Vert_{L^\infty}}{\vert \alpha \vert^{\frac{3}{2}}} \left(\int \frac{\Vert s_{\eta} f_{x} \Vert^2_{L^\infty}}{\eta^2} \ d\eta\right)^{\frac{1}{2}}  \ d\alpha \\
&\lesssim& \Vert f \Vert_{\dot H^3} \Vert f \Vert^2_{\dot H^2} \Vert f \Vert_{\dot B^{\frac{3}{2}}_{\infty,1}}
 \end{eqnarray*}   
 
 Using again the fact that $\dot  B^\frac{3}{2}_{\infty,1}=\left[\dot B^1_{\infty, \infty}, \dot B^2_{\infty, \infty}\right]_{\frac{1}{2},\frac{1}{2}}$
 One finds that,
 \begin{eqnarray*}
T_{4,1,2,1} &\lesssim& \Vert f \Vert_{\dot H^3} \Vert f \Vert^2_{\dot H^2}
 \Vert f \Vert^{\frac{1}{2}}_{\dot B^{1}_{\infty,\infty}}   \Vert f \Vert^{\frac{1}{2}}_{\dot B^{2}_{\infty,\infty}} \\
 &\lesssim&  \Vert f \Vert_{\dot H^3} \Vert f \Vert^2_{\dot H^2} \Vert f \Vert^{\frac{1}{2}}_{\dot H^{\frac{3}{2}}}\Vert f \Vert^{\frac{1}{2}}_{\dot H^{\frac{5}{2}}}
\end{eqnarray*}     
Since,
$$
\Vert f \Vert^{\frac{1}{2}}_{\dot H^\frac{5}{2}} \leq \Vert f \Vert^{\frac{1}{3}}_{\dot H^3} \Vert f \Vert^{\frac{1}{6}}_{\dot H^\frac{3}{2} } 
$$
and that $H^2=\left[\dot H^{\frac{3}{2}}, \dot H^{3} \right]_{\frac{2}{3},\frac{1}{3}}$ in particular
$$
\Vert f \Vert^2_{\dot H^2} \leq \Vert f \Vert^{\frac{2}{3}}_{\dot H^3} \Vert f \Vert^{\frac{4}{3}}_{\dot H^\frac{3}{2} } 
$$
Finally, 
\begin{eqnarray*}
T_{4,1,2,1} &\lesssim& \Vert f \Vert^2_{\dot H^3}  \Vert f \Vert^2_{\dot H^\frac{3}{2} }
\end{eqnarray*}
$\bullet$ {{Estimate of $T_{4,1,2,2}$}} \\

Recall that

\begin{eqnarray*}
T_{4,1,2,2}&=&-\frac{1}{2\pi} \int \int\int_{0}^{\infty}\int_{0}^{\infty} \ \gamma e^{-\gamma-\sigma} \ \Lambda^{3} f \ \frac{s_\alpha f_x}{\alpha^2} \cos(\frac{\gamma}{2}D_\alpha f)\cos(\frac{\gamma}{2}S_\alpha f)\frac{1}{\alpha^2} \int_{0}^{\alpha} s_{\eta} f_{x}    \ d\eta \\
&& \  \times\ {  \delta_\alpha f_x } \cos(\sigma \tau_{\alpha} f_x )\cos(\arctan(\tau_{\alpha} f_{x}))   \ d\gamma \ d\sigma\ d\alpha \ dx \\
&\lesssim& \Vert f \Vert_{\dot H^2} \int \frac{\Vert s_\alpha f_x \Vert_{L^{\infty}}}{\alpha^4} \Vert \delta_\alpha f_x \Vert_{L^{\infty}} \left( \int_{0}^{\alpha} \frac{\Vert s_{\eta} f_{x}\Vert^2_{L^2}}{\eta^2}    \ d\eta\right)^{\frac{1}{2}} \left( \int_{0}^{\alpha}\eta^2 \ d\eta \right)^{\frac{1}{2}} \ d\alpha\\
&\lesssim& \Vert f \Vert_{\dot H^2} \int \frac{\Vert s_\alpha f_x \Vert_{L^{\infty}}}{\vert\alpha\vert^{\frac{5}{2}}} \Vert \delta_\alpha f_x \Vert_{L^{\infty}} \left( \int \frac{\Vert s_{\eta} f_{x}\Vert^2_{L^2}}{\eta^2}    \ d\eta\right)^{\frac{1}{2}} \ d\alpha \\
&\lesssim& \Vert f \Vert_{\dot H^2} \Vert f \Vert_{\dot H^{\frac{3}{2}}} \left( \int \frac{\Vert s_\alpha f_x \Vert^2_{L^{\infty}}}{\vert\alpha \vert^{\frac{5}{2}}} \ d\alpha \right)^{\frac{1}{2}} \left( \int \frac{\Vert \delta_\alpha f_x \Vert^2_{L^{\infty}}}{\vert\alpha\vert^{\frac{5}{2}}} \ d\alpha \right)^{\frac{1}{2}}\\
&\lesssim& \Vert f \Vert_{\dot H^2} \Vert f \Vert_{\dot H^{\frac{3}{2}}} \Vert f \Vert^2_{\dot B^{\frac{7}{4}}_{\infty,2}}
\end{eqnarray*}
Since $\left[\dot H^{\frac{3}{2}}, \dot H^3 \right]_{\frac{1}{2},\frac{1}{2}}=\dot H^{\frac{9}{4}} \hookrightarrow \dot B^{\frac{7}{4}}_{\infty,2}$ we find that
\begin{eqnarray*}
T_{4,1,2,2}&\lesssim& \Vert f \Vert^2_{\dot H^3}  \Vert f \Vert^2_{\dot H^\frac{3}{2} }
\end{eqnarray*}

$\bullet$ {{Estimate of $T_{4,1,2,3}$}} \\

We have

\begin{eqnarray*}
T_{4,1,2,3}&=&-\frac{1}{2\pi} \int \int\int_{0}^{\infty}\int_{0}^{\infty} \int_0^\alpha  \int_0^\alpha  \ \gamma e^{-\gamma-\sigma} \ \Lambda^{3} f \ \frac{1}{\alpha^5} \ s_{\kappa}f_x \ s_{\eta} f_{x}  \  \cos(\frac{\gamma}{2}D_\alpha f)\cos(\frac{\gamma}{2}S_\alpha f)     \\
&& \  \times\ {  \delta_\alpha f_x } \cos(\sigma \tau_{\alpha} f_x )\cos(\arctan(\tau_{\alpha} f_{x}))\ d\kappa \ d\eta  \ d\gamma \ d\sigma  \ d\alpha \ dx \\
&\lesssim& \Vert f \Vert_{\dot H^3} \int \frac{\Vert \delta_\alpha f_x \Vert_{L^{\infty}}}{\vert\alpha\vert^5} \left(\int_0^\alpha \frac{ \Vert s_{\kappa}f_x \Vert^2_{L^{4}} }{\vert\kappa\vert^{\frac{5}{2}}} d\kappa \right)^{\frac{1}{2}} \left(\int_0^\alpha \frac{\Vert s_{\eta}f_x \Vert^2_{L^{4}}}{\vert \eta \vert^{\frac{5}{2}}}  d\eta \right)^{\frac{1}{2}} \\
&& \times \left( \int_{0}^{\alpha}\vert\eta\vert^{\frac{5}{2}} \ d\eta \right)^{\frac{1}{2}}\left( \int_{0}^{\alpha}\vert\kappa\vert^{\frac{5}{2}} \ d\kappa \right)^{\frac{1}{2}} \  d\alpha \\
&\lesssim& \Vert f \Vert_{\dot H^3} \int \frac{\Vert \delta_\alpha f_x \Vert_{L^{\infty}}}{\vert\alpha\vert^{\frac{3}{2}}} \left(\int \frac{ \Vert s_{\kappa}f_x \Vert^2_{L^{4}} }{\vert\kappa\vert^{\frac{5}{2}}} d\kappa \right)^{\frac{1}{2}} \left(\int \frac{\Vert s_{\eta}f_x \Vert^2_{L^{4}}}{\vert \eta \vert^{\frac{5}{2}}}  d\eta \right)^{\frac{1}{2}}\  d\alpha \\
&\lesssim& \Vert f \Vert_{\dot H^3} \int \frac{\Vert \delta_\alpha f_x \Vert_{L^{\infty}}}{\vert\alpha\vert^\frac{3}{2}} \left(\int \frac{ \Vert s_{\kappa}f_x \Vert^2_{L^{4}} }{\vert\kappa\vert^{\frac{5}{2}}} d\kappa \right)^{\frac{1}{2}} \left(\int \frac{\Vert s_{\eta}f_x \Vert^2_{L^{4}}}{\vert \eta \vert^{\frac{5}{2}}}  d\eta \right)^{\frac{1}{2}}\  d\alpha \\
&\lesssim& \Vert f \Vert_{\dot H^3} \Vert f \Vert_{\dot B^\frac{3}{2}_{\infty,1}}\Vert f \Vert^2_{\dot B^\frac{7}{4}_{4,2}} \\
&\lesssim& \Vert f \Vert_{\dot H^3} \Vert f \Vert^{\frac{2}{3}}_{\dot B^1_{\infty,\infty}} \Vert f \Vert^{\frac{1}{3}}_{\dot B^{\frac{5}{2}}_{\infty,\infty}} \Vert f \Vert^2_{\dot H^2}\\
&\lesssim& \Vert f \Vert_{\dot H^3} \Vert f \Vert^{\frac{2}{3}}_{\dot H^{\frac{3}{2}}} \Vert f \Vert^{\frac{1}{3}}_{\dot H^3} \Vert f \Vert^{\frac{2}{3}}_{\dot H^3} \Vert f \Vert^{\frac{4}{3}}_{\dot H^\frac{3}{2} } 
 \end{eqnarray*}
 We conclude as $T_{4,1,1}$, that is
 \begin{eqnarray*}
T_{4,1,2,3}&\lesssim& \Vert f \Vert^2_{\dot H^3} \Vert f \Vert^2_{\dot H^\frac{3}{2} }
 \end{eqnarray*}

$\bullet$ {Estimate of $T_{4,1,3}$} \\

Recall that
\begin{eqnarray*}
T_{4,1,3}&=&\frac{1}{8\pi} \int \int\int_{0}^{\infty}\int_{0}^{\infty}\int_{0}^{\alpha} \ \gamma^2 e^{-\gamma-\sigma} \ \Lambda^{3} f \ \frac{1}{\alpha^2}  s_{\eta} f_{x} \ (\partial_\alpha D_\alpha f)^2 \sin(\frac{\gamma}{2}D_\alpha f)\cos(\frac{\gamma}{2}S_\alpha f)     \\
&& \  \times\ {  \delta_\alpha f_x } \cos(\sigma \tau_{\alpha} f_x )\cos(\arctan(\tau_{\alpha} f_{x})) \ d\eta  \ d\gamma \ d\sigma  \ d\alpha \ dx. \\
 \end{eqnarray*}
Since,
\begin{eqnarray*}
(\partial_{\alpha}D_\alpha f)^2&=&\left(  \frac{s_\alpha f_x}{\alpha}+\frac{1}{\alpha^2} \ {\int_0^\alpha s_\kappa  f _x\ d\kappa} \right)^2
\end{eqnarray*}
we find,
\begin{eqnarray*}
T_{4,1,3}&=&\frac{1}{8\pi} \int \int\int_{0}^{\infty}\int_{0}^{\infty}\int_{0}^{\alpha} \ \gamma^2 e^{-\gamma-\sigma} \ \Lambda^{3} f \ \frac{s_{\eta} f_{x} \ s^2_\alpha f_x}{\alpha^4}   \sin(\frac{\gamma}{2}D_\alpha f)\cos(\frac{\gamma}{2}S_\alpha f) \\
&& \  \times\ {  \delta_\alpha f_x } \cos(\sigma \tau_{\alpha} f_x )\cos(\arctan(\tau_{\alpha} f_{x})) \ d\eta  \ d\gamma \ d\sigma  \ d\alpha\ dx \\
&+&\frac{1}{4\pi} \int \int\int_{0}^{\infty}\int_{0}^{\infty}\int_{0}^{\alpha}\int_0^\alpha \ \gamma^2 e^{-\gamma-\sigma} \ \Lambda^{3} f \ \frac{s_\alpha f_x \ s_\kappa f_x \ s_{\eta} f_{x}}{\alpha^5}  \ \sin(\frac{\gamma}{2}D_\alpha f)\cos(\frac{\gamma}{2}S_\alpha f) \\
&& \  \times\ {  \delta_\alpha f_x } \cos(\sigma \tau_{\alpha} f_x )\cos(\arctan(\tau_{\alpha} f_{x})) \ d\kappa \ d\eta  \ d\alpha \ d\gamma \ d\sigma \ dx \\
&+&\frac{1}{8\pi} \int \int\int_{0}^{\infty}\int_{0}^{\infty}\int_{0}^{\alpha} \ \gamma^2 e^{-\gamma-\sigma} \ \Lambda^{3} f \ \frac{s_{\eta} f_{x}}{\alpha^5}   \ \left(\int_0^\alpha s_\kappa f_x \ d\kappa\right)^2 \sin(\frac{\gamma}{2}D_\alpha f)\cos(\frac{\gamma}{2}S_\alpha f)     \\
&& \  \times\ {  \delta_\alpha f_x } \cos(\sigma \tau_{\alpha} f_x )\cos(\arctan(\tau_{\alpha} f_{x})) \ d\eta  \ d\gamma \ d\sigma \ d\alpha \ dx \\
&:=&T_{4,1,3,1}+T_{4,1,3,2}+T_{4,1,3,3}
 \end{eqnarray*}
 
 $\bullet$ {{Estimate of $T_{4,1,3,1}$}} \\
 \begin{eqnarray*}
 T_{4,1,3,1}&=&\frac{1}{8\pi} \int \int\int_{0}^{\infty}\int_{0}^{\infty}\int_{0}^{\alpha} \ \gamma^2 e^{-\gamma-\sigma} \ \Lambda^{3} f \ \frac{s_{\eta} f_{x} \ s^2_\alpha f_x}{\alpha^4}   \sin(\frac{\gamma}{2}D_\alpha f)\cos(\frac{\gamma}{2}S_\alpha f) \\
&& \  \times\ {  \delta_\alpha f_x } \cos(\sigma \tau_{\alpha} f_x )\cos(\arctan(\tau_{\alpha} f_{x})) \ d\eta  \ d\gamma \ d\sigma  \ d\alpha \ dx \\
&\lesssim& \Vert f \Vert_{\dot H^3} \int \frac{ \Vert \delta_\alpha f_x \Vert_{L^{\infty}}  \  \Vert s^2_\alpha f_x \Vert_{L^{\infty}} }{\alpha^4} \left(\int_{0}^{\alpha} \frac{\Vert s_{\eta} f_{x} \Vert^{2}_{L^{2}}}{\eta^{2}} \ d\eta \right)^{\frac{1}{2}} \left(\int_{0}^{\alpha} \eta^{2} \ d\eta \right)^{\frac{1}{2}} \ d\alpha \\
&\lesssim& \Vert f \Vert_{\dot H^3} \Vert f \Vert_{\dot H^{\frac{3}{2}}} \int \frac{\Vert \delta_\alpha f_x \Vert_{L^{\infty}}  \ \Vert s_\alpha f_x \Vert^2_{L^{\infty}}}{\vert\alpha\vert^{\frac{5}{2}}}  \ d\alpha \\
&\lesssim& \Vert f \Vert_{\dot H^3} \Vert f \Vert_{\dot H^{\frac{3}{2}}} 
\left( \int \frac{\Vert \delta_\alpha f_x \Vert^2_{L^{\infty}}}{\vert\alpha\vert^{2}}  \ d\alpha \right)^{\frac{1}{2}}
\left( \int \frac{\Vert s_\alpha f_x \Vert^4_{L^{\infty}}}{\vert\alpha\vert^{3}}  \ d\alpha \right)^{\frac{1}{2}} \\
&\lesssim& \Vert f \Vert_{\dot H^3} \Vert f \Vert_{\dot H^{\frac{3}{2}}} \Vert f \Vert^3_{\dot B^{\frac{3}{2}}_{\infty,2}} \\
&\lesssim& \Vert f \Vert_{\dot H^3} \Vert f \Vert_{\dot H^{\frac{3}{2}}} \Vert f \Vert^3_{\dot H^2} \\
&\lesssim& \Vert f \Vert^2_{\dot H^3} \Vert f \Vert^2_{\dot H^{\frac{3}{2}}} 
\end{eqnarray*}
 $\bullet$ {{Estimate of $T_{4,1,3,2}$}} \\
  \begin{eqnarray*}
T_{4,1,3,2}&=&\frac{1}{4\pi} \int \int\int_{0}^{\infty}\int_{0}^{\infty}\int_{0}^{\alpha}\int_0^\alpha \ \gamma^2 e^{-\gamma-\sigma} \ \Lambda^{3} f \ \frac{s_\alpha f_x \ s_\kappa f_x \ s_{\eta} f_{x}}{\alpha^5}  \ \sin(\frac{\gamma}{2}D_\alpha f)\cos(\frac{\gamma}{2}S_\alpha f)     \\
&& \  \times\ {  \delta_\alpha f_x } \cos(\sigma \tau_{\alpha} f_x )\cos(\arctan(\tau_{\alpha} f_{x})) \ d\kappa \ d\eta   \ d\gamma \ d\sigma \ d\alpha \ dx \\
&\lesssim&\Vert f \Vert_{\dot H^3} \int  \frac{\Vert s_\alpha f_x\Vert_{L^{\infty}}\Vert \delta_\alpha f_x\Vert_{L^{\infty}}}{\vert\alpha\vert^5} \left(\int_0^\alpha \frac{\Vert s_\kappa f_x \Vert^2_{L^{4}}}{\vert \kappa \vert^{\frac{3}{2}}} \ d\kappa \right)^{\frac{1}{2}} \left(\int_0^\alpha \frac{ \Vert s_\eta f_x \Vert^2_{L^{4}}}{\vert \eta \vert^{\frac{3}{2}}} \ d\eta \right)^{\frac{1}{2}} \\
&& \left(\int_0^\alpha \vert\kappa\vert^{\frac{3}{2}} d\kappa\right)^{\frac{1}{2}} \left(\int_0^\alpha \vert\eta\vert^{\frac{3}{2}} d\eta\right)^{\frac{1}{2}}\ d\alpha \\ 
&\lesssim&\Vert f \Vert_{\dot H^3} \Vert f_x \Vert^2_{\dot B^{\frac{1}{4}}_{4,2}} \int  \frac{\Vert s_\alpha f_x\Vert_{L^{\infty}}\Vert \delta_\alpha f_x\Vert_{L^{\infty}}}{\vert\alpha\vert^{\frac{5}{2}}} \ d\alpha \\
&\lesssim&\Vert f \Vert_{\dot H^3} \Vert f \Vert^2_{\dot H^{\frac{3}{2}}} 
\left( \int \frac{\Vert \delta_\alpha f_x \Vert^2_{L^{\infty}}}{\vert\alpha\vert^{\frac{5}{2}}}  \ d\alpha \right)^{\frac{1}{2}}
\left( \int \frac{\Vert s_\alpha f_x \Vert^2_{L^{\infty}}}{\vert\alpha\vert^{\frac{5}{2}}}  \ d\alpha \right)^{\frac{1}{2}} \\
&\lesssim& \Vert f \Vert_{\dot H^3} \Vert f \Vert^2_{\dot H^{\frac{3}{2}}} 
\Vert f \Vert^2_{\dot B^{\frac{7}{4}}_{\infty,2}} \\
&\lesssim& \Vert f \Vert_{\dot H^3} \Vert f \Vert^2_{\dot H^{\frac{3}{2}}} 
\Vert f \Vert^2_{\dot H^{\frac{9}{4}}} \\
&\lesssim& \Vert f \Vert^2_{\dot H^3} \Vert f \Vert^3_{\dot H^{\frac{3}{2}}}
\end{eqnarray*}
 $\bullet$ {{Estimate of $T_{4,1,3,3}$}} \\
 \begin{eqnarray*}
T_{4,1,3,3}&=&\frac{1}{8\pi} \int \int\int_{0}^{\infty}\int_{0}^{\infty}\int_{0}^{\alpha} \ \gamma^2 e^{-\gamma-\sigma} \ \Lambda^{3} f \ \frac{s_{\eta} f_{x}}{\alpha^6}   \ \left(\int_0^\alpha s_\kappa f_x \ d\kappa\right)^2     \sin(\frac{\gamma}{2}D_\alpha f)  \\
&& \  \times\  \cos(\frac{\gamma}{2}S_\alpha f) \delta_\alpha f_x   \cos(\sigma \tau_{\alpha} f_x )\cos(\arctan(\tau_{\alpha} f_{x})) \ d\eta  \ d\gamma \ d\sigma \ d\alpha \ dx \\
&\lesssim&\Vert f \Vert_{\dot H^3} \int \frac{\Vert \delta_\alpha f_x \Vert_{L^\infty}}{\vert \alpha\vert^6}  \left(\int  \frac{\Vert s_\eta f_x \Vert^2_{L^{2}}}{\vert \eta\vert^2} \ d\eta\right)^{\frac{1}{2}}  \left(\int  \frac{\Vert s_\kappa f_x \Vert^2_{L^{\infty}}}{\vert \kappa\vert^{\frac{9}{4}}} \ d\eta\right) \\
&&\ \times  \int_0^\alpha \vert\kappa\vert^{\frac{9}{4}} d\kappa \  \left(\int_0^\alpha \vert\eta\vert^{2} d\kappa\right)^{\frac{1}{2}}d\alpha \\
&\lesssim& \Vert f \Vert_{\dot H^3} \Vert f \Vert_{\dot H^{\frac{3}{2}}} \Vert f_x \Vert^2_{\dot B^{\frac{5}{8}}_{\infty,2}} \int \frac{\Vert \delta_\alpha f_x \Vert_{L^\infty}}{\vert \alpha\vert^{\frac{5}{4}}} \ d\alpha \\
&\lesssim& \Vert f \Vert_{\dot H^3} \Vert f \Vert_{\dot H^{\frac{3}{2}}} \Vert f \Vert^2_{\dot H^{\frac{17}{8}}} \int \frac{\Vert \delta_\alpha f_x \Vert_{L^\infty}}{\vert \alpha\vert^{\frac{5}{4}}} \ d\alpha
\end{eqnarray*}
Since $\dot H^{\frac{17}{8}}=\left[\dot H^3, \dot H^\frac{3}{2} \right]_{\frac{5}{12},\frac{7}{12}}$ and $\dot B^{\frac{5}{4}}_{\infty,2}=\left[ \dot B^{1}_{\infty,\infty}, \dot B^{\frac{3}{2}}_{\infty,\infty}\right]_{\frac{1}{2},\frac{1}{2}}$ we find that

\begin{eqnarray*}
T_{4,1,3,3}&\lesssim& \Vert f \Vert_{\dot H^3} \Vert f \Vert_{\dot H^{\frac{3}{2}}} \Vert f \Vert^{\frac{5}{6}}_{\dot H^3}\Vert f \Vert^{\frac{7}{6}}_{\dot H^\frac{3}{2}} \Vert f \Vert^{\frac{1}{2}}_{\dot B^{1}_{\infty,\infty}} \Vert f \Vert^{\frac{1}{2}}_{\dot B^{\frac{3}{2}}_{\infty,\infty}}
\end{eqnarray*}
Then, using the fact that $\left[\dot H^3, \dot H^\frac{3}{2} \right]_{\frac{1}{3},\frac{2}{3}}=\dot H^{2} \hookrightarrow \dot B^{\frac{3}{2}}_{\infty,\infty} $
and $\dot H^{\frac{3}{2}} \hookrightarrow \dot B^{1}_{\infty,\infty} $ we find
 \begin{eqnarray*}
T_{4,1,3,3}&\lesssim&\Vert f \Vert_{\dot H^3} \Vert f \Vert_{\dot H^{\frac{3}{2}}} \Vert f \Vert^{\frac{5}{6}}_{\dot H^3}\Vert f \Vert^{\frac{7}{6}}_{\dot H^\frac{3}{2}}\Vert f \Vert^{\frac{1}{2}}_{\dot H^2} \Vert f \Vert^{\frac{1}{2}}_{\dot H^\frac{3}{2}} \\
&\lesssim&\Vert f \Vert_{\dot H^3} \Vert f \Vert_{\dot H^{\frac{3}{2}}} \Vert f \Vert^{\frac{5}{6}}_{\dot H^3}\Vert f \Vert^{\frac{7}{6}}_{\dot H^\frac{3}{2}}
\Vert f \Vert^{\frac{1}{6}}_{\dot H^3} \Vert f \Vert^{\frac{1}{3}}_{\dot H^\frac{3}{2}}
\Vert f \Vert^{\frac{1}{2}}_{\dot H^\frac{3}{2}} \\
\end{eqnarray*}
Finally, we find the following estimate
 \begin{eqnarray*}
T_{4,1,3,3}&\lesssim& \Vert f \Vert^2_{\dot H^3} \Vert f \Vert^3_{\dot H^{\frac{3}{2}}}
\end{eqnarray*}

$\bullet$ {Estimate of $T_{4,1,4}$} \\

We want to estimate 
\begin{eqnarray*}
T_{4,1,4}&=&\frac{1}{4\pi} \int \int\int_{0}^{\infty}\int_{0}^{\infty} \ \gamma e^{-\gamma-\sigma} \ \Lambda^{3} f \ \partial^2_\alpha S_\alpha f \sin(\frac{\gamma}{2}D_\alpha f)\sin(\frac{\gamma}{2}S_\alpha f)\frac{1}{\alpha^2} \int_{0}^{\alpha} s_{\eta} f_{x}    \ d\eta \\
&& \  \times\ {  \delta_\alpha f_x } \cos(\sigma \tau_{\alpha} f_x )\cos(\arctan(\tau_{\alpha} f_{x}))  \ d\gamma \ d\sigma  \ d\alpha \ dx \\
\end{eqnarray*}
Using the identity 
\begin{eqnarray*}
\partial^2_{\alpha}S_\alpha f &=&\frac{s_{\alpha} f_{xx}}{\alpha}-\frac{\int_0^\alpha s_\eta f_{xx} \ d \eta}{\alpha^2}+\frac{d_{\alpha} f_{x}}{\alpha^2}+\frac{s_{\alpha} f}{\alpha^3}
\end{eqnarray*}
We find the following decomposition
\begin{eqnarray*}
T_{4,1,4}&=&\frac{1}{4\pi} \int \int\int_{0}^{\infty}\int_{0}^{\infty} \ \gamma e^{-\gamma-\sigma} \ \Lambda^{3} f \ {s_{\alpha} f_{xx}} \sin(\frac{\gamma}{2}D_\alpha f)\sin(\frac{\gamma}{2}S_\alpha f)\frac{1}{\alpha^3} \int_{0}^{\alpha} s_{\eta} f_{x}    \ d\eta \\
&& \  \times\ {  \delta_\alpha f_x } \cos(\sigma \tau_{\alpha} f_x )\cos(\arctan(\tau_{\alpha} f_{x})) \ d\eta \ d\gamma \ d\sigma  \ d\alpha \ dx  \\
&-&\frac{1}{4\pi} \int \int\int_{0}^{\infty}\int_{0}^{\infty}\int_0^\alpha\int_0^\alpha \ \gamma e^{-\gamma-\sigma} \ \Lambda^{3} f \ { s_\kappa f_{xx} \ d \kappa}\sin(\frac{\gamma}{2}D_\alpha f)\sin(\frac{\gamma}{2}S_\alpha f)\frac{1}{\alpha^4}  \ s_{\eta} f_{x}     \\
&& \  \times\ {  \delta_\alpha f_x } \cos(\sigma \tau_{\alpha} f_x )\cos(\arctan(\tau_{\alpha} f_{x})) \ d\kappa \ d\eta \ d\gamma \ d\sigma  \ d\alpha \ dx  \\
&+&\frac{1}{4\pi} \int \int\int_{0}^{\infty}\int_{0}^{\infty}\int_{0}^{\alpha} \ \gamma e^{-\gamma-\sigma} \ \Lambda^{3} f \ {D_{\alpha} f_x} \sin(\frac{\gamma}{2}D_\alpha f)\sin(\frac{\gamma}{2}S_\alpha f)\frac{1}{\alpha^3}  s_{\eta} f_{x}    \\
&& \  \times\ {  \delta_\alpha f_x } \cos(\sigma \tau_{\alpha} f_x )\cos(\arctan(\tau_{\alpha} f_{x})) \ d\eta \ d\gamma \ d\sigma \ d\alpha \ dx  \\
&+&\frac{1}{4\pi} \int \int\int_{0}^{\infty}\int_{0}^{\infty}\int_{0}^{\alpha} \ \gamma e^{-\gamma-\sigma} \ \Lambda^{3} f \ {s_{\alpha} f} \sin(\frac{\gamma}{2}D_\alpha f)\sin(\frac{\gamma}{2}S_\alpha f)\frac{1}{\alpha^4}  s_{\eta} f_{x}    \\
&& \  \times\ {  \delta_\alpha f_x } \cos(\sigma \tau_{\alpha} f_x )\cos(\arctan(\tau_{\alpha} f_{x})) \ d\eta \ d\gamma \ d\sigma \ d\alpha \ dx  \\
&:=& T_{4,1,4,1} + T_{4,1,4,2} + T_{4,1,4,3} + T_{4,1,4,4}
\end{eqnarray*}
 $\bullet$ {{Estimate of $T_{4,1,4,1}$}} \\

We have to balance the regularity and this is done through the artificial weight, more precisely, we have 
\begin{eqnarray*}
T_{4,1,4,1}&=&\frac{1}{4\pi} \int \int\int_{0}^{\infty}\int_{0}^{\infty} \ \gamma e^{-\gamma-\sigma} \ \Lambda^{3} f \ {s_{\alpha} f_{xx}} \sin(\frac{\gamma}{2}D_\alpha f)\sin(\frac{\gamma}{2}S_\alpha f) \\
&& \  \times \frac{1}{\alpha^3} \left(\int_{0}^{\alpha} s_{\eta} f_{x}    \ d\eta\right) \ {\delta_\alpha f_x } \cos(\sigma \tau_{\alpha} f_x )\cos(\arctan(\tau_{\alpha} f_{x}))  \ d\gamma \ d\sigma  \ d\alpha \ dx  \\
&\lesssim&\Vert f \Vert_{\dot H^3} \int \frac{\Vert s_{\alpha} f_{xx} \Vert_{L^{\infty}} \Vert\delta_\alpha f_x\Vert_{L^{\infty}}}{\vert \alpha\vert^{3}}  \left(\int_{0}^{\alpha} \frac{\Vert s_{\eta} f_{x} \Vert_{L^2}}{\vert \eta \vert^2}    \ d\eta \right)^{\frac{1}{2}} \left(\int_{0}^{\alpha} {\vert \eta \vert^2}    \ d\eta \right)^{\frac{1}{2}} \ d\alpha \\
&\lesssim&\Vert f \Vert_{\dot H^3} \Vert f \Vert_{\dot H^\frac{3}{2}}\int \frac{\Vert s_{\alpha} f_{xx} \Vert_{L^{\infty}} \Vert\delta_\alpha f_x\Vert_{L^{\infty}}}{\vert \alpha\vert^{\frac{3}{2}}}  \ d\alpha \\
&\lesssim&\Vert f \Vert_{\dot H^3} \Vert f \Vert_{\dot H^\frac{3}{2}}
\Vert f \Vert_{\dot B^{\frac{9}{4}}_{\infty,2}}\Vert f \Vert_{\dot B^{\frac{5}{4}}_{\infty,2}} \\
&\lesssim&\Vert f \Vert_{\dot H^3} \Vert f \Vert_{\dot H^\frac{3}{2}} 
\Vert f \Vert_{\dot H^\frac{11}{4}} \Vert f \Vert_{\dot H^\frac{7}{4}}
\end{eqnarray*}
Using that $\dot H^{\frac{11}{4}}=\left[\dot H^3, \dot H^\frac{3}{2} \right]_{\frac{5}{6},\frac{1}{6}}$ and $\dot H^{\frac{7}{4}}=\left[\dot H^3, \dot H^\frac{3}{2} \right]_{\frac{1}{6},\frac{5}{6}}$, we find that
\begin{eqnarray*}
T_{4,1,4,1}&\lesssim& \Vert f \Vert^2_{\dot H^3} \Vert f \Vert^2_{\dot H^{\frac{3}{2}}}
\end{eqnarray*}

 $\bullet$ {{Estimate of $T_{4,1,4,2}$}} \\

For this term, we can notice that we have to use different artificial weights to milder the singularity in $\alpha$ and get critical estimates
\begin{eqnarray*}
T_{4,1,4,2}&=&-\frac{1}{4\pi} \int \int\int_{0}^{\infty}\int_{0}^{\infty}\int_0^\alpha\int_0^\alpha \ \gamma e^{-\gamma-\sigma} \ \Lambda^{3} f \ { s_\kappa f_{xx} }\sin(\frac{\gamma}{2}D_\alpha f)\sin(\frac{\gamma}{2}S_\alpha f)     \\
&& \  \times \frac{1}{\alpha^4}  \ s_{\eta} f_{x} {  \delta_\alpha f_x } \cos(\sigma \tau_{\alpha} f_x )\cos(\arctan(\tau_{\alpha} f_{x})) \ d\kappa \ d\eta \ d\gamma \ d\sigma  \ d\alpha \ dx  \\
&\lesssim& \Vert f \Vert_{\dot H^3} \int \frac{\Vert \delta_\alpha f_x \Vert_{L^\infty}}{\alpha^4} \left(\int_0^\alpha \frac{ \Vert s_\kappa f_{xx} \Vert^2_{L^2} }{\vert \kappa\vert^{\frac{3}{2}}} \ d\kappa\right)^{\frac{1}{2}}  \left(\int_0^\alpha \frac{ \Vert s_\eta f_{x} \Vert^2_{L^2} }{\vert \eta\vert^{\frac{3}{2}}} \ d\kappa\right)^{\frac{1}{2}} \\
&& \times \left(\int_0^\alpha \vert \kappa\vert^{\frac{3}{2}} \ d\kappa\right)^{\frac{1}{2}}  \left(\int_0^\alpha \vert \eta\vert^{\frac{3}{2}} \ d\eta\right)^{\frac{1}{2}} \ d\alpha \\
&\lesssim& \Vert f \Vert_{\dot H^3} \int \frac{\Vert \delta_\alpha f_x \Vert_{L^\infty}}{\vert\alpha\vert^{\frac{3}{2}}} \left(\int \frac{ \Vert s_\kappa f_{xx} \Vert^2_{L^2} }{\vert \kappa\vert^{\frac{3}{2}}} \ d\kappa\right)^{\frac{1}{2}}  \left(\int \frac{ \Vert s_\eta f_{x} \Vert^2_{L^\infty} }{\vert \eta\vert^{\frac{3}{2}}} \ d\kappa\right)^{\frac{1}{2}} \ d\alpha\\
&\lesssim& \Vert f \Vert_{\dot H^3}\Vert f \Vert_{\dot H^{\frac{9}{4}}}\Vert f \Vert_{\dot B^{\frac{5}{4}}_{\infty,2}} \Vert f \Vert_{\dot B^{\frac{3}{2}}_{\infty,1}} \\
&\lesssim& \Vert f \Vert_{\dot H^3}\Vert f \Vert_{\dot H^{\frac{9}{4}}}\Vert f \Vert_{\dot H^{\frac{7}{4}}} \Vert f \Vert^{\frac{1}{2}}_{\dot B^{1}_{\infty,\infty}} \Vert f \Vert^{\frac{1}{2}}_{\dot B^{2}_{\infty,1}} \\
&\lesssim& \Vert f \Vert_{\dot H^3}\Vert f \Vert^{\frac{1}{2}}_{\dot H^{\frac{3}{2}}} \Vert f \Vert^{\frac{1}{2}}_{\dot H^{3}} \Vert f \Vert^{\frac{1}{6}}_{\dot H^{3}} \Vert f \Vert^{\frac{5}{6}}_{\dot H^{\frac{3}{2}}}  \Vert f \Vert^{\frac{1}{2}}_{\dot H^{\frac{3}{2}}} \Vert f \Vert^{\frac{1}{2}}_{\dot H^{\frac{5}{2}}} 
\end{eqnarray*}
But $\dot H^{\frac{5}{2}}=\left[\dot H^3, \dot H^\frac{3}{2} \right]_{\frac{2}{3},\frac{1}{3}}$, hence
\begin{eqnarray*}
T_{4,1,4,2}&\lesssim&\Vert f \Vert_{\dot H^3}\Vert f \Vert^{\frac{1}{2}}_{\dot H^{\frac{3}{2}}} \Vert f \Vert^{\frac{1}{2}}_{\dot H^{3}} 
\Vert f \Vert^{\frac{1}{6}}_{\dot H^{3}} \Vert f \Vert^{\frac{5}{6}}_{\dot H^{\frac{3}{2}}}  \Vert f \Vert^{\frac{1}{2}}_{\dot H^{\frac{3}{2}}} 
\Vert f \Vert^{\frac{1}{3}}_{\dot H^{3}} \Vert f \Vert^{\frac{1}{6}}_{\dot H^{\frac{3}{2}}}\\
&\lesssim& \Vert f \Vert^2_{\dot H^3} \Vert f \Vert^3_{\dot H^{\frac{3}{2}}}
\end{eqnarray*}

 $\bullet$ {{Estimate of $T_{4,1,4,3}$}} \\
 
 We have to consider the integral in $\eta$ which requires the introduction of a weight. One also notes that the term $D_{\alpha} f_x$ is as regular as $\Delta_{\alpha} f_x$ (or $\bar \Delta_{\alpha} f_x$), hence  using the notation $\delta^{\pm}_\alpha f_x=f(x\pm\alpha)-f(x)$, we find
 
\begin{eqnarray*}
T_{4,1,4,3}=&&\frac{1}{4\pi} \int \int\int_{0}^{\infty}\int_{0}^{\infty}\int_{0}^{\alpha} \ \gamma e^{-\gamma-\sigma} \ \Lambda^{3} f \ {D_{\alpha} f_x} \sin(\frac{\gamma}{2}D_\alpha f)\sin(\frac{\gamma}{2}S_\alpha f)    \\
&& \  \times \frac{1}{\alpha^3}  s_{\eta} f_{x} {  \delta_\alpha f_x } \cos(\sigma \tau_{\alpha} f_x )\cos(\arctan(\tau_{\alpha} f_{x})) \ d\eta \ d\gamma \ d\sigma \ d\alpha \ dx \\
 &\lesssim& \Vert f \Vert_{\dot H^3} \int \frac{\Vert \delta^\pm_{\alpha} f_x \Vert_{L^\infty} \Vert \delta_{\alpha} f_x \Vert_{L^\infty}}{\alpha^4} \left(\int_0^\alpha \frac{\Vert s_{\eta}f_x \Vert^2_{L^2}}{\vert \eta\vert^2} \ d\eta\right)^{\frac{1}{2}} \left(\int_0^\alpha \vert \eta\vert^2 \ d\eta\right)^{\frac{1}{2}} \ d\alpha \\
 &\lesssim& \Vert f \Vert_{\dot H^3} \Vert f \Vert_{\dot H^{\frac{3}{2}}} \int \frac{\Vert \delta^\pm_{\alpha} f_x \Vert_{L^\infty} \Vert \delta_{\alpha} f_x \Vert_{L^\infty}}{\vert\alpha\vert^{\frac{5}{2}}}  \ d\alpha \\
 &\lesssim& \Vert f \Vert_{\dot H^3} \Vert f \Vert_{\dot H^{\frac{3}{2}}}
  \int \frac{\Vert \delta^\pm_{\alpha} f_x \Vert_{L^\infty} \Vert \delta_{\alpha} f_x \Vert_{L^\infty}}{\vert\alpha\vert^{\frac{5}{2}}}  \ d\alpha \\
  &\lesssim& \Vert f \Vert_{\dot H^3} \Vert f \Vert_{\dot H^{\frac{3}{2}}} \Vert f_x \Vert^2_{\dot B^{\frac{1}{4}}_{\infty,2}} \\
  &\lesssim& \Vert f \Vert_{\dot H^3} \Vert f \Vert_{\dot H^{\frac{3}{2}}} \Vert f_x \Vert^2_{\dot B^{\frac{3}{4}}_{\infty,2}} \\
  &\lesssim& \Vert f \Vert_{\dot H^3} \Vert f \Vert_{\dot H^{\frac{3}{2}}} \Vert f_x \Vert^2_{\dot H^{\frac{9}{4}}}\\
  &\lesssim& \Vert f \Vert^2_{\dot H^3} \Vert f \Vert^2_{\dot H^{\frac{3}{2}}}
\end{eqnarray*}

$\bullet$ {{Estimate of $T_{4,1,4,4}$}} \\

For this term, we may  equally balance the artificial weight, we find
\begin{eqnarray*}
T_{4,1,4,4}&=&\frac{1}{4\pi} \int \int\int_{0}^{\infty}\int_{0}^{\infty}\int_{0}^{\alpha} \ \gamma e^{-\gamma-\sigma} \ \Lambda^{3} f \ {s_{\alpha} f} \sin(\frac{\gamma}{2}D_\alpha f)\sin(\frac{\gamma}{2}S_\alpha f)\frac{1}{\alpha^5}  s_{\eta} f_{x}    \\
&& \  \times\ {  \delta_\alpha f_x } \cos(\sigma \tau_{\alpha} f_x )\cos(\arctan(\tau_{\alpha} f_{x})) \ d\eta \ d\gamma \ d\sigma \ d\alpha \ dx  \\
&\lesssim& \Vert f \Vert_{\dot H^3} \int \frac{\Vert \delta_\alpha f_x \Vert_{L^{\infty}} \Vert s_\alpha f \Vert_{L^{\infty}}}{\vert \alpha \vert^5} \left(\int_0^\alpha \frac{\Vert s_{\eta}f_x \Vert^2_{L^2}}{\vert \eta\vert^3} \ d\eta\right)^{\frac{1}{2}} \left(\int_0^\alpha \vert \eta\vert^3 \ d\eta\right)^{\frac{1}{2}} \ d\alpha \\
&\lesssim& \Vert f \Vert_{\dot H^3} \Vert f \Vert_{\dot H^{2}} \int \frac{\Vert \delta_\alpha f_x \Vert_{L^{\infty}} \Vert s_\alpha f \Vert_{L^{\infty}}}{\vert \alpha \vert^{3}}  \ d\alpha \\
&\lesssim& \Vert f \Vert_{\dot H^3} \Vert f \Vert_{\dot H^{2}} 
\left(\int \frac{\Vert \delta_\alpha f_x \Vert^2_{L^{\infty}}}{\vert \alpha \vert^2} \ d\alpha\right)^{\frac{1}{2}} \left(\int \frac{\Vert s_\alpha f \Vert^2_{L^{\infty}}}{\vert \alpha \vert^4} \ d\alpha\right)^{\frac{1}{2}}  \\
&\lesssim& \Vert f \Vert_{\dot H^3} \Vert f \Vert^3_{\dot H^{2}} \\
&\lesssim& \Vert f \Vert^2_{\dot H^3} \Vert f \Vert_{\dot H^{\frac{3}{2}}} 
\end{eqnarray*}

$\bullet$ {Estimate of $T_{4,1,5}$} \\

We have,

\begin{eqnarray*}
T_{4,1,5}&=&\frac{1}{8\pi} \int \int\int_{0}^{\infty}\int_{0}^{\infty} \ \gamma^2 e^{-\gamma-\sigma} \ \Lambda^{3} f \ (\partial_\alpha S_\alpha f)^2 \sin(\frac{\gamma}{2}D_\alpha f)\cos(\frac{\gamma}{2}S_\alpha f) \\
&& \  \times  \frac{1}{\alpha^2} \left(\int_{0}^{\alpha} s_{\eta} f_{x} \ d\eta \right) \ {\delta_\alpha f_x } \cos(\sigma \tau_{\alpha} f_x )\cos(\arctan(\tau_{\alpha} f_{x}))   \ d\gamma \ d\sigma \ d\alpha\ dx. \\
\end{eqnarray*}
Since,
\begin{eqnarray*}
\left(\partial_{\alpha}S_\alpha f\right)^2 &=& \left(\bar\Delta_\alpha f_x-\Delta_\alpha f_x-\frac{s_\alpha f}{\alpha^2}\right)^2
\end{eqnarray*}
As $\bar\Delta_\alpha f_x$ and $\Delta_\alpha f_x$ have the same properties of regularity then it suffices to consider 
\begin{eqnarray*}
\left(\partial_{\alpha}S_\alpha f\right)^2 &\approx& \Delta^2_\alpha f_x +2\Delta_\alpha f_x \ \frac{s_\alpha f}{\alpha^2} + \left( \frac{s_\alpha f}{\alpha^2}\right)^2
\end{eqnarray*}
Therefore, 
\begin{eqnarray*}
T_{4,1,5}&\lesssim& \Vert f \Vert_{\dot H^3} \int \frac{\Vert \delta_\alpha f_x\Vert^3_{L^\infty}}{\vert \alpha \vert^4} \left(\int_{0}^{\alpha} \frac{\Vert s_{\eta} f_{x} \Vert_{L^2}}{\vert \eta \vert^2}    \ d\eta \right)^{\frac{1}{2}}\left(\int_{0}^{\alpha} {\vert \eta \vert^2}    \ d\eta \right)^{\frac{1}{2}} \ d\alpha \\
&& \ + \ \Vert f \Vert_{\dot H^3}\int \frac{\Vert \delta_\alpha f_x\Vert^2_{L^\infty} \Vert s_\alpha f\Vert_{L^\infty}}{\vert \alpha \vert^5} \left(\int_{0}^{\alpha} \frac{\Vert s_{\eta} f_{x} \Vert_{L^2}}{\vert \eta \vert^2}    \ d\eta \right)^{\frac{1}{2}}\left(\int_{0}^{\alpha} {\vert \eta \vert^2}    \ d\eta \right)^{\frac{1}{2}} \ d\alpha \\
&& \ + \ \Vert f \Vert_{\dot H^3}\int \frac{\Vert \delta_\alpha f_x\Vert_{L^\infty} \Vert s_\alpha f\Vert^2_{L^\infty}}{\vert \alpha \vert^6} \left(\int_{0}^{\alpha} \frac{\Vert s_{\eta} f_{x} \Vert_{L^2}}{\vert \eta \vert^2}    \ d\eta \right)^{\frac{1}{2}}\left(\int_{0}^{\alpha} {\vert \eta \vert^2}    \ d\eta \right)^{\frac{1}{2}} \ d\alpha \\
&\lesssim&\ \Vert f \Vert_{\dot H^3} \Vert f \Vert_{\dot H^{\frac{3}{2}}} \int \frac{\Vert \delta_\alpha f_x\Vert^3_{L^\infty}}{\vert \alpha \vert^{\frac{5}{2}}} \ d\alpha + \ \Vert f \Vert_{\dot H^3} \Vert f \Vert_{\dot H^{\frac{3}{2}}} \int \frac{\Vert \delta_\alpha f_x\Vert^2_{L^\infty} \Vert s_\alpha f\Vert_{L^\infty}}{\vert \alpha \vert^{\frac{7}{2}}} \ d\alpha \\
&& \ + \  \Vert f \Vert_{\dot H^3} \Vert f \Vert_{\dot H^{\frac{3}{2}}}\int \frac{\Vert \delta_\alpha f_x\Vert_{L^\infty} \Vert s_\alpha f\Vert^2_{L^\infty}}{\vert \alpha \vert^\frac{9}{2}} \ d\alpha \\
&\lesssim&\ \Vert f \Vert_{\dot H^3} \Vert f \Vert_{\dot H^{\frac{3}{2}}} \Vert f \Vert^3_{\dot B^{\frac{3}{2}}_{\infty,3}} + \ \Vert f \Vert_{\dot H^3} \Vert f \Vert_{\dot H^{\frac{3}{2}}} \left(\sup_{\alpha \in \mathbb R} \frac{\Vert s_\alpha f \Vert_{L^{\infty}}}{\vert \alpha\vert} \right) \int \frac{\Vert \delta_\alpha f_x\Vert^2_{L^\infty} }{\vert \alpha \vert^{\frac{5}{2}}} \ d\alpha \\
&& \ + \  \Vert f \Vert_{\dot H^3} \Vert f \Vert_{\dot H^{\frac{3}{2}}}  \int \frac{\Vert \delta_\alpha f_x\Vert_{L^\infty} \Vert s_\alpha f\Vert^2_{L^\infty}}{\vert \alpha \vert^\frac{9}{2}} \ d\alpha \\
&\lesssim&\ \Vert f \Vert_{\dot H^3} \Vert f \Vert_{\dot H^{\frac{3}{2}}} \Vert f \Vert^3_{\dot B^{\frac{3}{2}}_{\infty,3}} + \ \Vert f \Vert_{\dot H^3} \Vert f \Vert_{\dot H^{\frac{3}{2}}} \Vert f \Vert_{\dot B^{1}_{\infty,\infty}} \Vert f \Vert^2_{\dot B^{\frac{7}{4}}_{\infty,2}}  \\
&& \ + \  \Vert f \Vert_{\dot H^3} \Vert f \Vert_{\dot H^{\frac{3}{2}}} \left( \int \frac{\Vert \delta_\alpha f_x\Vert^2_{L^\infty}}{\vert \alpha\vert^2} \ d\alpha\right)^{\frac{1}{2}} \left( \int \frac{\Vert s_\alpha f\Vert^4_{L^\infty}}{\vert \alpha\vert^7} \ d\alpha\right)^{\frac{1}{2}}   \\
&\lesssim&\ \Vert f \Vert_{\dot H^3} \Vert f \Vert_{\dot H^{\frac{3}{2}}} \Vert f \Vert^3_{\dot B^{\frac{3}{2}}_{\infty,3}} + \ \Vert f \Vert_{\dot H^3} \Vert f \Vert_{\dot H^{\frac{3}{2}}} \Vert f \Vert_{\dot B^{1}_{\infty,\infty}} \Vert f \Vert^2_{\dot B^{\frac{7}{4}}_{\infty,2}}  \\
&& \ + \  \Vert f \Vert_{\dot H^3} \Vert f \Vert_{\dot H^{\frac{3}{2}}} \Vert f_x \Vert_{\dot B^{\frac{1}{2}}_{\infty,2}} \Vert f \Vert^2_{\dot B^{\frac{3}{2}}_{\infty,2}}   \\
\end{eqnarray*}
Then, since $\dot H^2 \hookrightarrow \dot B^{\frac{3}{2}}_{\infty,3}$ and 
$\dot H^{\frac{3}{2}}\hookrightarrow\dot B^{1}_{\infty,\infty}$. In particular,
$$
\Vert f_x \Vert_{\dot B^{\frac{1}{2}}_{\infty,2}} \Vert f \Vert^2_{\dot B^{\frac{3}{2}}_{\infty,2}} \lesssim \Vert f \Vert^3_{\dot H^3}
$$
Hence, since $\left[\dot H^3, \dot H^{\frac{3}{2}}\right]_{\frac{1}{2},\frac{1}{2}}=\dot H^{\frac{9}{4}} \hookrightarrow \dot B^{\frac{7}{4}}_{\infty,2}$, and $\dot H^2 =\left[\dot H^3, \dot H^{\frac{3}{2}}\right]_{\frac{1}{3},\frac{2}{3}}$
we find that

\begin{eqnarray*}
T_{4,1,5}&\lesssim& \Vert f \Vert_{\dot H^3} \Vert f \Vert_{\dot H^{\frac{3}{2}}} \Vert f \Vert^3_{\dot H^2} + \Vert f \Vert_{\dot H^3} \Vert f \Vert^2_{\dot H^{\frac{3}{2}}}  \Vert f \Vert^2_{\dot H^{\frac{9}{4}}} \\
&\lesssim& \Vert f \Vert^2_{\dot H^3} \Vert f \Vert^2_{\dot H^{\frac{3}{2}}}  + \Vert f \Vert^2_{\dot H^3} \Vert f \Vert^3_{\dot H^{\frac{3}{2}}}  
\end{eqnarray*}

$\bullet$ {Estimate of $T_{4,1,6}$} \\

Using the identity
\begin{eqnarray*}
\partial^2_{\alpha}D_\alpha f &=&\frac{d_\alpha f_{xx}}{\alpha}+2 \frac{s_\alpha f_x}{\alpha^2} + \frac{2}{\alpha^3}\int_0^\alpha s_{\eta}f_x \ d\eta
\end{eqnarray*}
we may estimate $T_{4,1,6}$ as follows

\begin{eqnarray*}
T_{4,1,6}&=&-\frac{1}{2\pi} \int \int\int_{0}^{\infty}\int_{0}^{\infty} \ e^{-\gamma-\sigma} \ \Lambda^{3} f \ \left[\frac{d_\alpha f_{xx}}{\alpha}+2 \frac{s_\alpha f_{x}}{\alpha^2} + \frac{2}{\alpha^3}\int_0^\alpha s_{\eta}f_x \ d\eta \right] \frac{1}{\alpha}\\
&& \ \times \sin(\frac{\gamma}{2}D_\alpha f)\cos(\frac{\gamma}{2}S_\alpha f)  {  \delta_\alpha f_x } \cos(\sigma \tau_{\alpha} f_x )\cos(\arctan(\tau_{\alpha} f_{x}))   \ d\gamma \ d\sigma  \ d\alpha \ dx\\
&\lesssim&\Vert f \Vert_{\dot H^3} \left( \int \frac{\Vert \delta_\alpha f_{xx} \Vert_{L^2} \Vert \delta_\alpha f_x \Vert_{L^\infty}}{\alpha^2} \ d\alpha + \int \frac{\Vert \bar\delta_\alpha f_{xx} \Vert_{L^\infty}\Vert \delta_\alpha f_x \Vert_{L^2}}{\alpha^2} \ d\alpha \right) \\
&& \ + \ \Vert f \Vert_{\dot H^3}  \int \frac{\Vert s_\alpha f_{x} \Vert_{L^2} \Vert \delta_\alpha f_x \Vert_{L^\infty}}{\alpha^3} \ d\alpha \\
&& \ + \ \Vert f \Vert_{\dot H^3} \int  \frac{\Vert \delta_\alpha f_x \Vert_{L^{\infty}}}{\alpha^4} \left(\int_0^\alpha\frac{\Vert s_\eta f_x \Vert^2_{L^2}}{\eta^4} \ d\eta \ \right)^{\frac{1}{2}} \left(\int_0^\alpha \eta^4 \ d\eta \ \right)^{\frac{1}{2}} \ d\alpha \\
&\lesssim&\Vert f \Vert_{\dot H^3} \left( \int \frac{\Vert \delta_\alpha f_{xx} \Vert^2_{L^2} }{\alpha^2} \ d\alpha \right)^{\frac{1}{2}} \left( \left( \int \frac{\Vert \delta_\alpha f_{x} \Vert^2_{L^\infty}}{\alpha^2} \ d\alpha \right)^{\frac{1}{2}}+ \left( \int \frac{\Vert \bar\delta_\alpha f_{x} \Vert^2_{L^\infty}}{\alpha^2} \ d\alpha \right)^{\frac{1}{2}}\right) \\
&& \ + \ \Vert f \Vert_{\dot H^3}  \left(\int \frac{\Vert s_\alpha f_{x} \Vert^2_{L^2} }{\alpha^4} \ d\alpha \right)^{\frac{1}{2}}  \left(\int \frac{\Vert \delta_\alpha f_x \Vert^2_{L^\infty}}{\alpha^2} \ d\alpha \right)^{\frac{1}{2}} \\
&& \ + \ \Vert f \Vert_{\dot H^3} \Vert f \Vert_{\dot H^{\frac{5}{2}}} \int  \frac{\Vert \delta_\alpha f_x \Vert_{L^{\infty}}}{\vert\alpha\vert^{\frac{3}{2}}}   d\alpha \\
&\lesssim&\Vert f \Vert_{\dot H^3} \Vert f \Vert_{\dot H^{\frac{5}{2}}}\left(  
\Vert f \Vert_{\dot B^\frac{3}{2}_{\infty,2}}  + \Vert f \Vert_{\dot B^\frac{3}{2}_{\infty,1}}   \right) \\
&\lesssim&\Vert f \Vert_{\dot H^3} \Vert f \Vert_{\dot H^{\frac{5}{2}}}\left(  
\Vert f \Vert_{\dot H^2}  + \Vert f \Vert^{\frac{1}{2}}_{\dot B^1_{\infty,\infty}}\Vert f \Vert^{\frac{1}{2}}_{\dot B^{2}_{\infty,\infty}}   \right)
\end{eqnarray*}
Since $\dot H^2 =\left[\dot H^3, \dot H^{\frac{3}{2}}\right]_{\frac{1}{3},\frac{2}{3}}$ and $\dot H^{\frac{5}{2}}=\left[\dot H^3, \dot H^\frac{3}{2} \right]_{\frac{2}{3},\frac{1}{3}}$, then
 $\Vert f \Vert_{\dot H^{\frac{5}{2}}}\Vert f \Vert_{\dot H^2}\leq \Vert f \Vert_{\dot H^3} \Vert f \Vert_{\dot H^{\frac{3}{2}}}  $
We also have that $\dot H^{\frac{3}{2}} \hookrightarrow \dot B^1_{\infty,\infty}$ and $\dot H^{\frac{5}{2}} \hookrightarrow \dot B^2_{\infty,\infty}$ therefore
$$
 \Vert f \Vert_{\dot H^{\frac{5}{2}}}\Vert f \Vert^{\frac{1}{2}}_{\dot B^1_{\infty,\infty}}\Vert f \Vert^{\frac{1}{2}}_{\dot B^{2}_{\infty,\infty}} \lesssim\Vert f \Vert^{\frac{1}{2}}_{\dot H^{\frac{3}{2}}}\Vert f \Vert^{\frac{3}{2}}_{\dot H^{\frac{5}{2}}}
$$
Since $\Vert f \Vert^{\frac{3}{2}}_{\dot H^{\frac{5}{2}}}\lesssim \Vert f \Vert^{\frac{1}{2}}_{\dot H^{\frac{3}{2}}} \Vert f \Vert_{\dot H^3}$, then
$$
 \Vert f \Vert_{\dot H^{\frac{5}{2}}}\Vert f \Vert^{\frac{1}{2}}_{\dot B^1_{\infty,\infty}}\Vert f \Vert^{\frac{1}{2}}_{\dot B^{2}_{\infty,\infty}} \lesssim \Vert f \Vert_{\dot H^3} \Vert f \Vert_{\dot H^{\frac{3}{2}}}
$$
Therefore,
\begin{eqnarray*}
T_{4,1,6}\lesssim \Vert f \Vert^2_{\dot H^3} \Vert f \Vert_{\dot H^{\frac{3}{2}}}
\end{eqnarray*}

$\bullet$ {Estimate of $T_{4,1,7}$} \\

We have,

\begin{eqnarray*}
T_{4,1,7}&=&\frac{1}{2\pi} \int \int\int_{0}^{\infty}\int_{0}^{\infty}\int_{0}^{\alpha} \ \gamma e^{-\gamma-\sigma} \ \Lambda^{3} f \ \partial_{\alpha}D_\alpha f \cos(\frac{\gamma}{2}D_\alpha f)\cos(\frac{\gamma}{2}S_\alpha f)\frac{1}{\alpha^3}  s_{\eta} f_{x}     \\
&& \  \times\ {  \delta_\alpha f_x } \cos(\sigma \tau_{\alpha} f_x )\cos(\arctan(\tau_{\alpha} f_{x})) \ d\eta   \ d\gamma \ d\sigma \ d\alpha \ dx\\
\end{eqnarray*}
Since
\begin{eqnarray*}
\partial_{\alpha}D_\alpha f&=& - \frac{s_\alpha f_x}{\alpha}-\frac{1}{\alpha^2} \ {\int_0^\alpha s_\kappa f_x \ d\kappa}
\end{eqnarray*}
We find that
\begin{eqnarray*}
T_{4,1,7}&=&-\frac{1}{2\pi} \int \int\int_{0}^{\infty}\int_{0}^{\infty} \int_{0}^{\alpha} \ \gamma e^{-\gamma-\sigma} \ \Lambda^{3} f \ s_\alpha f_x \ \cos(\frac{\gamma}{2}D_\alpha f)\cos(\frac{\gamma}{2}S_\alpha f)\frac{1}{\alpha^4}  s_{\eta} f_{x}     \\
&& \  \times\ {  \delta_\alpha f_x } \cos(\sigma \tau_{\alpha} f_x )\cos(\arctan(\tau_{\alpha} f_{x})) \ d\eta  \ d\gamma \ d\sigma \ d\alpha \ dx \\
&-&\frac{1}{2\pi} \int \int\int_{0}^{\infty}\int_{0}^{\infty}\int_{0}^{\alpha}\int_{0}^{\alpha} \ \gamma e^{-\gamma-\sigma} \ \Lambda^{3} f \ s_\kappa f_x \ \cos(\frac{\gamma}{2}D_\alpha f)\cos(\frac{\gamma}{2}S_\alpha f)\frac{1}{\alpha^5}  s_{\eta} f_{x}    \ d\eta \\
&& \  \times\ {  \delta_\alpha f_x } \cos(\sigma \tau_{\alpha} f_x )\cos(\arctan(\tau_{\alpha} f_{x})) \ d\kappa \ d\eta  \ d\gamma \ d\sigma  \ d\alpha \ dx \\
&:=& T_{4,1,7,1} + T_{4,1,7,2} 
\end{eqnarray*}

$\bullet$ {{Estimate of $T_{4,1,7,1}$}} \\
\begin{eqnarray} \label{T_{4,1,7,1}}
T_{4,1,7,1}&\lesssim& \Vert f \Vert_{\dot H^3} \int    \ \frac{\Vert s_\alpha f_x \Vert_{L^\infty} \Vert \delta_\alpha f_x \Vert_{L^\infty} }{\alpha^4}  \left(\int_{0}^{\alpha} \frac{\Vert s_{\eta} f_{x} \Vert^2_{L^2}}{ \eta^2} \ d\eta \right)^{\frac{1}{2}} \ \left(\int_{0}^{\alpha}{ \eta^2} \ d\eta \right)^{\frac{1}{2}} \  d\alpha \nonumber\\
&\lesssim& \Vert f \Vert_{\dot H^3}  \Vert f \Vert_{\dot H^{\frac{3}{2}}} \int    \ \frac{\Vert s_\alpha f_x \Vert_{L^\infty} \Vert \delta_\alpha f_x \Vert_{L^\infty} }{\vert\alpha\vert^{\frac{5}{2}}}    \  d\alpha  \\
&\lesssim&\Vert f \Vert_{\dot H^3}  \Vert f \Vert_{\dot H^{\frac{3}{2}}} \left(\int    \ \frac{\Vert \delta_\alpha f_x \Vert^2_{L^\infty} }{\vert\alpha\vert^{\frac{5}{2}}}    \  d\alpha\right)^{\frac{1}{2}} \left(\int    \ \frac{\Vert s_\alpha f_x \Vert^2_{L^\infty}  }{\vert\alpha\vert^{\frac{5}{2}}}    \  d\alpha\right)^{\frac{1}{2}}\nonumber\\
&\lesssim&\Vert f \Vert_{\dot H^3}  \Vert f \Vert_{\dot H^{\frac{3}{2}}} \Vert f \Vert^2_{\dot B^{\frac{7}{4}}_{\infty,2}} \nonumber\\
&\lesssim&\Vert f \Vert_{\dot H^3}  \Vert f \Vert_{\dot H^{\frac{3}{2}}} \Vert f \Vert^2_{\dot H^{\frac{9}{4}}}\nonumber\\
&\lesssim&\Vert f \Vert^2_{\dot H^3}  \Vert f \Vert^2_{\dot H^{\frac{3}{2}}}\nonumber
\end{eqnarray}

$\bullet$ {{Estimate of $T_{4,1,7,2}$}} \\
\begin{eqnarray} \label{t4172}
T_{4,1,7,2}&\lesssim&\Vert f \Vert_{\dot H^3} \int  \frac{\Vert \delta_\alpha f_x \Vert_{L^{\infty}}}{\vert\alpha\vert^5} 
\left(\int_{0}^{\alpha} \frac{\Vert s_{\eta} f_{x} \Vert^2_{L^{4}}}{\vert \eta \vert^{\frac{5}{2}}} \ d\eta  \right)^{\frac{1}{2}} \left(\int_{0}^{\alpha} \frac{\Vert s_{\kappa} f_{x} \Vert^2_{L^{4}}}{\vert \kappa \vert^{\frac{5}{2}}} \ d\kappa  \right)^{\frac{1}{2}}  \nonumber\\
&& \times \left(\int_{0}^{\alpha} {\vert \eta \vert^{\frac{5}{2}}} \ d\eta  \right)^{\frac{1}{2}} \left(\int_{0}^{\alpha} {\vert \kappa \vert^{\frac{5}{2}}} \ d\kappa  \right)^{\frac{1}{2}} \ d\alpha \nonumber\\
&\lesssim&\Vert f \Vert_{\dot H^3} \Vert f \Vert^2_{\dot B^{\frac{7}{4}}_{4,2}} \int  \frac{\Vert \delta_\alpha f_x \Vert_{L^{\infty}}}{\vert\alpha\vert^{\frac{3}{2}}}  \ d\alpha \\
&\lesssim&\Vert f \Vert_{\dot H^3} \Vert f \Vert^2_{\dot H^{2}}\Vert f \Vert_{\dot B^{\frac{3}{2}}_{\infty,1}} \\
&\lesssim&\Vert f \Vert_{\dot H^3} \Vert f \Vert^2_{\dot H^{2}}\Vert f \Vert^{\frac{1}{2}}_{\dot B^{1}_{\infty,\infty}} \Vert f \Vert^{\frac{1}{2}}_{\dot B^{2}_{\infty,\infty}} \nonumber\\
&\lesssim&\Vert f \Vert_{\dot H^3} \Vert f \Vert^2_{\dot H^{2}}\Vert f \Vert^{\frac{1}{2}}_{\dot H^{\frac{3}{2}}} \Vert f \Vert^{\frac{1}{2}}_{\dot H^{\frac{5}{2}}} \nonumber
\end{eqnarray}
Then, since by interpolation $\Vert f \Vert_{\dot H^\frac{5}{2}} \lesssim \Vert f \Vert^{\frac{2}{3}}_{\dot H^3} \Vert f \Vert^{\frac{1}{3}}_{\dot H^{\frac{3}{2}}} $ then $\Vert f \Vert^{\frac{1}{2}}_{\dot H^\frac{5}{2}} \lesssim \Vert f \Vert^{\frac{1}{3}}_{\dot H^3} \Vert f \Vert^{\frac{1}{6}}_{\dot H^{\frac{3}{2}}} $ 
and $\Vert f \Vert_{\dot H^2} \lesssim \Vert f \Vert^{\frac{1}{3}}_{\dot H^3} \Vert f \Vert^{\frac{2}{3}}_{\dot H^{\frac{3}{2}}} $ therefore $\Vert f \Vert^2_{\dot H^2} \lesssim \Vert f \Vert^{\frac{2}{3}}_{\dot H^3} \Vert f \Vert^{\frac{4}{3}}_{\dot H^{\frac{3}{2}}} $ therefore the product gives $\Vert f \Vert^2_{\dot H^{2}}
 \Vert f \Vert^{\frac{1}{2}}_{\dot H^{\frac{5}{2}}} \lesssim  \Vert f \Vert_{\dot H^3}\Vert f \Vert^{\frac{3}{2}}_{\dot H^{\frac{3}{2}}}.$  Finally we get
\begin{eqnarray*}
T_{4,1,7,2}&\lesssim& \Vert f \Vert^2_{\dot H^3}\Vert f \Vert^{2}_{\dot H^{\frac{3}{2}}}
\end{eqnarray*}

$\bullet$ {Estimate of $T_{4,1,8}$}
Since
\begin{eqnarray*}
\partial_{\alpha}S_\alpha f &=& \bar\Delta_\alpha f_x-\Delta_\alpha f_x-\frac{s_\alpha f}{\alpha^2}
\end{eqnarray*}
We find that 
\begin{eqnarray*}
T_{4,1,8}&=&-\frac{1}{2\pi} \int \int\int_{0}^{\infty}\int_{0}^{\infty} \ \gamma e^{-\gamma-\sigma} \ \Lambda^{3} f \ \partial_{\alpha}S_\alpha f \sin(\frac{\gamma}{2}D_\alpha f)\sin(\frac{\gamma}{2}S_\alpha f)\frac{1}{\alpha^3} \int_{0}^{\alpha} s_{\eta} f_{x}    \ d\eta \\
&& \  \times\ {  \delta_\alpha f_x } \cos(\sigma \tau_{\alpha} f_x )\cos(\arctan(\tau_{\alpha} f_{x}))  \ d\gamma \ d\sigma  \ d\alpha\ dx  \\
&\lesssim&\Vert f \Vert_{\dot H^3} \int  \Vert \bar\delta_\alpha f_x \Vert_{L^\infty} \frac{1}{\vert\alpha\vert^4} \Vert \delta_\alpha f_x \Vert_{L^\infty} \int_{0}^{\alpha} \Vert s_{\eta} f_{x} \Vert_{L^2}    \ d\eta     \ d\alpha \\
&+&\Vert f \Vert_{\dot H^3} \int  \Vert \delta_\alpha f_x \Vert^2_{L^\infty} \frac{1}{\vert\alpha\vert^4} \int_{0}^{\alpha} \Vert s_{\eta} f_{x} \Vert_{L^2}    \ d\eta     \ d\alpha \\
&+&\Vert f \Vert_{\dot H^3} \int  \Vert s_\alpha f \Vert_{L^\infty}\frac{1}{\vert\alpha\vert^5} \Vert \delta_\alpha f_x \Vert_{L^\infty} \int_{0}^{\alpha} \Vert s_{\eta} f_{x} \Vert_{L^2}    \ d\eta     \ d\alpha \\
\end{eqnarray*}
It suffices to estimate one of the first two terms as they are of the same type.
\begin{eqnarray*}
T_{4,1,8} &\lesssim&\Vert f \Vert_{\dot H^3} \int  \Vert \delta_\alpha f_x \Vert^2_{L^\infty}\frac{1}{\vert\alpha\vert^4}  \left(\int_{0}^{\alpha} \frac{\Vert s_{\eta} f_{x} \Vert^2_{L^2}}{\eta^2}    \ d\eta\right)^{\frac{1}{2}} \left(\int_{0}^{\alpha} \eta^2    \ d\eta\right)^{\frac{1}{2}}     \ d\alpha \\
&+&\Vert f \Vert_{\dot H^3} \int  \Vert s_\alpha f \Vert_{L^\infty} \Vert \delta_\alpha f_x \Vert_{L^\infty} \frac{1}{\vert\alpha\vert^5} \left(\int_{0}^{\alpha} \frac{\Vert s_{\eta} f_{x} \Vert^2_{L^2}}{\eta^4}    \ d\eta\right)^{\frac{1}{2}} \left(\int_{0}^{\alpha} \eta^4    \ d\eta\right)^{\frac{1}{2}}     \ d\alpha  \\
&\lesssim&\Vert f \Vert_{\dot H^3}\left( \Vert f \Vert_{\dot H^{\frac{3}{2}}}  \int  \frac{\Vert \delta_\alpha f_x \Vert^2_{L^\infty}}{\vert\alpha\vert^{\frac{5}{2}}} \ d\alpha + \Vert f \Vert_{\dot H^{\frac{5}{2}}}\int  \frac{ \Vert s_\alpha f \Vert_{L^\infty}\Vert \delta_\alpha f_x \Vert_{L^\infty}}{\vert\alpha\vert^{\frac{5}{2}}} \ d\alpha  \right) \\
&\lesssim&\Vert f \Vert_{\dot H^3}\left(\Vert f \Vert_{\dot H^{\frac{3}{2}}} 
\ \Vert f \Vert^2_{\dot B^{\frac{7}{4}}_{\infty,2}} + \Vert f \Vert_{\dot H^{\frac{5}{2}}} \left( \int  \frac{ \Vert s_\alpha f \Vert^2_{L^\infty}}{\vert\alpha\vert^3} \ d\alpha \right)^{\frac{1}{2}} \left( \int  \frac{ \Vert \delta_\alpha f_x \Vert^2_{L^\infty}}{\vert\alpha\vert^2}\ d\alpha\right)^{\frac{1}{2}} \right)   \\
&\lesssim&\Vert f \Vert_{\dot H^3} \left(\Vert f \Vert_{\dot H^{\frac{3}{2}}} 
\Vert f \Vert^2_{\dot H^{\frac{9}{4}}} + \Vert f \Vert_{\dot H^{\frac{5}{2}}}\Vert f \Vert_{\dot B^{1}_{\infty,2}}\Vert f \Vert_{\dot H^{2}}\right)  \\
&\lesssim&\Vert f \Vert_{\dot H^3} \left(\Vert f \Vert_{\dot H^{\frac{3}{2}}} 
\Vert f \Vert^2_{\dot H^{\frac{9}{4}}} + \Vert f \Vert_{\dot H^{\frac{5}{2}}}\Vert f \Vert_{\dot B^{1}_{\infty,2}}\Vert f \Vert_{\dot H^{2}}\right)  \\
\end{eqnarray*}
Since $\Vert f \Vert^2_{\dot H^{\frac{9}{4}}} \leq \Vert f \Vert_{\dot H^{\frac{3}{2}}}\Vert f \Vert_{\dot H^3}$ and $\Vert f \Vert_{\dot H^\frac{5}{2}} \Vert f \Vert_{\dot H^2} \lesssim \Vert f \Vert^{\frac{1}{3}}_{\dot H^3} \Vert f \Vert^{\frac{2}{3}}_{\dot H^{\frac{3}{2}}} \Vert f \Vert^{\frac{2}{3}}_{\dot H^3} \Vert f \Vert^{\frac{1}{3}}_{\dot H^{\frac{3}{2}}} $, we conclude that
\begin{eqnarray*}
T_{4,1,8} &\lesssim& \Vert f \Vert^2_{\dot H^3} \Vert f \Vert^2_{\dot H^{\frac{3}{2}}}
\end{eqnarray*}

$\bullet$ {Estimate of $T_{4,1,9}$} \\

Using the fact that 

\begin{eqnarray*}
\partial_{\alpha}D_\alpha f &=&  \frac{s_\alpha f_x}{\alpha}-\frac{1}{\alpha^2} \ {\int_0^\alpha s_\eta f_x \ d\eta}
\end{eqnarray*}

One finds

\begin{eqnarray*}
T_{4,1,9}&=&\frac{1}{\pi} \int \int\int_{0}^{\infty}\int_{0}^{\infty} \ e^{-\gamma-\sigma} \ \Lambda^{3} f \ \partial_\alpha D_\alpha f \sin(\frac{\gamma}{2}D_\alpha f)\cos(\frac{\gamma}{2}S_\alpha f)\frac{1}{\alpha^2}  \\
&& \  \times\ {  \delta_\alpha f_x } \cos(\sigma \tau_{\alpha} f_x )\cos(\arctan(\tau_{\alpha} f_{x}))   \ d\gamma \ d\sigma   \ d\alpha \ dx\\
&\lesssim& \Vert f \Vert_{\dot H^3} \left( \int \frac{\Vert s_\alpha f_x \Vert_{L^{2}} \Vert \delta_\alpha f_x \Vert_{L^{\infty}} }{\vert\alpha\vert^3} \ d\alpha \right. \\
&&\left. + \int \frac{\Vert \delta_\alpha f_x \Vert_{L^{\infty}} }{\vert\alpha\vert^4} \left(\int_0^\alpha \frac{\Vert s_\eta f_x \Vert^2_{L^2}}{\vert \eta \vert^4} \ d\eta\right)^{\frac{1}{2}} \left(\int_0^\alpha \vert \eta \vert^4 \ d\eta\right)^{\frac{1}{2}} \ d\alpha \right) 
\end{eqnarray*}
It is important to note that the term $s_\alpha f_x$ is important to be able to control the first integral as it is too singular in $\alpha$. We cannot balance  the regularity in $\alpha$ equally in  each term of the product $\Vert s_\alpha f_x \Vert_{L^{2}} \Vert \delta_\alpha f_x \Vert_{L^{\infty}}$ as the second term would not be well-defined. We write,
\begin{eqnarray*}
T_{4,1,9}&\lesssim& \Vert f \Vert_{\dot H^3} \left( \left(\int \frac{\Vert s_\alpha f_x \Vert^2_{L^{2}}  }{\vert\alpha\vert^4} \ d\alpha\right)^{\frac{1}{2}} \left(\int \frac{\Vert \delta_\alpha f_x \Vert^2_{L^{\infty}}  }{\vert\alpha\vert^2} \ d\alpha\right)^{\frac{1}{2}}  + \Vert f \Vert_{\dot H^2}\int \frac{\Vert \delta_\alpha f_x \Vert_{L^{\infty}} }{\vert\alpha\vert^\frac{3}{2}}  \ d\alpha\right) \\
&\lesssim& \Vert f \Vert_{\dot H^3} \left( \Vert f \Vert_{\dot H^{\frac{5}{2}}}
\Vert f \Vert_{\dot B^{\frac{3}{2}}_{\infty,2}} + \Vert f \Vert_{\dot H^\frac{5}{2}}\Vert f \Vert_{\dot B^{\frac{3}{2}}_{\infty,1}}   \right) \\
&\lesssim& \Vert f \Vert_{\dot H^3}\left(\Vert f \Vert_{\dot H^{\frac{5}{2}}}\Vert f \Vert_{\dot H^2} + \Vert f \Vert_{\dot H^\frac{5}{2}}\Vert f \Vert^{\frac{1}{2}}_{\dot B^{1}_{\infty,\infty}} \Vert f \Vert^{\frac{1}{2}}_{\dot B^{2}_{\infty,\infty}}  \right) \\
&\lesssim& \Vert f \Vert_{\dot H^3}\left(\Vert f \Vert_{\dot H^{\frac{5}{2}}}\Vert f \Vert_{\dot H^2} + \Vert f \Vert^{\frac{3}{2}}_{\dot H^\frac{5}{2}}\Vert f \Vert^{\frac{1}{2}}_{\dot H^{\frac{3}{2}}}   \right)
\end{eqnarray*}
As we have $\Vert f \Vert^{\frac{3}{2}}_{\dot H^{\frac{5}{2}}}\Vert f \Vert^{\frac{1}{2}}_{\dot H^{\frac{3}{2}}}\lesssim \Vert f \Vert_{\dot H^{\frac{3}{2}}} \Vert f \Vert_{\dot H^3}$ and that $\Vert f \Vert_{\dot H^{\frac{5}{2}}}\Vert f \Vert_{\dot H^2}\lesssim \Vert f \Vert_{\dot H^{\frac{3}{2}}} \Vert f \Vert_{\dot H^3}$, we find
\begin{eqnarray*}
T_{4,1,9}&\lesssim& \Vert f \Vert^2_{\dot H^3}\Vert f \Vert_{\dot H^{\frac{3}{2}}}
\end{eqnarray*}

$\bullet$ {Estimate of $T_{4,1,10}$}

\begin{eqnarray*}
T_{4,1,10}&=&\frac{1}{4\pi} \int \int\int_{0}^{\infty}\int_{0}^{\infty} \ \gamma^2 e^{-\gamma-\sigma} \ \Lambda^{3} f \ \partial_{\alpha} D_\alpha f \ \partial_{\alpha} S_\alpha f \cos(\frac{\gamma}{2}D_\alpha f)\sin(\frac{\gamma}{2}S_\alpha f) \\
&& \  \times \frac{1}{\alpha^2} \left(\int_{0}^{\alpha} s_{\eta} f_{x}    \ d\eta\right) {  \delta_\alpha f_x } \cos(\sigma \tau_{\alpha} f_x )\cos(\arctan(\tau_{\alpha} f_{x}))   \ d\gamma \ d\sigma  \ d\alpha \ dx. 
\end{eqnarray*}
Using the fact that, 

\begin{eqnarray*}
\partial_{\alpha}D_\alpha f \ \partial_{\alpha} S_\alpha f  &=& - \left( \frac{s_\alpha f_x}{\alpha}+\frac{1}{\alpha^2} \ {\int_0^\alpha s_\kappa f_x \ d\kappa}\right) \left(\bar\Delta_\alpha f_x-\Delta_\alpha f_x-\frac{s_\alpha f}{\alpha^2}\right)
\end{eqnarray*}
As $\bar\Delta_\alpha f_x$ and $\Delta_\alpha f_x$ have exactly the same regularity and structure and they will both lead to equivalent semi-norms, it is therefore enough to consider 
\begin{eqnarray} \label{prod}
\partial_{\alpha}D_\alpha f \ \partial_{\alpha} S_\alpha f  &\approx&  \left( \frac{s_\alpha f_x}{\alpha}+\frac{1}{\alpha^2} \ {\int_0^\alpha s_\kappa f_x \ d\kappa}\right) \left(\Delta_\alpha f_x+\frac{s_\alpha f}{\alpha^2}\right)
\end{eqnarray}
Therefore,
\begin{eqnarray*}
T_{4,1,10}&\lesssim& \Vert f \Vert_{\dot H^3}\int   \ 
\frac{\Vert \delta_\alpha f_x \Vert^2_{L^{\infty}} \Vert s_\alpha f_x \Vert_{L^{\infty}}}{\alpha^4}
\left( \int_{0}^{\alpha} \frac{\Vert s_{\eta} f_{x}\Vert^2_{L^2}}{\eta^2}\ d\eta\right)^{\frac{1}{2}} \left( \int_{0}^{\alpha} \eta^2 \ d\eta\right)^{\frac{1}{2}}   \ d\alpha \\
&+&\Vert f \Vert_{\dot H^3}\int   \ 
\frac{\Vert \delta_\alpha f_x \Vert_{L^{\infty}} \Vert s_\alpha f \Vert_{L^{\infty}} \Vert s_\alpha f_x \Vert_{L^{\infty}}}{\vert\alpha\vert^5}
\left( \int_{0}^{\alpha} \frac{\Vert s_{\eta} f_{x}\Vert^2_{L^2}}{\eta^2}\ d\eta\right)^{\frac{1}{2}} \left( \int_{0}^{\alpha} \eta^2 \ d\eta\right)^{\frac{1}{2}}   \ d\alpha \\
&+&\Vert f \Vert_{\dot H^3}\int   \ 
\frac{\Vert \delta_\alpha f_x \Vert^2_{L^{\infty}}}{\vert\alpha\vert^5}
\left( \int_{0}^{\alpha} \frac{\Vert s_{\eta} f_{x}\Vert^2_{L^2}}{\eta^2}\ d\eta\right)^{\frac{1}{2}}  \left({\int_0^\alpha \frac{\Vert s_\kappa f_x \Vert_{L^\infty}}{\vert \kappa \vert^2} \ d\kappa} \right)^{\frac{1}{2}}  \\&& \times \left( \int_{0}^{\alpha} \eta^2 \ d\eta\right)^{\frac{1}{2}} \left( \int_{0}^{\alpha} \kappa^2 \ d\kappa\right)^{\frac{1}{2}}   \ d\alpha \\
&+&\Vert f \Vert_{\dot H^3}\int   \ 
\frac{\Vert \delta_\alpha f_x \Vert_{L^{\infty}} \Vert s_\alpha f \Vert_{L^{\infty}}}{\vert\alpha\vert^6}
\left( \int_{0}^{\alpha} \frac{\Vert s_{\eta} f_{x}\Vert^2_{L^2}}{\eta^2}\ d\eta\right)^{\frac{1}{2}}  \left({\int_0^\alpha \frac{\Vert s_\kappa f_x \Vert_{L^\infty}}{\vert \kappa \vert^2} \ d\kappa} \right)^{\frac{1}{2}}  \\&& \times \left( \int_{0}^{\alpha} \eta^2 \ d\eta\right)^{\frac{1}{2}} \left( \int_{0}^{\alpha} \kappa^2 \ d\kappa\right)^{\frac{1}{2}}   \ d\alpha \\
\end{eqnarray*}
Then
\begin{eqnarray*}
T_{4,1,10}&\lesssim& \Vert f \Vert_{\dot H^3} \Vert f \Vert_{\dot H^\frac{3}{2}} \int   \ 
\frac{\Vert \delta_\alpha f_x \Vert^2_{L^{\infty}} \Vert s_\alpha f_x \Vert_{L^{\infty}}}{\vert\alpha\vert^{\frac{5}{2}}}   \ d\alpha \\
&+&\Vert f \Vert_{\dot H^3} \Vert f \Vert_{\dot H^\frac{3}{2}}\int   \ 
\frac{\Vert \delta_\alpha f_x \Vert_{L^{\infty}} \Vert s_\alpha f \Vert_{L^{\infty}} \Vert s_\alpha f_x \Vert_{L^{\infty}}}{\vert\alpha\vert^\frac{7}{2}}
   \ d\alpha \\
&+&\Vert f \Vert_{\dot H^3} \Vert f \Vert_{\dot H^\frac{3}{2}} \Vert f \Vert_{\dot H^2}\int   \ 
\frac{\Vert \delta_\alpha f_x \Vert^2_{L^{\infty}}}{\vert\alpha\vert^2}  \ d\alpha \\
&+&\Vert f \Vert_{\dot H^3} \Vert f \Vert_{\dot H^\frac{3}{2}}  \Vert f \Vert_{\dot H^{2}}\int   \ 
\frac{\Vert \delta_\alpha f_x \Vert_{L^{\infty}} \Vert s_\alpha f \Vert_{L^{\infty}}}{\vert\alpha\vert^3}    \ d\alpha. \\
\end{eqnarray*}
We find,
\begin{eqnarray*}
T_{4,1,10}&\lesssim& \Vert f \Vert_{\dot H^3} \Vert f \Vert_{\dot H^\frac{3}{2}} \left(\sup_{\alpha \in \mathbb R} \frac{\Vert s_\alpha f_x \Vert_{L^{\infty}} }{\vert\alpha\vert}\right) \int   \ 
\frac{\Vert \delta_\alpha f_x \Vert^2_{L^{\infty}}}{\vert \alpha \vert^{\frac{3}{2}}} \ d\alpha   \\
&&\ + \ \Vert f \Vert_{\dot H^3} \Vert f \Vert_{\dot H^\frac{3}{2}} \left(\sup_{\alpha \in \mathbb R} \frac{\Vert s_\alpha f \Vert_{L^{\infty}} }{\vert\alpha\vert}\right)\int   \ 
\frac{\Vert \delta_\alpha f_x \Vert_{L^{\infty}} 
 \Vert s_\alpha f_x \Vert_{L^{\infty}}}{\vert\alpha\vert^\frac{5}{2}}
   \ d\alpha \\
   &&\ + \ \Vert f \Vert_{\dot H^3} \Vert f \Vert_{\dot H^\frac{3}{2}} \Vert f \Vert^3_{\dot H^2} \\
   &&\ + \ \Vert f \Vert_{\dot H^3} \Vert f \Vert_{\dot H^\frac{3}{2}}  \Vert f \Vert_{\dot H^2}\int   \ 
\frac{\Vert \delta_\alpha f_x \Vert_{L^{\infty}} \Vert s_\alpha f \Vert_{L^{\infty}}}{\vert\alpha\vert^3}    \ d\alpha \\
&\lesssim& \Vert f \Vert_{\dot H^3} \Vert f \Vert_{\dot H^\frac{3}{2}} 
\Vert f \Vert_{\dot B^2_{\infty,\infty}} \Vert f \Vert^2_{\dot B^{\frac{5}{4}}_{\infty,2}}    \\
&&\ + \ \Vert f \Vert_{\dot H^3} \Vert f \Vert_{\dot H^\frac{3}{2}}\Vert f \Vert_{\dot B^1_{\infty,\infty}} \left( \int    \frac{\Vert \delta_\alpha f_x \Vert^2_{L^{\infty}} 
 }{\vert\alpha\vert^2} \ d\alpha \right)^{\frac{1}{2}} \left( \int    \frac{\Vert s_\alpha f_x \Vert^2_{L^{\infty}} 
 }{\vert\alpha\vert^3} \ d\alpha \right)^{\frac{1}{2}}\\
   &&\ + \ \Vert f \Vert_{\dot H^3} \Vert f \Vert_{\dot H^\frac{3}{2}} \Vert f \Vert^3_{\dot H^2} \\
   &&\ + \ \Vert f \Vert_{\dot H^3} \Vert f \Vert_{\dot H^\frac{3}{2}}  \Vert f \Vert_{\dot H^2} \left(\int \frac{\Vert \delta_\alpha f_x \Vert^2_{L^{\infty}}}{\vert \alpha\vert^2}    \ d\alpha \right)^{\frac{1}{2}}  \left(\int \frac{\Vert s_\alpha f \Vert^2_{L^{\infty}}}{\vert \alpha\vert^4}    \ d\alpha \right)^{\frac{1}{2}} \\
  \end{eqnarray*}
  Therefore,
\begin{eqnarray*}
   T_{4,1,10}&\lesssim& \Vert f \Vert_{\dot H^3} \Vert f \Vert_{\dot H^\frac{3}{2}} 
\Vert f \Vert_{\dot B^2_{\infty,\infty}} \Vert f \Vert^2_{\dot B^{\frac{5}{4}}_{\infty,2}} + \Vert f \Vert_{\dot H^3} \Vert f \Vert_{\dot H^\frac{3}{2}}\Vert f \Vert_{\dot B^1_{\infty,\infty}} \Vert f \Vert_{\dot B^{\frac{3}{2}}_{\infty,2}} \Vert f \Vert_{\dot B^{2}_{\infty,2}}\\
   && + \ \Vert f \Vert_{\dot H^3} \Vert f \Vert_{\dot H^\frac{3}{2}} \Vert f \Vert^3_{\dot H^2}  + \Vert f \Vert_{\dot H^3} \Vert f \Vert_{\dot H^\frac{3}{2}}  \Vert f \Vert_{\dot H^2} \Vert f \Vert^2_{\dot B^{\frac{3}{2}}_{\infty,2}}   \\
   &\lesssim& \Vert f \Vert_{\dot H^3} \Vert f \Vert_{\dot H^\frac{3}{2}} \Vert f \Vert_{\dot H^\frac{5}{2}} 
 \Vert f \Vert^2_{\dot H^{\frac{7}{4}}} + \Vert f \Vert_{\dot H^3} \Vert f \Vert^2_{\dot H^\frac{3}{2}} \Vert f \Vert_{\dot H^2} \Vert f \Vert_{\dot H^{\frac{5}{2}}}\\
   && \ + \ \Vert f \Vert_{\dot H^3} \Vert f \Vert_{\dot H^\frac{3}{2}} \Vert f \Vert^3_{\dot H^2} 
\end{eqnarray*}
Since  $\Vert f \Vert_{\dot H^\frac{5}{2}} \lesssim \Vert f \Vert^{\frac{2}{3}}_{\dot H^3} \Vert f \Vert^{\frac{1}{3}}_{\dot H^{\frac{3}{2}}} $ and
$\Vert f \Vert^2_{\dot H^\frac{7}{4}} \lesssim \Vert f \Vert^{\frac{1}{3}}_{\dot H^3} \Vert f \Vert^{\frac{2}{3}}_{\dot H^{\frac{3}{2}}} $ then we have that $$\Vert f \Vert_{\dot H^\frac{5}{2}} \Vert f \Vert^2_{\dot H^{\frac{7}{4}}} \lesssim \Vert f \Vert_{\dot H^3} \Vert f \Vert_{\dot H^{\frac{3}{2}}} $$ One may also check that $\Vert f \Vert^3_{\dot H^2} \lesssim \Vert f \Vert_{\dot H^3} \Vert f \Vert_{\dot H^{\frac{3}{2}}} $ and that $ \Vert f \Vert_{\dot H^2} \Vert f \Vert_{\dot H^{\frac{5}{2}}}\lesssim \Vert f \Vert_{\dot H^3} \Vert f \Vert_{\dot H^{\frac{3}{2}}} $. Finally, we have obtained that
\begin{eqnarray*}
   T_{4,1,10}&\lesssim&\Vert f \Vert^2_{\dot H^3}\left(\Vert f \Vert^2_{\dot H^{\frac{3}{2}}}+\Vert f \Vert_{\dot H^{\frac{3}{2}}} \right)
\end{eqnarray*}

$\bullet$ {Estimate of $T_{4,1,11}$}\\

We have 

\begin{eqnarray*}
T_{4,1,11}&=&-\frac{1}{2\pi} \int \int\int_{0}^{\infty}\int_{0}^{\infty} \ \gamma e^{-\gamma-\sigma} \ \Lambda^{3} f \ \frac{1}{\alpha} \ (\partial_{\alpha} D_\alpha f)^2 \  \cos(\frac{\gamma}{2}D_\alpha f)\cos(\frac{\gamma}{2}S_\alpha f)\\
&& \  \times\ {  \delta_\alpha f_x } \cos(\sigma \tau_{\alpha} f_x )\cos(\arctan(\tau_{\alpha} f_{x}))  \ d\gamma \ d\sigma   \ d\alpha\ dx \\
\end{eqnarray*}
Since,
\begin{eqnarray*}
\left(\partial_{\alpha}D_\alpha f \right)^2  &=&  \left( \frac{s_\alpha f_x}{\alpha}+\frac{1}{\alpha^2} \ {\int_0^\alpha s_\kappa f_x \ d\kappa}\right)^2 
\end{eqnarray*}
Then, we obtain
\begin{eqnarray*}
T_{4,1,11}&=&-\frac{1}{2\pi} \int \int\int_{0}^{\infty}\int_{0}^{\infty} \ \gamma e^{-\gamma-\sigma} \ \Lambda^{3} f  \ \frac{s^2_{\alpha} f_x}{\alpha^3} \  \cos(\frac{\gamma}{2}D_\alpha f)\cos(\frac{\gamma}{2}S_\alpha f)\\
&& \  \times\ {  \delta_\alpha f_x } \cos(\sigma \tau_{\alpha} f_x )\cos(\arctan(\tau_{\alpha} f_{x}))  \ d\gamma \ d\sigma   \ d\alpha \ dx\\
&-&\frac{1}{\pi} \int \int\int_{0}^{\infty}\int_{0}^{\infty} \ \gamma e^{-\gamma-\sigma} \ \Lambda^{3} f  \ \frac{s_{\alpha} f_x}{\alpha^4}  \  \left({\int_0^\alpha s_\kappa f_x \ d\kappa}\right)\cos(\frac{\gamma}{2}D_\alpha f)\cos(\frac{\gamma}{2}S_\alpha f)\\
&& \  \times\ {  \delta_\alpha f_x } \cos(\sigma \tau_{\alpha} f_x )\cos(\arctan(\tau_{\alpha} f_{x}))  \ d\gamma \ d\sigma  \ d\alpha\ dx  \\
&-&\frac{1}{2\pi} \int \int\int_{0}^{\infty}\int_{0}^{\infty} \ \gamma e^{-\gamma-\sigma} \ \Lambda^{3} f  \ \frac{1}{\alpha^5}  \  \left({\int_0^\alpha s_\kappa f_x \ d\kappa}\right)^2 \cos(\frac{\gamma}{2}D_\alpha f)\cos(\frac{\gamma}{2}S_\alpha f)\\
&& \  \times\ {  \delta_\alpha f_x } \cos(\sigma \tau_{\alpha} f_x )\cos(\arctan(\tau_{\alpha} f_{x}))  \ d\gamma \ d\sigma   \ d\alpha \ dx\\
&\lesssim& \Vert f \Vert_{\dot H^3}\left(\sup_{\alpha \in \mathbb R} \frac{\Vert s_\alpha f_x \Vert_{L^{\infty}}}{\vert\alpha\vert}\right) \int \frac{\Vert s_\alpha f_x  \Vert_{L^\infty} \Vert \delta_\alpha f_x  \Vert_{L^2}}{\vert \alpha \vert^2} \ d\alpha \\
&+&\Vert f \Vert_{\dot H^3} \int \frac{\Vert s_\alpha f_x  \Vert_{L^\infty} \Vert \delta_\alpha f_x  \Vert_{L^\infty}}{\vert \alpha \vert^4} \left(\int_0^\alpha \frac{\Vert s_\kappa f_x \Vert^2_{L^2}}{\vert\kappa\vert^2} \ d\kappa\right)^{\frac{1}{2}} \left(\int_0^\alpha {\vert\kappa\vert^2 }\ d\kappa\right)^{\frac{1}{2}} \ d\alpha \\
&+&\Vert f \Vert_{\dot H^3} \int \frac{\Vert \delta_\alpha f_x  \Vert_{L^\infty}}{\vert \alpha \vert^5} \left(\int_0^\alpha \frac{\Vert s_\kappa f_x \Vert^2_{L^4}}{\vert\kappa\vert^\frac{5}{2} \ }d\kappa\right) \left(\int_0^\alpha {\vert\kappa\vert^\frac{5}{2} }\ d\kappa\right) \ d\alpha \\
&\lesssim& \Vert f \Vert_{\dot H^3} \Vert f \Vert_{\dot B^2_{\infty,\infty}} \left(\int \frac{\Vert s_\alpha f_x  \Vert^2_{L^\infty}}{\vert \alpha \vert^2} \ d\alpha \right)^{\frac{1}{2}} \left(\int \frac{\Vert \delta_\alpha f_x  \Vert^2_{L^2}}{\vert \alpha \vert^2} \ d\alpha \right)^{\frac{1}{2}} \\
&& \ + \ \Vert f \Vert_{\dot H^3} \Vert f \Vert_{\dot H^{\frac{3}{2}}} \int \frac{\Vert s_\alpha f_x  \Vert_{L^\infty} \Vert \delta_\alpha f_x  \Vert_{L^\infty}}{\vert \alpha \vert^{\frac{5}{2}}}  \ d\alpha \\
&& \ + \ \Vert f \Vert_{\dot H^3} \Vert f \Vert^2_{\dot B^{\frac{7}{4}}_{4,2}} \int \frac{\Vert \delta_\alpha f_x  \Vert_{L^\infty}}{\vert \alpha \vert^\frac{3}{2}}  d\alpha \\
&\lesssim& \Vert f \Vert_{\dot H^3} \Vert f \Vert_{\dot H^\frac{5}{2}}\Vert f \Vert_{\dot H^2}\Vert f \Vert_{\dot H^{\frac{3}{2}}} +\Vert f \Vert^2_{\dot H^3}  \Vert f \Vert^2_{\dot H^{\frac{3}{2}}} +\Vert f \Vert_{\dot H^3} \Vert f \Vert^2_{\dot H^{2}}  \Vert f \Vert_{\dot B^{\frac{3}{2}}_{\infty,1}} \\
&\lesssim& \Vert f \Vert^2_{\dot H^3} \Vert f \Vert^2_{\dot H^\frac{3}{2}} + \Vert f \Vert_{\dot H^3} \Vert f \Vert^2_{\dot B^{\frac{7}{4}}_{4,2}} \int \frac{\Vert \delta_\alpha f_x  \Vert_{L^\infty}}{\vert \alpha \vert^\frac{3}{2}}  d\alpha, \\
\end{eqnarray*}
where to control the first term we used  that $\Vert f \Vert_{\dot H^\frac{5}{2}}\Vert f \Vert_{\dot H^2} \lesssim \Vert f \Vert_{\dot H^3}  \Vert f \Vert_{\dot H^{\frac{3}{2}}}$. The second term  can be directly controlled by $\Vert f \Vert^2_{\dot H^3}  \Vert f \Vert^2_{\dot H^{\frac{3}{2}}}$  as it is similar to $T_{4,1,7,1}$ (see \eqref{T_{4,1,7,1}}). For the third term, it is similar to the term  $T_{4,1,7,2}$ (see \eqref{t4172}). Finally, we get 

\begin{eqnarray*}
T_{4,1,11} \lesssim \Vert f \Vert^2_{\dot H^3}  \Vert f \Vert^2_{\dot H^{\frac{3}{2}}}
\end{eqnarray*}

$\bullet$ {Estimate of $T_{4,1,12}$}\\

It remains to estimate the following term 
\begin{eqnarray*}
T_{4,1,12}&=&\frac{1}{2\pi} \int \int\int_{0}^{\infty}\int_{0}^{\infty} \ \gamma e^{-\gamma-\sigma} \ \Lambda^{3} f \ \frac{1}{\alpha}\ \partial_{\alpha} D_\alpha f \ \partial_{\alpha} S_\alpha f \sin(\frac{\gamma}{2}D_\alpha f)\sin(\frac{\gamma}{2}S_\alpha f) \\
&& \  \times\ {  \delta_\alpha f_x } \cos(\sigma \tau_{\alpha} f_x )\cos(\arctan(\tau_{\alpha} f_{x}))   \ d\gamma \ d\sigma \ d\alpha \ dx. \\
\end{eqnarray*}

Using the fact the observation made in \eqref{prod}, that is 
\begin{eqnarray*} 
\partial_{\alpha}D_\alpha f \ \partial_{\alpha} S_\alpha f  &\approx&  \left( \frac{s_\alpha f_x}{\alpha}+\frac{1}{\alpha^2} \ {\int_0^\alpha s_\kappa f_x \ d\kappa}\right) 
\left(\Delta_\alpha f_x+\frac{s_\alpha f}{\alpha^2}\right)
\end{eqnarray*}
we have the following 4 terms to estimate

\begin{eqnarray*}
T_{4,1,12}&=&\frac{1}{2\pi} \int \int\int_{0}^{\infty}\int_{0}^{\infty} \ \gamma e^{-\gamma-\sigma} \ \Lambda^{3} f \ \
 \frac{s_\alpha f_x}{\alpha^2}\Delta_\alpha f_x
\sin(\frac{\gamma}{2}D_\alpha f)\sin(\frac{\gamma}{2}S_\alpha f) \\
&& \  \times\ {  \delta_\alpha f_x } \cos(\sigma \tau_{\alpha} f_x )\cos(\arctan(\tau_{\alpha} f_{x}))   \ d\gamma \ d\sigma  \ d\alpha\ dx \\
&+&\frac{1}{2\pi} \int \int\int_{0}^{\infty}\int_{0}^{\infty} \ \gamma e^{-\gamma-\sigma} \ \Lambda^{3} f \ \
 \frac{s_\alpha f_x \ s_\alpha f}{\alpha^4}
\sin(\frac{\gamma}{2}D_\alpha f)\sin(\frac{\gamma}{2}S_\alpha f) \\
&& \  \times\ {  \delta_\alpha f_x } \cos(\sigma \tau_{\alpha} f_x )\cos(\arctan(\tau_{\alpha} f_{x}))   \ d\gamma \ d\sigma  \ d\alpha\ dx \\
&+&\frac{1}{2\pi} \int \int\int_{0}^{\infty}\int_{0}^{\infty} \ \int_0^\alpha \gamma e^{-\gamma-\sigma} \ \Lambda^{3} f \
 \frac{1}{\alpha^3} \Delta_\alpha f_x \  s_\kappa f_x  
\sin(\frac{\gamma}{2}D_\alpha f)\sin(\frac{\gamma}{2}S_\alpha f) \\
&& \  \times\ {  \delta_\alpha f_x } \cos(\sigma \tau_{\alpha} f_x )\cos(\arctan(\tau_{\alpha} f_{x})) \ d\kappa  \ d\gamma \ d\sigma  \ d\alpha \ dx\\
&+&\frac{1}{2\pi} \int \int\int_{0}^{\infty}\int_{0}^{\infty} \int_0^\alpha\ \gamma e^{-\gamma-\sigma} \ \Lambda^{3} f  \ \frac{s_\alpha f}{\alpha^5} \ { s_\kappa f_x \ } 
\sin(\frac{\gamma}{2}D_\alpha f)\sin(\frac{\gamma}{2}S_\alpha f) \\
&& \  \times\ {  \delta_\alpha f_x } \cos(\sigma \tau_{\alpha} f_x )\cos(\arctan(\tau_{\alpha} f_{x})) \ d\kappa \ d\gamma \ d\sigma  \ d\alpha \ dx\\
&\lesssim& \Vert f \Vert_{\dot H^3}\int \frac{\Vert s_{\alpha} f_x \Vert_{L^{\infty}}\Vert \delta_{\alpha} f_x \Vert^2_{L^{4}}}{\vert \alpha \vert^3} \ d\alpha +  \Vert f \Vert_{\dot H^3}\int \frac{\Vert s_{\alpha} f_x \Vert_{L^{\infty}} \Vert s_{\alpha} f \Vert_{L^{\infty}}\Vert \delta_{\alpha} f_x \Vert_{L^{2}}}{\vert \alpha \vert^4} \ d\alpha \\
&\lesssim& \Vert f \Vert_{\dot H^3}\left(\sup_{\alpha \in \mathbb R}\frac{\Vert s_{\alpha} f_x \Vert_{L^{\infty}}}{\vert \alpha \vert}\right) \left(\int \frac{\Vert \delta_{\alpha} f_x \Vert_{L^{4}}}{\vert \alpha \vert^2} \ d\alpha + \int \frac{ \Vert s_{\alpha} f \Vert_{L^{\infty}}\Vert \delta_{\alpha} f_x \Vert_{L^{2}}}{\vert \alpha \vert^3} \ d\alpha \right)\\
&& +\ \Vert f \Vert_{\dot H^3}  \int    \
 \frac{\Vert \delta_\alpha f_x \Vert^2_{L^{\infty}} }{\vert\alpha\vert^4}  \left(  \int_0^\alpha \frac{\Vert s_\kappa f_x \Vert^2_{L^2}}{\vert \kappa \vert^2} \ d\kappa \right)^{\frac{1}{2}} \left(  \int_0^\alpha \vert \kappa \vert^2 \ d\kappa \right)^{\frac{1}{2}}   \ d\alpha \\
 && +\ \Vert f \Vert_{\dot H^3}  \int    \
 \frac{\Vert \delta_\alpha f_x \Vert_{L^{\infty}} \Vert s_\alpha f \Vert_{L^{\infty}} }{\vert\alpha\vert^5}  \left(  \int_0^\alpha \frac{\Vert s_\kappa f_x \Vert^2_{L^2}}{\vert \kappa \vert^2} \ d\kappa \right)^{\frac{1}{2}} \left(  \int_0^\alpha \vert \kappa \vert^2 \ d\kappa \right)^{\frac{1}{2}}   \ d\alpha \\
&\lesssim& \Vert f \Vert_{\dot H^3}\Vert f \Vert_{\dot B^{2}_{\infty,\infty}} \left(\Vert f \Vert^2_{\dot B^{\frac{3}{2}}_{4,2}} + \left(\int \frac{\Vert s_\alpha f \Vert^2_{L^\infty}}{\vert \alpha \vert^4} \ d\alpha\right)^{\frac{1}{2}} \left(\int \frac{\Vert \delta_\alpha f_x \Vert^2_{L^2}}{\vert \alpha \vert^2} \ d\alpha\right)^{\frac{1}{2}}  \right) \\
&& +\ \Vert f \Vert_{\dot H^3} \Vert f \Vert_{\dot H^{\frac{3}{2}}} \left(  \int    \
 \frac{\Vert \delta_\alpha f_x \Vert^2_{L^{\infty}} }{\vert\alpha\vert^{\frac{5}{2}}} \ d\alpha   + \int    \
 \frac{\Vert \delta_\alpha f_x \Vert_{L^{\infty}} \Vert s_\alpha f \Vert_{L^{\infty}} }{\vert\alpha\vert^{\frac{7}{2}}} \ d\alpha \right) \\
\end{eqnarray*}
Then, in order to carefully share the derivatives in the last term, we split the weight  $\frac{1}{\vert \alpha\vert^{\frac{7}{2}} }$ into $\frac{1}{\vert \alpha\vert^{\frac{5}{4}} }\frac{1}{\vert \alpha\vert^{\frac{9}{4}} }$, we obtain

\begin{eqnarray*}
T_{4,1,12}&\lesssim& \Vert f \Vert_{\dot H^3}\Vert f \Vert_{\dot B^{2}_{\infty,\infty}} \left(\Vert f \Vert^2_{\dot H^{\frac{7}{4}}} + \Vert f \Vert_{\dot B^{\frac{3}{2}}_{\infty,2}} \Vert f \Vert_{\dot H^{\frac{3}{2}}}   \right) \\
&& +\ \Vert f \Vert_{\dot H^3} \Vert f \Vert_{\dot H^{\frac{3}{2}}} \left( \Vert f \Vert^2_{\dot B^{\frac{7}{4}}_{\infty,2}}    + \int    \
 \frac{\Vert \delta_\alpha f_x \Vert_{L^{\infty}} \Vert s_\alpha f \Vert_{L^{\infty}} }{{\vert \alpha\vert^{\frac{5}{4}} }{\vert \alpha\vert^{\frac{9}{4}} }} \ d\alpha \right) \\
 &\lesssim& \Vert f \Vert_{\dot H^3}\Vert f \Vert_{\dot H^{\frac{5}{2}}} \left(\Vert f \Vert_{\dot H^{\frac{3}{2}}} \Vert f \Vert_{\dot H^{2}} + \Vert f \Vert_{\dot H^{2}} \Vert f \Vert_{\dot H^{\frac{3}{2}}}   \right) \\
&& +\ \Vert f \Vert_{\dot H^3} \Vert f \Vert_{\dot H^{\frac{3}{2}}} \left( \Vert f \Vert^2_{\dot H^{\frac{9}{4}}}    + \left(\int    \
 \frac{\Vert \delta_\alpha f_x \Vert^2_{L^{\infty}}  }{{\vert \alpha\vert^{\frac{5}{2}} } } \ d\alpha\right)^{\frac{1}{2}} \left(\int    \
 \frac{ \Vert s_\alpha f \Vert^2_{L^{\infty}} }{{\vert \alpha\vert^{\frac{9}{2}} }} \ d\alpha\right)^{\frac{1}{2}} \right) \\
\end{eqnarray*}
Using the fact that $\Vert f \Vert^2_{\dot H^{\frac{9}{4}}} \lesssim  \Vert f \Vert_{\dot H^3}\Vert f \Vert_{\dot H^{\frac{3}{2}}}$ together with the fact that  $\Vert f \Vert_{\dot H^{\frac{5}{2}}}\Vert f \Vert_{\dot H^{2}}\lesssim  \Vert f \Vert_{\dot H^3}\Vert f \Vert_{\dot H^{\frac{3}{2}}} $ and the following inequalities 
$$
\left(\int    \
 \frac{\Vert \delta_\alpha f_x \Vert^2_{L^{\infty}}  }{{\vert \alpha\vert^{\frac{5}{2}} } } \ d\alpha\right)^{\frac{1}{2}} \left(\int    \
 \frac{ \Vert s_\alpha f \Vert^2_{L^{\infty}} }{{\vert \alpha\vert^{\frac{9}{2}} }} \ d\alpha\right)^{\frac{1}{2}} \lesssim \Vert f \Vert^2_{\dot B^{\frac{7}{4}}_{\infty,2}} \lesssim \Vert f \Vert^2_{\dot H^{\frac{9}{4}}} \lesssim  \Vert f \Vert_{\dot H^3}\Vert f \Vert_{\dot H^{\frac{3}{2}}}
 $$
we find that
\begin{eqnarray*}
T_{4,1,12}&\lesssim& \Vert f \Vert^2_{\dot H^3}\Vert f \Vert^2_{\dot H^{\frac{3}{2}}}
\end{eqnarray*}

Gathering all the estimates of the $T_{4,1,i}$ we have proved that 
\begin{eqnarray*}
T_{4,1}&\lesssim& \Vert f \Vert^2_{\dot H^3}\left(\Vert f \Vert_{\dot H^{\frac{3}{2}}}+\Vert f \Vert^2_{\dot H^{\frac{3}{2}}}+\Vert f \Vert^3_{\dot H^{\frac{3}{2}}}\right)
\end{eqnarray*}

This ends the proof of Lemma {\ref{t41c}}. 

\qed

   $\bullet$ {Estimate of $T_{4,2}$}\\

Recall that

\begin{eqnarray*}
T_{4,2}&=&-\frac{1}{2\pi} \int \int\int_{0}^{\infty}\int_{0}^{\infty} \ \sigma e^{-\gamma-\sigma} \ \Lambda^{3} f \ \partial_{\alpha}\left[\frac{1}{\alpha}\sin(\frac{\gamma}{2}D_\alpha f)\cos(\frac{\gamma}{2}S_\alpha f)\frac{1}{\alpha} \int_{0}^{\alpha} s_{\eta} f_{x}    \ d\eta \right]\\
&& \  \times\ {  \delta_\alpha f_x } \left(\partial_{\alpha}\tau_{\alpha}f_x\right) \sin(\sigma \tau_{\alpha} f_x )\cos(\arctan(\tau_{\alpha} f_{x})) \ d\eta  \ d\gamma \ d\sigma  \ d\alpha\ dx.\\
\end{eqnarray*}
Opening the derivative in $\alpha$ gives rise to 4 terms,  in particular using the fact that $ D_\alpha f-2f_x=\frac{1}{\alpha} \int_{0}^{\alpha} s_{\eta} f_{x} \ d\eta,$ we infer that $ \partial_{\alpha}D_\alpha f=\partial_{\alpha}\left(\frac{1}{\alpha} \int_{0}^{\alpha} s_{\eta} f_{x} \ d\eta\right)$
\begin{eqnarray*}
T_{4,2}&=&\frac{1}{2\pi} \int \int\int_{0}^{\infty}\int_{0}^{\infty} \int_{0}^{\alpha}\ \sigma e^{-\gamma-\sigma} \ \Lambda^{3} f \ 
\frac{1}{\alpha^3}\sin(\frac{\gamma}{2}D_\alpha f)\cos(\frac{\gamma}{2}S_\alpha f)  s_{\eta} f_{x}     \\
&& \  \times\ {  \delta_\alpha f_x } \left(\partial_{\alpha}\tau_{\alpha}f_x\right) \sin(\sigma \tau_{\alpha} f_x )\cos(\arctan(\tau_{\alpha} f_{x})) \ d\eta  \ d\gamma \ d\sigma  \ d\alpha\ dx\\
&-& \frac{1}{4\pi} \int \int\int_{0}^{\infty}\int_{0}^{\infty} \int_{0}^{\alpha}\ \sigma \gamma e^{-\gamma-\sigma} \ \Lambda^{3} f \ 
\frac{1}{\alpha^2}\left(\partial_{\alpha}D_\alpha f \right)\cos(\frac{\gamma}{2}D_\alpha f)\cos(\frac{\gamma}{2}S_\alpha f)  s_{\eta} f_{x}    \\
&& \  \times\ {  \delta_\alpha f_x } \left(\partial_{\alpha}\tau_{\alpha}f_x\right) \sin(\sigma \tau_{\alpha} f_x )\cos(\arctan(\tau_{\alpha} f_{x})) \ d\eta  \ d\gamma \ d\sigma  \ d\alpha\ dx\\
&+& \frac{1}{4\pi} \int \int\int_{0}^{\infty}\int_{0}^{\infty} \int_{0}^{\alpha} \sigma \gamma e^{-\gamma-\sigma} \ \Lambda^{3} f \ \frac{1}{\alpha^2}\left(\partial_{\alpha}S_\alpha f \right)\sin(\frac{\gamma}{2}D_\alpha f)\sin(\frac{\gamma}{2}S_\alpha f)  s_{\eta} f_{x}     \\
&& \  \times\ {  \delta_\alpha f_x } \left(\partial_{\alpha}\tau_{\alpha}f_x\right) \sin(\sigma \tau_{\alpha} f_x )\cos(\arctan(\tau_{\alpha} f_{x})) \ d\eta  \ d\gamma \ d\sigma  \ d\alpha \ dx \\
&-&\frac{1}{2\pi} \int \int\int_{0}^{\infty}\int_{0}^{\infty} \ \sigma e^{-\gamma-\sigma} \ \Lambda^{3} f \ \frac{1}{\alpha}\left(\partial_{\alpha}D_\alpha f \right)\sin(\frac{\gamma}{2}D_\alpha f)\cos(\frac{\gamma}{2}S_\alpha f)\\
&& \  \times\ {  \delta_\alpha f_x } \left(\partial_{\alpha}\tau_{\alpha}f_x\right) \sin(\sigma \tau_{\alpha} f_x )\cos(\arctan(\tau_{\alpha} f_{x}))   \ d\gamma \ d\sigma  \ d\alpha\ dx\\
&=&T_{4,2,1}+T_{4,2,2}+T_{4,2,3}+T_{4,2,4}
\end{eqnarray*}

$\bullet$ {Estimate of $T_{4,2,1}$} \\

Since $\partial_{\alpha}\tau_{\alpha}f_x=-\partial^2_{x} \tau_{\alpha}f$

\begin{eqnarray*}
T_{4,2,1} &\lesssim& \Vert f \Vert_{\dot H^3}  \Vert f \Vert_{\dot H^2} \int\frac{\Vert \delta_\alpha f_x \Vert_{L^{\infty}}}{\vert\alpha\vert^3} \int_0^\alpha  \Vert s_{\eta} f_{x}  \Vert_{L^{\infty}}   \ d\eta \ d\alpha \\
&\lesssim&  \Vert f \Vert_{\dot H^3}  \Vert f \Vert_{\dot H^2} \int\frac{\Vert \delta_\alpha f_x \Vert_{L^{\infty}}}{\vert\alpha\vert^3} \left(\int_0^\alpha 
 \frac{\Vert s_{\eta} f_{x}  \Vert^2_{L^{\infty}}}{\vert \eta\vert^2 }   \ d\eta\right)^{\frac{1}{2}}  \left(\int_0^\alpha 
 \vert \eta\vert^2   \ d\eta\right)^{\frac{1}{2}} \ d\alpha \\
 &\lesssim&  \Vert f \Vert_{\dot H^3}  \Vert f \Vert^2_{\dot H^2} \int\frac{\Vert \delta_\alpha f_x \Vert_{L^{\infty}}}{\vert\alpha\vert^{\frac{3}{2}}} \ d\alpha \\
 &\lesssim&  \Vert f \Vert_{\dot H^3}  \Vert f \Vert^2_{\dot H^2} \Vert f \Vert_{\dot B^{\frac{3}{2}}_{\infty,1}} \\
 &\lesssim& \Vert f \Vert^2_{\dot H^3}\Vert f \Vert^{2}_{\dot H^{\frac{3}{2}}}
\end{eqnarray*}
Where one may noticed that the estimate at the end are the same as the term $T_{4,1,7,2}$. \\

$\bullet$ {Estimate of $T_{4,2,2}$} \\

Since we have the identity $\partial_{\alpha}D_\alpha f= - \frac{s_\alpha f_x}{\alpha}-\frac{1}{\alpha^2} \ {\int_0^\alpha s_\kappa f_x \ d\kappa}$, then we have to estimate
\begin{eqnarray*}
 T_{4,2,2}&=&\frac{1}{4\pi} \int \int\int_{0}^{\infty}\int_{0}^{\infty}\int_{0}^{\alpha} \ \sigma \gamma e^{-\gamma-\sigma} \ \Lambda^{3} f \ 
\frac{s_\alpha f_x}{\alpha^3} \cos(\frac{\gamma}{2}D_\alpha f)\cos(\frac{\gamma}{2}S_\alpha f)  s_{\eta} f_{x}     \\
&& \  \times\ {  \delta_\alpha f_x } \left(\partial_{\alpha}\tau_{\alpha}f_x\right) \sin(\sigma \tau_{\alpha} f_x )\cos(\arctan(\tau_{\alpha} f_{x})) \ d\eta  \ d\gamma \ d\sigma  \ d\alpha\ dx\\
&+&\frac{1}{4\pi} \int \int\int_{0}^{\infty}\int_{0}^{\infty} \int_{0}^{\alpha} \ \sigma \gamma e^{-\gamma-\sigma} \ \Lambda^{3} f \ 
\frac{1}{\alpha^4} \ { s_\kappa f_x }\cos(\frac{\gamma}{2}D_\alpha f)\cos(\frac{\gamma}{2}S_\alpha f)  s_{\eta} f_{x}    \ \\
&& \  \times\ {  \delta_\alpha f_x } \left(\partial_{\alpha}\tau_{\alpha}f_x\right) \sin(\sigma \tau_{\alpha} f_x )\cos(\arctan(\tau_{\alpha} f_{x})) \ d\eta \ d\kappa  \ d\gamma \ d\sigma  \ d\alpha\ dx\\
\end{eqnarray*}
While the first term has a mild decay in $\alpha$, the second one is more singular in $\alpha$ but the presence of the two artificial weights allows to get a nice control. More precisely, we write
\begin{eqnarray*}
T_{4,2,2}&\lesssim&  \Vert f \Vert_{\dot H^3}  \Vert f \Vert_{\dot H^2} \int \frac{\Vert s_\alpha f_x \Vert_{L^{\infty}}}{\vert \alpha \vert^3} \left( \int_0^\alpha \frac{\Vert s_\eta f_x \Vert^2_{L^\infty}}{\vert \eta \vert^2} \ d\eta \right)^{\frac{1}{2}} \left( \int_0^\alpha {\vert \eta \vert^2} \ d\eta \right)^{\frac{1}{2}} \ d\alpha\\
&+&\Vert f \Vert_{\dot H^3}  \Vert f \Vert_{\dot H^2} \int \frac{\Vert \delta_{\alpha}f_x \Vert_{L^{\infty}}}{\vert \alpha \vert^4} \left(\int_0^\alpha \frac{\Vert s_\eta f_x \Vert^2_{L^{\infty}}}{\vert \eta \vert^{\frac{3}{2}}} \ d\eta\right)^{\frac{1}{2}}  \left(\int_0^\alpha \frac{\Vert s_\kappa f_x \Vert^2_{L^{\infty}}}{\vert \kappa \vert^{\frac{3}{2}}} \ d\kappa\right)^{\frac{1}{2}} \\
&\times& \left(\int_0^\alpha {\vert \kappa \vert^{\frac{3}{2}}} \ d\kappa\right)^{\frac{1}{2}} \ \left(\int_0^\alpha {\vert \eta \vert^{\frac{3}{2}}} \ d\eta\right)^{\frac{1}{2}} \ d\alpha \\
 &\lesssim&\Vert f \Vert_{\dot H^3}  \Vert f \Vert_{\dot H^2}  \Vert f \Vert_{\dot B^{\frac{3}{2}}_{\infty,2}} \int \frac{\Vert s_\alpha f_x \Vert_{L^{\infty}}}{\vert \alpha \vert^{\frac{3}{2}}} \ d\alpha + \Vert f \Vert_{\dot H^3}  \Vert f \Vert_{\dot H^2} \Vert f \Vert^2_{\dot B^{\frac{5}{4}}_{\infty,2}} \int \frac{\Vert \delta_{\alpha}f_x \Vert_{L^{\infty}}}{\vert \alpha \vert^{\frac{3}{2}}} \ d\alpha\\
 &\lesssim&\Vert f \Vert_{\dot H^3}  \Vert f \Vert^2_{\dot H^2} \Vert f \Vert_{\dot B^{\frac{3}{2}}_{\infty,1}}+ \Vert f \Vert_{\dot H^3}  \Vert f \Vert_{\dot H^2} \Vert f \Vert^2_{\dot H^{\frac{7}{4}}} \Vert f \Vert_{\dot B^{\frac{3}{2}}_{\infty,1}} \\
 &\lesssim&\Vert f \Vert_{\dot H^3}  \Vert f \Vert^2_{\dot H^2} \left(\Vert f \Vert^{\frac{1}{2}}_{\dot B^{1}_{\infty,\infty}}\Vert f \Vert^{\frac{1}{2}}_{\dot B^{2}_{\infty,\infty}}+ \Vert f \Vert_{\dot H^{\frac{3}{2}}} \Vert f \Vert^{\frac{1}{2}}_{\dot B^{1}_{\infty,\infty}}\Vert f \Vert^{\frac{1}{2}}_{\dot B^{2}_{\infty,\infty}} \right) \\
 &\lesssim&\Vert f \Vert_{\dot H^3}  \Vert f \Vert^2_{\dot H^2} \left( \Vert f \Vert^{\frac{1}{2}}_{\dot H^{\frac{3}{2}}}\Vert f \Vert^{\frac{1}{2}}_{\dot H^{\frac{5}{2}}}+ \Vert f \Vert^{\frac{3}{2}}_{\dot H^{\frac{3}{2}}}\Vert f \Vert^{\frac{1}{2}}_{\dot H^{\frac{5}{2}}} \right)
\end{eqnarray*}
Since
$\Vert f \Vert^{\frac{1}{2}}_{\dot H^\frac{5}{2}} \lesssim \Vert f \Vert^{\frac{1}{3}}_{\dot H^3} \Vert f \Vert^{\frac{1}{6}}_{\dot H^{\frac{3}{2}}} $ then $\Vert f \Vert^{\frac{1}{2}}_{\dot H^{\frac{3}{2}}}\Vert f \Vert^{\frac{1}{2}}_{\dot H^{\frac{5}{2}}} \lesssim \Vert f \Vert^{\frac{1}{3}}_{\dot H^3} \Vert f \Vert^{\frac{2}{3}}_{\dot H^{\frac{3}{2}}}$ and since we have $\Vert f \Vert^2_{\dot H^2} \lesssim \Vert f \Vert^{\frac{2}{3}}_{\dot H^3} \Vert f \Vert^{\frac{4}{3}}_{\dot H^{\frac{3}{2}}} $ we find that
\begin{eqnarray*}
 T_{4,2,2}&\lesssim&\Vert f \Vert_{\dot H^3}   \Vert f \Vert^{\frac{2}{3}}_{\dot H^3} \Vert f \Vert^{\frac{4}{3}}_{\dot H^{\frac{3}{2}}}\left( \Vert f \Vert^{\frac{1}{3}}_{\dot H^3} \Vert f \Vert^{\frac{2}{3}}_{\dot H^{\frac{3}{2}}}+ \Vert f \Vert^{\frac{3}{2}}_{\dot H^{\frac{3}{2}}}\Vert f \Vert^{\frac{1}{3}}_{\dot H^3} \Vert f \Vert^{\frac{1}{6}}_{\dot H^{\frac{3}{2}}} \right) \\
 &\lesssim&\Vert f \Vert^2_{\dot H^3}  \left( \Vert f \Vert^2_{\dot H^{\frac{3}{2}}} +  \Vert f \Vert^3_{\dot H^{\frac{3}{2}}} \right)
\end{eqnarray*}

$\bullet$ {Estimate of $T_{4,2,3}$} \\

We have

\begin{eqnarray*}
T_{4,2,3}&=& \frac{1}{4\pi} \int \int\int_{0}^{\infty}\int_{0}^{\infty} \ \sigma \gamma e^{-\gamma-\sigma} \ \Lambda^{3} f \ \frac{1}{\alpha}\left(\partial_{\alpha}S_\alpha f \right)\sin(\frac{\gamma}{2}D_\alpha f)\sin(\frac{\gamma}{2}S_\alpha f)\frac{1}{\alpha} \int_{0}^{\alpha} s_{\eta} f_{x}    \ d\eta \\
&& \  \times\ {  \delta_\alpha f_x } \left(\partial_{\alpha}\tau_{\alpha}f_x\right) \sin(\sigma \tau_{\alpha} f_x )\cos(\arctan(\tau_{\alpha} f_{x})) \ d\eta  \ d\gamma \ d\sigma  \ d\alpha\ dx\\
\end{eqnarray*}
Since
\begin{eqnarray*}
\partial_{\alpha}S_\alpha f &=& \bar\Delta_\alpha f_x-\Delta_\alpha f_x-\frac{s_\alpha f}{\alpha^2}
\end{eqnarray*}
Therefore,
\begin{eqnarray*}
T_{4,2,3}&=& \frac{1}{4\pi} \int \int\int_{0}^{\infty}\int_{0}^{\infty} \ \sigma \gamma e^{-\gamma-\sigma} \ \Lambda^{3} f \ \frac{\bar\delta_\alpha f_x}{\alpha^3}  \sin(\frac{\gamma}{2}D_\alpha f)\sin(\frac{\gamma}{2}S_\alpha f) \int_{0}^{\alpha} s_{\eta} f_{x}    \ d\eta \\
&& \  \times\ {  \delta_\alpha f_x } \left(\partial_{\alpha}\tau_{\alpha}f_x\right) \sin(\sigma \tau_{\alpha} f_x )\cos(\arctan(\tau_{\alpha} f_{x})) \ d\eta  \ d\gamma \ d\sigma  \ d\alpha\ dx\\
&-&  \frac{1}{4\pi} \int \int\int_{0}^{\infty}\int_{0}^{\infty} \ \sigma \gamma e^{-\gamma-\sigma} \ \Lambda^{3} f \ \frac{\delta^2_\alpha f_x}{\alpha^3}  \sin(\frac{\gamma}{2}D_\alpha f)\sin(\frac{\gamma}{2}S_\alpha f) \int_{0}^{\alpha} s_{\eta} f_{x}    \ d\eta \\
&& \  \times\  \left(\partial_{\alpha}\tau_{\alpha}f_x\right) \sin(\sigma \tau_{\alpha} f_x )\cos(\arctan(\tau_{\alpha} f_{x})) \ d\eta  \ d\gamma \ d\sigma  \ d\alpha\ dx\\
&-&\frac{1}{4\pi} \int \int\int_{0}^{\infty}\int_{0}^{\infty} \ \sigma \gamma e^{-\gamma-\sigma} \ \Lambda^{3} f \ \frac{s_\alpha f}{\alpha^4}  \sin(\frac{\gamma}{2}D_\alpha f)\sin(\frac{\gamma}{2}S_\alpha f) \int_{0}^{\alpha} s_{\eta} f_{x}    \ d\eta \\
&& \  \times\ {  \delta_\alpha f_x } \left(\partial_{\alpha}\tau_{\alpha}f_x\right) \sin(\sigma \tau_{\alpha} f_x )\cos(\arctan(\tau_{\alpha} f_{x})) \ d\eta  \ d\gamma \ d\sigma  \ d\alpha\ dx\\
&=&T_{4,2,3,1} + T_{4,2,3,2} + T_{4,2,3,3}
\end{eqnarray*}
As the first two terms can be controlled in the same way, it suffices to estimate one of them, for example $T_{4,2,3,2}$
   
\begin{eqnarray*}
T_{4,2,3,2} &=& -\frac{1}{4\pi} \int \int\int_{0}^{\infty}\int_{0}^{\infty} \ \sigma \gamma e^{-\gamma-\sigma} \ \Lambda^{3} f \ \frac{\delta^2_\alpha f_x}{\alpha^3}  \sin(\frac{\gamma}{2}D_\alpha f)\sin(\frac{\gamma}{2}S_\alpha f) \int_{0}^{\alpha} s_{\eta} f_{x}    \ d\eta \\
&& \  \times\  \left(\partial_{\alpha}\tau_{\alpha}f_x\right) \sin(\sigma \tau_{\alpha} f_x )\cos(\arctan(\tau_{\alpha} f_{x})) \ d\eta  \ d\gamma \ d\sigma  \ d\alpha\ dx\\
&\lesssim& \Vert f \Vert_{\dot H^3}  \Vert f \Vert_{\dot H^2}\int \frac{\Vert s_\alpha f_x \Vert^2_{L^{\infty}}}{\vert \alpha \vert^3} \left( \int_0^\alpha \frac{\Vert s_\eta f_x \Vert^2_{L^\infty}}{\vert \eta \vert^2} \ d\eta \right)^{\frac{1}{2}} \left( \int_0^\alpha {\vert \eta \vert^2} \ d\eta \right)^{\frac{1}{2}} \ d\alpha\\
&\lesssim& \Vert f \Vert_{\dot H^3}  \Vert f \Vert^2_{\dot H^2}\int \frac{\Vert s_\alpha f_x \Vert^2_{L^{\infty}}}{\vert \alpha \vert^{\frac{3}{2}}} \ d\alpha
 \\
 &\lesssim& \Vert f \Vert_{\dot H^3}  \Vert f \Vert^2_{\dot H^2}\Vert f \Vert^2_{\dot B^{\frac{5}{4}}_{\infty,2}}
\end{eqnarray*}
Since $\left[\dot H^2,\dot H^{\frac{3}{2}} \right]_{\frac{1}{2},\frac{1}{2}}=\dot H^{\frac{7}{4}} \hookrightarrow \dot B^{\frac{5}{4}}_{\infty,2} $ we find that, 

\begin{eqnarray*}
T_{4,2,3,2} &\lesssim&\Vert f \Vert_{\dot H^3}  \Vert f \Vert^3_{\dot H^2} \Vert f \Vert_{\dot H^{\frac{3}{2}}} \\
&\lesssim& \Vert f \Vert^2_{\dot H^3} \Vert f \Vert^2_{\dot H^{\frac{3}{2}}} 
\end{eqnarray*}
and the same estimate holds for $T_{3,2,3,1}$. It remains to estimate $T_{4,2,3,3}$.

\begin{eqnarray*}
T_{4,2,3,3}&=&\frac{1}{4\pi} \int \int\int_{0}^{\infty}\int_{0}^{\infty} \ \sigma \gamma e^{-\gamma-\sigma} \ \Lambda^{3} f \ \frac{s_\alpha f}{\alpha^4}  \sin(\frac{\gamma}{2}D_\alpha f)\sin(\frac{\gamma}{2}S_\alpha f) \int_{0}^{\alpha} s_{\eta} f_{x}    \ d\eta \\
&& \  \times\ {  \delta_\alpha f_x } \left(\partial_{\alpha}\tau_{\alpha}f_x\right) \sin(\sigma \tau_{\alpha} f_x )\cos(\arctan(\tau_{\alpha} f_{x})) \ d\eta  \ d\gamma \ d\sigma  \ d\alpha \ dx\\
&\lesssim& \Vert f \Vert_{\dot H^3}  \Vert f \Vert_{\dot H^2} \int \frac{\Vert s_\alpha f \Vert_{L^{\infty}} \Vert \delta_\alpha f_x \Vert_{L^{\infty}}}{\vert \alpha \vert^3} \left( \int_0^\alpha \frac{\Vert s_\eta f_x \Vert_{L^\infty}}{\vert \eta \vert^2} \ d\eta \right)^{\frac{1}{2}} \left( \int_0^\alpha {\vert \eta \vert^2} \ d\eta \right)^{\frac{1}{2}} \ d\alpha\\
&\lesssim& \Vert f \Vert_{\dot H^3}  \Vert f \Vert^2_{\dot H^2} \int \frac{\Vert s_\alpha f \Vert_{L^{\infty}} \Vert \delta_\alpha f_x \Vert_{L^{\infty}}}{\vert \alpha \vert^{\frac{5}{2}}}   \ d\alpha\\
&\lesssim& \Vert f \Vert_{\dot H^3}  \Vert f \Vert^2_{\dot H^2} \left(\sup_{\alpha \in \mathbb R} \frac{\Vert s_\alpha f \Vert_{L^{\infty}} }{\vert \alpha \vert} \right) \int \frac{ \Vert \delta_\alpha f_x \Vert_{L^{\infty}}}{\vert \alpha \vert^{\frac{3}{2}}}   \ d\alpha\\
&\lesssim& \Vert f \Vert_{\dot H^3}  \Vert f \Vert^2_{\dot H^2} \Vert f \Vert_{\dot B^{1}_{\infty,\infty}} \Vert f \Vert_{\dot B^{\frac{3}{2}}_{\infty,1}} \\
&\lesssim& \Vert f \Vert_{\dot H^3}  \Vert f \Vert^2_{\dot H^2} \Vert f \Vert^{\frac{3}{2}}_{\dot B^{1}_{\infty,\infty}} \Vert f \Vert^{\frac{1}{2}}_{\dot B^{2}_{\infty,\infty}} \\
&\lesssim& \Vert f \Vert_{\dot H^3}  \Vert f \Vert^2_{\dot H^2} \Vert f \Vert^{\frac{3}{2}}_{\dot H^{\frac{3}{2}}} \Vert f \Vert^{\frac{1}{2}}_{\dot H^{\frac{5}{2}}}
\end{eqnarray*}
Using again the fact that $\Vert f \Vert^{\frac{1}{2}}_{\dot H^\frac{5}{2}} \lesssim \Vert f \Vert^{\frac{1}{3}}_{\dot H^3} \Vert f \Vert^{\frac{1}{6}}_{\dot H^{\frac{3}{2}}} $ and $\Vert f \Vert^2_{\dot H^2} \lesssim \Vert f \Vert^{\frac{2}{3}}_{\dot H^3} \Vert f \Vert^{\frac{4}{3}}_{\dot H^{\frac{3}{2}}} $ we conclude that
\begin{eqnarray*}
T_{4,2,3,3} &\lesssim& \Vert f \Vert_{\dot H^3}\Vert f \Vert^{\frac{2}{3}}_{\dot H^3} \Vert f \Vert^{\frac{4}{3}}_{\dot H^{\frac{3}{2}}} \Vert f \Vert^{\frac{3}{2}}_{\dot H^{\frac{3}{2}}} \Vert f \Vert^{\frac{1}{3}}_{\dot H^3} \Vert f \Vert^{\frac{1}{6}}_{\dot H^{\frac{3}{2}}}\\
 &\lesssim&\Vert f \Vert^2_{\dot H^3}  \Vert f \Vert^3_{\dot H^{\frac{3}{2}}} 
\end{eqnarray*}
Hence,
\begin{eqnarray*}
T_{4,2,3} &\lesssim&\Vert f \Vert^2_{\dot H^3}  \left( \Vert f \Vert^2_{\dot H^{\frac{3}{2}}} +  \Vert f \Vert^3_{\dot H^{\frac{3}{2}}} \right)
 \end{eqnarray*}
 
 $\bullet$ {Estimate of $T_{4,2,4}$}
 
 \begin{eqnarray*}
 T_{4,2,4}&=&-\frac{1}{2\pi} \int \int\int_{0}^{\infty}\int_{0}^{\infty} \ \sigma e^{-\gamma-\sigma} \ \Lambda^{3} f \ \frac{1}{\alpha}\left(\partial_{\alpha}D_\alpha f \right)\sin(\frac{\gamma}{2}D_\alpha f)\cos(\frac{\gamma}{2}S_\alpha f)\\
&& \  \times\ {  \delta_\alpha f_x } \left(\partial_{\alpha}\tau_{\alpha}f_x\right) \sin(\sigma \tau_{\alpha} f_x )\cos(\arctan(\tau_{\alpha} f_{x}))   \ d\gamma \ d\sigma  \ d\alpha\ dx.\\
\end{eqnarray*}
We use the identity 
\begin{eqnarray*}
\partial_{\alpha}D_\alpha f   &=&  - \frac{s_\alpha f_x}{\alpha}-\frac{1}{\alpha^2} \ {\int_0^\alpha s_\kappa f_x \ d\kappa} 
\end{eqnarray*}
Then, we write
\begin{eqnarray*}
 T_{4,2,4}&=&\frac{1}{2\pi} \int \int\int_{0}^{\infty}\int_{0}^{\infty} \ \sigma e^{-\gamma-\sigma} \ \Lambda^{3} f \ 
 \frac{s_\alpha f_x}{\alpha^2}
 \sin(\frac{\gamma}{2}D_\alpha f)\cos(\frac{\gamma}{2}S_\alpha f)\\
&& \  \times\ {  \delta_\alpha f_x } \left(\partial_{\alpha}\tau_{\alpha}f_x\right) \sin(\sigma \tau_{\alpha} f_x )\cos(\arctan(\tau_{\alpha} f_{x}))   \ d\gamma \ d\sigma \ dx \ d\alpha\\
&+&\frac{1}{2\pi} \int \int\int_{0}^{\infty}\int_{0}^{\infty} \ \sigma e^{-\gamma-\sigma} \ \Lambda^{3} f \
\frac{1}{\alpha^3} \ \left({\int_0^\alpha s_\kappa f_x \ d\kappa}\right)  \sin(\frac{\gamma}{2}D_\alpha f)\cos(\frac{\gamma}{2}S_\alpha f)\\
&& \  \times\ {  \delta_\alpha f_x } \left(\partial_{\alpha}\tau_{\alpha}f_x\right) \sin(\sigma \tau_{\alpha} f_x )\cos(\arctan(\tau_{\alpha} f_{x}))   \ d\gamma \ d\sigma  \ d\alpha\ dx\\
&\lesssim& \Vert f \Vert_{\dot H^3}  \Vert f \Vert_{\dot H^2} \int \frac{\Vert s_\alpha f_x \Vert_{L^{\infty}} \Vert \delta_\alpha f_x \Vert_{L^{\infty}}}{\vert \alpha \vert^{2}}   \ d\alpha\\
&\lesssim& \Vert f \Vert_{\dot H^3}  \Vert f \Vert_{\dot H^2} \left(\int \frac{ \Vert \delta_\alpha f_x \Vert^2_{L^{\infty}}}{\vert \alpha \vert^{2}}   \ d\alpha\right)^{\frac{1}{2}} \left(\int \frac{ \Vert s_\alpha f_x \Vert^2_{L^{\infty}}}{\vert \alpha \vert^{2}}   \ d\alpha\right)^{\frac{1}{2}}\\
&&  + \ \Vert f \Vert_{\dot H^3}  \Vert f \Vert_{\dot H^2}  \int \frac{\Vert s_\alpha f_x \Vert_{L^{\infty}}}{\vert \alpha \vert^3} \left( \int_0^\alpha \frac{\Vert s_\kappa f_x \Vert^2_{L^\infty}}{\vert \kappa \vert^2} \ d\eta \right)^{\frac{1}{2}} \left( \int_0^\alpha {\vert \kappa \vert^2} \ d\eta \right)^{\frac{1}{2}} \ d\alpha\\
&\lesssim& \Vert f \Vert_{\dot H^3}  \Vert f \Vert_{\dot H^2}\Vert f \Vert^2_{\dot B^{\frac{3}{2}}_{\infty,2}} + \Vert f \Vert_{\dot H^3}  \Vert f \Vert^2_{\dot H^2}\int \frac{\Vert s_\alpha f_x \Vert_{L^{\infty}}}{\vert \alpha \vert^{\frac{3}{2}}} \ d\alpha  \\
&\lesssim& \Vert f \Vert_{\dot H^3}  \Vert f \Vert^3_{\dot H^2} + \Vert f \Vert_{\dot H^3}  \Vert f \Vert^2_{\dot H^2}\Vert f \Vert^{\frac{1}{2}}_{\dot B^{1}_{\infty,\infty}}\Vert f \Vert^{\frac{1}{2}}_{\dot B^{2}_{\infty,\infty}} \\
&\lesssim& \Vert f \Vert^2_{\dot H^3}  \Vert f \Vert^2_{\dot H^{\frac{3}{2}}} + \Vert f \Vert_{\dot H^3}  \Vert f \Vert^2_{\dot H^2}\Vert f \Vert^{\frac{1}{2}}_{\dot H^{\frac{3}{2}}}\Vert f \Vert^{\frac{1}{2}}_{\dot H^{\frac{5}{2}}} \\
&\lesssim& \Vert f \Vert^2_{\dot H^3} \left( \Vert f \Vert^2_{\dot H^{\frac{3}{2}}}+ \Vert f \Vert^3_{\dot H^{\frac{3}{2}}} \right)
\end{eqnarray*} 
Finally,
\begin{eqnarray*}
 T_{4,2} \lesssim \Vert f \Vert^2_{\dot H^3} \left( \Vert f \Vert^2_{\dot H^{\frac{3}{2}}}+ \Vert f \Vert^3_{\dot H^{\frac{3}{2}}} \right)
 \end{eqnarray*}
 Gathering the estimate of $T_{4,1}$, $T_{4,2}$ and $T_{4,3}$ we have proved that 
 \begin{eqnarray*}
 T_{4} \lesssim \Vert f \Vert^2_{\dot H^3} \left(\Vert f \Vert_{\dot H^{\frac{3}{2}}} + \Vert f \Vert^2_{\dot H^{\frac{3}{2}}}+ \Vert f \Vert^3_{\dot H^{\frac{3}{2}}} \right)
 \end{eqnarray*}
 Hence Lemma is proved.
 
 \qed
 
 \subsection{Estimate of $T_{5}$} 
 
 $T_{5}$ is similar to $T_{4}$, up to interchanging the role of sin and cos as they were bounded by 1, hence
 
  \begin{eqnarray} \label{t5}
 T_{5} \lesssim \Vert f \Vert^2_{\dot H^3} \left( \Vert f \Vert^2_{\dot H^{\frac{3}{2}}}+ \Vert f \Vert^3_{\dot H^{\frac{3}{2}}} \right)
 \end{eqnarray}
 
  \subsection{Estimate of $T_6$}
  
  Recall that,

\begin{eqnarray*}
T_6 &=&- \frac{1}{4\pi} \int \int\int_{0}^{\infty}\int_{0}^{\infty} \ e^{-\gamma-\sigma} \ \Lambda^{3} f \ S_{\alpha} f \left(\sin(\gamma\Delta_\alpha f)-\sin(\gamma\bar\Delta_\alpha f)\right)   \\
&& \  \times\ \partial_x \left(\frac{\partial_x^2 \tau_{\alpha} f}{\alpha} \cos(\sigma \tau_{\alpha} f_x )\cos(\arctan(\tau_{\alpha} f_{x}))\right)  \ d\alpha \ d\gamma \ d\sigma \ dx \\
\end{eqnarray*}

  We decompose $T_6$ as  follows $T_6=T_{6,1}+T_{6,2}+T_{6,3}$ where
  
 \begin{eqnarray*}
T_{6,1}&=&-\frac{1}{2\pi} \int \int\int_{0}^{\infty}\int_{0}^{\infty} \ e^{-\gamma-\sigma} \ \Lambda^{3} f \ S_{\alpha} f \cos(\frac{\gamma}{2} S_\alpha f) \sin(\frac{\gamma}{2} D_\alpha f)  \\
&& \  \times\ \frac{ \partial_\alpha \delta_{\alpha} f_{xx}}{\alpha} \cos(\sigma \tau_{\alpha} f_x )\cos(\arctan(\tau_{\alpha} f_{x}))  \ d\alpha \ d\gamma \ d\sigma \ dx
\end{eqnarray*}

\begin{eqnarray*}
T_{6,2}&=&-\frac{1}{2\pi} \int \int\int_{0}^{\infty}\int_{0}^{\infty} \ \sigma e^{-\gamma-\sigma} \ \Lambda^{3} f \ S_{\alpha} f \cos(\frac{\gamma}{2} S_\alpha f) \sin(\frac{\gamma}{2} D_\alpha f)  \\
&& \  \times\ \frac{ \delta_{\alpha} f_{xx}}{\alpha} \tau_{\alpha} f_{xx}\sin(\sigma \tau_{\alpha} f_x )\cos(\arctan(\tau_{\alpha} f_{x}))  \ d\alpha \ d\gamma \ d\sigma \ dx
\end{eqnarray*} 

and,

\begin{eqnarray*}
T_{6,3}&=&-\frac{1}{2\pi} \int \int\int_{0}^{\infty}\int_{0}^{\infty} \ e^{-\gamma-\sigma} \ \Lambda^{3} f \ S_{\alpha} f \cos(\frac{\gamma}{2} S_\alpha f) \sin(\frac{\gamma}{2} D_\alpha f)  \\
&& \  \times\ \frac{ \partial_\alpha \delta_{\alpha} f_{xx}}{\alpha} \cos(\sigma \tau_{\alpha} f_x )\cos(\arctan(\tau_{\alpha} f_{x}))  \ d\alpha \ d\gamma \ d\sigma \ dx
\end{eqnarray*}

  $\bullet$ {Estimate of $T_{6,1}$ } \\

$T_{6,1}$ is clearly the most singular term among the $T_{6,i}$, for this term, we shall prove the following Lemma.
\begin{lemma} \label{te61}
The following estimate holds
\begin{equation} \label{t61}
T_{6,1} \lesssim \Vert f \Vert^2_{\dot H^3} \left( \Vert f \Vert_{\dot H^{3/2}}+\Vert f \Vert^2_{\dot H^{3/2}}  +\Vert f \Vert_{\dot H^{3/2}}  \Vert  f \Vert_{\dot B^1_{\infty,1}}\right)
\end{equation}
\end{lemma}
\noindent {\bf{Proof of Lemma \ref{te61}.}} One needs to first start balancing the derivatives by noticing that we have  $\partial_x^3 \tau_{\alpha} f= -\partial_\alpha \delta_{\alpha} f_{xx}$, this allows one to integrate by parts with respect to $\alpha$.

\begin{eqnarray*}
T_{6,1}&=&-\frac{1}{2\pi} \int \int\int_{0}^{\infty}\int_{0}^{\infty} \ e^{-\gamma-\sigma} \ \Lambda^{3} f \ S_{\alpha} f \cos(\frac{\gamma}{2} S_\alpha f) \sin(\frac{\gamma}{2} D_\alpha f)  \\
&& \  \times\ \frac{ \partial_\alpha \delta_{\alpha} f_{xx}}{\alpha} \cos(\sigma \tau_{\alpha} f_x )\cos(\arctan(\tau_{\alpha} f_{x}))  \ d\alpha \ d\gamma \ d\sigma \ dx \\
&=&-\frac{1}{2\pi} \int \int\int_{0}^{\infty}\int_{0}^{\infty} \ e^{-\gamma-\sigma} \ \Lambda^{3} f \ S_{\alpha} f  \sin(\frac{\gamma}{2} D_\alpha f)  \\
&& \  \times\  \frac{ \partial_\alpha \delta_{\alpha} f_{xx}}{\alpha} \cos(\sigma \tau_{\alpha} f_x )\cos(\arctan(\tau_{\alpha} f_{x}))  \ d\alpha \ d\gamma \ d\sigma \ dx \\
&+&\frac{1}{\pi} \int \int\int_{0}^{\infty}\int_{0}^{\infty} \ e^{-\gamma-\sigma} \ \Lambda^{3} f \ S_{\alpha} f \sin^2(\frac{\gamma}{4} S_\alpha f) \sin(\frac{\gamma}{2} D_\alpha f)  \\
&& \  \times\ \frac{ \partial_\alpha \delta_{\alpha} f_{xx}}{\alpha} \cos(\sigma \tau_{\alpha} f_x )\cos(\arctan(\tau_{\alpha} f_{x}))  \ d\alpha \ d\gamma \ d\sigma \ dx \\
&=& T_{6,1,1}+T_{6,1,2}
\end{eqnarray*}
In order to estimate $T_{6,1,1}$ we first integrate by parts in $\alpha$,

\begin{eqnarray*}
T_{6,1,1}&=&-\frac{1}{2\pi} \int \int\int_{0}^{\infty}\int_{0}^{\infty} \ e^{-\gamma-\sigma} \ \Lambda^{3} f \ S_{\alpha} f  \sin(\frac{\gamma}{2} D_\alpha f)  \\
&& \times\ \frac{  \delta_{\alpha} f_{xx}}{\alpha^2} \cos(\sigma \tau_{\alpha} f_x )\cos(\arctan(\tau_{\alpha} f_{x}))  \ d\alpha \ d\gamma \ d\sigma \ dx \\
 &+&\frac{1}{2\pi} \int \int\int_{0}^{\infty}\int_{0}^{\infty} \ e^{-\gamma-\sigma} \ \Lambda^{3} f \ \partial_\alpha(S_{\alpha} f)  \sin(\frac{\gamma}{2} D_\alpha f)\\
&& \times\ \frac{  \delta_{\alpha} f_{xx}}{\alpha} \cos(\sigma \tau_{\alpha} f_x )\cos(\arctan(\tau_{\alpha} f_{x}))  \ d\alpha \ d\gamma \ d\sigma \ dx \\
&+&\frac{1}{4\pi} \int \int\int_{0}^{\infty}\int_{0}^{\infty} \ \gamma e^{-\gamma-\sigma} \ \Lambda^{3} f \ \partial_\alpha(D_{\alpha} f)  \cos(\frac{\gamma}{2} D_\alpha f)\\
&& \times\ \frac{\delta_{\alpha} f_{xx}}{\alpha} \cos(\sigma \tau_{\alpha} f_x )\cos(\arctan(\tau_{\alpha} f_{x}))  \ d\alpha \ d\gamma \ d\sigma \ dx \\
&-&\frac{1}{2\pi} \int \int\int_{0}^{\infty}\int_{0}^{\infty} \ \sigma e^{-\gamma-\sigma} \ \Lambda^{3} f \ S_{\alpha} f  \sin(\frac{\gamma}{2} D_\alpha f)\\
&& \times\ \frac{  \delta_{\alpha} f_{xx}}{\alpha}\ \tau_{\alpha}f_{xx}  \sin(\sigma \tau_{\alpha} f_x )\cos(\arctan(\tau_{\alpha} f_{x}))  \ d\alpha \ d\gamma \ d\sigma \ dx \\
&-&\frac{1}{2\pi} \int \int\int_{0}^{\infty}\int_{0}^{\infty} \ e^{-\gamma-\sigma} \ \Lambda^{3} f \ S_{\alpha} f  \sin(\frac{\gamma}{2} D_\alpha f)\\
&& \times\ \frac{  \delta_{\alpha} f_{xx}}{\alpha}\ \frac{\tau_{\alpha}f_{xx}}{1+(\tau_{\alpha}f_{x})^2}  \cos(\sigma \tau_{\alpha} f_x )\sin(\arctan(\tau_{\alpha} f_{x}))  \ d\alpha \ d\gamma \ d\sigma \ dx \\
&:=&\sum_{i=1}^5 T_{6,1,1,i}
\end{eqnarray*}

$\bullet$ {Estimates of the $T_{6,1,1,i}$} \\

Since,
\begin{eqnarray*}
T_{6,1,1,1}&=&-\frac{1}{2\pi} \int \int\int_{0}^{\infty}\int_{0}^{\infty} \ e^{-\gamma-\sigma} \ \Lambda^{3} f \ S_{\alpha} f  \sin(\frac{\gamma}{2} D_\alpha f)  \\
&& \times\ \frac{  \delta_{\alpha} f_{xx}}{\alpha^2} \cos(\sigma \tau_{\alpha} f_x )\cos(\arctan(\tau_{\alpha} f_{x}))  \ d\alpha \ d\gamma \ d\sigma \ dx. 
\end{eqnarray*}
We find
\begin{eqnarray*}
T_{6,1,1,1}&\lesssim& \Vert f \Vert_{\dot H^3} \int \frac{\Vert s_{\alpha}f \Vert_{L^{2}}}{\vert\alpha \vert^2} \frac{\Vert \delta_{\alpha}f_{xx} \Vert_{L^{\infty}}}{\vert\alpha \vert}\ d\alpha  \\
&\lesssim& \Vert f \Vert_{\dot H^3} \left(\int \frac{\Vert s_{\alpha}f \Vert^2_{L^{2}}}{\vert\alpha \vert^4} \ d\alpha \right)^{1/2} \left(\frac{\Vert \delta_{\alpha}f_{xx} \Vert^2_{L^{\infty}}}{\vert\alpha \vert^2}\ d\alpha \right)^{1/2} \\
&\lesssim& \Vert f \Vert_{\dot H^3} \Vert f \Vert_{\dot H^{3/2}} \Vert f_{xx} \Vert_{\dot B^{1/2}_{\infty,2}} \\
&\lesssim& \Vert f \Vert_{\dot H^3} \Vert f \Vert_{\dot H^{3/2}} \Vert f_{xx} \Vert_{\dot B^{1/2}_{\infty,2}} \\
&\lesssim& \Vert f \Vert^2_{\dot H^3} \Vert f \Vert_{\dot H^{3/2}}    
\end{eqnarray*}
where we used $\dot H^{1} \hookrightarrow \dot B^{1/2}_{\infty,2}.$\\

$\bullet$ {{Estimate of $T_{6,1,1,2}$}} \\

By expanding $\partial_\alpha(S_{\alpha} f)$, we find the following decomposition

\begin{eqnarray*}
T_{6,1,1,2}&=&\frac{1}{2\pi} \int \int\int_{0}^{\infty}\int_{0}^{\infty} \ e^{-\gamma-\sigma} \ \Lambda^{3} f \ \partial_\alpha(S_{\alpha} f)  \sin(\frac{\gamma}{2} D_\alpha f)\\
&& \times\ \frac{  \delta_{\alpha} f_{xx}}{\alpha} \cos(\sigma \tau_{\alpha} f_x )\cos(\arctan(\tau_{\alpha} f_{x}))  \ d\alpha \ d\gamma \ d\sigma \ dx \\
&=&\frac{1}{2\pi} \int \int\int_{0}^{\infty}\int_{0}^{\infty} \ e^{-\gamma-\sigma} \ \Lambda^{3} f \    \bar\Delta_ \alpha f_{x} \sin(\frac{\gamma}{2} D_\alpha f)\\
&& \times\ \frac{  \delta_{\alpha} f_{xx}}{\alpha} \cos(\sigma \tau_{\alpha} f_x )\cos(\arctan(\tau_{\alpha} f_{x}))  \ d\alpha \ d\gamma \ d\sigma \ dx \\
&-&\frac{1}{2\pi} \int \int\int_{0}^{\infty}\int_{0}^{\infty} \ e^{-\gamma-\sigma} \ \Lambda^{3} f \   \Delta_ \alpha f_{x}  \sin(\frac{\gamma}{2} D_\alpha f)\\
&& \times\ \frac{  \delta_{\alpha} f_{xx}}{\alpha} \cos(\sigma \tau_{\alpha} f_x )\cos(\arctan(\tau_{\alpha} f_{x}))  \ d\alpha \ d\gamma \ d\sigma \ dx \\
&+&\frac{1}{2\pi} \int \int\int_{0}^{\infty}\int_{0}^{\infty} \ e^{-\gamma-\sigma} \ \Lambda^{3} f \   \frac{s_{\alpha}f}{\alpha^2} \sin(\frac{\gamma}{2} D_\alpha f)\\
&& \times\ \frac{\delta_{\alpha} f_{xx}}{\alpha} \cos(\sigma \tau_{\alpha} f_x )\cos(\arctan(\tau_{\alpha} f_{x}))  \ d\alpha \ d\gamma \ d\sigma \ dx \\
&=&T_{6,1,1,2,1}+T_{6,1,1,2,2}+T_{6,1,1,2,3}
\end{eqnarray*}
The first term may estimated as follows
\begin{eqnarray*}
T_{6,1,1,2,1} &\lesssim& \Vert f \Vert_{\dot H^3}  \int \frac{\Vert \bar \delta_{\alpha} f_x \Vert_{L^{\infty}}}{\vert\alpha\vert} \frac{\Vert \delta_{\alpha} f_{xx} \Vert_{L^{2}}}{\vert\alpha\vert}\ d\alpha \\
&\lesssim& \Vert f \Vert_{\dot H^3}  \Vert f_x \Vert_{\dot B^{1/2}_{\infty,2}} \Vert f_{xx} \Vert_{\dot B^{1/2}_{2,2}} \\
&\lesssim& \Vert f \Vert_{\dot H^3}  \Vert f \Vert_{\dot H^{2}} \Vert f \Vert_{\dot H^{5/2}} \\
&\lesssim& \Vert f \Vert^2_{\dot H^3}   \Vert f \Vert_{\dot H^{3/2}}
\end{eqnarray*}
The same holds replacing the operator $\bar \delta_\alpha$ by $\delta_\alpha$ as it leads to equivalent semi-norm, hence,
\begin{eqnarray*}
T_{6,1,1,2,2} &\lesssim& \Vert f \Vert^2_{\dot H^3}   \Vert f \Vert_{\dot H^{3/2}}
\end{eqnarray*}
As for $T_{6,1,1,2,3}$, one notices that it is similar to $T_{6,1,1}$, consquently, one immediately obtains that

\begin{eqnarray*}
T_{6,1,1,2,3} &=&\frac{1}{2\pi} \int \int\int_{0}^{\infty}\int_{0}^{\infty} \ e^{-\gamma-\sigma} \ \Lambda^{3} f \   \frac{s_{\alpha}f}{\alpha^2} \sin(\frac{\gamma}{2} D_\alpha f)\\
&& \times\ \frac{\delta_{\alpha} f_{xx}}{\alpha} \cos(\sigma \tau_{\alpha} f_x )\cos(\arctan(\tau_{\alpha} f_{x}))  \ d\alpha \ d\gamma \ d\sigma \ dx \\
&\lesssim&\Vert f \Vert_{\dot H^3} \int \frac{\Vert s_{\alpha}f \Vert_{L^{\infty}}}{\vert\alpha \vert^2} \frac{\Vert \delta_{\alpha}f_{xx} \Vert_{L^{\infty}}}{\vert\alpha \vert}\ d\alpha  \\
&\lesssim& \Vert f \Vert^2_{\dot H^3} \Vert f \Vert_{\dot H^{3/2}}    
\end{eqnarray*}
Hence,
\begin{eqnarray*}
T_{6,1,1,2}&\lesssim& \Vert f \Vert^2_{\dot H^3} \Vert f \Vert_{\dot H^{3/2}}    
\end{eqnarray*}

$\bullet$ {{Estimate of $T_{6,1,1,3}$}}\\

We have that 

\begin{eqnarray*}
T_{6,1,1,3}&=&\frac{1}{4\pi} \int \int\int_{0}^{\infty}\int_{0}^{\infty} \ \gamma e^{-\gamma-\sigma} \ \Lambda^{4} \mathcal{H} f \ \partial_\alpha(D_{\alpha} f)  \cos(\frac{\gamma}{2} D_\alpha f)\\
&& \times\ \frac{\delta_{\alpha} f_{xx}}{\alpha} \cos(\sigma \tau_{\alpha} f_x )\cos(\arctan(\tau_{\alpha} f_{x}))  \ d\alpha \ d\gamma \ d\sigma \ dx \\
&\lesssim& \Vert \Lambda^{3} f \Vert_{L^2} \int \frac{\Vert s_{\alpha} f_x \Vert_{L^{\infty}}}{\vert \alpha \vert} \frac{\Vert \delta_{\alpha} \partial^3_x f \Vert_{L^{2}}}{\vert \alpha \vert} \ d\alpha +
 \Vert \Lambda^{3} f \Vert_{L^2} \int \frac{1}{\vert\alpha\vert^2}\\
 &&\int_0^\alpha \frac{\Vert s_{\kappa} f_x \Vert_{L^{\infty}}}{\vert \kappa \vert^\eta} \vert \kappa \vert^\eta\frac{\Vert \delta_{\alpha} f_{xx} \Vert_{L^{2}}}{\vert \alpha \vert} \ \ d\kappa \ d\alpha  \\
&\lesssim&\Vert f \Vert_{\dot H^3}  \Vert f_x \Vert_{\dot B^{1/2}_{\infty,2}} \Vert f \Vert_{\dot B^{5/2}_{2,2}} + \Vert \Lambda^{3} f \Vert_{L^2} \int \frac{1}{\vert\alpha\vert^2}
\left(\int_0^\alpha \frac{\Vert s_{\kappa} f_x \Vert^{q}_{L^{\infty}}}{\vert \kappa \vert^{q\eta}} \ d\kappa\right)^{1/q} \\
&& \ \ \times \vert \alpha \vert^{\eta+\frac{1}{p}}\frac{\Vert \delta_{\alpha} f_{xx} \Vert_{L^{2}}}{\vert \alpha \vert}  \ d\alpha \\
&\lesssim&\Vert f \Vert_{\dot H^3}  \Vert f \Vert_{\dot B^{3/2}_{\infty,2}} \Vert f \Vert_{\dot H^{5/2}} + \Vert \Lambda^{3} f \Vert_{L^2} \\
&& \times \ \left(\int_0^\infty \frac{\Vert s_{\kappa} f_x \Vert^{q}_{L^{\infty}}}{\vert \kappa \vert^{q\eta}} \ d\kappa\right)^{1/q} \int    \frac{\Vert \delta_{\alpha} f_{xx} \Vert_{L^{2}}}{\vert \alpha \vert^{3-\eta-\frac{1}{p}}}  \ d\alpha \\
&\lesssim&\Vert f \Vert_{\dot H^3}  \Vert f \Vert_{\dot H^{2}} \Vert f \Vert_{\dot H^{5/2}} + \Vert \Lambda^{3} f \Vert_{L^2}  \Vert f_x \Vert_{\dot B^{\eta-\frac{1}{q}}_{\infty,q}} \Vert f_{xx}  \Vert_{\dot B^{2-\eta-\frac{1}{p}}_{2,1}} 
\end{eqnarray*}
By choosing  $\eta=5/4$ and $p=q=2$ note that for such a values we importantly have $2-\eta-\frac{1}{p} \in (0,1)$.  One finds, by using the real interpolation in homogeneous Besov spaces $\dot B^{\frac{7}{4}}_{\infty,2}=[\dot B^{\frac{3}{2}}_{2,\infty}, \dot B^{3}_{2,\infty}]_{\frac{1}{2}}$ together with the fact that 
\begin{eqnarray*}
T_{6,1,1,3}&\lesssim&\Vert f \Vert_{\dot H^3}  \Vert f \Vert_{\dot H^{2}} \Vert f \Vert_{\dot H^{5/2}} + \Vert f \Vert_{\dot H^3}  \Vert f \Vert_{\dot B^{\frac{7}{4}}_{\infty,2}} \Vert f 
\Vert_{\dot B^{\frac{9}{4}}_{2,1}}  \\
&\lesssim&\Vert f \Vert^2_{\dot H^3}  \Vert f \Vert_{\dot H^{\frac{3}{2}}} + \Vert f \Vert_{\dot H^3}  \Vert f \Vert_{\dot B^{\frac{9}{4}}_{2,2}} \Vert f 
\Vert^{\frac{1}{2}}_{\dot B^{\frac{3}{2}}_{2,\infty}} \Vert f \Vert^{\frac{1}{2}}_{\dot B^{3}_{2,\infty}} \\
&\lesssim&\Vert f \Vert^2_{\dot H^3}  \Vert f \Vert_{\dot H^{\frac{3}{2}}}  + \Vert f \Vert_{\dot H^3}  \Vert f \Vert^{\frac{1}{2}}_{\dot H^{\frac{3}{2}}} \Vert f \Vert^{\frac{1}{2}}_{\dot H^{3}}
\Vert f \Vert^{\frac{1}{2}}_{\dot H^{\frac{3}{2}}} \Vert f \Vert^{\frac{1}{2}}_{\dot H^{3}} \\
&\lesssim&\Vert f \Vert^2_{\dot H^3}  \Vert f \Vert_{\dot H^{\frac{3}{2}}} \\
\end{eqnarray*}
$\bullet$ {{Estimate of $T_{6,1,1, 4}$}} \\

Recall that,
\begin{eqnarray*}
T_{6,1,1, 4} &=&- \frac{1}{2\pi} \int \int\int_{0}^{\infty}\int_{0}^{\infty} \ \sigma e^{-\gamma-\sigma} \ \Lambda^{3} f \ S_{\alpha} f  \sin(\frac{\gamma}{2} D_\alpha f)\\
&& \times\ \frac{  \delta_{\alpha} f_{xx}}{\alpha}\ \tau_{\alpha}f_{xx}  \sin(\sigma \tau_{\alpha} f_x )\cos(\arctan(\tau_{\alpha} f_{x}))  \ d\alpha \ d\gamma \ d\sigma \ dx \\
\end{eqnarray*}

This term may be controlled as follows

\begin{eqnarray*}
T_{6,1,1, 4} &\lesssim&  \Vert \Lambda^{3} f \Vert_{L^2} \ \int \frac{\Vert s_{\alpha} f \Vert_{L^{\infty}}}{\vert \alpha \vert^2}
 \Vert \delta_{\alpha} f_{xx} \Vert_{L^4}  \Vert \tau_{\alpha}f_{xx} \Vert_{L^4}   \ d\alpha \\
 &\lesssim& \Vert \Lambda^{3} f \Vert_{L^2} \ \Vert f_{xx} \Vert^2_{L^4}  \int \frac{\Vert s_{\alpha} f \Vert_{L^{\infty}}}{\vert \alpha \vert^2}    \ d\alpha \\
 &\lesssim& \Vert \Lambda^{3} f \Vert_{L^2} \ \Vert f \Vert^2_{\dot H^{\frac{9}{4}}}  \int \frac{\Vert s_{\alpha} f \Vert_{L^{\infty}}}{\vert \alpha \vert^2}    \ d\alpha \\
 &\lesssim& \Vert f \Vert^2_{\dot H^3} \ \Vert f \Vert_{\dot H^{\frac{3}{2}}}  \Vert  f \Vert_{\dot B^1_{\infty,1}}   \\
\end{eqnarray*}

It is easy to see that the same holds for $T_{6,1,1, 5}$ as it is the same as $T_{6,1,1, 4}$ up to some bounded terms which are less or equal to 1. 

For the estimate of $T_{6,1,2}$, we do not need to integrate by parts, indeed we may use the fact that
\begin{eqnarray*}
T_{6,1,2}&=&\frac{1}{\pi} \int \int\int_{0}^{\infty}\int_{0}^{\infty} \ e^{-\gamma-\sigma} \ \Lambda^{3} f \ S_{\alpha} f \sin^2(\frac{\gamma}{4} S_\alpha f) \sin(\frac{\gamma}{2} D_\alpha f)  \\
&& \  \times\ \frac{ \partial_\alpha \delta_{\alpha} f_{xx}}{\alpha} \cos(\sigma \tau_{\alpha} f_x )\cos(\arctan(\tau_{\alpha} f_{x}))  \ d\alpha \ d\gamma \ d\sigma \ dx \\
&\lesssim& \Vert f \Vert^2_{\dot H^3} \int \frac{\Vert s_\alpha f \Vert^2_{L^\infty}}{\vert \alpha \vert^3} \ d\alpha \\
&\lesssim& \Vert f \Vert^2_{\dot H^3} \Vert f \Vert^2_{\dot H^{\frac{3}{2}}}
\end{eqnarray*}

This ends the estimate of $T_{6,1}$ and therefore the proof of Lemma \ref{te61} is done. \\

\qed

$\bullet$ {Estimate of $T_{6,2}$ and $T_{6, 3}$ } \\

The estimates of $T_{6,2}$  and $T_{6, 3}$  are similar to $T_{6,1,1, 5}$ or $T_{6,1,1, 4}$, more precisely

\begin{eqnarray} \label{t6243}
T_{6,2} + T_{6, 3} &\lesssim&  \Vert \Lambda^{3} f \Vert_{L^2} \ \int \frac{\Vert s_{\alpha} f \Vert_{L^{\infty}}}{\vert \alpha \vert^2}
 \Vert \delta_{\alpha} f_{xx} \Vert_{L^4}  \Vert \tau_{\alpha}f_{xx} \Vert_{L^4}   \ d\alpha \nonumber \\
  &\lesssim& \Vert \Lambda^{3} f \Vert_{L^2} \ \Vert f \Vert^2_{\dot H^{\frac{9}{4}}}  \int \frac{\Vert s_{\alpha} f \Vert_{L^{\infty}}}{\vert \alpha \vert^2}    \ d\alpha \nonumber \\
 &\lesssim& \Vert f \Vert^2_{\dot H^3} \ \Vert f \Vert_{\dot H^{\frac{3}{2}}}  \Vert  f \Vert_{\dot B^1_{\infty,1}}   
\end{eqnarray}

Combining the estimates  \eqref{t61} and \eqref{t6243} we conclude that

\begin{equation} \label{t6}
T_{6} \lesssim \Vert f \Vert^2_{\dot H^3} \left( \Vert f \Vert_{\dot H^{3/2}}+\Vert f \Vert^2_{\dot H^{3/2}}  +\Vert f \Vert_{\dot H^{3/2}}  \Vert  f \Vert_{\dot B^1_{\infty,1}}\right)
\end{equation}

Hence, gathering all the estimates of the $T_i, i=2,..,6$, we have obtained the following control
\begin{equation*} 
\frac{1}{2} \partial_t \Vert f \Vert^2_{\dot H^{3/2}} + \int   \ \frac{\vert\Lambda^3 f \vert^2}{(1+ \vert f_x \vert^2)^{\frac{3}{2}}} \ dx \lesssim \Vert f \Vert^2_{\dot H^3} \left( \Vert f \Vert_{\dot H^{3/2}}+\Vert f \Vert^2_{\dot H^{3/2}}  +\Vert f \Vert^3_{\dot H^{3/2}}+\Vert f \Vert_{\dot H^{3/2}}  \Vert  f \Vert_{\dot B^1_{\infty,1}}\right)
\end{equation*}
Integrating in time gives
 \begin{eqnarray*}
\Vert f(T) \Vert^2_{\dot H^{3/2}}+ \int_0^T\int   \ \frac{\vert\Lambda^3 f \vert^2}{(1+ \vert f_x \vert^2)^{\frac{3}{2}}} \ dx \ ds &\lesssim& \Vert f_0 \Vert^2_{\dot H^{3/2}} \\
&+& \int_0^T\Vert f \Vert^2_{\dot H^3} \left( \mathcal{P}(\Vert f \Vert_{\dot H^{3/2}})+\Vert f \Vert_{\dot H^{3/2}}  \Vert  f \Vert_{\dot B^1_{\infty,1}}\right) \ ds,
\end{eqnarray*}
where $\mathcal{P}(X)=X+X^2+X^3$. Then, since $0\leq a \mapsto \frac{1}{(1+ a^2)^{\frac{3}{2}}}$ is decreasing, we infer that
$$
\frac{1}{(1+ \vert f_x\vert^2)^{\frac{3}{2}}} \geq \frac{1}{(1+  \Vert f_x \Vert_{L^\infty}^2)^{\frac{3}{2}}}
$$
Hence, we may write
\begin{equation*} 
\Vert f(T) \Vert^2_{\dot H^{3/2}}+    \ \displaystyle\int_0^T \frac{\int\vert\Lambda^3 f \vert^2  \ dx}{(1+  \displaystyle  \Vert f_x \Vert_{L^\infty_{x}}^2)^{\frac{3}{2}}} \ ds \lesssim \Vert f_0 \Vert^2_{\dot H^{3/2}} +{\displaystyle\int_0^T\Vert f \Vert^2_{\dot H^3} \left( \mathcal{P}(\Vert f \Vert_{\dot H^{3/2}})+\Vert f \Vert_{\dot H^{3/2}}  \Vert  f \Vert_{\dot B^1_{\infty,1}}\right) \ ds}
\end{equation*}
Set $L:= \Vert f_x \Vert_{L^\infty_{x}}$. Then, the {\it{a priori}} estimates for the full equation, that is $(\mathcal{M}_\sigma)$ (see \eqref{CP0}) is

\begin{eqnarray*} 
&&\Vert f(T) \Vert^2_{\dot H^{3/2}}+  \ \int_0^T \frac{\displaystyle \Vert f \Vert_{\dot H^2}^2}{1+  L^2} \ ds +   \ \int_0^T \frac{\displaystyle \Vert f \Vert_{\dot H^3}^2}{(1+  L^2)^{\frac{3}{2}}} \ ds \\
&& \lesssim \Vert f_0 \Vert^2_{\dot H^{3/2}} +{\displaystyle\int_0^T\Vert f \Vert^2_{\dot H^3} \left( \mathcal{P}(\Vert f \Vert_{\dot H^{3/2}})+\Vert f \Vert_{\dot H^{3/2}}  \Vert  f \Vert_{\dot B^1_{\infty,1}}\right) \ ds}+{\displaystyle\int_0^T \mathcal{Q}(\Vert f \Vert_{\dot H^{3/2}}) \Vert f \Vert_{\dot H^{3/2}} \  \ ds},
\end{eqnarray*}
where $\mathcal{Q}(X)=X+X^2$, 
this ends the proof of Theorem \ref{H32}
\qed

\section{{\it{A priori}} estimates for the $L^2$ norm}

Interestingly, eventhough the linearized Muskat equation with surface tension is a nice  third order parabolic equation, one does not have a nice structure as in the Muskat equation without surface tension. Because of the lack of structure in the Muskat equation with surface tension, it  seems   impossible to get an $L^2$-maximum principle or any global control under the assumption that the data lies in $L^2$ only. The best we can do is a control of the $L^2$ under the control of higher order norms. More precisely, in this section we shall prove the following Proposition.
\begin{proposition} \label{L2-est}
Let $f_0 \in H^3$, then we have the following {\it{a priori estimates}}
\begin{eqnarray*}
 \frac{1}{2}\partial_t \Vert f \Vert^2_{L^2}+\displaystyle\int   \ \frac{\vert\Lambda^{\frac{3}{2}} f \vert^2}{(1+ \vert f_x \vert^2)^{\frac{3}{2}}} \ dx &\lesssim&   \Vert f \Vert^2_{H^4}  \left( \Vert f \Vert_{\dot H^{3/2}}+\Vert f \Vert^4_{\dot H^{3/2}}+\Vert f \Vert_{\dot H^{3/2}}  \Vert  f \Vert_{\dot B^1_{\infty,1}}\right) \\
  \end{eqnarray*}
  \end{proposition}
\noindent{{\bf{Proof of Proposition \ref{L2-est}}}.} We compute the dot product of $f$ and $\partial_t f$ with respect to the $L^2$ norm, one finds

\begin{eqnarray*}
\int f\partial_t f \ dx &=& -\int f \ \partial_x \left[\mathcal{H},\int_{0}^{\infty} \int_{0}^{\infty} e^{-\sigma}\cos(\sigma f_x )\cos(\arctan(f_{x}))   \ d \sigma \right] f_{xx} \ dx   \\
&-&  \int f\ \Lambda^{3} f \ \int_{0}^{\infty}  e^{-\sigma}\cos(\sigma f_x )\cos(\tau_{\alpha}\arctan(f_{x}))   \ d \sigma \ dx     \\
&+&\int f \ \mathcal{H}f_{xx} \ \partial_{x}\left(\int_{0}^{\infty} e^{-\sigma}\cos(\sigma f_x )\cos(\arctan(f_{x}))  \ d \sigma \right) \ dx \\
&+& \frac{1}{4\pi}\int f \int\int_{0}^{\infty}\int_{0}^{\infty}\int_{0}^{\alpha} \ e^{-\gamma-\sigma} \ \frac{1}{\alpha} s_{\eta} f_{x}  \left(\sin(\gamma\Delta_\alpha f)-\sin(\gamma\bar\Delta_\alpha f)\right) \\  
&& \ \times \ \partial_x \left(\frac{\partial_x^2 \tau_{\alpha} f}{\alpha} \cos(\sigma \tau_{\alpha} f_x )\cos(\arctan(\tau_{\alpha} f_{x}))\right) \ d\eta \ d\alpha \ d\gamma \ d\sigma \ dx  \\
&+& \frac{1}{4\pi}\int f \ \int\int_{0}^{\infty}\int_{0}^{\infty}\int_{0}^{\alpha} \ e^{-\gamma-\sigma} \ \frac{1}{\alpha} s_{\eta} f_{x}  \left(\sin(\gamma\Delta_\alpha f)+\sin(\gamma\bar\Delta_\alpha f)\right) \\  
&& \ \times \ \partial_x \left(\frac{\partial_x^2 \tau_{\alpha} f}{\alpha} \cos(\sigma \tau_{\alpha} f_x )\cos(\arctan(\tau_{\alpha} f_{x}))\right)  \ d\eta \ d\alpha \ d\gamma \ d\sigma \ dx \\
&+& \frac{1}{4\pi}\int f \ \int\int_{0}^{\infty}\int_{0}^{\infty} \ e^{-\gamma-\sigma} \ S_{\alpha} f \left(\sin(\gamma\Delta_\alpha f)-\sin(\gamma\bar\Delta_\alpha f)\right) \\
&& \ \times \  \partial_x \left(\frac{\partial_x^2 \tau_{\alpha} f}{\alpha} \cos(\sigma \tau_{\alpha} f_x )\cos(\arctan(\tau_{\alpha}f_{x}))\right) \ d\alpha \ d\gamma \ d\sigma \ dx  \\
&+& \frac{1}{4\pi}\int f \ \int\int_{0}^{\infty}\int_{0}^{\infty} \ e^{-\gamma-\sigma} \ S_{\alpha} f \left(\sin(\gamma\Delta_\alpha f)+\sin(\gamma\bar\Delta_\alpha f)\right) \\
&& \ \times \  \partial_x \left(\frac{\partial_x^2 \tau_{\alpha} f}{\alpha} \cos(\sigma \tau_{\alpha} f_x )\cos(\arctan(\tau_{\alpha}f_{x}))\right) \ d\alpha \ d\gamma \ d\sigma \ dx  \\
\end{eqnarray*} 

Unlike the second term which will give rise (up to some remainder terms) to a weighted dissipation of order
\begin{equation}
\mathcal W (f)= -\displaystyle\int   \ \frac{\vert\Lambda^{\frac{3}{2}} f \vert^2}{(1+ \vert f_x \vert^2)^{\frac{3}{2}}} \ dx,
\end{equation}
all the other terms  may be controlled exactly as done in the $\dot H^{\frac{3}{3}}$ estimates. Therefore, it suffices to focus on the term
$$
\mathcal{L}(f)=-  \int f\ \Lambda^{3} f \  \int_{0}^{\infty} e^{-\sigma}\cos(\sigma f_x )\cos(\arctan(f_{x}))   \ d \sigma \ dx $$
We shall prove the following estimate for the estimate of  $\mathcal{L}(f)$
\begin{lemma} \label{Lc} One has the following control of the term $\mathcal{L}(f)$
 \begin{eqnarray*}
  \mathcal{L}(f)  \lesssim -\displaystyle\int   \ \frac{\vert\Lambda^{\frac{3}{2}} f \vert^2}{(1+ \vert f_x \vert^2)^{\frac{3}{2}}} \ dx + \Vert f \Vert^2_{H^4} \left(\Vert f \Vert_{\dot H^{\frac{3}{2}}}+\Vert f \Vert^2_{\dot H^{\frac{3}{2}}}\right) 
  \end{eqnarray*} 
\end{lemma}
\noindent {\bf{Proof of Lemma \ref{Lc}.}}
In order to make appear the expected weighted dissipation one has to integrate by parts, more precisely, we have that
\begin{eqnarray*} 
\mathcal{L}(f)&=&-  \int  \vert \Lambda^{\frac{3}{2}} f \vert^2 \ \int_{0}^{\infty} \int_{0}^{\infty} e^{-\gamma-\sigma}\cos(\sigma f_x )\cos(\arctan(f_{x})) \ d \gamma  \ d \sigma \ dx \\
&-&  \int   f  \ \Lambda^{\frac{3}{2}} f \int_{0}^{\infty} \int_{0}^{\infty} e^{-\gamma-\sigma}\Lambda^{\frac{3}{2}}\left(\cos(\sigma f_x )\cos(\tau_{\alpha}\arctan(f_{x})) \right) \ d \gamma  \ d \sigma \ dx \\
 &+& C\int \Lambda^{\frac{3}{2}} f \int \frac{\delta_\beta f \ \delta_\beta\left( \int_{0}^{\infty} \int_{0}^{\infty} e^{-\gamma-\sigma}\left(\cos(\sigma f_x )\cos(\arctan(f_{x})) \right) \ d \gamma  \ d \sigma\right)}{\vert \beta \vert^{\frac{5}{2}}} \ d\beta \ dx \\
 &:=&\mathcal{L}_1(f) + \mathcal{L}_2(f) + \mathcal{L}_3(f)
\end{eqnarray*}
The $\mathcal{L}_1(f)$ is nothing but the weighted dissipation. Indeed, as 
$$ \int_{0}^{\infty} e^{-\sigma}\cos(\sigma f_x )\cos(\arctan(f_{x}))   \ d \sigma=\frac{1}{\left(1+(f_x)^2\right)^{\frac{3}{2}}}.$$ One finds that
\begin{eqnarray*}
\mathcal{L}_1(f)=-\displaystyle\int   \ \frac{\vert\Lambda^{\frac{3}{2}} f \vert^2}{(1+ \vert f_x \vert^2)^{\frac{3}{2}}} \ dx.
\end{eqnarray*}
For the second term $\mathcal{L}_2(f)$, we have 
\begin{eqnarray*}
\mathcal{L}_2(f)&\lesssim&     \Vert f \Vert_{L^2}  \ \Vert \Lambda^{\frac{3}{2}} f \Vert_{L^4} \left\Vert  \int_{0}^{\infty} e^{-\sigma}\Lambda^{\frac{3}{2}}\left(\cos(\sigma f_x )\cos(\tau_{\alpha}\arctan(f_{x})) \right)   \ d \sigma \right\Vert_{L^4} 
\end{eqnarray*}
The main effort will be devoted to estimate  the $L^4$ norm of 
$$
\widetilde{\mathcal{L}_2}(f):= \Lambda^{\frac{3}{2}}\left(\cos(\sigma f_x )\cos(\arctan(f_{x})) \right)  
$$
We shall prove the following lemma.
\begin{lemma} \label{Ltil}  We have the following control of 
$\Vert {\widetilde{\mathcal{L}_2}}(f)\Vert_{L^4}$
\begin{eqnarray*}
\Vert\widetilde{\mathcal{L}_2}(f)\Vert_{L^4}&\lesssim& \sigma \Vert f \Vert^{\frac{5}{6}}_{\dot H^{3}} \Vert f \Vert^{\frac{1}{6}}_{\dot H^{\frac{3}{2}}}+\sigma^2 \Vert f \Vert^{\frac{2}{3}}_{\dot H^{3}} \Vert f \Vert^{\frac{4}{3}}_{\dot H^{\frac{3}{2}}}+\sigma (k+1) \Vert f \Vert^{\frac{5}{6}}_{\dot H^{3}} \Vert f \Vert^{\frac{7}{6}}_{\dot H^{\frac{3}{2}}}
\end{eqnarray*}
Consequently, by integrating with respect to the positive measure $e^{-\sigma} d\sigma$, we find that $\mathcal{L}_2(f)$ may be estimated as follows
\begin{eqnarray*}
\mathcal{L}_2 (f)&\lesssim&  \Vert f \Vert_{L^2} \Vert f \Vert_{\dot H^3} \left(\Vert f \Vert_{\dot H^{\frac{3}{2}}}+\Vert f \Vert^2_{\dot H^{\frac{3}{2}}}\right)
\end{eqnarray*}
Consequently, 
\begin{eqnarray*}
\mathcal{L}_2 (f)&\lesssim&  \Vert f \Vert^2_{H^4} \left(\Vert f \Vert_{\dot H^{\frac{3}{2}}}+\Vert f \Vert^2_{\dot H^{\frac{3}{2}}}\right)
\end{eqnarray*}
\end{lemma}
\noindent {\bf{Proof of Lemma \ref{Ltil}}.}
One first need to lower the regularity of the fractional Laplacian as our goal will be to use the fact that low Sobolev regularity (at most 1) are preserved by composition with a Lipschitz outer functions. More precisely, we may write
\begin{eqnarray*}
-\widetilde{\mathcal{L}_2}(f)&=&\Lambda^{\frac{1}{2}}\left( \sigma f_{xx}\sin(\sigma f_x )\cos(\arctan(f_{x}))\right)+\Lambda^{\frac{1}{2}}\left(\frac{f_{xx}}{1+(f_x)^2}\cos(\sigma f_x )\sin(\arctan(f_{x})) \right)\\
&=&\mathcal{H}\Lambda^{\frac{1}{2}}\left( \sigma f_{xx}\sin(\sigma f_x )\cos(\tau_{\alpha}\arctan(f_{x}))\right)\\
&& \ + \ \mathcal{H}\Lambda^{\frac{1}{2}}\left(\int_0^\infty e^{-k} {f_{xx}}\cos(k f_x )\cos(\sigma f_x )\sin(\tau_{\alpha}\arctan(f_{x})) \ dk\right)  \\
&=&\sigma\mathcal{H}\Lambda^{\frac{1}{2}}(  f_{xx})\sin(\sigma f_x )\cos(\tau_{\alpha}\arctan(f_{x}))+\sigma f_{xx} \mathcal{H}\Lambda^{\frac{1}{2}}(\sin(\sigma f_x )\cos(\tau_{\alpha}\arctan(f_{x}))) \\
&&-c\ \mathcal{H} \left(\int \frac{\delta_y f_{xx} \ \delta_y(\sin(\sigma f_x )\cos(\arctan(f_{x}))) }{\vert y \vert^{\frac{3}{2}}} \ dy\right) \\
&& \ + \ \Lambda^{\frac{1}{2}}\left(\int_0^\infty e^{-k} {f_{xx}}\cos(k f_x )\cos(\sigma f_x )\sin(\arctan(f_{x})) \ dk\right)
\end{eqnarray*}
As the functions $\phi_\sigma (\cdot)=\sin(\sigma \cdot )\cos(\arctan(\cdot))$ and $\psi_{\sigma,k}(\cdot)=\cos(k \cdot )\cos(\sigma \cdot )\sin(\arctan(\cdot))$ are respectively $\sigma$-Lipschitz and $k\sigma$-Lipschitz we then have nice sublinear estimate in $\dot H^s$, $s \in (0,1)$, in particular $\Vert \phi_\sigma f_x \Vert_{\dot H^{\frac{3}{4}}}\leq \sigma \Vert  f_x \Vert_{\dot H^{\frac{3}{4}}}$. Then, we find using Sobolev embedding and the continuity of $\mathcal{H}$ on $L^p$, $p\in (1,\infty)$ that  
  \begin{eqnarray*}
\Vert\widetilde{\mathcal{L}_2}(f)\Vert_{L^4}&\lesssim&\sigma \Vert f \Vert_{\dot H^{\frac{11}{4}}} +\sigma^2\Vert f \Vert_{\dot H^{\frac{9}{4}}} \Vert  f_x \Vert_{\dot H^{\frac{3}{4}}}
+\int \frac{\Vert\delta_y f_{xx} \Vert_{L^4} \ \Vert\delta_y
\phi_{\sigma} f_{x}\Vert_{L^\infty} }{\vert y \vert^{\frac{3}{2}}} \ dy \\
&& \ + \ \int_0^\infty e^{-k} \Vert \Lambda^{\frac{1}{2}}({f_{xx}}\ \psi_{\sigma,k} (f_x)) \Vert_{L^4} \ dk \\
&\lesssim&\sigma \Vert f \Vert_{\dot H^{\frac{11}{4}}} +\sigma^2\Vert f \Vert_{\dot H^{\frac{9}{4}}} \Vert  f_x \Vert_{\dot H^{\frac{3}{4}}}
+\left(\int \frac{\Vert\delta_y f_{xx} \Vert^2_{L^4} }{\vert y \vert^{\frac{3}{2}}} \ dy\right)^{\frac{1}{2}} \left(\int \frac{ \Vert\delta_y\phi_{\sigma} f_{x}\Vert^2_{L^\infty} }{\vert y \vert^{\frac{3}{2}}} \ dy\right)^{\frac{1}{2}} \\
&& \ + \ \int_0^\infty e^{-k} \Vert \Lambda^{\frac{1}{2}}({f_{xx}}\ \psi_{\sigma,k} (f_x)) \Vert_{L^4} \ dk \\
&\lesssim&\sigma \Vert f \Vert_{\dot H^{\frac{11}{4}}} +\sigma^2\Vert f \Vert_{\dot H^{\frac{9}{4}}} \Vert  f_x \Vert_{\dot H^{\frac{3}{4}}}
+\sigma\Vert f  \Vert_{\dot H^{\frac{5}{2}}} \Vert f  \Vert_{\dot B^{\frac{5}{4}}_{\infty,2}} + \ \int_0^\infty e^{-k} \Vert \Lambda^{\frac{1}{2}}({f_{xx}}\ \psi_{\sigma,k} (f_x)) \Vert_{L^4} \ dk \\\\
&\lesssim&\sigma \Vert f \Vert_{\dot H^{\frac{11}{4}}} +\sigma^2\Vert f \Vert_{\dot H^{\frac{9}{4}}} \Vert  f_x \Vert_{\dot H^{\frac{3}{4}}}
+\sigma\Vert f  \Vert_{\dot H^{\frac{5}{2}}} \Vert f  \Vert_{\dot H^{\frac{7}{4}}} + \ \int_0^\infty e^{-k} \Vert \Lambda^{\frac{1}{2}}({f_{xx}}\ \psi_{\sigma,k} (f_x)) \Vert_{L^4} \ dk \\ \\
\end{eqnarray*}
Now it remains to estimate $\int_0^\infty e^{-k} \Vert \Lambda^{\frac{1}{2}}({f_{xx}}\ \psi_{\sigma,k} (f_x)) \Vert_{L^4} \ dk$. To do so, we write, as $\psi_{\sigma,k} (f_x))$ is uniformly bounded and using the Sobolev embedding $\dot H^{\frac{3}{8}} \hookrightarrow L^8$, that
\begin{eqnarray*}
 \Vert \Lambda^{\frac{1}{2}}({f_{xx}}\ \psi_{\sigma,k} (f_x)) \Vert_{L^4}
 &\lesssim&\Vert \Lambda^{\frac{1}{2}}{f_{xx}}\Vert_{L^4} + \Vert \Lambda^{\frac{1}{2}}\psi_{\sigma,k}(f_x)  \Vert_{L^8}\Vert {f_{xx}}\Vert_{L^8} \\
 && + \ \int_0^{\infty} \int e^{-k}  \frac{\Vert \delta_y f_{xx}\Vert_{L^8} \Vert \delta_y \psi_{\sigma,k}(f_x) \Vert_{L^8}}{\vert y \vert^{\frac{3}{2}}} \ dy\ dk \\
 &\lesssim&\Vert f\Vert_{\dot H^\frac{11}{4}} + \Vert \Lambda^{\frac{1}{2}}\psi_{\sigma,k}(f_x)  \Vert_{L^8}\Vert {f_{xx}}\Vert_{L^8} \\
 && + \ \int_0^{\infty} \int e^{-k}  \frac{\Vert \delta_y f_{xx}\Vert_{L^8} \Vert \delta_y \psi_{\sigma,k}(f_x) \Vert_{L^8}}{\vert y \vert^{\frac{3}{2}}} \ dy\ dk \\
 &\lesssim&\Vert f\Vert_{\dot H^\frac{11}{4}} + \Vert\psi_{\sigma,k}(f_x)  \Vert_{\dot H^{\frac{7}{8}}}\Vert {f_{xx}}\Vert_{L^8} \\
 && + \ \int_0^{\infty}  e^{-k} \int \frac{\Vert \delta_y f_{xx}\Vert_{L^4} \Vert \delta_y \psi_{\sigma,k}(f_x) \Vert_{L^\infty}}{\vert y \vert^{\frac{3}{2}}} \ dy\ dk \\
\end{eqnarray*}
Since $7/8<1$ and because $\psi_{\sigma,k}$ is a bounded $k\sigma$-Lipschitz function,  one has not only nice sublinear estimate in $\dot H^{\frac{7}{8}}$ but also a nice estimate on $\Vert \delta_y \psi_{\sigma,k}(f_x) \Vert_{L^\infty}$. More precisely, since we have that $\Vert\psi_{\sigma,k}(f_x)  \Vert_{\dot H^{\frac{7}{8}}} \lesssim k\sigma\Vert f_x \Vert_{\dot H^{\frac{7}{8}}}$ and $\Vert \delta_y \psi_{\sigma,k}(f_x) \Vert_{L^\infty}\lesssim k\sigma\Vert \delta_y f_x \Vert_{L^\infty},$ consequently,  we obtain
\begin{eqnarray*}
    \Vert \Lambda^{\frac{1}{2}}({f_{xx}}\ \psi_{\sigma,k} (f_x)) \Vert_{L^4} &\lesssim&\Vert f\Vert_{\dot H^\frac{11}{4}} + k\sigma  \Vert f  \Vert_{\dot H^{\frac{15}{8}}}\Vert f \Vert_{\dot H^\frac{19}{8}} \\
 && + \ k\sigma \  \int \frac{\Vert \delta_y f_{xx}\Vert_{L^4} \Vert \delta_y f_x \Vert_{L^\infty}}{\vert y \vert^{\frac{3}{2}}} \ dy \\
 &\lesssim&\Vert f\Vert_{\dot H^\frac{11}{4}} + k\sigma  \Vert f  \Vert_{\dot H^{\frac{15}{8}}}\Vert f \Vert_{\dot H^\frac{19}{8}} \\
 && + \ k\sigma   \left( \int \frac{\Vert \delta_y f_{xx}\Vert^2_{L^4}}{\vert y \vert^{\frac{3}{2}}} \ dy\right)^{\frac{1}{2}} \left( \int \frac{ \Vert \delta_y f_x \Vert^2_{L^\infty}}{\vert y \vert^{\frac{3}{2}}} \ dy\right)^{\frac{1}{2}}  \\
 &\lesssim&\Vert f\Vert_{\dot H^\frac{11}{4}} + k\sigma  \Vert f  \Vert_{\dot H^{\frac{15}{8}}}\Vert f \Vert_{\dot H^\frac{19}{8}} + k\sigma\Vert f  \Vert_{\dot H^{\frac{5}{2}}} \Vert f  \Vert_{\dot B^{\frac{5}{4}}_{\infty,2}} \\
 &\lesssim&\Vert f\Vert_{\dot H^\frac{11}{4}} + k\sigma  \Vert f  \Vert_{\dot H^{\frac{15}{8}}}\Vert f \Vert_{\dot H^\frac{19}{8}} + k\sigma\Vert f  \Vert_{\dot H^{\frac{5}{2}}} \Vert f  \Vert_{\dot H^{\frac{7}{4}}}
\end{eqnarray*}
where we used   the embedding $\dot H^{\frac{7}{4}} \hookrightarrow \dot B^{\frac{5}{4}}_{\infty,2}$. Gathering the estimates,
we have obtained that
\begin{eqnarray*}
\Vert\widetilde{\mathcal{L}_2}(f)\Vert_{L^4}&\lesssim& \sigma \Vert f \Vert_{\dot H^{\frac{11}{4}}} +\sigma^2\Vert f \Vert_{\dot H^{\frac{9}{4}}} \Vert  f\Vert_{\dot H^{\frac{7}{4}}}
+\sigma\Vert f  \Vert_{\dot H^{\frac{5}{2}}} \Vert f  \Vert_{\dot H^{\frac{7}{4}}} \\
&& \ +  k\sigma  \Vert f  \Vert_{\dot H^{\frac{15}{8}}}\Vert f \Vert_{\dot H^\frac{19}{8}} \ + \ k\sigma\Vert f  \Vert_{\dot H^{\frac{5}{2}}} \Vert f  \Vert_{\dot H^{\frac{7}{4}}}
\end{eqnarray*}
To conclude, we use Sobolev interpolation and find that
\begin{eqnarray*}
\Vert\widetilde{\mathcal{L}_2}(f)\Vert_{L^4}&\lesssim& \sigma \Vert f \Vert^{\frac{5}{6}}_{\dot H^{3}} \Vert f \Vert^{\frac{1}{6}}_{\dot H^{\frac{3}{2}}} +\sigma^2 \Vert f \Vert^{\frac{1}{2}}_{\dot H^{3}} \Vert f \Vert^{\frac{1}{2}}_{\dot H^{\frac{3}{2}}} \Vert f \Vert^{\frac{1}{6}}_{\dot H^{3}} \Vert f \Vert^{\frac{5}{6}}_{\dot H^{\frac{3}{2}}}
+\sigma\Vert f \Vert^{\frac{2}{3}}_{\dot H^{3}} \Vert f \Vert^{\frac{1}{3}}_{\dot H^{\frac{3}{2}}}  \Vert f \Vert^{\frac{1}{6}}_{\dot H^{3}} \Vert f \Vert^{\frac{5}{6}}_{\dot H^{\frac{3}{2}}} \\
&+&   k\sigma  \Vert f \Vert^{\frac{1}{4}}_{\dot H^{3}} \Vert f \Vert^{\frac{3}{4}}_{\dot H^{\frac{3}{2}}} \Vert f \Vert^{\frac{7}{12}}_{\dot H^{3}} \Vert f \Vert^{\frac{5}{12}}_{\dot H^{\frac{3}{2}}} \ + \ k\sigma
 \Vert f \Vert^{\frac{2}{3}}_{\dot H^{3}} \Vert f \Vert^{\frac{1}{3}}_{\dot H^{\frac{3}{2}}}  \Vert f \Vert^{\frac{1}{6}}_{\dot H^{3}} \Vert f \Vert^{\frac{5}{6}}_{\dot H^{\frac{3}{2}}}
\end{eqnarray*}
Hence,
\begin{eqnarray*}
\Vert\widetilde{\mathcal{L}_2}(f)\Vert_{L^4}&\lesssim& \sigma \Vert f \Vert^{\frac{5}{6}}_{\dot H^{3}} \Vert f \Vert^{\frac{1}{6}}_{\dot H^{\frac{3}{2}}}+\sigma^2 \Vert f \Vert^{\frac{2}{3}}_{\dot H^{3}} \Vert f \Vert^{\frac{4}{3}}_{\dot H^{\frac{3}{2}}}+\sigma (k+1) \Vert f \Vert^{\frac{5}{6}}_{\dot H^{3}} \Vert f \Vert^{\frac{7}{6}}_{\dot H^{\frac{3}{2}}}
\end{eqnarray*}

Since, by Sobolev embedding,

\begin{eqnarray*}
\mathcal{L}_2(f)&\lesssim&     \Vert f \Vert_{L^2}  \ \Vert \Lambda^{\frac{3}{2}} f \Vert_{L^4}  \int_{0}^{\infty}  e^{-\sigma}\left\Vert\Lambda^{\frac{3}{2}}\left(\cos(\sigma f_x )\cos(\tau_{\alpha}\arctan(f_{x})) \right) \right\Vert_{L^4}   \ d \sigma  \\
&\lesssim&     \Vert f \Vert_{L^2}  \ \Vert f \Vert^{\frac{1}{6}}_{\dot H^{3}} \Vert f \Vert^{\frac{5}{6}}_{\dot H^{\frac{3}{2}}} \left\Vert \widetilde{\mathcal{L}_2}(f) \right\Vert_{L^4} \\
\end{eqnarray*}

Therefore, we find 

\begin{eqnarray*}
\mathcal{L}_2(f)&\lesssim& \sigma \Vert f \Vert_{L^2}  \ \Vert f \Vert^{\frac{1}{6}}_{\dot H^{3}} \Vert f \Vert^{\frac{5}{6}}_{\dot H^{\frac{3}{2}}}\Vert f \Vert^{\frac{5}{6}}_{\dot H^{3}} \Vert f \Vert^{\frac{1}{6}}_{\dot H^{\frac{3}{2}}}+\sigma^2 \Vert f \Vert_{L^2}  \ \Vert f \Vert^{\frac{1}{6}}_{\dot H^{3}} \Vert f \Vert^{\frac{5}{6}}_{\dot H^{\frac{3}{2}}}\Vert f \Vert^{\frac{2}{3}}_{\dot H^{3}} \Vert f \Vert^{\frac{4}{3}}_{\dot H^{\frac{3}{2}}}\\
&&+\sigma (k+1) \Vert f \Vert_{L^2}  \ \Vert f \Vert^{\frac{1}{6}}_{\dot H^{3}} \Vert f \Vert^{\frac{5}{6}}_{\dot H^{\frac{3}{2}}}\Vert f \Vert^{\frac{5}{6}}_{\dot H^{3}} \Vert f \Vert^{\frac{7}{6}}_{\dot H^{\frac{3}{2}}} 
\end{eqnarray*}
Therefore,
\begin{eqnarray*}
\mathcal{L}_2(f)&\lesssim& (\sigma^2+1)(k+1) \Vert f \Vert_{L^2} \Vert f \Vert_{\dot H^3} \left(\Vert f \Vert_{\dot H^{\frac{3}{2}}}+\Vert f \Vert^2_{\dot H^{\frac{3}{2}}}\right)
\end{eqnarray*}
Finally,
\begin{eqnarray*}
\mathcal{L}_2(f)&\lesssim& (\sigma^2+1)(k+1) \Vert f \Vert^2_{H^4} \left(\Vert f \Vert_{\dot H^{\frac{3}{2}}}+\Vert f \Vert^2_{\dot H^{\frac{3}{2}}}\right)
\end{eqnarray*}

It remains to estimate $\mathcal{L}_3(f)$, as $\lambda_\sigma(\cdot):=\cos(\sigma \cdot )\cos(\arctan(\cdot))$ is a bounded $\sigma$-Lipschitz function, we obtain that
\begin{eqnarray*}
\mathcal{L}_3(f)&=& C\int \Lambda^{\frac{3}{2}} f \int \frac{\delta_\beta f \ \delta_\beta\left( \int_{0}^{\infty}  e^{-\sigma}\left(\cos(\sigma f_x )\cos(\arctan(f_{x})) \right)   \ d \sigma\right)}{\vert \beta \vert^{\frac{5}{2}}} \ d\beta \ dx \\
 &\lesssim& \Vert \Lambda^{\frac{3}{2}} f \Vert_{L^4} \int_{0}^{\infty} e^{-\sigma}\int \frac{\Vert \delta_\beta f \Vert_{L^2} \ \left\Vert \delta_\beta\lambda_{\sigma}(f_x)  \right\Vert_{L^4} }{\vert \beta \vert^{\frac{5}{2}}} \ d\beta \ d\sigma  \\
 &\lesssim& \Vert f \Vert_{\dot H{^\frac{7}{4}}}\int_{0}^{\infty} \sigma e^{-\sigma}\int \frac{\Vert \delta_\beta f \Vert_{L^2} \ \left\Vert \delta_\beta f_x  \right\Vert_{L^4} }{\vert \beta \vert^{\frac{5}{2}}}  \ d\beta \ d\sigma  \\
 &\lesssim& \Vert f \Vert_{\dot H{^\frac{7}{4}}} \int_{0}^{\infty} \sigma e^{-\sigma}\left(\int \frac{\Vert \delta_\beta f \Vert^2_{L^2} }{\vert \beta \vert^{\frac{5}{2}}}  \ d\beta \right)^{\frac{1}{2}} \left(\int \frac{\Vert \delta_\beta f_x \Vert^2_{L^4}  }{\vert \beta \vert^{\frac{5}{2}}}  \ d\beta \right)^{\frac{1}{2}} \ d\sigma  \\
 &\lesssim& \Vert f \Vert_{\dot H{^\frac{7}{4}}}  \Vert f \Vert_{\dot H{^\frac{3}{4}}} \Vert f \Vert_{\dot B^{\frac{7}{4}}_{4,2}}\\
 &\lesssim& \Vert f \Vert_{\dot H{^\frac{7}{4}}}  \Vert f \Vert_{\dot H{^\frac{3}{4}}} \Vert f \Vert_{\dot H^{2}}\\
 &\lesssim& \Vert f \Vert^{\frac{1}{6}}_{\dot H^{3}} \Vert f \Vert^{\frac{5}{6}}_{\dot H^{\frac{3}{2}}} \Vert f \Vert^{\frac{1}{2}}_{L^2} \Vert f \Vert^{\frac{1}{2}}_{\dot H^{\frac{3}{2}}}\Vert f \Vert^{\frac{1}{3}}_{\dot H^{3}} \Vert f \Vert^{\frac{2}{3}}_{\dot H^{\frac{3}{2}}} \\
 &\lesssim& \Vert f \Vert^{\frac{1}{2}}_{\dot H^{3}} \Vert f \Vert^{2}_{\dot H^{\frac{3}{2}}} \Vert f \Vert^{\frac{1}{2}}_{L^2} \\
 &\lesssim& \Vert f \Vert^2_{H^{4}} \Vert f \Vert^{2}_{\dot H^{\frac{3}{2}}}  
 \end{eqnarray*}
 where we used Sobolev embedding. 
 
 \qed
 
 Therefore, gathering all the estimates of the $\mathcal{L}_i(f)$ for $i=1,2,3$, we have proved that,
  \begin{eqnarray*}
  \mathcal{L}(f)  \leq -\displaystyle\int   \ \frac{\vert\Lambda^{\frac{3}{2}} f \vert^2}{(1+ \vert f_x \vert^2)^{\frac{3}{2}}} \ dx + C \Vert f \Vert^2_{\dot H^4} \left(\Vert f \Vert_{\dot H^{\frac{3}{2}}}+\Vert f \Vert^2_{\dot H^{\frac{3}{2}}}\right) 
  \end{eqnarray*}
 
 As we already metionned, all the other are similar as it suffices to interchange the role of one of the $\Lambda^3 f$ and $f$. So we have obtained the following {\it{a priori}} estimate for the $L^2$ norm
   \begin{eqnarray*}
 \frac{1}{2}\partial_t \Vert f \Vert^2_{L^2}+\displaystyle\int   \ \frac{\vert\Lambda^{\frac{3}{2}} f \vert^2}{(1+ \vert f_x \vert^2)^{\frac{3}{2}}} \ dx &\lesssim&   \Vert f \Vert^2_{H^4}  \left( \Vert f \Vert_{\dot H^{3/2}}+\Vert f \Vert^2_{\dot H^{3/2}}  +\Vert f \Vert^3_{\dot H^{3/2}}+\Vert f \Vert_{\dot H^{3/2}}  \Vert  f \Vert_{\dot B^1_{\infty,1}}\right) \\
  \end{eqnarray*}
  Which ends the proof of the $L^2$ estimates
 \qed


\section{Control of  the $\dot B^1_{\infty,1}$ semi-norm}

 To control the $\dot B^1_{\infty,1}$ semi-norm it suffices for example to control the $H^{\frac{5}{2}}$ norm, this will also allows us in particular to prevent the Lipschitz semi-norm from blowing-up. More precisely, the shall prove the main result of this paper, namely, the following theorem.
 
\begin{theorem} \label{main}
There exists a universal constant $C>0,$ such that for any $f_0 \in  H^{\frac{5}{2}}$ satisfying 
\begin{equation}
\left(C\Vert f_0 \Vert_{\dot H^{\frac{3}{2}}}+ C\Vert f_0 \Vert^4_{\dot H^{\frac{3}{2}}} + \Vert f_0 \Vert_{\dot H^{\frac{3}{2}}} \Vert f_0 \Vert_{\dot B^{1}_{\infty,1}}\right)\left(1+\Vert f_0\Vert^2_{\dot W^{1,\infty}}\right)^{\frac{3}{2}}<1,
\end{equation}  
there exists a unique global solution $f$ to the Muskat problem with surface tension which verifies
$$
f \in L^\infty([0,\infty], H^{\frac{5}{2}}) \cap L^2([0,\infty], \dot H^4 \cap \dot H^{3}).
$$
\end{theorem}

\noindent {\bf{Proof of Theorem \ref{main}}.} We study  the time evolution of the  $\dot H^{\frac{5}{2}}$ semi-norm, by doing an integration by parts, we find
 
 \begin{eqnarray*}
\int \Lambda^{4}f \ \partial_t \Lambda f \ dx &=&\int \Lambda^{4}f  \ \Lambda\partial_{x} \left[\mathcal{H}, \int_{0}^{\infty} e^{-\sigma}\cos(\sigma f_x )\cos(\arctan(f_{x}))   \ d \sigma \right]f_{xx}  \ dx  \\
&-&  \int \Lambda^{4}f  \ \Lambda\left( \Lambda^{3} f \  \int_{0}^{\infty} e^{-\sigma}\cos(\sigma f_x )\cos(\arctan(f_{x}))\ d \sigma \right) \ dx    \\
&+&\int \Lambda^{4}f \ \Lambda\mathcal{H}f_{xx} \ \partial_{x}\left(\int_{0}^{\infty} e^{-\sigma}\cos(\sigma f_x )\cos(\arctan(f_{x}))  \ d \sigma \right) \ dx \\
&-& \frac{1}{4\pi}\int \Lambda^{4}\mathcal{H}f \ \partial_x \int\int_{0}^{\infty}\int_{0}^{\infty} \ e^{-\gamma-\sigma} \ S_{\alpha} f \left(\sin(\gamma\Delta_\alpha f)+\sin(\gamma\bar\Delta_\alpha f)\right) \\
&& \ \times \  \partial_x \left(\frac{\partial_x^2 \tau_{\alpha} f}{\alpha} \cos(\sigma \tau_{\alpha} f_x )\cos(\arctan(\tau_{\alpha}f_{x}))\right) \ d\alpha \ d\gamma \ d\sigma \ dx \\
&-& \frac{1}{4\pi}\int \Lambda^{4}\mathcal{H}f \ \partial_x \int_{0}^{\infty}\int_{0}^{\infty}\int_{0}^{\alpha} \ e^{-\gamma-\sigma} \ \frac{1}{\alpha} s_{\eta} f_{x}  \left(\sin(\gamma\Delta_\alpha f)-\sin(\gamma\bar\Delta_\alpha f)\right) \\  
&& \ \times \ \partial_x \left(\frac{\partial_x^2 \tau_{\alpha} f}{\alpha} \cos(\sigma \tau_{\alpha} f_x )\cos(\arctan(\tau_{\alpha} f_{x}))\right) \ d\eta \ d\alpha \ d\gamma \ d\sigma\ dx \\
&-& \frac{1}{4\pi}\int \Lambda^{4}\mathcal{H}f \ \partial_x\int\int_{0}^{\infty}\int_{0}^{\infty}\int_{0}^{\alpha} \ e^{-\gamma-\sigma} \ \frac{1}{\alpha}  s_{\eta} f_{x}  \left(\sin(\gamma\Delta_\alpha f)+\sin(\gamma\bar\Delta_\alpha f)\right)  \\  
&& \times \ \partial_x \left(\frac{\partial_x^2 \tau_{\alpha} f}{\alpha} \cos(\sigma \tau_{\alpha} f_x )\cos(\arctan(\tau_{\alpha} f_{x}))\right)  \ d\eta \ d\alpha \ d\gamma \ d\sigma\ dx  \\
&-& \frac{1}{4\pi}\int \Lambda^{4}\mathcal{H}f \ \partial_x \int\int_{0}^{\infty}\int_{0}^{\infty} \ e^{-\gamma-\sigma} \ S_{\alpha} f \left(\sin(\gamma\Delta_\alpha f)-\sin(\gamma\bar\Delta_\alpha f)\right) \\
&& \ \times \  \partial_x \left(\frac{\partial_x^2 \tau_{\alpha} f}{\alpha} \cos(\sigma \tau_{\alpha} f_x )\cos(\arctan(\tau_{\alpha}f_{x}))\right) \ d\alpha \ d\gamma \ d\sigma \ dx \\
&=&\mathcal{U}_{1} + \mathcal{U}_{2,dissip} + \mathcal{U}_{2,rem} +\sum_{i=3}^6 \mathcal{U}_i
\end{eqnarray*} 

Next, we shall estimates  the $\mathcal{U}_i$, $i=1,...,6$. We shall prove the following Lemma
\begin{lemma} The following controls holds
 \begin{eqnarray} \label{U1}
\mathcal{U}_1 &\lesssim&\Vert f \Vert^2_{\dot H^4}\left(\Vert f \Vert^2_{\dot H^{\frac{3}{2}}}+\Vert f \Vert_{\dot H^{\frac{3}{2}}}\right)
  \end{eqnarray}
 \begin{eqnarray}  \label{U2}
\mathcal{U}_{2,dissip} &\lesssim&- \ \int \frac{\vert \Lambda^{4}f \vert^2}{(1+(f_x)^2)^{\frac{3}{2}}}    \ dx + \Vert f \Vert^2_{\dot H^4}\left(\Vert f \Vert^2_{\dot H^{\frac{3}{2}}}+\Vert f \Vert_{\dot H^{\frac{3}{2}}}\right)
  \end{eqnarray}
  \begin{eqnarray}  \label{U2rem}
  \mathcal{U}_{2,rem} &\lesssim& \Vert f \Vert^2_{\dot H^4} \Vert f \Vert_{\dot H^{\frac{3}{2}}}
  \end{eqnarray}
  \begin{eqnarray} \label{U3}
    \mathcal{U}_{3}&\lesssim& \Vert f \Vert^2_{\dot H^4} \left( \Vert f \Vert^2_{\dot H^{\frac{3}{2}}}  + \Vert f \Vert^3_{\dot H^{\frac{3}{2}}} +  \Vert f \Vert^4_{\dot H^{\frac{3}{2}}} \right)
\end{eqnarray}
\begin{eqnarray}
\mathcal{U}_{4}\lesssim  \Vert f \Vert^2_{\dot H^4} \left(\Vert f \Vert_{\dot H^{\frac{3}{2}}}+\Vert f \Vert^2_{\dot H^{\frac{3}{2}}}+\Vert f \Vert^3_{\dot H^{\frac{3}{2}}}+\Vert f \Vert^4_{\dot H^{\frac{3}{2}}}\right)
\end{eqnarray}
\begin{eqnarray}
\mathcal{U}_{5}\lesssim  \Vert f \Vert^2_{\dot H^4} \left(\Vert f \Vert_{\dot H^{\frac{3}{2}}}+\Vert f \Vert^2_{\dot H^{\frac{3}{2}}}+\Vert f \Vert^3_{\dot H^{\frac{3}{2}}}+\Vert f \Vert^4_{\dot H^{\frac{3}{2}}}\right)
\end{eqnarray}
\begin{eqnarray} \label{U6}
\mathcal{U}_{6} \lesssim \Vert f \Vert^2_{\dot H^4} \left( \Vert f \Vert_{\dot H^{3/2}}+\Vert f \Vert^2_{\dot H^{3/2}}  +\Vert f \Vert_{\dot H^{3/2}}  \Vert  f \Vert_{\dot B^1_{\infty,1}}\right)
\end{eqnarray}
\end{lemma}

\subsection{Estimate of  $\mathcal{U}_1$}

Using the fact that $\mathcal{H}$ is skew-symmetric and that it is a continuous operator on $L^p$
\begin{eqnarray*}
\mathcal{U}_1&=&\int \Lambda^{4}f  \ \Lambda\partial_{x} \left[\mathcal{H},\int_{0}^{\infty}  e^{-\sigma}\cos(\sigma f_x )\cos(\arctan(f_{x}))   \ d \sigma \right]f_{xx}  \ dx \\
&=&-\int \Lambda^{4}\mathcal{H}f  \  \partial_{xx} \left[\mathcal{H},\int_{0}^{\infty} \int_{0}^{\infty} e^{-\sigma}\cos(\sigma f_x )\cos(\arctan(f_{x}))   \ d \sigma \right]f_{xx}  \ dx \\
&\lesssim& \Vert f \Vert_{\dot H^4}  \int_{0}^{\infty} e^{-\sigma} \left\Vert \partial_{xx}\left(\cos(\sigma f_x )\cos(\arctan(f_{x}))\right)\right\Vert_{BMO} \Vert f \Vert_{\dot H^2} \ d\sigma
 \end{eqnarray*}
 
As we do not aim to use any symmetry, it is more convenient to linearize in order to compute the second derivative, using that
 $$
\cos(\sigma f_x )\cos(\arctan(f_{x}))=\frac{1}{2}\left(\cos(\sigma f_x +\arctan(f_{x}) ) + \cos(\sigma f_x - \arctan(f_{x}) )\right)
 $$
 It is easy to see that, up to some bounded functions, one has
 \begin{equation} \label{dsec}
 \partial^2_{x}\left(\cos(\sigma f_x )\cos(\arctan(f_{x})) \approx (\sigma+\sigma^2)\left((\partial^2_{x}f)^2+ \partial^3_{x}f \right)\right)
 \end{equation}
 Therefore,
 \begin{eqnarray*}
\mathcal{U}_1 &\lesssim& \Vert f \Vert_{\dot H^4}  \int_{0}^{\infty} e^{-\sigma} \left\Vert \partial^2_{x}\left(\cos(\sigma f_x )\cos(\arctan(f_{x}))\right)\right\Vert_{BMO} \Vert f \Vert_{\dot H^2} \ d\sigma \\
&\lesssim&\Vert f \Vert_{\dot H^4} \Vert f \Vert_{\dot H^2}  \int_{0}^{\infty} (\sigma+\sigma^2) e^{-\sigma} \left(\left\Vert (\partial^2_{x}f)^2 \right\Vert_{\dot H^{\frac{1}{2}}}  +\left\Vert\partial^3_{x}f\right\Vert_{\dot H^{\frac{1}{2}}}  \right) \ d\sigma \\
&\lesssim&\Vert f \Vert_{\dot H^4} \Vert f \Vert_{\dot H^2} \int_{0}^{\infty} (\sigma+\sigma^2) e^{-\sigma} \left(\left\Vert \Lambda^{\frac{1}{2}}\left((\partial^2_{x}f)^2\right) \right\Vert_{L^2}  +\left\Vert f \right\Vert_{\dot H^{\frac{7}{2}}}  \right) \ d\sigma 
 \end{eqnarray*}
 In order to avoid the use of $L^\infty$ in the estimate of the $L^2$ norm of $\Lambda^{\frac{1}{2}}((\partial^2_{x}f))^2$, we may use  the fact that
 $$
 \Lambda^{\frac{1}{2}}((\partial^2_{x}f))^2=2\partial^2_{x}f\Lambda^{\frac{1}{2}}(\partial^2_{x}f)-C\int \frac{\left(\delta_{\beta} \partial^2_{x}f\right)^{2}}{\vert \beta \vert^{\frac{3}{2}}} \ d\beta
 $$
 Therefore, 
 \begin{eqnarray*}
 \Vert \Lambda^{\frac{1}{2}}((\partial^2_{x}f))^2\Vert_{L^2}&\lesssim& \Vert \partial^2_{x}f \Vert_{L^4} \ \Vert\Lambda^{\frac{1}{2}}(\partial^2_{x}f) \Vert_{L^4}+\int \frac{\Vert\delta_{\beta} \partial^2_{x}f\Vert_{L^4}^{2}}{\vert \beta \vert^{\frac{3}{2}}} \ d\beta \\
 &\lesssim& \Vert f \Vert_{\dot H^{\frac{9}{4}}} \ \Vert f \Vert_{\dot H^{\frac{11}{4}}}+  \Vert f \Vert^2_{\dot B^{\frac{9}{4}}_{4,2}}\\
 &\lesssim& \Vert f \Vert_{\dot H^{\frac{9}{4}}} \ \Vert f \Vert_{\dot H^{\frac{11}{4}}}+  \Vert f \Vert^2_{\dot H^{\frac{5}{2}}}\\
 \end{eqnarray*}
 Hence, using this estimate we infer
 \begin{eqnarray*}
\mathcal{U}_1 &\lesssim&\Vert f \Vert_{\dot H^4}\Vert f \Vert_{\dot H^2}\left(\Vert f \Vert_{\dot H^{\frac{9}{4}}} \ \Vert f \Vert_{\dot H^{\frac{11}{4}}}+  \Vert f \Vert^2_{\dot H^{\frac{5}{2}}}+\left\Vert f \right\Vert_{\dot H^{\frac{7}{2}}}\right)
  \end{eqnarray*}
  Then, by interpolation
  \begin{eqnarray*}
\mathcal{U}_1 &\lesssim&\Vert f \Vert_{\dot H^4}\Vert f \Vert^{\frac{1}{5}}_{\dot H^4}\Vert f \Vert^{\frac{4}{5}}_{\dot H^{\frac{3}{2}}}\left(\Vert f \Vert^{\frac{3}{10}}_{\dot H^4}\Vert f \Vert^{\frac{7}{10}}_{\dot H^{\frac{3}{2}}} \ \Vert f \Vert^{\frac{1}{2}}_{\dot H^4}\Vert f \Vert^{\frac{1}{2}}_{\dot H^{\frac{3}{2}}}+  \Vert f \Vert^{\frac{4}{5}}_{\dot H^{4}}
\Vert f \Vert^{\frac{6}{5}}_{\dot H^{\frac{3}{2}}}+\Vert f \Vert^{\frac{4}{5}}_{\dot H^{4}} \Vert f \Vert^{\frac{1}{5}}_{\dot H^{\frac{3}{2}}}\right)
  \end{eqnarray*}
  
  Finally, one finds
  
   \begin{eqnarray*}
\mathcal{U}_1 &\lesssim&\Vert f \Vert^2_{\dot H^4}\left(\Vert f \Vert^2_{\dot H^{\frac{3}{2}}}+\Vert f \Vert_{\dot H^{\frac{3}{2}}}\right)
  \end{eqnarray*}
  
  \subsection{Estimate of  $\mathcal{U}_{2,dissip}$}
This term is important as it will give the elliptic component plus some controlled remainders. More precisely, we have

\begin{eqnarray*}
  \mathcal{U}_2&=&\ - \ \int \Lambda^{4}f  \ \Lambda\left( \Lambda^{3} f \  \int_{0}^{\infty} e^{-\sigma}\cos(\sigma f_x )\cos(\arctan(f_{x}))\ d \sigma \right)  \ dx    \\
  \end{eqnarray*}
 Set $\Psi(f)(x):=\int_{0}^{\infty} e^{-\sigma}\cos(\sigma f_x )\cos(\arctan(f_{x}))\ d \sigma.$ Then,
 \begin{eqnarray*}
  \mathcal{U}_2&=&\ - \ \int \vert \Lambda^{4}f \vert^2  \  \Psi(f) \ dx - \int  \Lambda^{4}f \ \Lambda^{3} f \  \Lambda \Psi(f)  \ dx + c \int \int\Lambda^{4}f  \frac{\delta_{\beta}\Lambda^{3} f \ \delta_{\beta}\Psi(f)}{\vert \beta \vert^2} \ d\beta \ dx    \\
  \end{eqnarray*}
  
 On notices that the first term is the elliptic component, namely we have that 
   \begin{eqnarray*}
\mathcal{V}(f)= - \ \int \vert \Lambda^{4}f \vert^2  \  \Psi(f) \ dx=- \ \int \frac{\vert \Lambda^{4}f \vert^2}{(1+(f_x)^2)^{\frac{3}{2}}}  \   \ dx
  \end{eqnarray*}
  The second term of $\mathcal{U}_2$ is easy to estimate as one may write, using Holder and then Sobolev embedding
     \begin{eqnarray*}
  \int  \Lambda^{4}f \ \Lambda^{3} f \  \mathcal{H} \partial_x\Psi(f)  \ dx &\lesssim& \Vert f \Vert_{\dot H^4} \Vert f \Vert_{\dot H^{\frac{13}{4}}} \Vert \partial_x\Psi(f) \Vert_{L^4} \\
  &\lesssim& \Vert f \Vert_{\dot H^4} \Vert f \Vert_{\dot H^{\frac{13}{4}}} \Vert f_{xx} \Vert_{L^4} \\
  &\lesssim& \Vert f \Vert_{\dot H^4} \Vert f \Vert_{\dot H^{\frac{13}{4}}} \Vert f \Vert_{\dot H^{\frac{9}{4}}} 
  \end{eqnarray*}
 Now, it suffices to interpolate 
  \begin{eqnarray*}
  \int  \Lambda^{4}f \ \Lambda^{3} f \  \mathcal{H}\partial_x\Psi(f)  \ dx &\lesssim&\Vert f \Vert_{\dot H^4} \Vert f \Vert^{\frac{7}{10}}_{\dot H^4}\Vert f \Vert^{\frac{3}{10}}_{\dot H^{\frac{3}{2}}} \Vert f \Vert^{\frac{3}{10}}_{\dot H^4}\Vert f \Vert^{\frac{7}{10}}_{\dot H^{\frac{3}{2}}} \\
  &\lesssim&\Vert f \Vert^2_{\dot H^4}\Vert f \Vert_{\dot H^{\frac{3}{2}}}
   \end{eqnarray*}
   It remains to estimate the third and last term of $\mathcal{U}_2$, that is
     \begin{eqnarray*}
   \int \int\Lambda^{4}f  \frac{\delta_{\beta}\Lambda^{3} f \ \delta_{\beta}\Psi(f)}{\vert \beta \vert^2} \ d\beta \ dx &\lesssim&\Vert f \Vert_{\dot H^4} \left(\int   \frac{\Vert \delta_{\beta}\Lambda^{3} f \Vert^2_{L^4}}{\vert \beta \vert^2} \ d\beta\right)^{\frac{1}{2}} \left(\int   \frac{ \Vert\delta_{\beta}\Psi(f) \Vert^2_{L^4}}{\vert \beta \vert^2} \ d\beta\right)^{\frac{1}{2}}
      \end{eqnarray*}
 As $\Psi$ is a bounded Lipschitz function, one finds
\begin{eqnarray*}
   \int \int\Lambda^{4}f  \frac{\delta_{\beta}\Lambda^{3} f \ \delta_{\beta}\Psi(f)}{\vert \beta \vert^2} \ d\beta \ dx &\lesssim&\Vert f \Vert_{\dot H^4} \left(\int   \frac{\Vert \delta_{\beta}\Lambda^{3} f \Vert^2_{L^4}}{\vert \beta \vert^2} \ d\beta\right)^{\frac{1}{2}} \left(\int   \frac{ \Vert\delta_{\beta}f_x\Vert^2_{L^4}}{\vert \beta \vert^2} \ d\beta\right)^{\frac{1}{2}}\\
   &\lesssim&\Vert f \Vert_{\dot H^4} \left(\int   \frac{\Vert \delta_{\beta}\Lambda^{3} f \Vert^2_{L^4}}{\vert \beta \vert^2} \ d\beta\right)^{\frac{1}{2}} \left(\int   \frac{ \Vert\delta_{\beta}f_x\Vert^2_{L^4}}{\vert \beta \vert^2} \ d\beta\right)^{\frac{1}{2}}\\
   &\lesssim&\Vert f \Vert_{\dot H^4} \Vert f \Vert_{\dot B^{\frac{7}{2}}_{4,2}} \Vert f \Vert_{\dot B^{\frac{3}{2}}_{4,2}} \\
   &\lesssim&\Vert f \Vert_{\dot H^4}\Vert f \Vert_{\dot H^{\frac{15}{4}}} \Vert f \Vert_{\dot H^{\frac{7}{4}}}
      \end{eqnarray*}
By interpolation,
\begin{eqnarray*}
   \int \int\Lambda^{4}f  \frac{\delta_{\beta}\Lambda^{3} f \ \delta_{\beta}\Psi(f)}{\vert \beta \vert^2} \ d\beta \ dx &\lesssim&\Vert f \Vert_{\dot H^4}  \Vert f \Vert^{\frac{9}{10}}_{\dot H^4}\Vert f \Vert^{\frac{1}{10}}_{\dot H^{\frac{3}{2}}} \Vert f \Vert^{\frac{1}{10}}_{\dot H^4}\Vert f \Vert^{\frac{9}{10}}_{\dot H^{\frac{3}{2}}} \\
   &\lesssim&\Vert f \Vert^2_{\dot H^4}\Vert f \Vert_{\dot H^{\frac{3}{2}}}
\end{eqnarray*}
Finally, we have proved that
\begin{eqnarray*}
\mathcal{U}_2 &\lesssim&- \ \int \frac{\vert \Lambda^{4}f \vert^2}{(1+(f_x)^2)^{\frac{3}{2}}}    \ dx + \Vert f \Vert^2_{\dot H^4}\left(\Vert f \Vert^2_{\dot H^{\frac{3}{2}}}+\Vert f \Vert_{\dot H^{\frac{3}{2}}}\right)
  \end{eqnarray*}
  
   \subsection{Estimate of  $\mathcal{U}_{2,rem}$}
   \begin{eqnarray*}
  \mathcal{U}_{2,rem}&=&\int \Lambda^{4}f \ \Lambda\mathcal{H}f_{xx} \ \partial_{x}\left(\int_{0}^{\infty} e^{-\sigma}\cos(\sigma f_x )\cos(\arctan(f_{x}))  \ d \sigma \right) \ dx \\
  &\lesssim& \Vert f \Vert_{\dot H^4} \Vert \partial^3_x f \Vert_{L^4} 
  \Vert \partial^2_x f \Vert_{L^4}
  \end{eqnarray*}
  where we used that,
  $$
  \left \Vert \int_{0}^{\infty} e^{-\sigma}\cos(\sigma f_x )\cos(\arctan(f_{x}))  \ d \sigma \right\Vert_{\dot W^{1,4}} \lesssim \Vert  \partial^2_x f  \Vert_{L^4}.
  $$
  Therefore, using Sobolev embedding and interpolation we find that
  \begin{eqnarray*}
  \mathcal{U}_{2,rem}&\lesssim& \Vert f \Vert_{\dot H^4} \Vert f \Vert_{\dot H^\frac{13}{4}} \Vert f \Vert_{\dot H^\frac{9}{4}} \\
  &\lesssim& \Vert f \Vert_{\dot H^4} \Vert f \Vert^{\frac{7}{10}}_{\dot H^4} \Vert f \Vert^{\frac{3}{10}}_{\dot H^{\frac{3}{2}}} \Vert f \Vert^{\frac{3}{10}}_{\dot H^4} \Vert f \Vert^{\frac{7}{10}}_{\dot H^{\frac{3}{2}}}  \\
  &\lesssim& \Vert f \Vert^2_{\dot H^4} \Vert f \Vert_{\dot H^{\frac{3}{2}}}
    \end{eqnarray*}
    
    \subsection{Estimate of  $\mathcal{U}_3$}
     \begin{eqnarray*}
\mathcal{U}_3&=&- \frac{1}{4\pi}\int \Lambda^{4}\mathcal{H}f \ \partial_x \int\int_{0}^{\infty}\int_{0}^{\infty} \ e^{-\gamma-\sigma} \ S_{\alpha} f \left(\sin(\gamma\Delta_\alpha f)+\sin(\gamma\bar\Delta_\alpha f)\right) \\
&& \ \times \  \partial_x \left(\frac{\partial_x^2 \tau_{\alpha} f}{\alpha} \cos(\sigma \tau_{\alpha} f_x )\cos(\arctan(\tau_{\alpha}f_{x}))\right) \ d\alpha \ d\gamma \ d\sigma \ dx \\
&=&- \frac{1}{4\pi}\int \Lambda^{4}\mathcal{H}f \  \int\int_{0}^{\infty}\int_{0}^{\infty} \ e^{-\gamma-\sigma} \ S_{\alpha} f \ \ S_{\alpha} f_x \ \cos(\frac{1}{2}\gamma S_\alpha f)\cos(\frac{1}{2}\gamma D_\alpha f) \\
&& \ \times \  \partial_x \left(\frac{\partial_x^2 \tau_{\alpha} f}{\alpha} \cos(\sigma \tau_{\alpha} f_x )\cos(\arctan(\tau_{\alpha}f_{x}))\right) \ d\alpha \ d\gamma \ d\sigma \ dx \\
&+& \frac{1}{4\pi}\int \Lambda^{4}\mathcal{H}f \  \int\int_{0}^{\infty}\int_{0}^{\infty} \ e^{-\gamma-\sigma} \ S_{\alpha} f \ \ D_{\alpha} f_x \ \sin(\frac{1}{2}\gamma S_\alpha f)\sin(\frac{1}{2}\gamma D_\alpha f) \\
&& \ \times \  \partial_x \left(\frac{\partial_x^2 \tau_{\alpha} f}{\alpha} \cos(\sigma \tau_{\alpha} f_x )\cos(\arctan(\tau_{\alpha}f_{x}))\right) \ d\alpha \ d\gamma \ d\sigma \ dx \\
&-& \frac{1}{4\pi}\int \Lambda^{4}\mathcal{H}f \  \int\int_{0}^{\infty}\int_{0}^{\infty} \ e^{-\gamma-\sigma} \ S_{\alpha} f \left(\sin(\gamma\Delta_\alpha f)+\sin(\gamma\bar\Delta_\alpha f)\right) \\
&& \ \times \  \partial^2_x \left(\frac{\partial_x^2 \tau_{\alpha} f}{\alpha} \cos(\sigma \tau_{\alpha} f_x )\cos(\arctan(\tau_{\alpha}f_{x}))\right) \ d\alpha \ d\gamma \ d\sigma \ dx \\
&:=& \mathcal{U}_{3,1}+\mathcal{U}_{3,2}+\mathcal{U}_{3,3}
 \end{eqnarray*}
 
 $\bullet$ Estimate of $\mathcal{U}_{3,1}$
 
 \begin{eqnarray*}
\mathcal{U}_{3,1}&\lesssim&\Vert f \Vert_{\dot H^4} \left( \Vert f \Vert_{\dot H^3} + \Vert f_{xx} \Vert^2_{L^4} \right)  \int \Vert S_{\alpha} f \Vert_{L^\infty} \ \Vert S_{\alpha} f_x \Vert_{L^\infty} \frac{1}{\vert \alpha \vert} \ d\alpha \\
&\lesssim&\Vert f \Vert_{\dot H^4} \left( \Vert f \Vert_{\dot H^3} + \Vert f \Vert^2_{\dot H^{\frac{9}{4}}} \right) \sup_{\alpha \in \mathbb{R}}\Vert S_{\alpha} f \Vert_{L^\infty}  \int  \frac{\Vert s_{\alpha} f_x \Vert_{L^\infty}}{\vert\alpha\vert} \frac{1}{\vert \alpha \vert} \ d\alpha \\
&\lesssim&\Vert f \Vert_{\dot H^4} \left( \Vert f \Vert_{\dot H^3} +\Vert f \Vert^2_{\dot H^{\frac{9}{4}}} \right) \Vert  f \Vert_{\dot B^1_{\infty,\infty}}  \Vert  f \Vert_{\dot B^2_{\infty,1}} \\
&\lesssim&\Vert f \Vert_{\dot H^4} \Vert f \Vert_{\dot H^{\frac{3}{2}}} \left( \Vert f \Vert_{\dot H^3} + \Vert f \Vert^2_{\dot H^{\frac{9}{4}}} \right)   \Vert  f \Vert_{\dot B^2_{\infty,1}} \\
 \end{eqnarray*}
where we used that $\dot H^{\frac{3}{2}} \hookrightarrow \dot B^1_{\infty,\infty}$ and interpolation. Then, since we have
$$
  \Vert  f \Vert_{\dot B^2_{\infty,1}} \lesssim \Vert  f\Vert^{\frac{3}{5}}_{\dot B^{1}_{\infty,\infty}}   \Vert  f \Vert^{\frac{2}{5}}_{\dot B^{7/2}_{\infty,\infty}}\lesssim \Vert  f\Vert^{\frac{3}{5}}_{\dot H^{\frac{3}{2}}}   \Vert  f \Vert^{\frac{2}{5}}_{\dot H^4}
$$
Therefore, 
\begin{eqnarray*}
\mathcal{U}_{3,1}&\lesssim&\Vert f \Vert_{\dot H^4} \Vert f \Vert_{\dot H^{\frac{3}{2}}} \left( \Vert f \Vert^{\frac{3}{5}}_{\dot H^4} \Vert f \Vert^{\frac{2}{5}}_{\dot H^\frac{3}{2}} + \Vert f \Vert^{\frac{3}{5}}_{\dot H^4}\Vert f \Vert^{\frac{7}{5}}_{\dot H^{\frac{3}{2}}} \right) \Vert  f \Vert^{\frac{2}{5}}_{\dot H^4}  \Vert  f\Vert^{\frac{3}{5}}_{\dot H^{\frac{3}{2}}} \\
&\lesssim& \Vert f \Vert^2_{\dot H^4} \left(  \Vert f \Vert^2_{\dot H^\frac{3}{2}} + \Vert f \Vert^3_{\dot H^\frac{3}{2}} \right)
 \end{eqnarray*}

  $\bullet$ Estimate of $\mathcal{U}_{3,2}$\\

  To estimate this term, we use need to make appear a quadratic term in \(S\) and then balance the regularity

  \begin{eqnarray*}
  \mathcal{U}_{3,2}&\lesssim& \frac{1}{4\pi}\int \Lambda^{4}\mathcal{H}f \  \int\int_{0}^{\infty}\int_{0}^{\infty} \ e^{-\gamma-\sigma} \ S_{\alpha} f \ \ D_{\alpha} f_x \ \sin(\frac{1}{2}\gamma S_\alpha f)\sin(\frac{1}{2}\gamma D_\alpha f) \\
&& \ \times \  \partial_x \left(\frac{\partial_x^2 \tau_{\alpha} f}{\alpha} \cos(\sigma \tau_{\alpha} f_x )\cos(\arctan(\tau_{\alpha}f_{x}))\right) \ d\alpha \ d\gamma \ d\sigma \ dx \\
&\lesssim& \Vert f \Vert_{\dot H^4} \left(\Vert f \Vert_{\dot H^3}+\Vert f \Vert^2_{\dot H^2}\right) \int \frac{\Vert s_\alpha f \Vert^2_{L^{\infty}} \Vert d_\alpha f_x \Vert_{L^{\infty}}}{\vert\alpha \vert^4} \ d\alpha \\
&\lesssim& \Vert f \Vert_{\dot H^4} \left(\Vert f \Vert_{\dot H^3} + \Vert f_{xx} \Vert^2_{L^4}\right) \left(\int \frac{\Vert s_\alpha f \Vert^4_{L^{\infty}} }{\vert\alpha \vert^6} \ d\alpha \right)^{\frac{1}{2}} \left(\int \frac{\Vert d_\alpha f_x \Vert^2_{L^{\infty}} }{\vert\alpha \vert^2} \ d\alpha \right)^{\frac{1}{2}} \\
&\lesssim& \Vert f \Vert_{\dot H^4} \left(\Vert f \Vert_{\dot H^3}+\Vert f \Vert^2_{\dot H^\frac{9}{4}}\right) \Vert f \Vert^2_{\dot B^{\frac{5}{4}}_{\infty,4}} \Vert f \Vert_{\dot B^{\frac{3}{2}}_{\infty,2}} \\
&\lesssim& \Vert f \Vert_{\dot H^4} \left(\Vert f \Vert_{\dot H^3}+\Vert f \Vert^2_{\dot H^\frac{9}{4}}\right) \Vert f \Vert^2_{\dot H^{\frac{7}{4}}} \Vert f \Vert_{\dot H^2} \\
&\lesssim& \Vert f \Vert_{\dot H^4} \left( \Vert f \Vert^{\frac{3}{5}}_{\dot H^4} \Vert f \Vert^{\frac{2}{5}}_{\dot H^\frac{3}{2}} + \Vert f \Vert^{\frac{3}{5}}_{\dot H^4}\Vert f \Vert^{\frac{7}{5}}_{\dot H^{\frac{3}{2}}} \right) \Vert f \Vert^2_{\dot H^{\frac{7}{4}}} \Vert f \Vert^{\frac{1}{5}}_{\dot H^4}\Vert f \Vert^{\frac{4}{5}}_{\dot H^{\frac{3}{2}}}
\end{eqnarray*}
Therefore, since $\Vert f \Vert^2_{\dot H^{\frac{7}{4}}}\lesssim \Vert f \Vert^{\frac{9}{5}}_{\dot H^{\frac{3}{2}}} \Vert f \Vert^{\frac{1}{5}}_{\dot H^4} $

\begin{eqnarray*}
  \mathcal{U}_{3,2}&\lesssim&\Vert f \Vert_{\dot H^4} \left( \Vert f \Vert^{\frac{3}{5}}_{\dot H^4} \Vert f \Vert^{\frac{2}{5}}_{\dot H^\frac{3}{2}} + \Vert f \Vert^{\frac{3}{5}}_{\dot H^4}\Vert f \Vert^{\frac{7}{5}}_{\dot H^{\frac{3}{2}}} \right)  \Vert f \Vert^{\frac{13}{5}}_{\dot H^{\frac{3}{2}}} \Vert f \Vert^{\frac{2}{5}}_{\dot H^4} \\
  &\lesssim&\Vert f \Vert_{\dot H^4} \left(\Vert f \Vert^3_{\dot H^{\frac{3}{2}}} + \Vert f \Vert^4_{\dot H^{\frac{3}{2}}} \right)
  \end{eqnarray*}
  
    $\bullet$ Estimate of $\mathcal{U}_{3,3}$
    
   \begin{eqnarray*} 
    \mathcal{U}_{3,3}&=&- \frac{1}{4\pi}\int \Lambda^{4}\mathcal{H}f \  \int\int_{0}^{\infty}\int_{0}^{\infty} \ e^{-\gamma-\sigma} \ S_{\alpha} f \left(\sin(\gamma\Delta_\alpha f)+\sin(\gamma\bar\Delta_\alpha f)\right)\frac{1}{\alpha} \\
&& \ \times \  \partial^2_x \left({\partial_x^2 \tau_{\alpha} f} \cos(\sigma \tau_{\alpha} f_x )\cos(\arctan(\tau_{\alpha}f_{x}))\right) \ d\alpha \ d\gamma \ d\sigma \ dx \\
\end{eqnarray*}
It not difficult to check that,
\begin{eqnarray*}
\left\vert \partial^2_x \left(\partial_x^2 \tau_{\alpha} f \cos(\sigma \tau_{\alpha} f_x )\cos(\arctan(\tau_{\alpha}f_{x}))\right) \right\vert &\lesssim&  \left\vert{\partial_x^4 \tau_{\alpha} f} +  U(\sigma){\partial_x^2 \tau_{\alpha} f} \left((\partial^2_{x}f)^2+ \partial^3_{x}f \right) \right\vert
\end{eqnarray*}
where $U(\sigma):=1+\sigma+\sigma^2.$ Therefore, 
 \begin{eqnarray*} 
    \mathcal{U}_{3,3}&\lesssim&\Vert f \Vert_{\dot H^4}\ \left(\Vert f \Vert_{\dot H^4} + \Vert \partial^2_{x}f\Vert^3_{L^6}+\Vert \partial^2_{x}f\Vert_{L^4}\Vert \partial^3_{x}f\Vert_{L^4}\right)   \\
&& \ \times  \int\int_{0}^{\infty}\int_{0}^{\infty} \ \gamma U(\sigma) e^{-\gamma-\sigma} \ \frac{\Vert s_{\alpha} f \Vert^2_{L^{\infty}}}{\vert \alpha \vert^3}\  \ d\alpha \ d\gamma \ d\sigma  \\
&\lesssim&\Vert f \Vert^2_{\dot H^4}\ \Vert f \Vert^2_{\dot B^1_{\infty,2}}
  + \Vert f \Vert_{\dot H^4}\Vert \partial^2_{x}f\Vert^3_{L^6}\Vert f \Vert^2_{\dot B^1_{\infty,2}}
+\Vert f \Vert_{\dot H^4}\Vert \partial^2_{x}f\Vert_{L^4}\Vert \partial^3_{x}f\Vert_{L^4} \Vert f \Vert^2_{\dot B^1_{\infty,2}}
\end{eqnarray*}
As $\dot H^{\frac{1}{4}} \hookrightarrow L^4$ and $\dot H^{\frac{1}{3}} \hookrightarrow L^6$ and $\dot H^{\frac{3}{2}} \hookrightarrow \dot B^1_{\infty,2}$, we find

 \begin{eqnarray*} 
    \mathcal{U}_{3,3}&\lesssim& \Vert f \Vert^2_{\dot H^4} \ \Vert f \Vert^2_{\dot H^{\frac{3}{2}}} + \Vert f \Vert_{\dot H^4}  \Vert f \Vert^3_{\dot H^{\frac{7}{3}}}\Vert f \Vert^2_{\dot H^{\frac{3}{2}}} + \Vert f \Vert_{\dot H^4}\Vert f \Vert_{\dot H^{\frac{9}{4}}}\Vert f \Vert_{\dot H^{\frac{13}{4}}}\Vert f \Vert^2_{\dot H^{\frac{3}{2}}}
\end{eqnarray*}
Since  $\Vert f \Vert^3_{\dot H^{\frac{7}{3}}}\lesssim \Vert f \Vert_{\dot H^4}  \Vert f \Vert^2_{\dot H^{\frac{3}{2}}}$ and $\Vert f \Vert_{\dot H^{\frac{9}{4}}}\Vert f \Vert_{\dot H^{\frac{13}{4}}}\lesssim \Vert f \Vert_{\dot H^4}  \Vert f \Vert_{\dot H^{\frac{3}{2}}}$, then
\begin{eqnarray*} 
    \mathcal{U}_{3,3}&\lesssim& \Vert f \Vert^2_{\dot H^4} \ \Vert f \Vert^2_{\dot H^{\frac{3}{2}}} + \Vert f \Vert^2_{\dot H^4}  \Vert f \Vert^4_{\dot H^{\frac{3}{2}}} + \Vert f \Vert^2_{\dot H^4}\Vert f \Vert^3_{\dot H^{\frac{3}{2}}}
\end{eqnarray*}
 
Garthering all the estimates of the $ \mathcal{U}_{3,i}, i=1,2,3$ we have obtained
\begin{eqnarray*} 
    \mathcal{U}_{3}&\lesssim& \Vert f \Vert^2_{\dot H^4} \left( \Vert f \Vert^2_{\dot H^{\frac{3}{2}}}  + \Vert f \Vert^3_{\dot H^{\frac{3}{2}}} +  \Vert f \Vert^4_{\dot H^{\frac{3}{2}}} \right)
\end{eqnarray*}

\subsection{Estimate of $\mathcal{U}_{4}$}

We want to estimate 
\begin{eqnarray*}
\mathcal{U}_{4}&=&- \frac{1}{4\pi}\int \Lambda^{4}\mathcal{H}f \ \partial_x \int_{0}^{\infty}\int_{0}^{\infty}\int_{0}^{\alpha} \ e^{-\gamma-\sigma} \ \frac{1}{\alpha} s_{\eta} f_{x}  \left(\sin(\gamma\Delta_\alpha f)-\sin(\gamma\bar\Delta_\alpha f)\right) \\  
&& \ \times \ \partial_x \left(\frac{\partial_x^2 \tau_{\alpha} f}{\alpha} \cos(\sigma \tau_{\alpha} f_x )\cos(\arctan(\tau_{\alpha} f_{x}))\right) \ d\eta \ d\alpha \ d\gamma \ d\sigma\ dx. \\
\end{eqnarray*}
We notice that this term is somehow analogous to $T_4$, indeed, we may integrate by parts in the space variable, we get

\begin{eqnarray*}
\mathcal{U}_{4}&=&\frac{1}{4\pi}\int \Lambda^{4}\mathcal{H}f_x \  \int_{0}^{\infty}\int_{0}^{\infty}\int_{0}^{\alpha} \ e^{-\gamma-\sigma} \ \frac{1}{\alpha} s_{\eta} f_{x}  \left(\sin(\gamma\Delta_\alpha f)-\sin(\gamma\bar\Delta_\alpha f)\right) \\  
&& \ \times \ \partial_x \left(\frac{\partial_x^2 \tau_{\alpha} f}{\alpha} \cos(\sigma \tau_{\alpha} f_x )\cos(\arctan(\tau_{\alpha} f_{x}))\right) \ d\eta \ d\alpha \ d\gamma \ d\sigma\ dx \\
\end{eqnarray*}

Up to interchanging the role of $\Lambda^3 f$ and $\Lambda^{4}\mathcal{H}f_x$, we deduce that the identity obtained for $T_4$ can be easily extended to $\mathcal{U}_{4}$ as $\Lambda^3 f$ did not play any specific role in the proof. We just systematically put  this term in $L^2$ in order to use the maximal expected regularity from the dissipation. Therefore, we readily deduce  the following decomposition for $\mathcal{U}_{4}$

 \begin{eqnarray*}
  \mathcal{U}_{4}&=&\frac{1}{2\pi} \int \int\int_{0}^{\infty}\int_{0}^{\infty} \ e^{-\gamma-\sigma} \ \Lambda^{4}\mathcal{H}f_x \ \partial^2_{\alpha}\left[\frac{1}{\alpha}\sin(\frac{\gamma}{2}D_\alpha f)\cos(\frac{\gamma}{2}S_\alpha f)\frac{1}{\alpha} \int_{0}^{\alpha} s_{\eta} f_{x}    \ d\eta \right]\\
&& \  \times\ {  \delta_\alpha f_x } \cos(\sigma \tau_{\alpha} f_x )\cos(\arctan(\tau_{\alpha} f_{x})) \ d\eta  \ d\gamma \ d\sigma \ dx \ d\alpha \\
&+&\frac{1}{2\pi} \int \int\int_{0}^{\infty}\int_{0}^{\infty} \ \sigma e^{-\gamma-\sigma} \ \Lambda^{4}\mathcal{H}f_x \ \partial_{\alpha}\left[\frac{1}{\alpha}\sin(\frac{\gamma}{2}D_\alpha f)\cos(\frac{\gamma}{2}S_\alpha f)\frac{1}{\alpha} \int_{0}^{\alpha} s_{\eta} f_{x}    \ d\eta \right]\\
&& \  \times\ {  \delta_\alpha f_x } \left(\partial_{\alpha}\tau_{\alpha}f_x\right) \sin(\sigma \tau_{\alpha} f_x )\cos(\arctan(\tau_{\alpha} f_{x})) \ d\eta  \ d\gamma \ d\sigma \ dx \ d\alpha\\
&+&\frac{1}{2\pi} \int \int\int_{0}^{\infty}\int_{0}^{\infty} \ \sigma e^{-\gamma-\sigma} \ \Lambda^{4}\mathcal{H}f_x \ \partial_{\alpha}\left[\frac{1}{\alpha}\sin(\frac{\gamma}{2}D_\alpha f)\cos(\frac{\gamma}{2}S_\alpha f)\frac{1}{\alpha} \int_{0}^{\alpha} s_{\eta} f_{x}    \ d\eta \right]\\
&& \  \times\ {  \delta_\alpha f_x } \frac{\partial_{\alpha}\tau_{\alpha}f_x}{1+(\tau_{\alpha} f_{x})^2} \sin(\sigma \tau_{\alpha} f_x )\cos(\arctan(\tau_{\alpha} f_{x})) \ d\eta  \ d\gamma \ d\sigma \ dx \ d\alpha \\
&:=&\mathcal{U}_{4,1}+\mathcal{U}_{4,2}+\mathcal{U}_{4,3}
   \end{eqnarray*}
   Then, we integrate by parts as the term $\Lambda^{4}\mathcal{H}f_x$ is too singular so we need to balance the regularity in $x$. Remark that it is clearly not enough to integrate by parts in $x$.
    \begin{eqnarray*}
  \mathcal{U}_{4}&=&\frac{1}{2\pi} \int \int\int_{0}^{\infty}\int_{0}^{\infty} \ e^{-\gamma-\sigma} \ \Lambda^{4}\mathcal{H}f \ \partial_{x}\partial^2_{\alpha}\left[\frac{1}{\alpha}\sin(\frac{\gamma}{2}D_\alpha f)\cos(\frac{\gamma}{2}S_\alpha f)\frac{1}{\alpha} \int_{0}^{\alpha} s_{\eta} f_{x}    \ d\eta \right]\\
&& \  \times\ {  \delta_\alpha f_x } \cos(\sigma \tau_{\alpha} f_x )\cos(\arctan(\tau_{\alpha} f_{x})) \ d\eta  \ d\gamma \ d\sigma \ dx \ d\alpha \\
&+&\frac{1}{2\pi} \int \int\int_{0}^{\infty}\int_{0}^{\infty} \ e^{-\gamma-\sigma} \ \Lambda^{4}\mathcal{H}f \ \partial^2_{\alpha}\left[\frac{1}{\alpha}\sin(\frac{\gamma}{2}D_\alpha f)\cos(\frac{\gamma}{2}S_\alpha f)\frac{1}{\alpha} \int_{0}^{\alpha} s_{\eta} f_{x}    \ d\eta \right]\\
&& \  \times\ \partial_x \left( {  \delta_\alpha f_x } \cos(\sigma \tau_{\alpha} f_x )\cos(\arctan(\tau_{\alpha} f_{x})) )\right) \ d\eta  \ d\gamma \ d\sigma \ dx \ d\alpha \\
&+&\frac{1}{2\pi} \int \int\int_{0}^{\infty}\int_{0}^{\infty} \ \sigma e^{-\gamma-\sigma} \ \Lambda^{4}\mathcal{H}f \ \partial_{x}\partial_{\alpha}\left[\frac{1}{\alpha}\sin(\frac{\gamma}{2}D_\alpha f)\cos(\frac{\gamma}{2}S_\alpha f)\frac{1}{\alpha} \int_{0}^{\alpha} s_{\eta} f_{x}    \ d\eta \right]\\
&& \  \times\ {  \delta_\alpha f_x } \left(\partial_{\alpha}\tau_{\alpha}f_x\right) \sin(\sigma \tau_{\alpha} f_x )\cos(\arctan(\tau_{\alpha} f_{x})) \ d\eta  \ d\gamma \ d\sigma \ dx \ d\alpha\\
&+&\frac{1}{2\pi} \int \int\int_{0}^{\infty}\int_{0}^{\infty} \ \sigma e^{-\gamma-\sigma} \ \Lambda^{4}\mathcal{H}f \ \partial_{\alpha}\left[\frac{1}{\alpha}\sin(\frac{\gamma}{2}D_\alpha f)\cos(\frac{\gamma}{2}S_\alpha f)\frac{1}{\alpha} \int_{0}^{\alpha} s_{\eta} f_{x}    \ d\eta \right]\\
&& \  \times\  \partial_{x}\left({  \delta_\alpha f_x } \left(\partial_{\alpha}\tau_{\alpha}f_x\right) \sin(\sigma \tau_{\alpha} f_x )\cos(\arctan(\tau_{\alpha} f_{x})) \right) \ d\eta  \ d\gamma \ d\sigma \ dx \ d\alpha\\
&+&\frac{1}{2\pi} \int \int\int_{0}^{\infty}\int_{0}^{\infty} \ \sigma e^{-\gamma-\sigma} \ \Lambda^{4}\mathcal{H}f \ \partial_{x}\partial_{\alpha}\left[\frac{1}{\alpha}\sin(\frac{\gamma}{2}D_\alpha f)\cos(\frac{\gamma}{2}S_\alpha f)\frac{1}{\alpha} \int_{0}^{\alpha} s_{\eta} f_{x}    \ d\eta \right]\\
&& \  \times\ {  \delta_\alpha f_x } \frac{\partial_{\alpha}\tau_{\alpha}f_x}{1+(\tau_{\alpha} f_{x})^2} \sin(\sigma \tau_{\alpha} f_x )\cos(\arctan(\tau_{\alpha} f_{x})) \ d\eta  \ d\gamma \ d\sigma \ dx \ d\alpha \\
&+&\frac{1}{2\pi} \int \int\int_{0}^{\infty}\int_{0}^{\infty} \ \sigma e^{-\gamma-\sigma} \ \Lambda^{4}\mathcal{H}f \ \partial_{\alpha}\left[\frac{1}{\alpha}\sin(\frac{\gamma}{2}D_\alpha f)\cos(\frac{\gamma}{2}S_\alpha f)\frac{1}{\alpha} \int_{0}^{\alpha} s_{\eta} f_{x}    \ d\eta \right]\\
&& \  \times\ \partial_{x}({  \delta_\alpha f_x } \frac{\partial_{\alpha}\tau_{\alpha}f_x}{1+(\tau_{\alpha} f_{x})^2} \sin(\sigma \tau_{\alpha} f_x )\cos(\arctan(\tau_{\alpha} f_{x}))) \ d\eta  \ d\gamma \ d\sigma \ dx \ d\alpha \\
&:=&\sum_{i=1}^6 \mathcal{U}_{4,i}
   \end{eqnarray*}
   In other words, this term is somehow similar to $\mathcal{T}_4$, the only difference is that one has to compute one derivative more in space in the nonlinearity. Unlike the derivatives in $\alpha$, the derivative in $x$ does preserve the symmetric structure and in general the finite difference. Hence, the only changes will be in the use of homogeneous Sobolev or Besov embeddings but as the maximal regularity is $\dot H^4$, one will have to interpolate between $\dot H^{\frac{3}{2}}$ and $\dot H^4$. It is therefore not difficult to adapt the estimates we obtained for   ${T}_4$ to the case of $\mathcal{U}_4$. 
   \begin{eqnarray*}
 \mathcal{U}_{4,1}&=&\frac{1}{\pi} \int \int\int_{0}^{\infty}\int_{0}^{\infty} \ e^{-\gamma-\sigma} \ \Lambda^{4}\mathcal{H}f \ \frac{1}{\alpha^3} \partial_x \left(\sin(\frac{\gamma}{2}D_\alpha f)\cos(\frac{\gamma}{2}S_\alpha f)\frac{1}{\alpha} \int_{0}^{\alpha} s_{\eta} f_{x}    \ d\eta\right) \\
&& \  \times\ {  \delta_\alpha f_x } \cos(\sigma \tau_{\alpha} f_x )\cos(\arctan(\tau_{\alpha} f_{x}))   \ d\gamma \ d\sigma \ dx \ d\alpha \\
&+&\frac{1}{4\pi} \int \int\int_{0}^{\infty}\int_{0}^{\infty} \ \gamma e^{-\gamma-\sigma} \ \Lambda^{4}\mathcal{H}f \ \partial_x \left(\partial^2_\alpha D_\alpha f \cos(\frac{\gamma}{2}D_\alpha f)\cos(\frac{\gamma}{2}S_\alpha f)\frac{1}{\alpha^2} \int_{0}^{\alpha} s_{\eta} f_{x}    \ d\eta\right) \\
&& \  \times\ {  \delta_\alpha f_x } \cos(\sigma \tau_{\alpha} f_x )\cos(\arctan(\tau_{\alpha} f_{x}))   \ d\gamma \ d\sigma \ dx \ d\alpha \\
&-&\frac{1}{8\pi} \int \int\int_{0}^{\infty}\int_{0}^{\infty} \ \gamma^2 e^{-\gamma-\sigma} \ \Lambda^{4}\mathcal{H}f \ \partial_x \left((\partial_\alpha D_\alpha f)^2 \sin(\frac{\gamma}{2}D_\alpha f)\cos(\frac{\gamma}{2}S_\alpha f)\frac{1}{\alpha^2} \int_{0}^{\alpha} s_{\eta} f_{x}    \ d\eta\right) \\
&& \  \times\ {  \delta_\alpha f_x } \cos(\sigma \tau_{\alpha} f_x )\cos(\arctan(\tau_{\alpha} f_{x}))   \ d\gamma \ d\sigma \ dx \ d\alpha \\
&-&\frac{1}{4\pi} \int \int\int_{0}^{\infty}\int_{0}^{\infty} \ \gamma e^{-\gamma-\sigma} \ \Lambda^{4}\mathcal{H}f \ \partial_x \left(\partial^2_\alpha S_\alpha f \sin(\frac{\gamma}{2}D_\alpha f)\sin(\frac{\gamma}{2}S_\alpha f)\frac{1}{\alpha^2} \int_{0}^{\alpha} s_{\eta} f_{x}    \ d\eta\right) \\
&& \  \times\ {  \delta_\alpha f_x } \cos(\sigma \tau_{\alpha} f_x )\cos(\arctan(\tau_{\alpha} f_{x})) \ d\eta \ d\gamma \ d\sigma \ dx \ d\alpha  \\
&-&\frac{1}{8\pi} \int \int\int_{0}^{\infty}\int_{0}^{\infty} \ \gamma^2 e^{-\gamma-\sigma} \ \Lambda^{4}\mathcal{H}f \ \partial_x \left((\partial_\alpha S_\alpha f)^2 \sin(\frac{\gamma}{2}D_\alpha f)\cos(\frac{\gamma}{2}S_\alpha f)\frac{1}{\alpha^2} \int_{0}^{\alpha} s_{\eta} f_{x}    \ d\eta\right) \\
&& \  \times\ {  \delta_\alpha f_x } \cos(\sigma \tau_{\alpha} f_x )\cos(\arctan(\tau_{\alpha} f_{x}))   \ d\gamma \ d\sigma \ dx\ d\alpha \\
&+&\frac{1}{2\pi} \int \int\int_{0}^{\infty}\int_{0}^{\infty} \ e^{-\gamma-\sigma} \ \Lambda^{4}\mathcal{H}f \ \partial_x\left(\partial^2_{\alpha}D_\alpha f \ \frac{1}{\alpha}\sin(\frac{\gamma}{2}D_\alpha f)\cos(\frac{\gamma}{2}S_\alpha f)\right)\\
&& \  \times\ {  \delta_\alpha f_x } \cos(\sigma \tau_{\alpha} f_x )\cos(\arctan(\tau_{\alpha} f_{x}))   \ d\gamma \ d\sigma \ dx \ d\alpha \\
&-&\frac{1}{2\pi} \int \int\int_{0}^{\infty}\int_{0}^{\infty} \ \gamma e^{-\gamma-\sigma} \ \Lambda^{4}\mathcal{H}f \ \partial_x \left(\partial_{\alpha}D_\alpha f \cos(\frac{\gamma}{2}D_\alpha f)\cos(\frac{\gamma}{2}S_\alpha f)\frac{1}{\alpha^3} \int_{0}^{\alpha} s_{\eta} f_{x}    \ d\eta\right) \\
&& \  \times\ {  \delta_\alpha f_x } \cos(\sigma \tau_{\alpha} f_x )\cos(\arctan(\tau_{\alpha} f_{x}))  \ d\gamma \ d\sigma \ dx \ d\alpha\\
&+&\frac{1}{2\pi} \int \int\int_{0}^{\infty}\int_{0}^{\infty} \ \gamma e^{-\gamma-\sigma} \ \Lambda^{4}\mathcal{H}f \ \partial_x \left(\partial_{\alpha}S_\alpha f \sin(\frac{\gamma}{2}D_\alpha f)\sin(\frac{\gamma}{2}S_\alpha f)\frac{1}{\alpha^3} \int_{0}^{\alpha} s_{\eta} f_{x}    \ d\eta\right) \\
&& \  \times\ {  \delta_\alpha f_x } \cos(\sigma \tau_{\alpha} f_x )\cos(\arctan(\tau_{\alpha} f_{x}))  \ d\gamma \ d\sigma \ dx  \ d\alpha \\
&-&\frac{1}{\pi} \int \int\int_{0}^{\infty}\int_{0}^{\infty} \ e^{-\gamma-\sigma} \ \Lambda^{4}\mathcal{H}f \ \partial_x \left(\partial_\alpha D_\alpha f \sin(\frac{\gamma}{2}D_\alpha f)\cos(\frac{\gamma}{2}S_\alpha f)\frac{1}{\alpha^2}\right)  \\
&& \  \times\ {  \delta_\alpha f_x } \cos(\sigma \tau_{\alpha} f_x )\cos(\arctan(\tau_{\alpha} f_{x}))   \ d\gamma \ d\sigma \ dx  \ d\alpha \\
&-&\frac{1}{4\pi} \int \int\int_{0}^{\infty}\int_{0}^{\infty} \ \gamma^2 e^{-\gamma-\sigma} \ \Lambda^{4}\mathcal{H}f \ \partial_x \left(\partial_{\alpha} D_\alpha f \ \partial_{\alpha} S_\alpha f \cos(\frac{\gamma}{2}D_\alpha f)\sin(\frac{\gamma}{2}S_\alpha f)\frac{1}{\alpha^2} \int_{0}^{\alpha} s_{\eta} f_{x}    \ d\eta\right) \\
&& \  \times\ {  \delta_\alpha f_x } \cos(\sigma \tau_{\alpha} f_x )\cos(\arctan(\tau_{\alpha} f_{x}))   \ d\gamma \ d\sigma \ dx  \ d\alpha \\
&+&\frac{1}{2\pi} \int \int\int_{0}^{\infty}\int_{0}^{\infty} \ \gamma e^{-\gamma-\sigma} \ \Lambda^{4}\mathcal{H}f \ \frac{1}{\alpha} \ \partial_x \left((\partial_{\alpha} D_\alpha f)^2 \  \cos(\frac{\gamma}{2}D_\alpha f)\cos(\frac{\gamma}{2}S_\alpha f)\right)\\
&& \  \times\ {  \delta_\alpha f_x } \cos(\sigma \tau_{\alpha} f_x )\cos(\arctan(\tau_{\alpha} f_{x}))  \ d\gamma \ d\sigma \ dx  \ d\alpha \\
&-&\frac{1}{2\pi} \int \int\int_{0}^{\infty}\int_{0}^{\infty} \ \gamma e^{-\gamma-\sigma} \ \Lambda^{4}\mathcal{H}f \ \frac{1}{\alpha}\ \partial_x \left(\partial_{\alpha} D_\alpha f \ \partial_{\alpha} S_\alpha f \sin(\frac{\gamma}{2}D_\alpha f)\sin(\frac{\gamma}{2}S_\alpha f)\right) \\
&& \  \times\ {  \delta_\alpha f_x } \cos(\sigma \tau_{\alpha} f_x )\cos(\arctan(\tau_{\alpha} f_{x}))   \ d\gamma \ d\sigma \ dx \ d\alpha \\
&:=& \sum_{i=1}^{12}  \mathcal{U}_{4,1,i}
\end{eqnarray*}

In the next part, we shall prove the following Lemma 

\begin{lemma} \label{u41} For any $i=1,...,12$, we have the following estimate
\begin{eqnarray*}
\mathcal{U}_{4,1,i}&\lesssim& \Vert f \Vert^2_{\dot H^{4}} \left(\Vert f \Vert_{\dot H^{\frac{3}{2}}}+\Vert f \Vert^{2}_{\dot H^{\frac{3}{2}}}+\Vert f \Vert^{3}_{\dot H^{\frac{3}{2}}}\right).
  \end{eqnarray*}
\end{lemma}
\noindent{{\bf{Proof of Lemma}}{ \ref{u41}} \\

$\bullet$ \ {{Estimate of $\mathcal{U}_{4,1,1}$}} \\

Recall that,

 \begin{eqnarray*}
 \mathcal{U}_{4,1,1}&=&\frac{1}{\pi} \int \int\int_{0}^{\infty}\int_{0}^{\infty} \ e^{-\gamma-\sigma} \ \Lambda^{4}\mathcal{H}f \ \frac{1}{\alpha^3}\partial_x \left(\sin(\frac{\gamma}{2}D_\alpha f)\cos(\frac{\gamma}{2}S_\alpha f)\frac{1}{\alpha} \int_{0}^{\alpha} s_{\eta} f_{x}    \ d\eta\right) \\
&& \  \times\ {  \delta_\alpha f_x } \cos(\sigma \tau_{\alpha} f_x )\cos(\arctan(\tau_{\alpha} f_{x}))   \ d\gamma \ d\sigma \ dx \ d\alpha. \\
\end{eqnarray*}
To estimate $\mathcal{U}_{4,1,1},$ we see that if the derivative falls onto one of the oscillatory term it would give rise to (up to some functions uniformly bounded by 1) either $\delta_\alpha f_x$ or $\bar\delta_\alpha f_x$  so it suffices for example to treat the case $\delta_\alpha f_x$ (the other case being analogous by using $L^2$-$L^2$ in $\alpha$ and the fact that one gets equivalent homogeneous Besov semi-norm). Otherwise, the derivative falls on the term $\int_{0}^{\alpha} s_{\eta} f_{x}    \ d\eta$ so that we gets $\int_{0}^{\alpha} s_{\eta} f_{xx}    \ d\eta$ and one may easily follow the same steps as ${T}_{4,1,1}$ by replacing  $\int_{0}^{\alpha} s_{\eta} f_{x}    \ d\eta$ by $\int_{0}^{\alpha} s_{\eta} f_{xx}    \ d\eta$ (and of course $\dot H^3$ by $\dot H^4$ ) and using the same estimates we find that it is controlled by $\Vert f \Vert_{\dot H^4} \Vert f \Vert_{\dot H^\frac{5}{2}} \Vert f \Vert_{\dot B^\frac{5}{2}_{\infty,1}}.$ By noticing that, $\dot B^\frac{5}{2}_{\infty,1}=\left[\dot B^1_{\infty,\infty},\dot B^\frac{7}{2}_{\infty,\infty}\right]_{\frac{2}{5},\frac{3}{5}}$, we find
\begin{eqnarray*}
\mathcal{U}_{4,1,1} &\lesssim& \Vert f \Vert_{\dot H^4} \int \frac{\Vert \delta_\alpha f_x \Vert^2_{L^\infty}}{\vert \alpha \vert^5} \left(\int_{0}^{\alpha} \frac{\Vert s_{\eta} f_{x} \Vert^2_{L^2}}{\vert \eta\vert^4}    \ d\eta\right)^{\frac{1}{2}} \vert \alpha \vert^{\frac{5}{2}} \ d\alpha \\
 &+& \Vert f \Vert_{\dot H^4} \Vert f \Vert_{\dot H^\frac{5}{2}} \Vert f \Vert_{\dot B^\frac{5}{2}_{\infty,1}} \\
&\lesssim& \Vert f \Vert_{\dot H^4} \int \frac{\Vert \delta_\alpha f_x \Vert^2_{L^\infty}}{\vert \alpha \vert^\frac{5}{2}} \left(\int \frac{\Vert s_{\eta} f_{x} \Vert^2_{L^2}}{\vert \eta\vert^4}    \ d\eta\right)^{\frac{1}{2}}  \ d\alpha + \Vert f \Vert_{\dot H^4} \Vert f \Vert_{\dot H^\frac{5}{2}} \Vert f \Vert^\frac{2}{5}_{\dot B^1_{\infty,\infty}} \Vert f \Vert^\frac{3}{5}_{\dot B^\frac{7}{2}_{\infty,\infty}} \\
&\lesssim& \Vert f \Vert_{\dot H^4} \Vert f \Vert_{\dot H^\frac{5}{2}} \left( \Vert f \Vert^2_{\dot B^\frac{7}{4}_{\infty,2}} +  \Vert f \Vert^\frac{2}{5}_{\dot B^1_{\infty,\infty}} \Vert f \Vert^\frac{3}{5}_{\dot B^\frac{7}{2}_{\infty,\infty}} \right) \\
&\lesssim& \Vert f \Vert_{\dot H^4}  \Vert f \Vert_{\dot H^\frac{5}{2}} \left(\Vert f \Vert^2_{\dot H^\frac{9}{4}} + \Vert f \Vert^\frac{2}{5}_{\dot H^\frac{3}{2}} \Vert f \Vert^\frac{3}{5}_{\dot H^4} \right)
\end{eqnarray*}
Since, $\Vert f \Vert_{\dot H^\frac{9}{4}}\lesssim  \Vert f \Vert^\frac{7}{10}_{\dot H^\frac{3}{2}}\Vert f \Vert^\frac{3}{10}_{\dot H^4}$ and $\Vert f \Vert_{\dot H^\frac{5}{2}}\lesssim  \Vert f \Vert^\frac{3}{5}_{\dot H^\frac{3}{2}}\Vert f \Vert^\frac{2}{5}_{\dot H^4}$ we find
\begin{eqnarray*}
\mathcal{U}_{4,1,1} &\lesssim& \Vert f \Vert_{\dot H^4} \Vert f \Vert^\frac{3}{5}_{\dot H^\frac{3}{2}}\Vert f \Vert^\frac{2}{5}_{\dot H^4} \left(\Vert f \Vert^\frac{7}{5}_{\dot H^\frac{3}{2}}\Vert f \Vert^\frac{3}{5}_{\dot H^4} + \Vert f \Vert^\frac{2}{5}_{\dot H^\frac{3}{2}} \Vert f \Vert^\frac{3}{5}_{\dot H^4}\right) \\
&\lesssim& \Vert f \Vert^2_{\dot H^4} \Vert f \Vert^2_{\dot H^\frac{3}{2}}
\end{eqnarray*}

$\bullet$ \ {{Estimate of $\mathcal{U}_{4,1,2}$}} \\

Using the decomposition obtained in ${T}_{4,1,2}$, one immediatly finds that

\begin{eqnarray*}
\mathcal{U}_{4,1,2}&=&\frac{1}{4\pi} \int \int\int_{0}^{\infty}\int_{0}^{\infty} \ \gamma e^{-\gamma-\sigma} \ \Lambda^{4}\mathcal{H} f \ \partial_{x}\left(\frac{d_\alpha f_{xx}}{\alpha} \cos(\frac{\gamma}{2}D_\alpha f)\cos(\frac{\gamma}{2}S_\alpha f)\frac{1}{\alpha^2} \int_{0}^{\alpha} s_{\eta} f_{x}    \ d\eta \right) \\
&& \  \times\ {  \delta_\alpha f_x } \cos(\sigma \tau_{\alpha} f_x )\cos(\arctan(\tau_{\alpha} f_{x}))   \ d\gamma \ d\sigma \ dx \ d\alpha \\
&+&\frac{1}{2\pi} \int \int\int_{0}^{\infty}\int_{0}^{\infty} \ \gamma e^{-\gamma-\sigma} \ \Lambda^{4}\mathcal{H} f \ \partial_{x}\left(\frac{s_\alpha f_x}{\alpha^2} \cos(\frac{\gamma}{2}D_\alpha f)\cos(\frac{\gamma}{2}S_\alpha f)\frac{1}{\alpha^2} \int_{0}^{\alpha} s_{\eta} f_{x}    \ d\eta \right) \\
&& \  \times\ {  \delta_\alpha f_x } \cos(\sigma \tau_{\alpha} f_x )\cos(\arctan(\tau_{\alpha} f_{x}))   \ d\gamma \ d\sigma \ dx \ d\alpha \\
&+&\frac{1}{2\pi} \int \int\int_{0}^{\infty}\int_{0}^{\infty} \ \gamma e^{-\gamma-\sigma} \ \Lambda^{4}\mathcal{H} f \ \frac{1}{\alpha^3} \partial_{x}\left(\int_0^\alpha s_{\kappa}f_x \ d\kappa \cos(\frac{\gamma}{2}D_\alpha f)\cos(\frac{\gamma}{2}S_\alpha f)\frac{1}{\alpha^2} \int_{0}^{\alpha} s_{\eta} f_{x}    \ d\eta \right)\\
&& \  \times\ {  \delta_\alpha f_x } \cos(\sigma \tau_{\alpha} f_x )\cos(\arctan(\tau_{\alpha} f_{x}))   \ d\gamma \ d\sigma \ dx \ d\alpha \\
&=& \sum_{i=1}^3 \mathcal{U}_{4,1,2,i}
\end{eqnarray*}

We start by estimating $\mathcal{U}_{4,1,2,1}$. The most singular term is when the derivative falls onto $d_{\alpha} f_{xx}$, in this case we find
\begin{eqnarray*}
\mathcal{U}_{4,1,2,1,sing}&\lesssim& \Vert f \Vert_{\dot H^4} \Vert f \Vert_{\dot H^3} \int \frac{\Vert \delta_\alpha f_x \Vert_{L^\infty}}{\vert \alpha \vert^3} \left(\int_{0}^{\alpha} \frac{\Vert s_{\eta} f_{x} \Vert^2}{\eta^2} \ d\eta\right)^{\frac{1}{2}} \left(\int_{0}^{\alpha} \eta^2  \ d\eta\right)^{\frac{1}{2}} \ d\alpha \\
&\lesssim& \Vert f \Vert_{\dot H^4} \Vert f \Vert_{\dot H^3} \int \frac{\Vert \delta_\alpha f_x \Vert_{L^\infty}}{\vert \alpha \vert^{\frac{3}{2}}} \left(\int \frac{\Vert s_{\eta} f_{x} \Vert^2_{L^\infty}}{\eta^2} \ d\eta\right)^{\frac{1}{2}}  \ d\alpha \\
&\lesssim& \Vert f \Vert_{\dot H^4}  \Vert f \Vert_{\dot H^3} \Vert f \Vert_{\dot H^2} \Vert f \Vert_{\dot B^{\frac{3}{2}}_{\infty,1}}
 \end{eqnarray*}   
 Since,
 \begin{equation}
 \Vert f \Vert_{\dot B^{\frac{3}{2}}_{\infty,1}} \lesssim \Vert f \Vert^{\frac{4}{5}}_{\dot B^{1}_{\infty,\infty}} \Vert f \Vert^{\frac{1}{5}}_{\dot B^{\frac{7}{2}}_{\infty,\infty}} \lesssim \Vert f \Vert^{\frac{4}{5}}_{\dot H^{\frac{3}{2}}} \Vert f \Vert^{\frac{1}{5}}_{\dot H^{4}}
\end{equation}
and 
\begin{equation}
 \Vert f \Vert_{\dot H^3} \lesssim \Vert f \Vert^{\frac{2}{5}}_{\dot H^{\frac{3}{2}}} \Vert f \Vert^{\frac{3}{5}}_{\dot H^{4}}
\end{equation}
\begin{equation}
 \Vert f \Vert_{\dot H^2} \lesssim \Vert f \Vert^{\frac{4}{5}}_{\dot H^{\frac{3}{2}}} \Vert f \Vert^{\frac{1}{5}}_{\dot H^{4}}
\end{equation}
we get
\begin{eqnarray*}
\mathcal{U}_{4,1,2,1,1}&\lesssim&  \Vert f \Vert^{2}_{\dot H^{4}}\Vert f \Vert^{2}_{\dot H^{\frac{3}{2}}}
 \end{eqnarray*}

Otherwise, the derivatives falls onto the oscillatory term or the integral term in $\eta$. In this case, we find 
\begin{eqnarray*}
\mathcal{U}_{4,1,2,1,rem1}&\lesssim& \Vert f \Vert_{\dot H^{4}} \Vert f \Vert_{\dot H^{2}} \int \frac{\Vert\delta_\alpha f_x \Vert^2_{L^\infty}}{\vert \alpha\vert^\frac{5}{2}} \ \left(\int_0^\alpha \frac{\Vert s_\eta f_x \Vert^2_{L^\infty} }{\vert\eta\vert^2} \ d\eta\right)^{\frac{1}{2}} \ d\alpha \\
&+&  \Vert f \Vert_{\dot H^{4}} \Vert f \Vert_{\dot H^2} \int \frac{\Vert \delta_\alpha f_x \Vert_{L^\infty}}{\vert \alpha \vert^{\frac{3}{2}}} \left(\int_0^\alpha \frac{\Vert s_{\eta} f_{xx} \Vert^2_{L^\infty}}{\vert\eta\vert^2} \ d\eta\right)^{\frac{1}{2}}  \ d\alpha \\
 \end{eqnarray*}
 where we have neglected again the case $\Vert\delta_\alpha f_x \Vert_{L^\infty} \Vert\bar\delta_\alpha f_x \Vert_{L^\infty}$ as it is similar to $\Vert\delta_\alpha f_x \Vert^2_{L^\infty}$ so we only treat the latter for brevity. Therefore, we get
 \begin{eqnarray*}
\mathcal{U}_{4,1,2,1,rem1}&\lesssim&\Vert f \Vert_{\dot H^{4}} \Vert f \Vert_{\dot H^{2}} \left(\Vert f\Vert^2_{\dot B^{\frac{7}{4}}_{\infty,2}}\Vert f\Vert_{\dot B^{\frac{3}{2}}_{\infty,2}}+\Vert f\Vert_{\dot B^{\frac{3}{2}}_{\infty,1}}\Vert f\Vert_{\dot B^{\frac{5}{2}}_{\infty,2}}\right)\\
&\lesssim&\Vert f \Vert_{\dot H^{4}} \Vert f \Vert_{\dot H^{2}} \left(\Vert f\Vert^2_{\dot H^{\frac{9}{4}}}\Vert f\Vert_{\dot H^{2}}+\Vert f \Vert^{\frac{4}{5}}_{\dot B^{1}_{\infty,\infty}} \Vert f \Vert^{\frac{1}{5}}_{\dot B^{\frac{7}{2}}_{\infty,\infty}}\Vert f\Vert_{\dot H^3}\right)\\
\end{eqnarray*}
since $$\Vert f\Vert^2_{\dot H^{\frac{9}{4}}} \lesssim \Vert f \Vert^\frac{7}{5}_{\dot H^\frac{3}{2}}\Vert f \Vert^\frac{3}{5}_{\dot H^4}$$ we find
\begin{eqnarray*}
\mathcal{U}_{4,1,2,1,rem1}&\lesssim&\Vert f \Vert_{\dot H^{4}}\left(\Vert f \Vert^{\frac{8}{5}}_{\dot H^{\frac{3}{2}}} \Vert f \Vert^{\frac{2}{5}}_{\dot H^{4}}\Vert f \Vert^\frac{7}{5}_{\dot H^\frac{3}{2}}\Vert f \Vert^\frac{3}{5}_{\dot H^4} +\Vert f \Vert^{\frac{4}{5}}_{\dot H^{\frac{3}{2}}} \Vert f \Vert^{\frac{1}{5}}_{\dot H^{4}}
 \Vert f \Vert^{\frac{4}{5}}_{\dot H^{\frac{3}{2}}}  \Vert f \Vert^{\frac{1}{5}}_{\dot H^{4}}\Vert f \Vert^{\frac{2}{5}}_{\dot H^{\frac{3}{2}}} \Vert f \Vert^{\frac{3}{5}}_{\dot H^{4}}\right)
\end{eqnarray*}
Therefore
\begin{eqnarray*}
\mathcal{U}_{4,1,2,1,rem1} \lesssim \Vert f \Vert^2_{\dot H^{4}} \left(\Vert f \Vert^{2}_{\dot H^{\frac{3}{2}}}+\Vert f \Vert^{3}_{\dot H^{\frac{3}{2}}}\right)
\end{eqnarray*}
Finally,
\begin{eqnarray*}
\mathcal{U}_{4,1,2,1} \lesssim \Vert f \Vert^2_{\dot H^{4}} \left(\Vert f \Vert_{\dot H^{\frac{3}{2}}}+\Vert f \Vert^{2}_{\dot H^{\frac{3}{2}}}+\Vert f \Vert^{3}_{\dot H^{\frac{3}{2}}}\right)
\end{eqnarray*}

Let us estimate $\mathcal{U}_{4,1,2,2}$. The most singular term is when the derivative hits the term $s_{\alpha} f_{xx}$. In this case we have
\begin{eqnarray*}
\mathcal{U}_{4,1,2,2,sing}  &\lesssim& \Vert f \Vert_{\dot H^4} \int \frac{\Vert s_\alpha f_{xx} \Vert_{L^{\infty}}}{\alpha^4} \Vert \delta_\alpha f_x \Vert_{L^{\infty}} \left( \int_{0}^{\alpha} \frac{\Vert s_{\eta} f_{x}\Vert^2_{L^2}}{\eta^2}    \ d\eta\right)^{\frac{1}{2}} \left( \int_{0}^{\alpha}\eta^2 \ d\eta \right)^{\frac{1}{2}} \ d\alpha\\
&\lesssim& \Vert f \Vert_{\dot H^4} \int \frac{\Vert s_\alpha f_{xx} \Vert_{L^{\infty}}}{\vert\alpha\vert^{\frac{5}{2}}} \Vert \delta_\alpha f_x \Vert_{L^{\infty}} \left( \int \frac{\Vert s_{\eta} f_{x}\Vert^2_{L^2}}{\eta^2}    \ d\eta\right)^{\frac{1}{2}} \ d\alpha \\
&\lesssim& \Vert f \Vert_{\dot H^4} \Vert f \Vert_{\dot H^{\frac{3}{2}}} \left( \int \frac{\Vert s_\alpha f_{xx} \Vert^2_{L^{\infty}}}{\vert\alpha \vert^{\frac{5}{2}}} \ d\alpha \right)^{\frac{1}{2}} \left( \int \frac{\Vert \delta_\alpha f_x \Vert^2_{L^{\infty}}}{\vert\alpha\vert^{\frac{5}{2}}} \ d\alpha \right)^{\frac{1}{2}}\\
&\lesssim& \Vert f \Vert_{\dot H^4} \Vert f \Vert_{\dot H^{\frac{3}{2}}} \Vert f \Vert^2_{\dot B^{\frac{11}{4}}_{\infty,2}}
\end{eqnarray*}
Since $\left[\dot H^{\frac{3}{2}}, \dot H^4 \right]_{\frac{1}{2},\frac{1}{2}}=\dot H^{\frac{13}{4}} \hookrightarrow \dot B^{\frac{11}{4}}_{\infty,2}$ we find that
\begin{eqnarray*}
\mathcal{U}_{4,1,2,2,sing}&\lesssim& \Vert f \Vert^2_{\dot H^3}  \Vert f \Vert^2_{\dot H^\frac{3}{2} }
\end{eqnarray*}
Otherwise, the derivative hits the remaining terms which are the oscillatory terms and the integral in $\eta$. In this case, we have, 
\begin{eqnarray*}
\mathcal{U}_{4,1,2,2,rem1}  &\lesssim& \Vert f \Vert_{\dot H^4} \sup_{\alpha \in \mathbb R} \frac{\Vert s_{\alpha}f_{x} \Vert_{L^\infty}}{\vert \alpha \vert} \int \frac{\Vert\delta_{\alpha} f_x \Vert^2_{L^\infty}}{\vert \alpha \vert^4} \left( \int_{0}^{\alpha} \frac{\Vert s_{\eta} f_{x}\Vert^2_{L^2}}{\eta^2}    \ d\eta\right)^{\frac{1}{2}} \left( \int_{0}^{\alpha}\eta^2 \ d\eta \right)^{\frac{1}{2}} \ d\alpha\\
&\lesssim& \Vert f \Vert_{\dot H^4}  \Vert f \Vert_{\dot B^{2}_{\infty,\infty}}
\Vert f \Vert_{\dot H^\frac{3}{2}}
\int \frac{\Vert\delta_{\alpha} f_x \Vert^2_{L^\infty}}{\vert \alpha \vert^{\frac{5}{2}}} \ d\alpha \\
&\lesssim& \Vert f \Vert_{\dot H^4}  \Vert f \Vert_{\dot H^{\frac{5}{2}}}
\Vert f \Vert_{\dot H^\frac{3}{2}} \Vert f \Vert^2_{\dot B^{\frac{7}{4}}_{\infty,2}} \\
&\lesssim& \Vert f \Vert_{\dot H^4}  \Vert f \Vert_{\dot H^{\frac{5}{2}}}
\Vert f \Vert_{\dot H^\frac{3}{2}} \Vert f \Vert^2_{\dot H^{\frac{9}{4}}}
\end{eqnarray*}
By interpolation,
\begin{eqnarray*}
\mathcal{U}_{4,1,2,2,rem1}&\lesssim&\Vert f \Vert^2_{\dot H^4} \Vert f \Vert^2_{\dot H^\frac{3}{2}}
\end{eqnarray*}
For the second remainder term corresponding to the differentiation of the integral in $\eta$, one finds 
\begin{eqnarray*}
\mathcal{U}_{4,1,2,2,rem2}&\lesssim& \Vert f \Vert_{\dot H^4} \int \frac{\Vert s_\alpha f_x \Vert_{L^{\infty}}}{\alpha^4} \Vert \delta_\alpha f_x \Vert_{L^{\infty}} \left( \int_{0}^{\alpha} \frac{\Vert s_{\eta} f_{xx}\Vert^2_{L^2}}{\eta^2}    \ d\eta\right)^{\frac{1}{2}} \left( \int_{0}^{\alpha}\eta^2 \ d\eta \right)^{\frac{1}{2}} \ d\alpha\\
&\lesssim& \Vert f \Vert_{\dot H^4} \int \frac{\Vert s_\alpha f_x \Vert_{L^{\infty}}}{\vert\alpha\vert^{\frac{5}{2}}} \Vert \delta_\alpha f_x \Vert_{L^{\infty}} \left( \int \frac{\Vert s_{\eta} f_{xx}\Vert^2_{L^2}}{\eta^2}    \ d\eta\right)^{\frac{1}{2}} \ d\alpha \\
&\lesssim& \Vert f \Vert_{\dot H^4} \Vert f \Vert_{\dot H^{\frac{5}{2}}} \left( \int \frac{\Vert s_\alpha f_{x} \Vert^2_{L^{\infty}}}{\vert\alpha \vert^{\frac{5}{2}}} \ d\alpha \right)^{\frac{1}{2}} \left( \int \frac{\Vert \delta_\alpha f_x \Vert^2_{L^{\infty}}}{\vert\alpha\vert^{\frac{5}{2}}} \ d\alpha \right)^{\frac{1}{2}}\\
&\lesssim& \Vert f \Vert_{\dot H^4} \Vert f \Vert_{\dot H^{\frac{5}{2}}} \Vert f \Vert^2_{\dot B^{\frac{7}{4}}_{\infty,2}} \\
&\lesssim& \Vert f \Vert_{\dot H^4} \Vert f \Vert_{\dot H^{\frac{5}{2}}} \Vert f \Vert^2_{\dot H^{\frac{9}{4}}}
\end{eqnarray*}
Hence, interpolation gives
\begin{eqnarray*}
\mathcal{U}_{4,1,2,2,rem2}&\lesssim& \Vert f \Vert^2_{\dot H^4}  \Vert f \Vert^2_{\dot H^\frac{3}{2} }
\end{eqnarray*}
Therefore, collection the singular and remainders estimates we have obtained
\begin{eqnarray*}
\mathcal{U}_{4,1,2,2}&\lesssim& \Vert f \Vert^2_{\dot H^4}  \Vert f \Vert^2_{\dot H^\frac{3}{2}}
\end{eqnarray*}

Now it remains to estimate $\mathcal{U}_{4,1,2,2}$. The integrals in $\eta$ or $\kappa$ are obviously analogous, hence, we may only control the term where the derivative hits in the integral in $\kappa$. This term (that we call $\mathcal{U}_{4,1,2,2,1}$) can be estimated as follows
\begin{eqnarray*}
\mathcal{U}_{4,1,2,2,1}&\lesssim& \Vert f \Vert_{\dot H^4} \int \frac{\Vert \delta_\alpha f_x \Vert_{L^{\infty}}}{\vert\alpha\vert^5} \left(\int_0^\alpha \frac{ \Vert s_{\kappa}f_{xx} \Vert^2_{L^{4}} }{\vert\kappa\vert^{\frac{5}{2}}} d\kappa \right)^{\frac{1}{2}} \left(\int_0^\alpha \frac{\Vert s_{\eta}f_x \Vert^2_{L^{4}}}{\vert \eta \vert^{\frac{5}{2}}}  d\eta \right)^{\frac{1}{2}} \\
&& \times \left( \int_{0}^{\alpha}\vert\eta\vert^{\frac{5}{2}} \ d\eta \right)^{\frac{1}{2}}\left( \int_{0}^{\alpha}\vert\kappa\vert^{\frac{5}{2}} \ d\kappa \right)^{\frac{1}{2}} \  d\alpha \\
&\lesssim& \Vert f \Vert_{\dot H^4} \int \frac{\Vert \delta_\alpha f_x \Vert_{L^{\infty}}}{\vert\alpha\vert^{\frac{3}{2}}} \left(\int \frac{ \Vert s_{\kappa}f_{xx} \Vert^2_{L^{4}} }{\vert\kappa\vert^{\frac{5}{2}}} d\kappa \right)^{\frac{1}{2}} \left(\int \frac{\Vert s_{\eta}f_x \Vert^2_{L^{4}}}{\vert \eta \vert^{\frac{5}{2}}}  d\eta \right)^{\frac{1}{2}}\  d\alpha \\
&\lesssim& \Vert f \Vert_{\dot H^4} \int \frac{\Vert \delta_\alpha f_x \Vert_{L^{\infty}}}{\vert\alpha\vert^\frac{3}{2}} \left(\int \frac{ \Vert s_{\kappa}f_{xx} \Vert^2_{L^{4}} }{\vert\kappa\vert^{\frac{5}{2}}} d\kappa \right)^{\frac{1}{2}} \left(\int \frac{\Vert s_{\eta}f_x \Vert^2_{L^{4}}}{\vert \eta \vert^{\frac{5}{2}}}  d\eta \right)^{\frac{1}{2}}\  d\alpha \\
&\lesssim& \Vert f \Vert_{\dot H^4} \Vert f \Vert_{\dot B^\frac{3}{2}_{\infty,1}}\Vert f \Vert_{\dot B^\frac{11}{4}_{4,2}}\Vert f \Vert_{\dot B^\frac{7}{4}_{4,2}} \\
&\lesssim& \Vert f \Vert_{\dot H^4} \Vert f \Vert^{\frac{4}{5}}_{\dot B^1_{\infty,\infty}} \Vert f \Vert^{\frac{1}{5}}_{\dot B^{\frac{7}{2}}_{\infty,\infty}} \Vert f \Vert_{\dot H^3} \Vert f \Vert_{\dot H^2}\\
&\lesssim& \Vert f \Vert_{\dot H^4} \Vert f \Vert^{\frac{4}{5}}_{\dot H^{\frac{3}{2}}} \Vert f \Vert^{\frac{1}{5}}_{\dot H^4}
\Vert f \Vert^{\frac{2}{5}}_{\dot H^{\frac{3}{2}}} \Vert f \Vert^{\frac{3}{5}}_{\dot H^{4}}
\Vert f \Vert^{\frac{4}{5}}_{\dot H^{\frac{3}{2}}} \Vert f \Vert^{\frac{1}{5}}_{\dot H^{4}}
 \end{eqnarray*}
 so
 \begin{eqnarray*}
\mathcal{U}_{4,1,2,2,1}&\lesssim& \Vert f \Vert^2_{\dot H^4} \Vert f \Vert^2_{\dot H^\frac{3}{2} }
 \end{eqnarray*}
 Therefore, 
 \begin{eqnarray*}
\mathcal{U}_{4,1,2}&\lesssim& \Vert f \Vert^2_{\dot H^{4}} \left(\Vert f \Vert_{\dot H^{\frac{3}{2}}}+\Vert f \Vert^{2}_{\dot H^{\frac{3}{2}}}+\Vert f \Vert^{3}_{\dot H^{\frac{3}{2}}}\right)
\end{eqnarray*}

$\bullet$ \ {{Estimate of $\mathcal{U}_{4,1,3}$}} \\

We have 

\begin{eqnarray*}
\mathcal{U}_{4,1,3}&=&-\frac{1}{8\pi} \int \int\int_{0}^{\infty}\int_{0}^{\infty} \ \gamma^2 e^{-\gamma-\sigma} \ \Lambda^{4}\mathcal{H}f \ \partial_x \left((\partial_\alpha D_\alpha f)^2 \sin(\frac{\gamma}{2}D_\alpha f)\cos(\frac{\gamma}{2}S_\alpha f)\frac{1}{\alpha^2} \int_{0}^{\alpha} s_{\eta} f_{x}    \ d\eta\right) \\
&& \  \times\ {  \delta_\alpha f_x } \cos(\sigma \tau_{\alpha} f_x )\cos(\arctan(\tau_{\alpha} f_{x}))   \ d\gamma \ d\sigma \ dx \ d\alpha \\
\end{eqnarray*}

The most singular term corresponds to the term where the derivative hits the term $(\partial_\alpha D_\alpha f)^2$ that is $\partial_\alpha D_\alpha f \ \partial_\alpha D_\alpha f_x$. Since

\begin{eqnarray*}
\partial_{\alpha}D_\alpha f&=& - \frac{s_\alpha f_x}{\alpha}-\frac{1}{\alpha^2} \ {\int_0^\alpha s_\kappa f_x \ d\kappa}
\end{eqnarray*}

Therefore,

\begin{eqnarray*}
\partial_{\alpha}D_\alpha f \ \partial_{\alpha}D_\alpha f_x&=&  \left(\frac{s_\alpha f_x}{\alpha}+\frac{1}{\alpha^2} \ {\int_0^\alpha s_\kappa f_x \ d\kappa}\right) \left(\frac{s_\alpha f_{xx}}{\alpha}+\frac{1}{\alpha^2} \ {\int_0^\alpha s_\mu f_{xx} \ d\mu}\right)
\end{eqnarray*}
The crossing terms are analogous to  $\mathcal{U}_{4,1,2,2,rem2}$ and $\mathcal{U}_{4,1,2,2,1}$ except the term
\begin{eqnarray*}
\mathcal{U}_{4,1,3,sing}&=&-\frac{1}{8\pi} \int \int\int_{0}^{\infty}\int_{0}^{\infty} \int_0^\alpha \int_0^\alpha \int_0^\alpha \ \gamma^2 e^{-\gamma-\sigma} \ \Lambda^{4}\mathcal{H}f \ \frac{\delta_\alpha f_x}{\alpha^6} \ s_\mu f_{xx} \  s_\kappa f_x \ s_{\eta} f_{x}  \\
&& \sin(\frac{\gamma}{2}D_\alpha f)\cos(\frac{\gamma}{2}S_\alpha f)  \cos(\sigma \tau_{\alpha} f_x )\cos(\arctan(\tau_{\alpha} f_{x})) \ d\mu \ d\kappa \ d\eta   \ d\gamma \ d\sigma \ dx \ d\alpha \\
\end{eqnarray*}
The idea is to make appear some weights in order to milder the strong decay in $\alpha$. To do so, we write

\begin{eqnarray*}
\mathcal{U}_{4,1,3,sing}&\lesssim&\Vert f \Vert_{\dot H^4} \int \frac{\Vert \delta_\alpha f_x \Vert_{L^{\infty}}}{\vert \alpha\vert^{\frac{3}{2}}}  \left(\int \frac{ \Vert s_{\mu}f_{xx} \Vert^2_{L^{6}} }{\vert\mu\vert^{2}} d\mu\right)^{\frac{1}{2}}  \left(\int \frac{ \Vert s_{\kappa}f_{x} \Vert^2_{L^{6}} }{\vert\kappa\vert^{2}} d\kappa \right)^{\frac{1}{2}}  \left(\int \frac{ \Vert s_{\eta}f_{x} \Vert^2_{L^{6}} }{\vert\eta\vert^{2}} d\eta \right)^{\frac{1}{2}} \ d\alpha \\
&\lesssim& \Vert f \Vert_{\dot H^4} \Vert f \Vert_{\dot B^{\frac{3}{2}}_{\infty,1}} \Vert f \Vert_{\dot B^{\frac{5}{2}}_{6,2}} \Vert f \Vert^2_{\dot B^{\frac{3}{2}}_{6,2}} \\
&\lesssim& \Vert f \Vert_{\dot H^4} \Vert f \Vert^{\frac{4}{5}}_{\dot B^{1}_{\infty,\infty}} \Vert f \Vert^{\frac{1}{5}}_{\dot B^{\frac{7}{2}}_{\infty,\infty}} \Vert f \Vert_{\dot H^{\frac{17}{6}}} \Vert f \Vert^2_{\dot H^{\frac{11}{6}}} \\
&\lesssim& \Vert f \Vert_{\dot H^4} \Vert f \Vert^{\frac{4}{5}}_{\dot H^{\frac{3}{2}}} \Vert f \Vert^{\frac{1}{5}}_{\dot H^{4}} \Vert f \Vert^\frac{8}{15}_{\dot H^{4}} \Vert f \Vert^\frac{7}{15}_{\dot H^{\frac{3}{2}}} \Vert f \Vert^{\frac{26}{15}}_{\dot H^{\frac{3}{2}}} \Vert f \Vert^{\frac{4}{15}}_{\dot H^4}
\end{eqnarray*}
Hence, we get 

\begin{eqnarray*}
\mathcal{U}_{4,1,3,sing}&\lesssim& \Vert f \Vert^2_{\dot H^4} \Vert f \Vert^{3}_{\dot H^{\frac{3}{2}}} 
\end{eqnarray*}

If the derivative hits the oscillatory term one observes that
\begin{eqnarray*}
\mathcal{U}_{4,1,3,rem1}&=&\frac{1}{8\pi} \int \int\int_{0}^{\infty}\int_{0}^{\infty}\int_{0}^{\alpha} \ \gamma^2 e^{-\gamma-\sigma} \ \Lambda^{4} \mathcal{H}f \ \frac{s_{\eta} f_{x} \ s^2_\alpha f_x}{\alpha^4}   \partial_x\left(\sin(\frac{\gamma}{2}D_\alpha f)\cos(\frac{\gamma}{2}S_\alpha f)\right) \\
&& \  \times\ {  \delta_\alpha f_x } \cos(\sigma \tau_{\alpha} f_x )\cos(\arctan(\tau_{\alpha} f_{x})) \ d\eta  \ d\gamma \ d\sigma \ dx \ d\alpha \\
&+&\frac{1}{4\pi} \int \int\int_{0}^{\infty}\int_{0}^{\infty}\int_{0}^{\alpha}\int_0^\alpha \ \gamma^2 e^{-\gamma-\sigma} \ \Lambda^{4} \mathcal{H}f  \ \frac{s_\alpha f_x \ s_\kappa f_x \ s_{\eta} f_{x}}{\alpha^5}  \\ 
&& \partial_x\left(\sin(\frac{\gamma}{2}D_\alpha f)\cos(\frac{\gamma}{2}S_\alpha f)\right)      {  \delta_\alpha f_x } \cos(\sigma \tau_{\alpha} f_x ) \cos(\arctan(\tau_{\alpha} f_{x})) d\kappa \ d\eta  \ d\alpha \ d\gamma \ d\sigma \ dx \\
&+&\frac{1}{8\pi} \int \int\int_{0}^{\infty}\int_{0}^{\infty}\int_{0}^{\alpha} \ \gamma^2 e^{-\gamma-\sigma} \ \Lambda^{4} \mathcal{H}f  \ \frac{s_{\eta} f_{x}}{\alpha^5}   \ \left(\int_0^\alpha s_\kappa f_x \ d\kappa\right)^2 \\
&&\partial_{x}\left(\sin(\frac{\gamma}{2}D_\alpha f)\cos(\frac{\gamma}{2}S_\alpha f)\right)   \\
&& \  \times\ {  \delta_\alpha f_x } \cos(\sigma \tau_{\alpha} f_x )\cos(\arctan(\tau_{\alpha} f_{x})) \ d\eta  \ d\gamma \ d\sigma \ dx\ d\alpha \\
&=&\mathcal{U}_{4,1,3,1,rem1}+\mathcal{U}_{4,1,3,2,rem1}+\mathcal{U}_{4,1,3,3,rem1}
 \end{eqnarray*}

Then, we  estimate  $\mathcal{U}_{4,1,3,1,rem1}$ \\
\begin{eqnarray*}
 \mathcal{U}_{4,1,3,1}&\lesssim& \Vert f \Vert_{\dot H^4} \int \frac{ \Vert \delta_\alpha f_x \Vert^2_{L^{\infty}}  \  \Vert s^2_\alpha f_x \Vert_{L^{\infty}} }{\vert\alpha\vert^5} \left(\int_{0}^{\alpha} \frac{\Vert s_{\eta} f_{x} \Vert^{2}_{L^{2}}}{\eta^{2}} \ d\eta \right)^{\frac{1}{2}} \left(\int_{0}^{\alpha} \eta^{2} \ d\eta \right)^{\frac{1}{2}} \ d\alpha \\
&\lesssim& \Vert f \Vert_{\dot H^4} \Vert f \Vert_{\dot H^{\frac{3}{2}}} \int \frac{\Vert \delta_\alpha f_x \Vert^2_{L^{\infty}}  \ \Vert s_\alpha f_x \Vert^2_{L^{\infty}}}{\vert\alpha\vert^{\frac{7}{2}}}  \ d\alpha \\
&\lesssim& \Vert f \Vert_{\dot H^4} \Vert f \Vert_{\dot H^{\frac{3}{2}}} 
\left( \int \frac{\Vert \delta_\alpha f_x \Vert^4_{L^{\infty}}}{\vert\alpha\vert^{4}}  \ d\alpha \right)^{\frac{1}{2}}
\left( \int \frac{\Vert s_\alpha f_x \Vert^4_{L^{\infty}}}{\vert\alpha\vert^{3}}  \ d\alpha \right)^{\frac{1}{2}} \\
&\lesssim& \Vert f \Vert_{\dot H^4} \Vert f \Vert_{\dot H^{\frac{3}{2}}} \Vert f \Vert^2_{\dot B^{\frac{7}{4}}_{\infty,4}} \Vert f \Vert_{\dot B^{\frac{3}{2}}_{\infty,2}} \\
&\lesssim& \Vert f \Vert_{\dot H^4} \Vert f \Vert_{\dot H^{\frac{3}{2}}} \Vert f \Vert^2_{\dot H^{\frac{9}{4}}} \Vert f \Vert^2_{\dot H^{2}} \\
&\lesssim& \Vert f \Vert^2_{\dot H^4} \Vert f \Vert_{\dot H^{\frac{3}{2}}} \Vert f \Vert^{\frac{8}{5}}_{\dot H^{\frac{3}{2}}} \Vert f \Vert^{\frac{2}{5}}_{\dot H^{4}}  \Vert f \Vert^\frac{7}{5}_{\dot H^\frac{3}{2}}\Vert f \Vert^\frac{3}{5}_{\dot H^4} \\
&\lesssim& \Vert f \Vert^2_{\dot H^4} \Vert f \Vert^4_{\dot H^{\frac{3}{2}}}
\end{eqnarray*}
As for the {{estimate of $\mathcal{U}_{4,1,3,2,rem1}$}}, we have that \\
  \begin{eqnarray*}
\mathcal{U}_{4,1,3,2,rem1}&\lesssim&\Vert f \Vert_{\dot H^4} \int  \frac{\Vert s_\alpha f_x\Vert_{L^{\infty}}\Vert \delta_\alpha f_x\Vert^2_{L^{\infty}}}{\vert\alpha\vert^6} \left(\int_0^\alpha \frac{\Vert s_\kappa f_x \Vert^2_{L^{4}}}{\vert \kappa \vert^{\frac{3}{2}}} \ d\kappa \right)^{\frac{1}{2}} \left(\int_0^\alpha \frac{ \Vert s_\eta f_x \Vert^2_{L^{4}}}{\vert \eta \vert^{\frac{3}{2}}} \ d\eta \right)^{\frac{1}{2}} \\
&& \left(\int_0^\alpha \vert\kappa\vert^{\frac{3}{2}} d\kappa\right)^{\frac{1}{2}} \left(\int_0^\alpha \vert\eta\vert^{\frac{3}{2}} d\eta\right)^{\frac{1}{2}}\ d\alpha \\ 
&\lesssim&\Vert f \Vert_{\dot H^4} \Vert f_x \Vert^2_{\dot B^{\frac{1}{4}}_{4,2}} \int  \frac{\Vert s_\alpha f_x\Vert_{L^{\infty}}\Vert \delta_\alpha f_x\Vert^2_{L^{\infty}}}{\vert\alpha\vert^{\frac{7}{2}}} \ d\alpha \\
&\lesssim&\Vert f \Vert_{\dot H^4} \Vert f \Vert^2_{\dot H^{\frac{3}{2}}} 
\left( \int \frac{\Vert \delta_\alpha f_x \Vert^4_{L^{\infty}}}{\vert\alpha\vert^{\frac{9}{2}}}  \ d\alpha \right)^{\frac{1}{2}}
\left( \int \frac{\Vert s_\alpha f_x \Vert^2_{L^{\infty}}}{\vert\alpha\vert^{\frac{5}{2}}}  \ d\alpha \right)^{\frac{1}{2}} \\
&\lesssim& \Vert f \Vert_{\dot H^4} \Vert f \Vert^2_{\dot H^{\frac{3}{2}}} 
\Vert f \Vert^2_{\dot B^{\frac{15}{8}}_{\infty,4}} \Vert f \Vert_{\dot B^{\frac{7}{4}}_{\infty,2}} \\
&\lesssim& \Vert f \Vert_{\dot H^4} \Vert f \Vert^2_{\dot H^{\frac{3}{2}}} 
\Vert f \Vert^2_{\dot H^{\frac{19}{8}}} \Vert f \Vert_{\dot H^{\frac{9}{4}}} \\
&\lesssim& \Vert f \Vert_{\dot H^4} \Vert f \Vert^2_{\dot H^{\frac{3}{2}}} 
\Vert f \Vert^\frac{7}{10}_{\dot H^{4}}\Vert f \Vert^\frac{13}{10}_{\dot H^{\frac{3}{2}}} \Vert f \Vert^\frac{7}{10}_{\dot H^\frac{3}{2}}\Vert f \Vert^\frac{3}{10}_{\dot H^4} \\
&\lesssim& \Vert f \Vert^2_{\dot H^4} \Vert f \Vert^3_{\dot H^{\frac{3}{2}}}
\end{eqnarray*}

It remains to control $\mathcal{U}_{4,1,3,3,rem1}$, we have that \\
 \begin{eqnarray*}
\mathcal{U}_{4,1,3,3}&=&-\frac{1}{8\pi} \int \int\int_{0}^{\infty}\int_{0}^{\infty}\int_{0}^{\alpha} \ \gamma^2 e^{-\gamma-\sigma} \ \Lambda^{4} \mathcal{H}f \ \frac{s_{\eta} f_{x}}{\alpha^6}   \ \left(\int_0^\alpha s_\kappa f_x \ d\kappa\right)^2  \\
&&    \partial_x\left(\sin(\frac{\gamma}{2}D_\alpha f)  
  \cos(\frac{\gamma}{2}S_\alpha f) \right) \delta_\alpha f_x   \cos(\sigma \tau_{\alpha} f_x )\cos(\arctan(\tau_{\alpha} f_{x})) \ d\eta  \ d\gamma \ d\sigma \ dx\ d\alpha. \\
&\lesssim&\Vert f \Vert_{\dot H^4} \int \frac{\Vert \delta_\alpha f_x \Vert^2_{L^\infty}}{\vert \alpha\vert^7}  \left(\int_0^\alpha  \frac{\Vert s_\eta f_x \Vert^2_{L^{2}}}{\vert \eta\vert^2} \ d\eta\right)^{\frac{1}{2}}  \left(\int_0^\alpha  \frac{\Vert s_\kappa f_x \Vert^2_{L^{\infty}}}{\vert \kappa\vert^{\frac{9}{4}}} \ d\eta\right)^{\frac{1}{2}} \\
&&\ \times  \int_0^\alpha \vert\kappa\vert^{\frac{9}{4}} d\kappa \  \left(\int_0^\alpha \vert\eta\vert^{2} d\kappa\right)^{\frac{1}{2}}d\alpha \\
&\lesssim& \Vert f \Vert_{\dot H^4} \Vert f \Vert_{\dot H^{\frac{3}{2}}} 
\Vert f_x \Vert^2_{\dot B^{\frac{5}{8}}_{\infty,2}} \int \frac{\Vert \delta_\alpha f_x \Vert^2_{L^\infty}}{\vert \alpha\vert^{\frac{9}{4}}} \ d\alpha \\
&\lesssim& \Vert f \Vert_{\dot H^4} \Vert f \Vert_{\dot H^{\frac{3}{2}}} \Vert f \Vert^2_{\dot H^{\frac{17}{8}}} \Vert f \Vert^2_{\dot B^{\frac{13}{8}}_{\infty,2}}\\
&\lesssim& \Vert f \Vert_{\dot H^4} \Vert f \Vert_{\dot H^{\frac{3}{2}}} \Vert f \Vert^4_{\dot H^{\frac{17}{8}}} 
\end{eqnarray*}
Since $$ \Vert f \Vert_{\dot H^{\frac{17}{8}}} \lesssim \Vert f \Vert^{\frac{1}{4}}_{\dot H^{4}} \Vert f \Vert^{\frac{3}{4}}_{\dot H^{\frac{3}{2}}}$$
We finally find,

\begin{eqnarray*}
\mathcal{U}_{4,1,3,3,rem1} \lesssim   \Vert f \Vert^2_{\dot H^4} \Vert f \Vert^4_{\dot H^{\frac{3}{2}}}
\end{eqnarray*}

There is a last reminder to study, it corresponds to the derivative of the integral in $\eta$ in the term

\begin{eqnarray*}
\mathcal{U}_{4,1,3,rem2}&=&-\frac{1}{8\pi} \int \int\int_{0}^{\infty}\int_{0}^{\infty} \ \gamma^2 e^{-\gamma-\sigma} \ \Lambda^{4}\mathcal{H}f \\
&& \partial_x \left((\partial_\alpha D_\alpha f)^2 \sin(\frac{\gamma}{2}D_\alpha f)\cos(\frac{\gamma}{2}S_\alpha f)\frac{1}{\alpha^2} \int_{0}^{\alpha} s_{\eta} f_{x}    \ d\eta\right) \\
&& \  \times\ {  \delta_\alpha f_x } \cos(\sigma \tau_{\alpha} f_x )\cos(\arctan(\tau_{\alpha} f_{x}))   \ d\gamma \ d\sigma \ dx \ d\alpha \\
\end{eqnarray*}
which gives
\begin{eqnarray*}
\mathcal{U}_{4,1,3,rem2}&=&-\frac{1}{8\pi} \int \int\int_{0}^{\infty}\int_{0}^{\infty} \ \gamma^2 e^{-\gamma-\sigma} \ \Lambda^{4}\mathcal{H}f \\
&& (\partial_\alpha D_\alpha f)^2 \sin(\frac{\gamma}{2}D_\alpha f)\cos(\frac{\gamma}{2}S_\alpha f)\frac{1}{\alpha^2} \int_{0}^{\alpha} s_{\eta} f_{xx}    \ d\eta \\
&& \  \times\ {  \delta_\alpha f_x } \cos(\sigma \tau_{\alpha} f_x )\cos(\arctan(\tau_{\alpha} f_{x}))   \ d\gamma \ d\sigma \ dx \ d\alpha \\
\end{eqnarray*}

This term may be decomposed as follows

\begin{eqnarray*}
\mathcal{U}_{4,1,3,1,rem2}&=&-\frac{1}{8\pi} \int \int\int_{0}^{\infty}\int_{0}^{\infty}\int_{0}^{\alpha} \ \gamma^2 e^{-\gamma-\sigma} \ \Lambda^{4} \mathcal{H}f \ \frac{1}{\alpha^2}  s_{\eta} f_{xx} \ (\partial_\alpha D_\alpha f)^2      \\
&& \sin(\frac{\gamma}{2}D_\alpha f)\cos(\frac{\gamma}{2}S_\alpha f) {  \delta_\alpha f_x } \cos(\sigma \tau_{\alpha} f_x )\cos(\arctan(\tau_{\alpha} f_{x})) \ d\eta  \ d\gamma \ d\sigma \ dx \ d\alpha \\
 \end{eqnarray*}
since
\begin{eqnarray*}
(\partial_{\alpha}D_\alpha f)^2&=&\left(  \frac{s_\alpha f_x}{\alpha}+\frac{1}{\alpha^2} \ {\int_0^\alpha s_\kappa f_x \ d\kappa} \right)^2
\end{eqnarray*}
we find,
\begin{eqnarray*}
\mathcal{U}_{4,1,3,1,rem2}&=&-\frac{1}{8\pi} \int \int\int_{0}^{\infty}\int_{0}^{\infty}\int_{0}^{\alpha} \ \gamma^2 e^{-\gamma-\sigma} \ \Lambda^{4} \mathcal{H} f \ \frac{s_{\eta} f_{xx} \ s^2_\alpha f_x}{\alpha^4}   \sin(\frac{\gamma}{2}D_\alpha f)\cos(\frac{\gamma}{2}S_\alpha f)
   \ d\eta \\
&& \  \times\ {  \delta_\alpha f_x } \cos(\sigma \tau_{\alpha} f_x )\cos(\arctan(\tau_{\alpha} f_{x})) \ d\eta  \ d\gamma \ d\sigma \ dx \ d\alpha \\
&-&\frac{1}{4\pi} \int \int\int_{0}^{\infty}\int_{0}^{\infty}\int_{0}^{\alpha}\int_0^\alpha \ \gamma^2 e^{-\gamma-\sigma} \ \Lambda^{4} \mathcal{H} f \ \frac{s_\alpha f_x \ s_\kappa f_x \ s_{\eta} f_{xx}}{\alpha^5}  \     \\
&& \sin(\frac{\gamma}{2}D_\alpha f)\cos(\frac{\gamma}{2}S_\alpha f)  {  \delta_\alpha f_x } \cos(\sigma \tau_{\alpha} f_x )\cos(\arctan(\tau_{\alpha} f_{x})) d\kappa \ d\eta  \ d\alpha \ d\gamma \ d\sigma \ dx \\
&-&\frac{1}{8\pi} \int \int\int_{0}^{\infty}\int_{0}^{\infty}\int_{0}^{\alpha} \ \gamma^2 e^{-\gamma-\sigma} \ \Lambda^{4} \mathcal{H} f \ \frac{s_{\eta} f_{xx}}{\alpha^5}   \ \left(\int_0^\alpha s_\kappa f_x \ d\kappa\right)^2     \\
&&\sin(\frac{\gamma}{2}D_\alpha f)\cos(\frac{\gamma}{2}S_\alpha f)  {  \delta_\alpha f_x } \cos(\sigma \tau_{\alpha} f_x )\cos(\arctan(\tau_{\alpha} f_{x})) \ d\eta  \ d\gamma \ d\sigma \ dx\ d\alpha \\
&=&\mathcal{U}_{4,1,3,1,1,rem}+\mathcal{U}_{4,1,3,1,2,rem}+\mathcal{U}_{4,1,3,1,3,rem}
 \end{eqnarray*}
 
Again, we estimate these 3 remainders, we find for $\mathcal{U}_{4,1,3,1,1,rem2}$ the following control \\
 \begin{eqnarray*}
 \mathcal{U}_{4,1,3,1,1,rem2}&\lesssim& \Vert f \Vert_{\dot H^4} \int \frac{ \Vert \delta_\alpha f_x \Vert_{L^{\infty}}  \  \Vert s^2_\alpha f_x \Vert_{L^{\infty}} }{\alpha^4} \left(\int_{0}^{\alpha} \frac{\Vert s_{\eta} f_{xx} \Vert^{2}_{L^{2}}}{\eta^{2}} \ d\eta \right)^{\frac{1}{2}} \left(\int_{0}^{\alpha} \eta^{2} \ d\eta \right)^{\frac{1}{2}} \ d\alpha \\
&\lesssim& \Vert f \Vert_{\dot H^4} \Vert f \Vert_{\dot H^{\frac{5}{2}}} \int \frac{\Vert \delta_\alpha f_x \Vert_{L^{\infty}}  \ \Vert s_\alpha f_x \Vert^2_{L^{\infty}}}{\vert\alpha\vert^{\frac{5}{2}}}  \ d\alpha \\
&\lesssim& \Vert f \Vert_{\dot H^3} \Vert f \Vert_{\dot H^{\frac{5}{2}}} 
\left( \int \frac{\Vert \delta_\alpha f_x \Vert^2_{L^{\infty}}}{\vert\alpha\vert^{2}}  \ d\alpha \right)^{\frac{1}{2}}
\left( \int \frac{\Vert s_\alpha f_x \Vert^4_{L^{\infty}}}{\vert\alpha\vert^{3}}  \ d\alpha \right)^{\frac{1}{2}} \\
&\lesssim& \Vert f \Vert_{\dot H^4} \Vert f \Vert_{\dot H^{\frac{5}{2}}} \Vert f \Vert^3_{\dot B^{\frac{3}{2}}_{\infty,2}} \\
&\lesssim& \Vert f \Vert_{\dot H^4} \Vert f \Vert_{\dot H^{\frac{5}{2}}} \Vert f \Vert^3_{\dot H^2} \\
&\lesssim& \Vert f \Vert_{\dot H^4} \Vert f \Vert^{\frac{3}{5}}_{\dot H^{\frac{3}{2}}} \Vert f \Vert^{\frac{2}{5}}_{\dot H^{4}}
\Vert f \Vert^{\frac{12}{5}}_{\dot H^{\frac{3}{2}}} \Vert f \Vert^{\frac{3}{5}}_{\dot H^{4}} \\
&\lesssim& \Vert f \Vert^2_{\dot H^4}  \Vert f \Vert^3_{\dot H^{\frac{3}{2}}}
\end{eqnarray*}
For the term $\mathcal{U}_{4,1,3,1,2,rem2}$ we observe that \\
  \begin{eqnarray*}
\mathcal{U}_{4,1,3,1,2,rem2}&\lesssim&\Vert f \Vert_{\dot H^4} \int  \frac{\Vert s_\alpha f_x\Vert_{L^{\infty}}\Vert \delta_\alpha f_x\Vert_{L^{\infty}}}{\vert\alpha\vert^5} \left(\int_0^\alpha \frac{\Vert s_\kappa f_x \Vert^2_{L^{4}}}{\vert \kappa \vert^{\frac{3}{2}}} \ d\kappa \right)^{\frac{1}{2}} \left(\int_0^\alpha \frac{ \Vert s_\eta f_{xx} \Vert^2_{L^{4}}}{\vert \eta \vert^{\frac{3}{2}}} \ d\eta \right)^{\frac{1}{2}} \\
&& \left(\int_0^\alpha \vert\kappa\vert^{\frac{3}{2}} d\kappa\right)^{\frac{1}{2}} \left(\int_0^\alpha \vert\eta\vert^{\frac{3}{2}} d\eta\right)^{\frac{1}{2}}\ d\alpha \\ 
&\lesssim&\Vert f \Vert_{\dot H^4} \Vert f \Vert^2_{\dot B^{\frac{9}{4}}_{4,2}} \int  \frac{\Vert s_\alpha f_{x}\Vert_{L^{\infty}}\Vert \delta_\alpha f_x\Vert_{L^{\infty}}}{\vert\alpha\vert^{\frac{5}{2}}} \ d\alpha \\
&\lesssim&\Vert f \Vert_{\dot H^4} \Vert f \Vert^2_{\dot H^{\frac{5}{2}}} 
\left( \int \frac{\Vert \delta_\alpha f_x \Vert^2_{L^{\infty}}}{\vert\alpha\vert^{\frac{5}{2}}}  \ d\alpha \right)^{\frac{1}{2}}
\left( \int \frac{\Vert s_\alpha f_x \Vert^2_{L^{\infty}}}{\vert\alpha\vert^{\frac{5}{2}}}  \ d\alpha \right)^{\frac{1}{2}} \\
&\lesssim& \Vert f \Vert_{\dot H^4} \Vert f \Vert^2_{\dot H^{\frac{5}{2}}} 
\Vert f \Vert^2_{\dot B^{\frac{7}{4}}_{\infty,2}} \\
&\lesssim& \Vert f \Vert_{\dot H^4} \Vert f \Vert^{\frac{6}{5}}_{\dot H^{\frac{3}{2}}} \Vert f \Vert^{\frac{4}{5}}_{\dot H^{4}}
\Vert f \Vert^\frac{7}{5}_{\dot H^\frac{3}{2}}\Vert f \Vert^\frac{3}{5}_{\dot H^4} \\
&\lesssim& \Vert f \Vert^2_{\dot H^4} \Vert f \Vert^3_{\dot H^{\frac{3}{2}}}
\end{eqnarray*}
Finally for $\mathcal{U}_{4,1,3,1,3,rem2}$, we use the following estimates \\
 \begin{eqnarray*}
    \mathcal{U}_{4,1,3,1,3,rem2} &\lesssim&\Vert f \Vert_{\dot H^4} \int \frac{\Vert \delta_\alpha f_x \Vert_{L^\infty}}{\vert \alpha\vert^6}  \left(\int  \frac{\Vert s_\eta f_{xx} \Vert^2_{L^{2}}}{\vert \eta\vert^2} \ d\eta\right)^{\frac{1}{2}}  \left(\int  \frac{\Vert s_\kappa f_x \Vert^2_{L^{\infty}}}{\vert \kappa\vert^{\frac{9}{4}}} \ d\eta\right) \\
&&\ \times  \int_0^\alpha \vert\kappa\vert^{\frac{9}{4}} d\kappa \  \left(\int_0^\alpha \vert\eta\vert^{2} d\kappa\right)^{\frac{1}{2}}d\alpha \\
&\lesssim& \Vert f \Vert_{\dot H^4} \Vert f \Vert_{\dot H^{\frac{5}{2}}} \Vert f_x \Vert^2_{\dot B^{\frac{5}{8}}_{\infty,2}} \int \frac{\Vert \delta_\alpha f_x \Vert_{L^\infty}}{\vert \alpha\vert^{\frac{5}{4}}} \ d\alpha \\
&\lesssim& \Vert f \Vert_{\dot H^4} \Vert f \Vert_{\dot H^{\frac{5}{2}}} \Vert f \Vert^2_{\dot H^{\frac{17}{8}}} \Vert f \Vert_{B^{\frac{5}{4}}_{\infty,1}}
\end{eqnarray*}
Since we have $\Vert f \Vert_{\dot H^{\frac{17}{8}}} \lesssim \Vert f \Vert^{\frac{1}{4}}_{\dot H^{4}} \Vert f \Vert^{\frac{3}{4}}_{\dot H^{\frac{3}{2}}}$ and $\Vert f \Vert_{B^{\frac{5}{4}}_{\infty,1}}\lesssim \Vert f \Vert^{\frac{9}{10}}_{\dot B^{1}_{\infty,\infty}}\Vert f \Vert^{\frac{1}{10}}_{\dot B^{\frac{7}{2}}_{\infty,\infty}}$ we conclude that
\begin{eqnarray*}
    \mathcal{U}_{4,1,3,1,3,rem2} &\lesssim&\Vert f \Vert_{\dot H^4}  \Vert f \Vert^{\frac{2}{5}}_{\dot H^{4}} \Vert f \Vert^{\frac{3}{5}}_{\dot H^{\frac{3}{2}}} \Vert f \Vert^{\frac{1}{2}}_{\dot H^{4}} \Vert f \Vert^{\frac{3}{2}}_{\dot H^{\frac{3}{2}}}\Vert f \Vert^{\frac{9}{10}}_{\dot H^{\frac{3}{2}}}\Vert f \Vert^{\frac{1}{10}}_{\dot H^4} \\
    &\lesssim&\Vert f \Vert^2_{\dot H^4}  \Vert f \Vert^{3}_{\dot H^{\frac{3}{2}}}
    \end{eqnarray*}

    Gathering the estimate of the singular and remainders, one concludes that
    \begin{eqnarray*}
    \mathcal{U}_{4,1,3} &\lesssim&\Vert f \Vert^2_{\dot H^4}  \left(\Vert f \Vert^{3}_{\dot H^{\frac{3}{2}}}+\Vert f \Vert^{4}_{\dot H^{\frac{3}{2}}}\right)
\end{eqnarray*}

$\bullet$  { {Estimate of $\mathcal{U}_{4,1,4}$}} \\

Recall that,

\begin{eqnarray*}
\mathcal{U}_{4,1,4}&=&-\frac{1}{4\pi} \int \int\int_{0}^{\infty}\int_{0}^{\infty} \ \gamma e^{-\gamma-\sigma} \ \Lambda^{4}\mathcal{H}f \ \partial_x \left(\partial^2_\alpha S_\alpha f \sin(\frac{\gamma}{2}D_\alpha f)\sin(\frac{\gamma}{2}S_\alpha f)\frac{1}{\alpha^2} \int_{0}^{\alpha} s_{\eta} f_{x}    \ d\eta\right) \\
&& \  \times\ {  \delta_\alpha f_x } \cos(\sigma \tau_{\alpha} f_x )\cos(\arctan(\tau_{\alpha} f_{x})) \ d\eta \ d\gamma \ d\sigma \ dx \ d\alpha  \\
\end{eqnarray*}

Since we have the identity 
\begin{eqnarray*}
\partial^2_{\alpha}S_\alpha f &=&\frac{s_{\alpha} f_{xx}}{\alpha}-\frac{\int_0^\alpha s_\eta f_{xx} \ d \eta}{\alpha^2}+\frac{d_{\alpha} f_{x}}{\alpha^2}+\frac{s_{\alpha} f}{\alpha^3}
\end{eqnarray*}

We clearly see that the most singular is when the derivatives hits $s_{\alpha} f_{xx}$  (note that $\int_0^\alpha s_\eta f_{xx} \ d \eta$ is of order 1). The most singular is therefore
\begin{eqnarray*}
\mathcal{U}_{4,1,4,sing}&=&-\frac{1}{4\pi} \int \int\int_{0}^{\infty}\int_{0}^{\infty} \ \gamma e^{-\gamma-\sigma} \ \Lambda^{3} f \ {s_{\alpha} f_{xxx}} \sin(\frac{\gamma}{2}D_\alpha f)\sin(\frac{\gamma}{2}S_\alpha f)\frac{1}{\alpha^3} \int_{0}^{\alpha} s_{\eta} f_{x}    \ d\eta \\
&& \  \times\ {  \delta_\alpha f_x } \cos(\sigma \tau_{\alpha} f_x )\cos(\arctan(\tau_{\alpha} f_{x})) \ d\eta \ d\gamma \ d\sigma \ dx \ d\alpha  \\ 
\end{eqnarray*}

We first estimate $\mathcal{U}_{4,1,4,sing}$, we have

\begin{eqnarray*}
\mathcal{U}_{4,1,4,sing}&=&\frac{1}{4\pi} \int \int\int_{0}^{\infty}\int_{0}^{\infty}\int_{0}^{\alpha} \ \gamma e^{-\gamma-\sigma} \ \Lambda^{4}\mathcal{H} f \ {s_{\alpha} f_{xxx}} \sin(\frac{\gamma}{2}D_\alpha f)\sin(\frac{\gamma}{2}S_\alpha f) \\
&& \  \times \frac{1}{\alpha^3}  s_{\eta} f_{x}    \ {  \delta_\alpha f_x } \cos(\sigma \tau_{\alpha} f_x )\cos(\arctan(\tau_{\alpha} f_{x})) \ d\eta \ d\gamma \ d\sigma \ dx \ d\alpha  \\
&\lesssim&\Vert f \Vert_{\dot H^4} \int \frac{\Vert s_{\alpha} f_{xxx} \Vert_{L^{\infty}} \Vert\delta_\alpha f_x\Vert_{L^{\infty}}}{\vert \alpha\vert^{3}}  \left(\int_{0}^{\alpha} \frac{\Vert s_{\eta} f_{x} \Vert_{L^2}}{\vert \eta \vert^2}    \ d\eta \right)^{\frac{1}{2}} \left(\int_{0}^{\alpha} {\vert \eta \vert^2}    \ d\eta \right)^{\frac{1}{2}} \\
&\lesssim&\Vert f \Vert_{\dot H^4} \Vert f \Vert_{\dot H^\frac{3}{2}}\int \frac{\Vert s_{\alpha} f_{xxx} \Vert_{L^{\infty}} \Vert\delta_\alpha f_x\Vert_{L^{\infty}}}{\vert \alpha\vert^{\frac{3}{2}}}  \ d\alpha \\
&\lesssim&\Vert f \Vert_{\dot H^4} \Vert f \Vert_{\dot H^\frac{3}{2}}
\Vert f \Vert_{\dot B^{\frac{13}{4}}_{\infty,2}}\Vert f \Vert_{\dot B^{\frac{5}{4}}_{\infty,2}} \\
&\lesssim&\Vert f \Vert_{\dot H^4} \Vert f \Vert_{\dot H^\frac{3}{2}} 
\Vert f \Vert_{\dot H^\frac{15}{4}} \Vert f \Vert_{\dot H^\frac{7}{4}}
\end{eqnarray*}

Since $\Vert f \Vert_{\dot H^\frac{15}{4}} \lesssim \Vert f \Vert^{\frac{9}{10}}_{\dot H^4} \Vert f \Vert^{\frac{1}{10}}_{\dot H^{\frac{3}{2}}}$ and $\Vert f \Vert_{\dot H^\frac{7}{4}} \lesssim \Vert f \Vert^{\frac{1}{10}}_{\dot H^4} \Vert f \Vert^{\frac{9}{10}}_{\dot H^{\frac{3}{2}}}$
 we find that 
 \begin{eqnarray*}
\mathcal{U}_{4,1,4,sing} \lesssim \Vert f \Vert^2_{\dot H^4} \Vert f \Vert^2_{\dot H^\frac{3}{2}} 
\end{eqnarray*}
 So that we get
\begin{eqnarray*}
\mathcal{U}_{4,1,4,sing}&\lesssim& \Vert f \Vert^2_{\dot H^3} \Vert f \Vert^2_{\dot H^{\frac{3}{2}}}
\end{eqnarray*}
All the other remainder terms may be estimated exactly as $\mathcal{T}_{4,1,4}$ by using the appropriate embeddings.  \\

$\bullet$ \ {{Estimate of $\mathcal{U}_{4,1,5}$}} \\

We have

\begin{eqnarray*}
\mathcal{U}_{4,1,5}&=&-\frac{1}{8\pi} \int \int\int_{0}^{\infty}\int_{0}^{\infty}\int_{0}^{\alpha}  \ \gamma^2 e^{-\gamma-\sigma} \ \Lambda^{4} \mathcal{H}f \ \partial_{x}\left((\partial_\alpha S_\alpha f)^2 \sin(\frac{\gamma}{2}D_\alpha f)\cos(\frac{\gamma}{2}S_\alpha f)\right) \\
&& \  \times \frac{1}{\alpha^2} s_{\eta} f_{x}   {  \delta_\alpha f_x } \cos(\sigma \tau_{\alpha} f_x )\cos(\arctan(\tau_{\alpha} f_{x})) d\eta  \ d\gamma \ d\sigma \ dx\ d\alpha \\
\end{eqnarray*}
Since
\begin{eqnarray*}
\partial_x\left(\partial_{\alpha}S_\alpha f\right)^2 &=& \partial_x\left(\bar\Delta_\alpha f_x-\Delta_\alpha f_x-\frac{s_\alpha f}{\alpha^2}\right)^2
\end{eqnarray*}
The most singular term is $\Delta_\alpha f_{xx}\Delta_\alpha f_{x}$ (or $\bar\Delta_\alpha f_{xx}\bar\Delta_\alpha f_{x}$ but they are the analogous with respect to the estimates), so the singular is
\begin{eqnarray*}
\mathcal{U}_{4,1,5,sing} &\lesssim& \Vert f \Vert_{\dot H^4} \int \frac{\Vert \delta_{\alpha} f_{xx} \Vert_{L^{\infty}}\Vert \delta_{\alpha} f_{x} \Vert^2_{L^{\infty}} }{\vert \alpha\vert^4} \left(\int_0^\alpha \frac{\Vert s_{\eta}f_x \Vert^2_{L^2}}{\vert \eta \vert^2}\ d\eta\right)^{\frac{1}{2}} \left(\int_0^\alpha {\vert \eta \vert^2}\ d\eta\right)^{\frac{1}{2}} \ d\alpha\\
&\lesssim& \Vert f \Vert_{\dot H^4} \Vert f \Vert_{\dot H^{\frac{3}{2}}}  \int \frac{\Vert \delta_{\alpha} f_{xx} \Vert_{L^{\infty}}\Vert \delta_{\alpha} f_{x} \Vert^2_{L^{\infty}} }{\vert \alpha\vert^{\frac{5}{2}}} \ d\alpha \\
&\lesssim& \Vert f \Vert_{\dot H^4} \Vert f \Vert_{\dot H^{\frac{3}{2}}} \left(\int \frac{\Vert \delta_{\alpha} f_{xx} \Vert^2_{L^{\infty}} }{\vert \alpha\vert^{2}} \ d\alpha\right)^{\frac{1}{2}} \left(\int \frac{\Vert \delta_{\alpha} f_{x} \Vert^4_{L^{\infty}} }{\vert \alpha\vert^{3}} \ d\alpha\right)^{\frac{1}{2}} \\
&\lesssim& \Vert f \Vert_{\dot H^4} \Vert f \Vert_{\dot H^{\frac{3}{2}}} \Vert f \Vert_{\dot H^{3}} \Vert f \Vert^2_{\dot H^{2}} \\
&\lesssim& \Vert f \Vert_{\dot H^4} \Vert f \Vert_{\dot H^{\frac{3}{2}}} 
\Vert f \Vert^{\frac{3}{5}}_{\dot H^4} \Vert f \Vert^{\frac{2}{5}}_{\dot H^{\frac{3}{2}}} \Vert f \Vert^{\frac{2}{5}}_{\dot H^4} \Vert f \Vert^{\frac{8}{5}}_{\dot H^{\frac{3}{2}}} \\
&\lesssim& \Vert f \Vert_{\dot H^4} \Vert f \Vert^{3}_{\dot H^{\frac{3}{2}}}
\end{eqnarray*}
Otherwise the derivatives are spread and one may estimated these terms following the same steps as the control of $\mathcal{T}_{4,1,5}$.\\

$\bullet$ \ {{Estimate of $\mathcal{U}_{4,1,5}$}} \\

One has,

\begin{eqnarray*}
\mathcal{U}_{4,1,6}&=&\frac{1}{2\pi} \int \int\int_{0}^{\infty}\int_{0}^{\infty} \ e^{-\gamma-\sigma} \ \Lambda^{4}\mathcal{H}f \ \partial_x\left(\partial^2_{\alpha}D_\alpha f \ \frac{1}{\alpha}\sin(\frac{\gamma}{2}D_\alpha f)\cos(\frac{\gamma}{2}S_\alpha f)\right)\\
&& \  \times\ {  \delta_\alpha f_x } \cos(\sigma \tau_{\alpha} f_x )\cos(\arctan(\tau_{\alpha} f_{x}))   \ d\gamma \ d\sigma \ dx \ d\alpha \\
\end{eqnarray*}

Using the identity
\begin{eqnarray*}
\partial_x\partial^2_{\alpha}D_\alpha f &=&\frac{d_\alpha f_{xxx}}{\alpha}+2 \frac{s_\alpha f_x}{\alpha^2} + \frac{2}{\alpha^3}\int_0^\alpha s_{\eta}f_x \ d\eta
\end{eqnarray*}
One observes that the most singular term is the one containing $d_\alpha f_{xxx}$, in order to estimate it, we write
\begin{eqnarray*}
\mathcal{U}_{4,1,6,sing}&\lesssim&\Vert f \Vert_{\dot H^4} \left( \int \frac{\Vert \delta_\alpha f_{xxx} \Vert_{L^2} \Vert \delta_\alpha f_x \Vert_{L^\infty}}{\alpha^2} \ d\alpha  \right) \\
&\lesssim&\Vert f \Vert_{\dot H^4} \Vert f \Vert_{\dot H^{\frac{7}{2}}} \Vert f \Vert_{\dot H^{2}} \\
&\lesssim&\Vert f \Vert_{\dot H^4} \Vert f \Vert^{\frac{4}{5}}_{\dot H^4}\Vert f \Vert^{\frac{1}{5}}_{\dot H^{\frac{3}{2}}} \Vert f \Vert^{\frac{1}{5}}_{\dot H^4}\Vert f \Vert^{\frac{4}{5}}_{\dot H^{\frac{3}{2}}} \\
&\lesssim&\Vert f \Vert^2_{\dot H^4} \Vert f \Vert_{\dot H^{\frac{3}{2}}}
\end{eqnarray*}

The term containing  $\frac{2}{\alpha^3}\int_0^\alpha s_{\eta}f_x \ d\eta$ is also singular in $\alpha$ but this singular term has been already treated in $T_{4,1,6}$ and the only change is the additional harmless derivative in $x$ which does not prevent one from adapting the estimates done for $T_{4,1,6}$ to the case of this special singular term. \\

$\bullet$ \ {{Estimate of $\mathcal{U}_{4,1,7}$}} \\

We want to estimate
\begin{eqnarray*} 
\mathcal{U}_{4,1,7}&=&\frac{1}{2\pi} \int \int\int_{0}^{\infty}\int_{0}^{\infty} \ \gamma e^{-\gamma-\sigma} \ \Lambda^{4}\mathcal{H}f \ \partial_x \left(\partial_{\alpha}D_\alpha f \cos(\frac{\gamma}{2}D_\alpha f)\cos(\frac{\gamma}{2}S_\alpha f)\frac{1}{\alpha^3} \int_{0}^{\alpha} s_{\eta} f_{x}    \ d\eta\right) \\
&& \  \times\ {  \delta_\alpha f_x } \cos(\sigma \tau_{\alpha} f_x )\cos(\arctan(\tau_{\alpha} f_{x}))  \ d\gamma \ d\sigma \ dx \ d\alpha\\
\end{eqnarray*}
Since, we have the identity
\begin{eqnarray*}
\partial_{x}\partial_{\alpha}D_\alpha f&=& - \frac{s_\alpha f_{xx}}{\alpha}-\frac{1}{\alpha^2} \ {\int_0^\alpha s_\kappa f_{xx} \ d\kappa}
\end{eqnarray*}
we see that the most singular term is $\frac{s_\alpha f_{xx}}{\alpha}$. The second integral in $\kappa$ and the integral in $\eta$ in $\mathcal{U}_{4,1,7}$ have the same regularity, therefore, is suffices to consider the case where the derivatives hits the integral in $\kappa$. More precisely, we have
\begin{eqnarray*}
\mathcal{U}_{4,1,7,sing}&=&\frac{1}{2\pi} \int \int\int_{0}^{\infty}\int_{0}^{\infty} \int_{0}^{\alpha} \ \gamma e^{-\gamma-\sigma} \ \Lambda^{4}\mathcal{H} f \ s_\alpha f_{xx} \ \cos(\frac{\gamma}{2}D_\alpha f)\cos(\frac{\gamma}{2}S_\alpha f)\frac{1}{\alpha^4}  s_{\eta} f_{x}     \\
&& \  \times\ {  \delta_\alpha f_x } \cos(\sigma \tau_{\alpha} f_x )\cos(\arctan(\tau_{\alpha} f_{x})) \ d\eta  \ d\gamma \ d\sigma \ d\alpha \ dx \\
&+&\frac{1}{2\pi} \int \int\int_{0}^{\infty}\int_{0}^{\infty}\int_{0}^{\alpha}\int_{0}^{\alpha} \ \gamma e^{-\gamma-\sigma} \ \Lambda^{4}\mathcal{H} f \ s_\kappa f_{xx} \ \cos(\frac{\gamma}{2}D_\alpha f)\cos(\frac{\gamma}{2}S_\alpha f)\frac{1}{\alpha^5}  s_{\eta} f_{x}    \ d\eta \\
&& \  \times\ {  \delta_\alpha f_x } \cos(\sigma \tau_{\alpha} f_x )\cos(\arctan(\tau_{\alpha} f_{x})) \ d\kappa \ d\eta  \ d\gamma \ d\sigma \ dx \ d\alpha\\
&=& \mathcal{U}_{4,1,7,1,sing} + \mathcal{U}_{4,1,7,2,rem} 
\end{eqnarray*}

Let us estimate  $\mathcal{U}_{4,1,7,1,sing}$, we have that 
\begin{eqnarray} \label{{4,1,7,1,sing}}
\mathcal{U}_{4,1,7,1,sing}&\lesssim& \Vert f \Vert_{\dot H^4} \int    \ \frac{\Vert s_\alpha f_{xx} \Vert_{L^\infty} \Vert \delta_\alpha f_x \Vert_{L^\infty} }{\alpha^4}  \left(\int_{0}^{\alpha} \frac{\Vert s_{\eta} f_{x} \Vert^2_{L^2}}{ \eta^2} \ d\eta \right)^{\frac{1}{2}} \ \left(\int_{0}^{\alpha}{ \eta^2} \ d\eta \right)^{\frac{1}{2}} \  d\alpha \nonumber\\
&\lesssim& \Vert f \Vert_{\dot H^3}  \Vert f \Vert_{\dot H^{\frac{3}{2}}} \int    \ \frac{\Vert s_\alpha f_{xx} \Vert_{L^\infty} \Vert \delta_\alpha f_x \Vert_{L^\infty} }{\vert\alpha\vert^{\frac{5}{2}}}    \  d\alpha  \\
&\lesssim&\Vert f \Vert_{\dot H^4}  \Vert f \Vert_{\dot H^{\frac{3}{2}}} \left(\int    \ \frac{\Vert \delta_\alpha f_x \Vert^2_{L^\infty} }{\vert\alpha\vert^{\frac{5}{2}}}    \  d\alpha\right)^{\frac{1}{2}} \left(\int    \ \frac{\Vert s_\alpha f_{xx} \Vert_{L^\infty}  }{\vert\alpha\vert^{\frac{5}{2}}}    \  d\alpha\right)^{\frac{1}{2}}\nonumber\\
&\lesssim&\Vert f \Vert_{\dot H^4}  \Vert f \Vert_{\dot H^{\frac{3}{2}}} \Vert f \Vert_{\dot B^{\frac{7}{4}}_{\infty,2}} \Vert f \Vert_{\dot B^{\frac{11}{4}}_{\infty,2}} \nonumber\\
&\lesssim&\Vert f \Vert_{\dot H^4}  \Vert f \Vert_{\dot H^{\frac{3}{2}}} \Vert f \Vert_{\dot H^{\frac{9}{4}}} \Vert f \Vert_{\dot H^{\frac{13}{4}}}\nonumber\\
&\lesssim&\Vert f \Vert^2_{\dot H^4}  \Vert f \Vert_{\dot H^{\frac{3}{2}}}
\Vert f \Vert^\frac{7}{10}_{\dot H^\frac{3}{2}}\Vert f \Vert^\frac{3}{10}_{\dot H^4}
 \Vert f \Vert^{\frac{7}{10}}_{\dot H^4} \Vert f \Vert^{\frac{3}{10}}_{\dot H^{\frac{3}{2}}} \\
 &\lesssim&\Vert f \Vert^2_{\dot H^4}  \Vert f \Vert^2_{\dot H^{\frac{3}{2}}}
\nonumber
\end{eqnarray}

As for the singular term $\mathcal{U}_{4,1,7,2,sing}$. The idea here is to split equally the regularity in $x$ and $\alpha$ and make sure that the artificial weighted decay is sufficiently strong so as to milder the singular kernel. In other words, we write

\begin{eqnarray} 
\mathcal{U}_{4,1,7,2,sing}&\lesssim&\Vert f \Vert_{\dot H^4} \int  \frac{\Vert \delta_\alpha f_x \Vert_{L^{\infty}}}{\vert\alpha\vert^5} 
\left(\int_{0}^{\alpha} \frac{\Vert s_{\eta} f_{x} \Vert^2_{L^{4}}}{\vert \eta \vert^{\frac{5}{2}}} \ d\eta  \right)^{\frac{1}{2}} \left(\int_{0}^{\alpha} \frac{\Vert s_{\kappa} f_{xx} \Vert^2_{L^{4}}}{\vert \kappa \vert^{\frac{5}{2}}} \ d\kappa  \right)^{\frac{1}{2}}  \nonumber\\
&& \times \left(\int_{0}^{\alpha} {\vert \eta \vert^{\frac{5}{2}}} \ d\eta  \right)^{\frac{1}{2}} \left(\int_{0}^{\alpha} {\vert \kappa \vert^{\frac{5}{2}}} \ d\kappa  \right)^{\frac{1}{2}} \ d\alpha \nonumber\\
&\lesssim&\Vert f \Vert_{\dot H^4}   \Vert f \Vert_{\dot B^{\frac{7}{4}}_{4,2}} \Vert f \Vert_{\dot B^{\frac{11}{4}}_{4,2}} \int  \frac{\Vert \delta_\alpha f_x \Vert_{L^{\infty}}}{\vert\alpha\vert^{\frac{3}{2}}}  \ d\alpha \\
&\lesssim&\Vert f \Vert_{\dot H^4} \Vert f \Vert_{\dot H^{2}} \Vert f \Vert_{\dot H^{3}} \Vert f \Vert_{\dot B^{\frac{3}{2}}_{\infty,1}} \\
&\lesssim&\Vert f \Vert_{\dot H^4} \Vert f \Vert_{\dot H^{2}} \Vert f \Vert_{\dot H^{3}}\Vert f \Vert^{\frac{4}{5}}_{\dot B^{1}_{\infty,\infty}} \Vert f \Vert^{\frac{1}{5}}_{\dot B^{\frac{7}{2}}_{\infty,\infty}} \nonumber\\
&\lesssim&\Vert f \Vert_{\dot H^4} \Vert f \Vert^{\frac{1}{5}}_{\dot H^4}\Vert f \Vert^{\frac{4}{5}}_{\dot H^{\frac{3}{2}}} \Vert f \Vert^{\frac{2}{5}}_{\dot H^{\frac{3}{2}}} \Vert f \Vert^{\frac{3}{5}}_{\dot H^{4}}\Vert f \Vert^{\frac{4}{5}}_{\dot H^{\frac{3}{2}}} \Vert f \Vert^{\frac{1}{5}}_{\dot H^{4}} \nonumber\\
&\lesssim&\Vert f \Vert^2_{\dot H^4}\Vert f \Vert^{2}_{\dot H^{\frac{3}{2}}}
\end{eqnarray}
If the derivative falls on the oscillatory terms we find that it is bounded in $L^\infty$ by $\vert \Delta_\alpha f_x\vert$, therefore we find that
\begin{eqnarray*}
\mathcal{U}_{4,1,7,rem}&=&\frac{1}{2\pi} \int \int\int_{0}^{\infty}\int_{0}^{\infty} \int_{0}^{\alpha} \ \gamma e^{-\gamma-\sigma} \ \Lambda^{4}\mathcal{H} f \ s_\alpha f_{x}  \ \partial_x\left( \cos(\frac{\gamma}{2}D_\alpha f)\cos(\frac{\gamma}{2}S_\alpha f)\right)\frac{1}{\alpha^4}  s_{\eta} f_{x}     \\
&& \  \times\ {  \delta_\alpha f_x } \cos(\sigma \tau_{\alpha} f_x )\cos(\arctan(\tau_{\alpha} f_{x})) \ d\eta  \ d\gamma \ d\sigma \ d\alpha \ dx \\
&+&\frac{1}{2\pi} \int \int\int_{0}^{\infty}\int_{0}^{\infty}\int_{0}^{\alpha}\int_{0}^{\alpha} \ \gamma e^{-\gamma-\sigma} \ \Lambda^{4}\mathcal{H} f \ s_\kappa f_x \ \partial_x\left( \cos(\frac{\gamma}{2}D_\alpha f)\cos(\frac{\gamma}{2}S_\alpha f)\right)\frac{1}{\alpha^5}  s_{\eta} f_{x}    \ d\eta \\
&& \  \times\ {  \delta_\alpha f_x } \cos(\sigma \tau_{\alpha} f_x )\cos(\arctan(\tau_{\alpha} f_{x})) \ d\kappa \ d\eta  \ d\gamma \ d\sigma \ dx \ d\alpha\\
&=& \mathcal{U}_{4,1,7,1,rem} + \mathcal{U}_{4,1,7,2,rem} 
\end{eqnarray*}

We first estimate  $\mathcal{U}_{4,1,7,1,rem} $ \\

As we have quadratic term, one does not need too much artificial decay and we may write
\begin{eqnarray} 
\mathcal{U}_{4,1,7,1,rem} &\lesssim& \Vert f \Vert_{\dot H^4} \int    \ \frac{\Vert s_\alpha f_x \Vert_{L^\infty} \Vert \delta_\alpha f_x \Vert^2_{L^\infty} }{\vert\alpha\vert^5}  \left(\int_{0}^{\alpha} \frac{\Vert s_{\eta} f_{x} \Vert^2_{L^2}}{ \eta^2} \ d\eta \right)^{\frac{1}{2}} \ \left(\int_{0}^{\alpha}{ \eta^2} \ d\eta \right)^{\frac{1}{2}} \  d\alpha \nonumber\\
&\lesssim& \Vert f \Vert_{\dot H^4}  \Vert f \Vert_{\dot H^{\frac{3}{2}}} \int    \ \frac{\Vert s_\alpha f_x \Vert_{L^\infty} \Vert \delta_\alpha f_x \Vert^2_{L^\infty} }{\vert\alpha\vert^{\frac{7}{2}}}    \  d\alpha  \\
&\lesssim&\Vert f \Vert_{\dot H^4}  \Vert f \Vert_{\dot H^{\frac{3}{2}}} \left(\int    \ \frac{\Vert \delta_\alpha f_x \Vert^2_{L^\infty} }{\vert\alpha\vert^{\frac{5}{2}}}    \  d\alpha\right)^{\frac{1}{2}} \left(\int    \ \frac{\Vert s_\alpha f_x \Vert^4_{L^\infty}  }{\vert\alpha\vert^{\frac{9}{2}}}    \  d\alpha\right)^{\frac{1}{2}}\nonumber\\
&\lesssim&\Vert f \Vert_{\dot H^4}  \Vert f \Vert_{\dot H^{\frac{3}{2}}} 
\Vert f \Vert_{\dot B^{\frac{7}{4}}_{\infty,2}} \Vert f \Vert^2_{\dot B^{\frac{15}{8}}_{\infty,2}} \nonumber\\
&\lesssim&\Vert f \Vert_{\dot H^3}  \Vert f \Vert_{\dot H^{\frac{3}{2}}} 
\Vert f \Vert_{\dot H^{\frac{9}{4}}} \Vert f \Vert^2_{\dot H^{\frac{19}{8}}}\nonumber\\
&\lesssim&\Vert f \Vert_{\dot H^3}  \Vert f \Vert_{\dot H^{\frac{3}{2}}} 
\Vert f \Vert^\frac{7}{10}_{\dot H^\frac{3}{2}}\Vert f \Vert^\frac{3}{10}_{\dot H^4} \Vert f \Vert^\frac{3}{10}_{\dot H^\frac{3}{2}}\Vert f \Vert^\frac{7}{10}_{\dot H^4} \nonumber\\
&\lesssim&\Vert f \Vert^2_{\dot H^4}  \Vert f \Vert^2_{\dot H^{\frac{3}{2}}}\nonumber
\end{eqnarray}

$\bullet$ {{Estimate of $\mathcal{U}_{4,1,7,2,rem} $}} \\
\begin{eqnarray*} 
\mathcal{U}_{4,1,7,2,rem}&\lesssim&\Vert f \Vert_{\dot H^4} \int  \frac{\Vert \delta_\alpha f_x \Vert^2_{L^{\infty}}}{\vert\alpha\vert^6} 
\left(\int_{0}^{\alpha} \frac{\Vert s_{\eta} f_{x} \Vert^2_{L^{4}}}{\vert \eta \vert^{\frac{5}{2}}} \ d\eta  \right)^{\frac{1}{2}} \left(\int_{0}^{\alpha} \frac{\Vert s_{\kappa} f_{x} \Vert^2_{L^{4}}}{\vert \kappa \vert^{\frac{5}{2}}} \ d\kappa  \right)^{\frac{1}{2}}  \nonumber\\
&& \times \left(\int_{0}^{\alpha} {\vert \eta \vert^{\frac{5}{2}}} \ d\eta  \right)^{\frac{1}{2}} \left(\int_{0}^{\alpha} {\vert \kappa \vert^{\frac{5}{2}}} \ d\kappa  \right)^{\frac{1}{2}} \ d\alpha \nonumber\\
&\lesssim&\Vert f \Vert_{\dot H^4} \Vert f \Vert^2_{\dot B^{\frac{7}{4}}_{4,2}} \int  \frac{\Vert \delta_\alpha f_x \Vert^2_{L^{\infty}}}{\vert\alpha\vert^{\frac{5}{2}}}  \ d\alpha \\
&\lesssim&\Vert f \Vert_{\dot H^4} \Vert f \Vert^2_{\dot H^{2}}\Vert f \Vert^2_{\dot B^{\frac{7}{4}}_{\infty,2}} \\
&\lesssim&\Vert f \Vert_{\dot H^4} \Vert f \Vert^2_{\dot H^{2}} \Vert f \Vert^{{2}}_{\dot H^{\frac{9}{4}}}
 \nonumber\\
&\lesssim&\Vert f \Vert_{\dot H^4}\Vert f \Vert^{\frac{2}{5}}_{\dot H^4} \Vert f \Vert^{\frac{8}{5}}_{\dot H^{\frac{3}{2}}} \Vert f \Vert^\frac{7}{5}_{\dot H^\frac{3}{2}}\Vert f \Vert^\frac{3}{5}_{\dot H^4} \\
&\lesssim&\Vert f \Vert^2_{\dot H^4} \Vert f \Vert^3_{\dot H^\frac{3}{2}}
\end{eqnarray*}

$\bullet$ \ {{Estimate of $\mathcal{U}_{4,1,8}$}} \\

Recall that, we have the following the identity
\begin{eqnarray*}
\mathcal{U}_{4,1,8}&=&\frac{1}{2\pi} \int \int\int_{0}^{\infty}\int_{0}^{\infty} \ \gamma e^{-\gamma-\sigma} \ \Lambda^{4}\mathcal{H}f \ \partial_x \left(\partial_{\alpha}S_\alpha f \sin(\frac{\gamma}{2}D_\alpha f)\sin(\frac{\gamma}{2}S_\alpha f)\frac{1}{\alpha^3} \int_{0}^{\alpha} s_{\eta} f_{x}    \ d\eta\right) \\
&& \  \times\ {  \delta_\alpha f_x } \cos(\sigma \tau_{\alpha} f_x )\cos(\arctan(\tau_{\alpha} f_{x}))  \ d\gamma \ d\sigma \ dx  \ d\alpha \\
\end{eqnarray*}
The most singular term is $\partial_{x}\partial_{\alpha}S_\alpha f$ which can be decomposed as follows. First, since one has
\begin{eqnarray*}
\partial_{x}\partial_{\alpha}S_\alpha f &=& \bar\Delta_\alpha f_{xx}-\Delta_\alpha f_{xx}-\frac{s_\alpha f_{x}}{\alpha^2},
\end{eqnarray*}
one finds that (by omitting to estimate  $\bar\Delta_\alpha f_{xx}$ as it will lead to equivalent semi-norm)
\begin{eqnarray*}
\mathcal{U}_{4,1,8,sing}&=&\frac{1}{2\pi} \int \int\int_{0}^{\infty}\int_{0}^{\infty} \ \gamma e^{-\gamma-\sigma} \ \Lambda^{3} f \ \partial_{\alpha}S_\alpha f_x \sin(\frac{\gamma}{2}D_\alpha f)\sin(\frac{\gamma}{2}S_\alpha f)\frac{1}{\alpha^3} \int_{0}^{\alpha} s_{\eta} f_{x}    \ d\eta \\
&& \  \times\ {  \delta_\alpha f_x } \cos(\sigma \tau_{\alpha} f_x )\cos(\arctan(\tau_{\alpha} f_{x}))  \ d\gamma \ d\sigma  \ d\alpha\ dx  \\
&\lesssim& \Vert f \Vert_{\dot H^3} \int  \Vert \delta_\alpha f_{xx} \Vert^2_{L^\infty} \frac{1}{\vert\alpha\vert^4} \int_{0}^{\alpha} \Vert s_{\eta} f_{x} \Vert_{L^2}    \ d\eta     \ d\alpha \\
&+&\Vert f \Vert_{\dot H^3} \int  \Vert s_\alpha f_x \Vert_{L^\infty}\frac{1}{\vert\alpha\vert^5} \Vert \delta_\alpha f_{x} \Vert_{L^\infty} \int_{0}^{\alpha} \Vert s_{\eta} f_{x} \Vert_{L^2}    \ d\eta     \ d\alpha \\
\end{eqnarray*}
It suffices to estimate one of the first two terms as they have the same control.
\begin{eqnarray*}
\mathcal{U}_{4,1,8,sing} &\lesssim&\Vert f \Vert_{\dot H^4} \int  \Vert \delta_\alpha f_{xx} \Vert^2_{L^\infty}\frac{1}{\vert\alpha\vert^4}  \left(\int_{0}^{\alpha} \frac{\Vert s_{\eta} f_{x} \Vert^2_{L^2}}{\eta^2}    \ d\eta\right)^{\frac{1}{2}} \left(\int_{0}^{\alpha} \eta^2    \ d\eta\right)^{\frac{1}{2}}     \ d\alpha \\
&+&\Vert f \Vert_{\dot H^4} \int  \Vert s_\alpha f_x \Vert_{L^\infty} \Vert \delta_\alpha f_{x} \Vert_{L^\infty} \frac{1}{\vert\alpha\vert^5} \left(\int_{0}^{\alpha} \frac{\Vert s_{\eta} f_{x} \Vert^2_{L^2}}{\eta^4}    \ d\eta\right)^{\frac{1}{2}} \left(\int_{0}^{\alpha} \eta^4    \ d\eta\right)^{\frac{1}{2}}     \ d\alpha  \\
&\lesssim&\Vert f \Vert_{\dot H^4}\left( \Vert f \Vert_{\dot H^{\frac{3}{2}}}  \int  \frac{\Vert \delta_\alpha f_x \Vert^2_{L^\infty}}{\vert\alpha\vert^{\frac{5}{2}}} \ d\alpha + \Vert f \Vert_{\dot H^{\frac{5}{2}}}\int  \frac{ \Vert s_\alpha f_x \Vert_{L^\infty}\Vert \delta_\alpha f_x \Vert_{L^\infty}}{\vert\alpha\vert^{\frac{5}{2}}} \ d\alpha  \right) \\
&\lesssim&\Vert f \Vert_{\dot H^4}\left(\Vert f \Vert_{\dot H^{\frac{3}{2}}} 
\ \Vert f \Vert^2_{\dot B^{\frac{9}{4}}_{\infty,2}} + \Vert f \Vert_{\dot H^{\frac{5}{2}}} \left( \int  \frac{ \Vert s_\alpha f_x \Vert^2_{L^\infty}}{\vert\alpha\vert^3} \ d\alpha \right)^{\frac{1}{2}} \left( \int  \frac{ \Vert \delta_\alpha f_x \Vert^2_{L^\infty}}{\vert\alpha\vert^2}\ d\alpha\right)^{\frac{1}{2}} \right)   \\
&\lesssim&\Vert f \Vert_{\dot H^4} \left(\Vert f \Vert_{\dot H^{\frac{3}{2}}} 
\Vert f \Vert^2_{\dot H^{\frac{9}{4}}} + \Vert f \Vert_{\dot H^{\frac{5}{2}}}\Vert f \Vert_{\dot B^{2}_{\infty,2}}\Vert f \Vert_{\dot H^{2}}\right)  \\
&\lesssim&\Vert f \Vert_{\dot H^4} \left(\Vert f \Vert_{\dot H^{\frac{3}{2}}} 
\Vert f \Vert^2_{\dot H^{\frac{11}{4}}} + \Vert f \Vert^2_{\dot H^{\frac{5}{2}}}\Vert f \Vert_{\dot H^{2}}\right)  \\
&\lesssim&\Vert f \Vert_{\dot H^4} \left(\Vert f \Vert_{\dot H^4}\Vert f \Vert^2_{\dot H^{\frac{3}{2}}}  + \Vert f \Vert^{\frac{4}{5}}_{\dot H^4} \Vert f \Vert^{\frac{6}{5}}_{\dot H^{\frac{3}{2}}} \Vert f \Vert^{\frac{1}{5}}_{\dot H^4} \Vert f \Vert^{\frac{4}{5}}_{\dot H^{\frac{3}{2}}}\right)  \\
&\lesssim&\Vert f \Vert^2_{\dot H^4}\Vert f \Vert^{2}_{\dot H^{\frac{3}{2}}}
\end{eqnarray*}
 
The remainders may be estimated following the same steps as above and one may checked that we obtain the same control as the singular term. \\

$\bullet$ \ {{Estimate of $\mathcal{U}_{4,1,9}$}} \\

We have
\begin{eqnarray*}
\mathcal{U}_{4,1,9}&=&-\frac{1}{\pi} \int \int\int_{0}^{\infty}\int_{0}^{\infty} \ e^{-\gamma-\sigma} \ \Lambda^{4}\mathcal{H}f \ \partial_x \left(\partial_\alpha D_\alpha f \sin(\frac{\gamma}{2}D_\alpha f)\cos(\frac{\gamma}{2}S_\alpha f)\frac{1}{\alpha^2}\right)  \\
&& \  \times\ {  \delta_\alpha f_x } \cos(\sigma \tau_{\alpha} f_x )\cos(\arctan(\tau_{\alpha} f_{x}))   \ d\gamma \ d\sigma \ dx  \ d\alpha \\
\end{eqnarray*}
We first use the fact that 

\begin{eqnarray*}
\partial_{x}\partial_{\alpha}D_\alpha f &=&  \frac{s_\alpha f_{xx}}{\alpha}-\frac{1}{\alpha^2} \ {\int_0^\alpha s_\eta f_{xx} \ d\eta}
\end{eqnarray*}

Then, one observes that the singular term may be controlled as follows

\begin{eqnarray*}
\mathcal{U}_{4,1,9,sing}&=&\frac{1}{\pi} \int \int\int_{0}^{\infty}\int_{0}^{\infty} \ e^{-\gamma-\sigma} \ \Lambda^{4} \mathcal{H} f \ \partial_x \left(\partial_\alpha D_\alpha f \sin(\frac{\gamma}{2}D_\alpha f)\cos(\frac{\gamma}{2}S_\alpha f)\right)\frac{1}{\alpha^2}  \\
&& \  \times\ {  \delta_\alpha f_x } \cos(\sigma \tau_{\alpha} f_x )\cos(\arctan(\tau_{\alpha} f_{x}))   \ d\gamma \ d\sigma \ dx  \ d\alpha \\
&\lesssim& \Vert f \Vert_{\dot H^3} \left( \int \frac{\Vert s_\alpha f_{xx} \Vert_{L^{2}} \Vert \delta_\alpha f_x \Vert_{L^{\infty}} }{\vert\alpha\vert^3} \ d\alpha \right. \\
&&\left. + \int \frac{\Vert \delta_\alpha f_x \Vert_{L^{\infty}} }{\vert\alpha\vert^4} \left(\int_0^\alpha \frac{\Vert s_\eta f_{xx} \Vert^2_{L^2}}{\vert \eta \vert^4} \ d\eta\right)^{\frac{1}{2}} \left(\int_0^\alpha \vert \eta \vert^4 \ d\eta\right)^{\frac{1}{2}} \ d\alpha \right) \\
&\lesssim& \Vert f \Vert_{\dot H^3} \left( \left(\int \frac{\Vert s_\alpha f_{xx} \Vert^2_{L^{2}}  }{\vert\alpha\vert^4} \ d\alpha\right)^{\frac{1}{2}} \left(\int \frac{ \Vert \delta_\alpha f_x \Vert^2_{L^{\infty}} }{\vert\alpha\vert^2} \ d\alpha\right)^{\frac{1}{2}} \right. \\
&&\left. + \Vert f \Vert_{\dot H^\frac{7}{2}}\int \frac{\Vert \delta_\alpha f_x \Vert_{L^{\infty}} }{\vert\alpha\vert^{\frac{3}{2}}}   \ d\alpha \right)\\
&\lesssim& \Vert f \Vert_{\dot H^3} \left(\Vert f \Vert_{\dot H^\frac{7}{2}} \Vert f \Vert_{\dot H^2} + \Vert f \Vert_{\dot H^\frac{7}{2}} \Vert f \Vert_{\dot B^{\frac{3}{2}}_{\infty,1}} \right) 
\end{eqnarray*}
Since
$$
\Vert f \Vert_{\dot B^{\frac{3}{2}}_{\infty,1}} \lesssim \Vert f \Vert^{\frac{4}{5}}_{\dot B^{1}_{\infty,\infty}} \Vert f \Vert^{\frac{1}{5}}_{\dot B^{\frac{7}{2}}_{\infty,\infty}} \lesssim \Vert f \Vert^{\frac{4}{5}}_{\dot H^{\frac{3}{2}}} \Vert f \Vert^{\frac{1}{5}}_{\dot H^{4}}
$$
and
$$
\Vert f \Vert_{\dot H^\frac{7}{2}}\lesssim \Vert f \Vert^{\frac{1}{5}}_{\dot H^\frac{3}{2}} \Vert f \Vert^{\frac{4}{5}}_{\dot H^4} , \ \ \Vert f \Vert_{\dot H^2}\lesssim   \Vert f \Vert^{\frac{4}{5}}_{\dot H^\frac{3}{2}} \Vert f \Vert^{\frac{1}{5}}_{\dot H^4}
$$
we conclude that,
\begin{eqnarray*}
\mathcal{U}_{4,1,9,sing} \lesssim \Vert f \Vert^2_{\dot H^4}\Vert f \Vert_{\dot H^\frac{3}{2}}
\end{eqnarray*}

As for the remainders, as $\left\vert\partial_x(\sin(\frac{\gamma}{2}D_\alpha f)\cos(\frac{\gamma}{2}S_\alpha f))\right\vert \lesssim \gamma \left( \left\vert\Delta_\alpha f_x \right\vert + \left\vert\bar\Delta_\alpha f_x \right\vert\right)$, one may write 
\begin{eqnarray*}
\mathcal{U}_{4,1,9,rem}&\lesssim& \Vert f \Vert_{\dot H^4} \left( \int \frac{\Vert s_\alpha f_x \Vert_{L^{2}} \Vert \delta_\alpha f_x \Vert^2_{L^{\infty}} }{\vert\alpha\vert^4} \ d\alpha \right. \\
&&\left. + \int \frac{\Vert \delta_\alpha f_x \Vert^2_{L^{\infty}} }{\vert\alpha\vert^5} \left(\int_0^\alpha \frac{\Vert s_\eta f_x \Vert^2_{L^2}}{\vert \eta \vert^4} \ d\eta\right)^{\frac{1}{2}} \left(\int_0^\alpha \vert \eta \vert^4 \ d\eta\right)^{\frac{1}{2}} \ d\alpha \right) \\
&\lesssim& \Vert f \Vert_{\dot H^4} \left(\left( \int \frac{\Vert s_\alpha f_x \Vert^2_{L^{2}}  }{\vert\alpha\vert^4} \ d\alpha  \right)^{\frac{1}{2}}\left( \int \frac{ \Vert \delta_\alpha f_x \Vert^4_{L^{\infty}} }{\vert\alpha\vert^4} \ d\alpha  \right)^{\frac{1}{2}}\right. \\
&&\left. + \ \Vert f \Vert_{\dot H^{\frac{5}{2}}}\int \frac{\Vert \delta_\alpha f_x \Vert^2_{L^{\infty}} }{\vert\alpha\vert^{\frac{5}{2}}} \ d\alpha  \right) \\
&\lesssim& \Vert f \Vert_{\dot H^4} \Vert f \Vert_{\dot H^{\frac{5}{2}}} \left( \Vert f \Vert^2_{\dot B^{\frac{7}{4}}_{\infty,4}} + \ \Vert f \Vert_{\dot B^{\frac{3}{2}}_{\infty,1}} \right) \\
&\lesssim& \Vert f \Vert_{\dot H^4} \Vert f \Vert_{\dot H^{\frac{5}{2}}} \left( \Vert f \Vert^2_{\dot H^{\frac{9}{4}}} + \Vert f \Vert^{\frac{4}{5}}_{\dot B^{1}_{\infty,\infty}} \Vert f \Vert^{\frac{1}{5}}_{\dot B^{\frac{7}{2}}_{\infty,\infty}} \right) \\
 &\lesssim& \Vert f \Vert_{\dot H^4}\Vert f \Vert^{\frac{2}{5}}_{\dot H^{4}} \Vert f \Vert^{\frac{3}{5}}_{\dot H^{\frac{3}{2}}}\left( \Vert f \Vert^\frac{7}{5}_{\dot H^\frac{3}{2}}\Vert f \Vert^\frac{3}{5}_{\dot H^4}
 + \Vert f \Vert^{\frac{4}{5}}_{\dot H^{\frac{3}{2}}} \Vert f \Vert^{\frac{1}{5}}_{\dot H^{4}} \right)
\end{eqnarray*}
Finally, we have obtained
\begin{eqnarray*}
\mathcal{U}_{4,1,9,rem}&\lesssim& \Vert f \Vert^2_{\dot H^4}  \Vert f \Vert^2_{\dot H^{\frac{3}{2}}}
\end{eqnarray*}

$\bullet$ \ {{Estimate of $\mathcal{U}_{4,1,10}$}} \\

We have

\begin{eqnarray*}
\mathcal{U}_{4,1,10}&=&-\frac{1}{4\pi} \int \int\int_{0}^{\infty}\int_{0}^{\infty}\int_{0}^{\alpha} \ \gamma^2 e^{-\gamma-\sigma} \ \Lambda^{4} \mathcal{H} f \ \partial_x\left(\partial_{\alpha} D_\alpha f \ \partial_{\alpha} S_\alpha f \cos(\frac{\gamma}{2}D_\alpha f)\sin(\frac{\gamma}{2}S_\alpha f)\right) \\
&& \  \times \frac{s_{\eta} f_{x} }{\alpha^2}      {  \delta_\alpha f_x } \cos(\sigma \tau_{\alpha} f_x )\cos(\arctan(\tau_{\alpha} f_{x}))    \ d\eta \ d\gamma \ d\sigma \ dx  \ d\alpha 
\end{eqnarray*}
We use  the fact that 
\begin{eqnarray*}
\partial_{\alpha}D_\alpha f \ \partial_{\alpha} S_\alpha f  &=& - \left( \frac{s_\alpha f_x}{\alpha}+\frac{1}{\alpha^2} \ {\int_0^\alpha s_\kappa f_x \ d\kappa}\right) \left(\bar\Delta_\alpha f_x-\Delta_\alpha f_x-\frac{s_\alpha f}{\alpha^2}\right)
\end{eqnarray*}
As $\bar\Delta_\alpha f_x$ and $\bar\Delta_\alpha f_x$ would give rise to equivalent homogeneous Besov semi-norm, it therefore suffices to consider  
\begin{eqnarray*} 
\partial_{\alpha}D_\alpha f \ \partial_{\alpha} S_\alpha f  &\approx&  \left( \frac{s_\alpha f_x}{\alpha}+\frac{1}{\alpha^2} \ {\int_0^\alpha s_\kappa f_x \ d\kappa}\right) \left(\Delta_\alpha f_x+\frac{s_\alpha f}{\alpha^2}\right)
\end{eqnarray*}
Hence, we find 
\begin{eqnarray*} 
\partial_{x} \left(\partial_{\alpha}D_\alpha f \ \partial_{\alpha} S_\alpha f\right)  &\approx&  \left( \frac{s_\alpha f_{xx}}{\alpha}+\frac{1}{\alpha^2} \ {\int_0^\alpha s_\kappa f_{xx} \ d\kappa}\right) \left(\Delta_\alpha f_{x}+\frac{s_\alpha f}{\alpha^2}\right) \\
&& \ + \left( \frac{s_\alpha f_{x}}{\alpha}+\frac{1}{\alpha^2} \ {\int_0^\alpha s_\kappa f_{x} \ d\kappa}\right)\left(\Delta_\alpha f_{xx}+\frac{s_\alpha f_{x}}{\alpha^2}\right)
\end{eqnarray*}
We shall estimate only the first 4 singular terms, the other one may be treated analogously by using the same balance in the regularity in $\alpha$, we find
\begin{eqnarray*}
\mathcal{U}_{4,1,10, sing}&\lesssim& \Vert f \Vert_{\dot H^4}\int   \ 
\frac{\Vert \delta_\alpha f_x \Vert^2_{L^{\infty}} \Vert s_\alpha f_{xx} \Vert_{L^{\infty}}}{\alpha^4}
\left( \int_{0}^{\alpha} \frac{\Vert s_{\eta} f_{x}\Vert^2_{L^2}}{\eta^2}\ d\eta\right)^{\frac{1}{2}} \left( \int_{0}^{\alpha} \eta^2 \ d\eta\right)^{\frac{1}{2}}   \ d\alpha \\
&+&\Vert f \Vert_{\dot H^4}\int   \ 
\frac{\Vert \delta_\alpha f_x \Vert_{L^{\infty}} \Vert s_\alpha f \Vert_{L^{\infty}} \Vert s_\alpha f_{xx} \Vert_{L^{\infty}}}{\vert\alpha\vert^5}
\left( \int_{0}^{\alpha} \frac{\Vert s_{\eta} f_{x}\Vert^2_{L^2}}{\eta^2}\ d\eta\right)^{\frac{1}{2}} \left( \int_{0}^{\alpha} \eta^2 \ d\eta\right)^{\frac{1}{2}}   \ d\alpha \\
&+&\Vert f \Vert_{\dot H^4}\int   \ 
\frac{\Vert \delta_\alpha f_x \Vert^2_{L^{\infty}}}{\vert\alpha\vert^5}
\left( \int_{0}^{\alpha} \frac{\Vert s_{\eta} f_{x}\Vert^2_{L^2}}{\eta^2}\ d\eta\right)^{\frac{1}{2}}  \left({\int_0^\alpha \frac{\Vert s_\kappa f_{xx} \Vert_{L^\infty}}{\vert \kappa \vert^2} \ d\kappa} \right)^{\frac{1}{2}}  \\&& \times \left( \int_{0}^{\alpha} \eta^2 \ d\eta\right)^{\frac{1}{2}} \left( \int_{0}^{\alpha} \kappa^2 \ d\kappa\right)^{\frac{1}{2}}   \ d\alpha \\
&+&\Vert f \Vert_{\dot H^4}\int   \ 
\frac{\Vert \delta_\alpha f_x \Vert_{L^{\infty}} \Vert s_\alpha f \Vert_{L^{\infty}}}{\vert\alpha\vert^6}
\left( \int_{0}^{\alpha} \frac{\Vert s_{\eta} f_{x}\Vert^2_{L^2}}{\eta^2}\ d\eta\right)^{\frac{1}{2}}  \left({\int_0^\alpha \frac{\Vert s_\kappa f_{xx} \Vert^2_{L^\infty}}{\vert \kappa \vert^2} \ d\kappa} \right)^{\frac{1}{2}}  \\&& \times \left( \int_{0}^{\alpha} \eta^2 \ d\eta\right)^{\frac{1}{2}} \left( \int_{0}^{\alpha} \kappa^2 \ d\kappa\right)^{\frac{1}{2}}   \ d\alpha \\
\end{eqnarray*}
Then

\begin{eqnarray*}
\mathcal{U}_{4,1,10, sing}&\lesssim& \Vert f \Vert_{\dot H^4} \Vert f \Vert_{\dot H^\frac{3}{2}} \int   \ 
\frac{\Vert \delta_\alpha f_x \Vert^2_{L^{\infty}} \Vert s_\alpha f_{xx} \Vert_{L^{\infty}}}{\vert\alpha\vert^{\frac{5}{2}}}   \ d\alpha \\
&+&\Vert f \Vert_{\dot H^4} \Vert f \Vert_{\dot H^\frac{3}{2}}\int   \ 
\frac{\Vert \delta_\alpha f_x \Vert_{L^{\infty}} \Vert s_\alpha f \Vert_{L^{\infty}} \Vert s_\alpha f_{xx} \Vert_{L^{\infty}}}{\vert\alpha\vert^\frac{7}{2}}
   \ d\alpha \\
&+&\Vert f \Vert_{\dot H^4} \Vert f \Vert_{\dot H^\frac{3}{2}} \Vert f \Vert_{\dot H^3}\int   \ 
\frac{\Vert \delta_\alpha f_x \Vert^2_{L^{\infty}}}{\vert\alpha\vert^2}  \ d\alpha \\
&+&\Vert f \Vert_{\dot H^4} \Vert f \Vert_{\dot H^\frac{3}{2}} 
 \Vert f \Vert_{\dot H^{3}}\int   \ 
\frac{\Vert \delta_\alpha f_x \Vert_{L^{\infty}} \Vert s_\alpha f \Vert_{L^{\infty}}}{\vert\alpha\vert^3}    \ d\alpha \\
\end{eqnarray*}
We find
\begin{eqnarray*}
\mathcal{U}_{4,1,10,sing}&\lesssim& \Vert f \Vert_{\dot H^4} \Vert f \Vert_{\dot H^\frac{3}{2}} \left(\sup_{\alpha \in \mathbb R} \frac{\Vert s_\alpha f_{xx} \Vert_{L^{\infty}} }{\vert\alpha\vert}\right) \int   \ 
\frac{\Vert \delta_\alpha f_x \Vert^2_{L^{\infty}}}{\vert \alpha \vert^{\frac{3}{2}}} \ d\alpha   \\
&+&\Vert f \Vert_{\dot H^4} \Vert f \Vert_{\dot H^\frac{3}{2}} \left(\sup_{\alpha \in \mathbb R} \frac{\Vert s_\alpha f \Vert_{L^{\infty}} }{\vert\alpha\vert}\right)\int   \ 
\frac{\Vert \delta_\alpha f_x \Vert_{L^{\infty}} 
 \Vert s_\alpha f_{xx} \Vert_{L^{\infty}}}{\vert\alpha\vert^\frac{5}{2}}
   \ d\alpha \\
   &+&\Vert f \Vert_{\dot H^4} \Vert f \Vert_{\dot H^\frac{3}{2}} \Vert f \Vert_{\dot H^3} \Vert f \Vert^2_{\dot H^2} \\
   &+&\Vert f \Vert_{\dot H^3} \Vert f \Vert_{\dot H^\frac{3}{2}}  \Vert f \Vert_{\dot H^3}\int   \ 
\frac{\Vert \delta_\alpha f_x \Vert_{L^{\infty}} \Vert s_\alpha f \Vert_{L^{\infty}}}{\vert\alpha\vert^3}    \ d\alpha \\
\end{eqnarray*}
Therefore,
\begin{eqnarray*}
\mathcal{U}_{4,1,10,sing}&\lesssim& \Vert f \Vert_{\dot H^4} \Vert f \Vert_{\dot H^\frac{3}{2}} 
\Vert f \Vert_{\dot B^3_{\infty,\infty}} \Vert f \Vert^2_{\dot B^{\frac{5}{4}}_{\infty,2}}    \\
&+&\Vert f \Vert_{\dot H^4} \Vert f \Vert_{\dot H^\frac{3}{2}}\Vert f \Vert_{\dot B^1_{\infty,\infty}} \left( \int    \frac{\Vert \delta_\alpha f_x \Vert^2_{L^{\infty}} 
 }{\vert\alpha\vert^2} \ d\alpha \right)^{\frac{1}{2}} \left( \int    \frac{\Vert s_\alpha f_{xx} \Vert^2_{L^{\infty}} 
 }{\vert\alpha\vert^3} \ d\alpha \right)^{\frac{1}{2}}\\
   &+&\Vert f \Vert_{\dot H^4} \Vert f \Vert_{\dot H^\frac{3}{2}} \Vert f \Vert_{\dot H^3} \Vert f \Vert^2_{\dot H^2} \\
   &+&\Vert f \Vert_{\dot H^4} \Vert f \Vert_{\dot H^\frac{3}{2}}  \Vert f \Vert_{\dot H^3} \left(\int \frac{\Vert \delta_\alpha f_x \Vert^2_{L^{\infty}}}{\vert \alpha\vert^2}    \ d\alpha \right)^{\frac{1}{2}}  \left(\int \frac{\Vert s_\alpha f \Vert^2_{L^{\infty}}}{\vert \alpha\vert^4}    \ d\alpha \right)^{\frac{1}{2}} \\
   &\lesssim& \Vert f \Vert_{\dot H^4} \Vert f \Vert_{\dot H^\frac{3}{2}} 
\Vert f \Vert^{\frac{4}{5}}_{\dot B^{\frac{7}{2}}_{\infty,\infty}} \Vert f \Vert^{\frac{1}{5}}_{\dot B^{1}_{\infty,\infty}} \Vert f \Vert^2_{\dot H^{\frac{7}{4}}}    \\
&+&\Vert f \Vert_{\dot H^4} \Vert f \Vert_{\dot H^\frac{3}{2}}\Vert f \Vert_{\dot B^1_{\infty,\infty}} \Vert f \Vert_{\dot H^2} \Vert f \Vert_{\dot H^\frac{7}{2}} \\
   &+&\Vert f \Vert_{\dot H^4} \Vert f \Vert_{\dot H^\frac{3}{2}} \Vert f \Vert_{\dot H^3} \Vert f \Vert^2_{\dot H^2} \\
  \end{eqnarray*}
 We find,
 \begin{eqnarray*}
\mathcal{U}_{4,1,10,sing}&\lesssim&  \Vert f \Vert_{\dot H^4} \Vert f \Vert_{\dot H^\frac{3}{2}} 
\Vert f \Vert^{\frac{4}{5}}_{\dot H^4} \Vert f \Vert^{\frac{1}{5}}_{\dot H^\frac{3}{2}} \Vert f \Vert^{\frac{1}{5}}_{\dot H^4} \Vert f \Vert^{\frac{4}{5}}_{\dot H^\frac{3}{2}} \\
&+& \Vert f \Vert_{\dot H^4} \Vert f \Vert^2_{\dot H^\frac{3}{2}}  \Vert f \Vert^{\frac{1}{5}}_{\dot H^4} \Vert f \Vert^{\frac{4}{5}}_{\dot H^\frac{3}{2}} \Vert f \Vert^{\frac{4}{5}}_{\dot H^4} \Vert f \Vert^{\frac{1}{5}}_{\dot H^\frac{3}{2}} \\
&+&\Vert f \Vert_{\dot H^4} \Vert f \Vert_{\dot H^\frac{3}{2}}
 \Vert f \Vert^{\frac{3}{5}}_{\dot H^4} \Vert f \Vert^{\frac{2}{5}}_{\dot H^\frac{3}{2}} \Vert f \Vert^{\frac{2}{5}}_{\dot H^4} \Vert f \Vert^{\frac{3}{5}}_{\dot H^\frac{3}{2}} \\
\end{eqnarray*}
Finally,
\begin{eqnarray*}
\mathcal{U}_{4,1,10,sing}&\lesssim&  \Vert f \Vert^2_{\dot H^4} \left( \Vert f \Vert^2_{\dot H^\frac{3}{2}}+\Vert f \Vert^3_{\dot H^\frac{3}{2}} \right)
 \end{eqnarray*}
 
  
$\bullet$ \ {{Estimate of $\mathcal{U}_{4,1,11}$}} \\

We have 

\begin{eqnarray*}
\mathcal{U}_{4,1,11}&=&\frac{1}{2\pi} \int \int\int_{0}^{\infty}\int_{0}^{\infty} \ \gamma e^{-\gamma-\sigma} \ \Lambda^{4} \mathcal{H}f \ \frac{1}{\alpha} \ \partial_x\left((\partial_{\alpha} D_\alpha f)^2 \  \cos(\frac{\gamma}{2}D_\alpha f)\cos(\frac{\gamma}{2}S_\alpha f)\right)\\
&& \  \times\ {  \delta_\alpha f_x } \cos(\sigma \tau_{\alpha} f_x )\cos(\arctan(\tau_{\alpha} f_{x}))  \ d\gamma \ d\sigma \ dx  \ d\alpha \\
\end{eqnarray*}
Since
\begin{eqnarray*}
\left(\partial_{\alpha}D_\alpha f \right)^2  &=&  \left( \frac{s_\alpha f_x}{\alpha}+\frac{1}{\alpha^2} \ {\int_0^\alpha s_\kappa f_x \ d\kappa}\right)^2 
\end{eqnarray*}
Then, we see that
\begin{eqnarray*}
\partial_x\left(\partial_{\alpha}D_\alpha f \right)^2  &\approx&  \left( \frac{s_\alpha f_{xx}}{\alpha}+\frac{1}{\alpha^2} \ {\int_0^\alpha s_\kappa f_{xx} \ d\kappa}\right)  \left( \frac{s_\alpha f_x}{\alpha}+\frac{1}{\alpha^2} \ {\int_0^\alpha s_\kappa f_x \ d\kappa}\right)
\end{eqnarray*}
Hence, the singular term gives rise to 4 terms which can be estimated as follows
\begin{eqnarray*}
\mathcal{U}_{4,1,11} &\lesssim& \Vert f \Vert_{\dot H^4} \int \frac{\Vert s_\alpha f_{xx}\Vert_{L^2} \Vert s_\alpha f_x \Vert_{L^\infty} \Vert \delta_\alpha f_x \Vert_{L^\infty} }{\vert \alpha\vert^3} \ d\alpha \\
&+&\Vert f \Vert_{\dot H^4} \int \frac{1}{\vert \alpha\vert^{\frac{5}{2}}}{\Vert s_\alpha f_{xx}\Vert_{L^\infty} \left(\int_0^\alpha \frac{\Vert s_\kappa f_x \Vert^2_{L^2}}{\vert\kappa \vert^2} \ d\kappa\right)^{\frac{1}{2}} \Vert \delta_\alpha f_x \Vert_{L^\infty} }{} \ d\alpha \\
&+&\Vert f \Vert_{\dot H^4} \int \frac{1}{\vert \alpha\vert^{\frac{5}{2}}}{\Vert s_\alpha f_{x}\Vert_{L^\infty} \left(\int_0^\alpha \frac{\Vert s_\kappa f_{xx} \Vert^2_{L^2}}{\vert\kappa \vert^2} \ d\kappa\right)^{\frac{1}{2}} \Vert \delta_\alpha f_x \Vert_{L^\infty} }{} \ d\alpha \\
&+&\Vert f \Vert_{\dot H^4} \int \frac{1}{\vert \alpha\vert^{\frac{3}{2}}}{\left(\int_0^\alpha \frac{\Vert s_\kappa f_{xx} \Vert^2_{L^2}}{\vert\kappa \vert^2} \ d\kappa\right)^{\frac{1}{2}} \left(\int_0^\alpha \frac{\Vert s_\kappa f_{x} \Vert^2_{L^2}}{\vert\kappa \vert^3} \ d\kappa\right)^{\frac{1}{2}} \Vert \delta_\alpha f_x \Vert_{L^\infty} }{} \ d\alpha \\
&\lesssim& \Vert f \Vert_{\dot H^4} \Vert f \Vert_{\dot B^2_{\infty,\infty}} \left(\int \frac{\Vert s_\alpha f_{xx}\Vert^2_{L^2}}{\alpha^2} \ d\alpha\right)^{\frac{1}{2}} \left(\int \frac{\Vert \delta_\alpha f_{x}\Vert^2_{L^\infty}}{\alpha^2} \ d\alpha\right)^{\frac{1}{2}} \\
&+& \Vert f \Vert_{\dot H^4} \Vert f \Vert_{\dot H^\frac{3}{2}} \left(\int \frac{\Vert s_\alpha f_{xx}\Vert^2_{L^\infty}}{\vert\alpha\vert^3} \ d\alpha\right)^{\frac{1}{2}} \left(\int \frac{\Vert \delta_\alpha f_{x}\Vert^2_{L^\infty}}{\alpha^2} \ d\alpha\right)^{\frac{1}{2}} \\
&+& \Vert f \Vert_{\dot H^4} \Vert f \Vert_{\dot H^\frac{5}{2}} \left(\int \frac{\Vert s_\alpha f_{x}\Vert^2_{L^\infty}}{\vert\alpha\vert^3} \ d\alpha\right)^{\frac{1}{2}} \left(\int \frac{\Vert \delta_\alpha f_{x}\Vert^2_{L^\infty}}{\alpha^2} \ d\alpha\right)^{\frac{1}{2}} \\
&+& \Vert f \Vert_{\dot H^4}\Vert f \Vert_{\dot H^\frac{5}{2}} \Vert f \Vert_{\dot H^2} \Vert f \Vert_{\dot B^{\frac{3}{2}}_{\infty,1}} \\
&\lesssim& \Vert f \Vert_{\dot H^4} \Vert f \Vert^2_{\dot H^\frac{5}{2}}\Vert f \Vert_{\dot H^2} + \Vert f \Vert_{\dot H^4}\Vert f \Vert_{\dot H^\frac{3}{2}} \Vert f \Vert_{\dot H^\frac{7}{2}}\Vert f \Vert_{\dot H^2} \\
&+&  \Vert f \Vert_{\dot H^4}\Vert f \Vert_{\dot H^\frac{5}{2}} \Vert f \Vert_{\dot H^2} \Vert f \Vert^{\frac{1}{5}}_{\dot H^{4}}\Vert f \Vert^{\frac{4}{5}}_{\dot H^{\frac{3}{2}}}   \\
\end{eqnarray*}
Hence, by interpolation

\begin{eqnarray*}
\mathcal{U}_{4,1,11} &\lesssim& \Vert f \Vert_{\dot H^4}\left(\Vert f \Vert_{\dot H^\frac{3}{2}} +\Vert f \Vert^2_{\dot H^\frac{3}{2}} \right)
\end{eqnarray*}

$\bullet$ \ {{Estimate of $\mathcal{U}_{4,1,12}$}} \\

Recall that
\begin{eqnarray*}
\mathcal{U}_{4,1,12}&=&-\frac{1}{2\pi} \int \int\int_{0}^{\infty}\int_{0}^{\infty} \ \gamma e^{-\gamma-\sigma} \ \Lambda^{4}\mathcal{H}f \ \frac{1}{\alpha}\ \partial_x \left(\partial_{\alpha} D_\alpha f \ \partial_{\alpha} S_\alpha f \sin(\frac{\gamma}{2}D_\alpha f)\sin(\frac{\gamma}{2}S_\alpha f)\right) \\
&& \  \times\ {  \delta_\alpha f_x } \cos(\sigma \tau_{\alpha} f_x )\cos(\arctan(\tau_{\alpha} f_{x}))   \ d\gamma \ d\sigma \ dx \ d\alpha \\
\end{eqnarray*}

Since
\begin{eqnarray*} 
\partial_{x} \left(\partial_{\alpha}D_\alpha f \ \partial_{\alpha} S_\alpha f\right)  &\approx&  \left( \frac{s_\alpha f_{xx}}{\alpha}+\frac{1}{\alpha^2} \ {\int_0^\alpha s_\kappa f_{xx} \ d\kappa}\right) \left(\Delta_\alpha f_{x}+\frac{s_\alpha f}{\alpha^2}\right) \\
&& \ + \left( \frac{s_\alpha f_{x}}{\alpha}+\frac{1}{\alpha^2} \ {\int_0^\alpha s_\kappa f_{x} \ d\kappa}\right)\left(\Delta_\alpha f_{xx}+\frac{s_\alpha f_{x}}{\alpha^2}\right)
\end{eqnarray*}
We find that the singular term may be controlled as follows
\begin{eqnarray*}
\mathcal{U}_{4,1,12,sing}&\lesssim& \Vert f \Vert_{\dot H^4} \int \frac{\Vert s_\alpha f_{xx}\Vert_{L^\infty}}{\vert\alpha\vert} \frac{\Vert \delta_\alpha f_x\Vert^2_{L^4}}{\vert\alpha\vert^2} \ d\alpha \\
&+&\Vert f \Vert_{\dot H^4} \int \frac{\Vert s_\alpha f_{xx}\Vert_{L^\infty} \Vert s_\alpha f\Vert_{L^\infty}}{\vert\alpha\vert^3} \frac{\Vert \delta_\alpha f_x\Vert_{L^2}}{\vert\alpha\vert} \ d\alpha\\
&+&\Vert f \Vert_{\dot H^4} \int  \frac{\Vert \delta_\alpha f_x\Vert^2_{L^\infty}}{\vert\alpha\vert^4} \ \int_0^\alpha \Vert s_\kappa f_{xx} \Vert_{L^2} \ d\kappa  \ d\alpha\\
&+&\Vert f \Vert_{\dot H^4} \int  \frac{\Vert s_\alpha f\Vert_{L^\infty}\Vert \delta_\alpha f_x\Vert_{L^\infty}}{\vert\alpha\vert^5} \ \int_0^\alpha \Vert s_\kappa f_{xx} \Vert_{L^2} \ d\kappa  \ d\alpha\\
&\lesssim& \Vert f \Vert_{\dot H^4} \ \sup_{\alpha \in \mathbb R} \frac{\Vert s_\alpha f_{xx}\Vert_{L^\infty}}{\vert\alpha\vert} \int  \frac{\Vert \delta_\alpha f_x\Vert^2_{L^4}}{\vert\alpha\vert^2} \ d\alpha \\
&+&\Vert f \Vert_{\dot H^4} \sup_{\alpha \in \mathbb R} \frac{\Vert s_\alpha f_{xx}\Vert_{L^\infty}}{\vert\alpha\vert} \int \frac{ \Vert s_\alpha f\Vert_{L^\infty}}{\vert\alpha\vert^2} \frac{\Vert \delta_\alpha f_x\Vert_{L^2}}{\vert\alpha\vert} \ d\alpha\\
&+&\Vert f \Vert_{\dot H^4} \int  \frac{\Vert \delta_\alpha f_x\Vert^2_{L^\infty}}{\vert\alpha\vert^\frac{5}{2}} \ \left(\int \frac{\Vert s_\kappa f_{xx} \Vert_{L^2}}{\kappa^2} \ d\kappa\right)^{\frac{1}{2}}  \ d\alpha\\
&+&\Vert f \Vert_{\dot H^4} \int  \frac{\Vert s_\alpha f\Vert_{L^\infty}\Vert \delta_\alpha f_x\Vert_{L^\infty}}{\vert\alpha\vert^{\frac{7}{2}}} \ \left(\int \frac{\Vert s_\kappa f_{xx} \Vert_{L^2}}{\kappa^2} \ d\kappa\right)^{\frac{1}{2}}  \ d\alpha\\
&\lesssim& \Vert f \Vert_{\dot H^4} \Vert f \Vert_{\dot B^3_{\infty,1}} \left(  \Vert f \Vert^2_{\dot H^\frac{7}{4}} + \left(\int \frac{\Vert s_\alpha f\Vert^2_{L^\infty}}{\alpha^4} \ d\alpha \right)^{\frac{1}{2}} \left(\int \frac{\Vert \delta_\alpha f_x\Vert^2_{L^2}}{\alpha^2} \ d\alpha \right)^{\frac{1}{2}} \right) \\
&+& \Vert f \Vert_{\dot H^4} \Vert f \Vert_{\dot H^{\frac{5}{2}}} \left( \Vert f \Vert^2_{\dot H^{\frac{9}{4}}} + \left(\int \frac{\Vert s_\alpha f\Vert^2_{L^\infty}}{\vert\alpha\vert^{\frac{9}{2}}} \ d\alpha \right)^{\frac{1}{2}} \left(\int \frac{\Vert \delta_\alpha f_x\Vert^2_{L^\infty}}{\vert\alpha\vert^{\frac{5}{2}}} \ d\alpha \right)^{\frac{1}{2}}\right) 
\end{eqnarray*}
Since  $\Vert f \Vert_{\dot B^3_{\infty,1}} \lesssim\Vert f \Vert^{\frac{4}{5}}_{\dot B^{\frac{7}{2}}_{\infty,\infty}} \Vert f \Vert^{\frac{1}{5}}_{\dot B^{1}_{\infty,\infty}}$ and  $\dot H^4 \hookrightarrow \dot B^{\frac{7}{2}}_{\infty,\infty}$, we find 
\begin{eqnarray*}
\mathcal{U}_{4,1,12,sing}&\lesssim& \Vert f \Vert_{\dot H^4} \Vert f \Vert^{\frac{4}{5}}_{\dot B^{\frac{7}{2}}_{\infty,\infty}} \Vert f \Vert^{\frac{1}{5}}_{\dot B^{1}_{\infty,\infty}} \left(  \Vert f \Vert^2_{\dot H^\frac{7}{4}} + \Vert f \Vert_{\dot H^2} \Vert f \Vert_{\dot H^\frac{3}{2}} \right) \\
&+& \Vert f \Vert_{\dot H^4} \Vert f \Vert_{\dot H^{\frac{5}{2}}} \Vert f \Vert^2_{\dot H^{\frac{9}{4}}} \\
&\lesssim& \Vert f \Vert_{\dot H^4} \Vert f \Vert^{\frac{4}{5}}_{\dot H^{4}} \Vert f \Vert^{\frac{1}{5}}_{\dot H^\frac{3}{2}} \left(  \Vert f \Vert^{\frac{1}{5}}_{\dot H^4} \Vert f \Vert^{\frac{4}{5}}_{\dot H^\frac{3}{2}} +  \Vert f \Vert^{\frac{1}{5}}_{\dot H^4} \Vert f \Vert^{\frac{4}{5}}_{\dot H^\frac{3}{2}} \Vert f \Vert_{\dot H^\frac{3}{2}}\right)  \\
&+& \Vert f \Vert_{\dot H^4} \Vert f \Vert^{\frac{2}{5}}_{\dot H^4} \Vert f \Vert^{\frac{3}{5}}_{\dot H^\frac{3}{2}}  \Vert f \Vert^{\frac{3}{5}}_{\dot H^4} \Vert f \Vert^{\frac{2}{5}}_{\dot H^\frac{3}{2}} \\
\end{eqnarray*}
Finally, we find that
\begin{eqnarray*}
\mathcal{U}_{4,1,12,sing}&\lesssim& \Vert f \Vert^2_{\dot H^4}\left(\Vert f \Vert_{\dot H^\frac{3}{2}}+\Vert f \Vert^2_{\dot H^\frac{3}{2}}\right)
  \end{eqnarray*}
Collecting all the $\mathcal{U}_{4,1,i}$ terms for $i=1,...12$, we have obtain the estimate claimed in Lemma \ref{u41}. \\

Then, we turn to the estimate of $\mathcal{U}_{4,2}$, more precisely, we want to estimate
\begin{eqnarray*}
\mathcal{U}_{4,2}&=&\frac{1}{2\pi} \int \int\int_{0}^{\infty}\int_{0}^{\infty} \ e^{-\gamma-\sigma} \ \Lambda^{4}\mathcal{H}f \ \partial^2_{\alpha}\left[\frac{1}{\alpha}\sin(\frac{\gamma}{2}D_\alpha f)\cos(\frac{\gamma}{2}S_\alpha f)\frac{1}{\alpha} \int_{0}^{\alpha} s_{\eta} f_{x}    \ d\eta \right]\\
&& \  \times\ \partial_x \left( {  \delta_\alpha f_x } \cos(\sigma \tau_{\alpha} f_x )\cos(\arctan(\tau_{\alpha} f_{x})) )\right) \ d\eta  \ d\gamma \ d\sigma \ dx \ d\alpha \\
\end{eqnarray*}
When the differentiation in $x$ hits either  $ \cos(\sigma \tau_{\alpha} f_x )$ or   $\cos(\arctan(\tau_{\alpha} f_{x}))$ both will be controlled by $\vert \tau_{\alpha} f_{xx}\vert$.  Otherwise, it hits . By following the estimates done for the term $T_{4,1}$ and sharing the $L^2$ (resp $L^4-L^4$) estimate for each term into $L^{4}-L^4$ (resp $L^{6}-L^6-L^6$) by putting $\vert \tau_{\alpha} f_{xx}\vert$ in $L^4$ (resp $L^{6}$). It is clear, by Sobolev embedding, the homogeneous Besov space would m that one obtains a Sobolev space having 1/4 derivative more resp 1/3 derivative more). Let 
$\mathcal{U}_{4,2,1}$ be the first term of the 12 terms in the decomposition of $\mathcal{U}_{4,2}$ (note that the decomposition is obtained in the same manner as the decomposition of $T_{4,1}$). For this first term, we have
\begin{eqnarray*}
\mathcal{U}_{4,2,1}&=&\frac{1}{\pi} \int \int\int_{0}^{\infty}\int_{0}^{\infty} \ e^{-\gamma-\sigma} \ \Lambda^{4} \mathcal{H}f \ \frac{1}{\alpha^3}\sin(\frac{\gamma}{2}D_\alpha f)\cos(\frac{\gamma}{2}S_\alpha f)\frac{1}{\alpha} \int_{0}^{\alpha} s_{\eta} f_{x}    \ d\eta \\
&& \  \times\ \delta_\alpha f_x \ \partial_x\left(   \cos(\sigma \tau_{\alpha} f_x )\cos(\arctan(\tau_{\alpha} f_{x})) \right)   \ d\gamma \ d\sigma \ d\alpha \ dx \\
&\lesssim& \Vert f \Vert_{\dot H^4}  \int \frac{1}{\alpha^4} \Vert \delta_\alpha f_x \Vert_{L^\infty} \int_{0}^{\alpha} \frac{\Vert s_{\eta} f_{x} \Vert_{L^4}}{\eta^2} \eta^2 \ \Vert \tau_{\alpha} f_{xx}\Vert_{L^4}    \ d\eta \ d\alpha \\
&\lesssim& \Vert f \Vert_{\dot H^4} \Vert f \Vert_{\dot H^{\frac{9}{4}}}\int \frac{1}{\alpha^4} \Vert \delta_\alpha f_x \Vert_{L^\infty}  \left(\int \frac{\Vert s_{\eta} f_{x} \Vert^2_{L^4}}{\eta^4}    \ d\eta \right)^{\frac{1}{2}} \left(\int_{0}^{\alpha} \eta^4 \ d\eta \right)^{\frac{1}{2}} \ d\alpha \\
&\lesssim& \Vert f \Vert_{\dot H^4} \Vert f \Vert_{\dot H^\frac{9}{4}} \Vert f \Vert_{\dot H^\frac{11}{4}} \int  \Vert \delta_\alpha f_x \Vert_{L^\infty} \frac{1}{\vert\alpha\vert^{\frac{3}{2}}} \ d\alpha \\
&\lesssim& \Vert f \Vert_{\dot H^4} \Vert f \Vert^{\frac{3}{10}}_{\dot H^4}\Vert f \Vert^{\frac{7}{10}}_{\dot H^\frac{3}{2}} \Vert f \Vert^{\frac{1}{2}}_{\dot H^4}\Vert f \Vert^{\frac{1}{2}}_{\dot H^\frac{3}{2}} \Vert f \Vert_{\dot B^\frac{3}{2}_{\infty,1}} \\
\end{eqnarray*}
To conclude, it suffices to use $\Vert f \Vert_{\dot B^\frac{3}{2}_{\infty,1}} \lesssim \Vert f \Vert^{\frac{4}{5}}_{\dot B^1_{\infty,\infty}} \Vert f \Vert^{\frac{1}{5}}_{\dot B^{7/2}_{\infty,\infty}}\lesssim \Vert f \Vert^{\frac{4}{5}}_{\dot H^{\frac{3}{2}}} \Vert f \Vert^{\frac{1}{5}}_{\dot H^4}.$
Therefore,
\begin{eqnarray*}
\mathcal{U}_{4,2,1}&\lesssim& \Vert f \Vert^2_{\dot H^4} \Vert f \Vert^2_{\dot H^\frac{3}{2}}
\end{eqnarray*}
Hence we get the same control as $T_{4,1,1}$ and this is expected as we only need to balance the regularity in $x$ and therefore it is easy to adapt the inequalities obtained in the estimates $T_{4,1}$ to the case of $\mathcal{U}_{4,2}$. Hence,
\begin{eqnarray*}
\mathcal{U}_{4,2}&\lesssim& \Vert f \Vert^2_{\dot H^4} \Vert f \Vert^2_{\dot H^\frac{3}{2}}
\end{eqnarray*}
Finally, 
\begin{eqnarray*}
 \mathcal{U}_{4,2} \lesssim \Vert f \Vert^2_{\dot H^3} \left( \Vert f \Vert^2_{\dot H^{\frac{3}{2}}}+ \Vert f \Vert^3_{\dot H^{\frac{3}{2}}} \right)
 \end{eqnarray*}
 $\bullet$ {Estimate of $\mathcal{U}_{4,3}$} \\
 
 To estimate $\mathcal{U}_{4,3}$, that is
 \begin{eqnarray*}
\mathcal{U}_{4,3}&=& \frac{1}{2\pi} \int \int\int_{0}^{\infty}\int_{0}^{\infty} \ \sigma e^{-\gamma-\sigma} \ \Lambda^{4}\mathcal{H}f \ \partial_{x}\partial_{\alpha}\left[\frac{1}{\alpha}\sin(\frac{\gamma}{2}D_\alpha f)\cos(\frac{\gamma}{2}S_\alpha f)\frac{1}{\alpha} \int_{0}^{\alpha} s_{\eta} f_{x}    \ d\eta \right]\\
&& \  \times\ {  \delta_\alpha f_x } \left(\partial_{\alpha}\tau_{\alpha}f_x\right) \sin(\sigma \tau_{\alpha} f_x )\cos(\arctan(\tau_{\alpha} f_{x})) \ d\eta  \ d\gamma \ d\sigma \ dx \ d\alpha 
\end{eqnarray*}
We first integrate by parts in $\alpha$, one gets
\begin{eqnarray*}
\mathcal{U}_{4,3}&=& \frac{1}{4\pi} \int \int\int_{0}^{\infty}\int_{0}^{\infty} \gamma \sigma e^{-\gamma-\sigma} \ \Lambda^{4}\mathcal{H}f \ \partial_{x}D_\alpha f \ \partial_{\alpha}\left[\frac{1}{\alpha}\cos(\frac{\gamma}{2}D_\alpha f)\cos(\frac{\gamma}{2}S_\alpha f)\frac{1}{\alpha} \int_{0}^{\alpha} s_{\eta} f_{x}    \ d\eta \right]\\
&& \  \times\ {  \delta_\alpha f_x } \left(\partial_{\alpha}\tau_{\alpha}f_x\right) \sin(\sigma \tau_{\alpha} f_x )\cos(\arctan(\tau_{\alpha} f_{x})) \ d\eta  \ d\gamma \ d\sigma \ dx \ d\alpha \\
&+& \frac{1}{4\pi} \int \int\int_{0}^{\infty}\int_{0}^{\infty} \gamma \sigma e^{-\gamma-\sigma} \ \Lambda^{4}\mathcal{H}f \ \partial_{\alpha}D_\alpha f_x \ \left[\frac{1}{\alpha}\cos(\frac{\gamma}{2}D_\alpha f)\cos(\frac{\gamma}{2}S_\alpha f)\frac{1}{\alpha} \int_{0}^{\alpha} s_{\eta} f_{x}    \ d\eta \right]\\
&& \  \times\ {  \delta_\alpha f_x } \left(\partial_{\alpha}\tau_{\alpha}f_x\right) \sin(\sigma \tau_{\alpha} f_x )\cos(\arctan(\tau_{\alpha} f_{x})) \ d\eta  \ d\gamma \ d\sigma \ dx \ d\alpha \\
&-& \frac{1}{4\pi} \int \int\int_{0}^{\infty}\int_{0}^{\infty} \gamma \sigma e^{-\gamma-\sigma} \ \Lambda^{4}\mathcal{H}f \ \partial_{x}S_\alpha f \ \partial_{\alpha}\left[\frac{1}{\alpha}\sin(\frac{\gamma}{2}D_\alpha f)\sin(\frac{\gamma}{2}S_\alpha f)\frac{1}{\alpha} \int_{0}^{\alpha} s_{\eta} f_{x}    \ d\eta \right]\\
&& \  \times\ {  \delta_\alpha f_x } \left(\partial_{\alpha}\tau_{\alpha}f_x\right) \sin(\sigma \tau_{\alpha} f_x )\cos(\arctan(\tau_{\alpha} f_{x})) \ d\eta  \ d\gamma \ d\sigma \ dx \ d\alpha \\
&-& \frac{1}{4\pi} \int \int\int_{0}^{\infty}\int_{0}^{\infty} \gamma \sigma e^{-\gamma-\sigma} \ \Lambda^{4}\mathcal{H}f \ \partial_{\alpha}S_\alpha f_x \ \left[\frac{1}{\alpha}\sin(\frac{\gamma}{2}D_\alpha f)\sin(\frac{\gamma}{2}S_\alpha f)\frac{1}{\alpha} \int_{0}^{\alpha} s_{\eta} f_{x}    \ d\eta \right]\\
&& \  \times\ {  \delta_\alpha f_x } \left(\partial_{\alpha}\tau_{\alpha}f_x\right) \sin(\sigma \tau_{\alpha} f_x )\cos(\arctan(\tau_{\alpha} f_{x})) \ d\eta  \ d\gamma \ d\sigma \ dx \ d\alpha \\
&+&\frac{1}{2\pi} \int \int\int_{0}^{\infty}\int_{0}^{\infty} \ \sigma e^{-\gamma-\sigma} \ \Lambda^{4}\mathcal{H}f \ \partial_{\alpha}\left[\frac{1}{\alpha}\sin(\frac{\gamma}{2}D_\alpha f)\cos(\frac{\gamma}{2}S_\alpha f)\frac{1}{\alpha} \int_{0}^{\alpha} s_{\eta} f_{xx}    \ d\eta \right]\\
&& \  \times\ {  \delta_\alpha f_x } \left(\partial_{\alpha}\tau_{\alpha}f_x\right) \sin(\sigma \tau_{\alpha} f_x )\cos(\arctan(\tau_{\alpha} f_{x})) \ d\eta  \ d\gamma \ d\sigma \ dx \ d\alpha \\
&=&\sum_{i=1}^5 \mathcal{U}_{4,3,i}
 \end{eqnarray*}
 
 We shall estimate the $\mathcal{U}_{4,3,i}$ for $i=1,...,5$. \\
 
$\bullet$ \ {Estimate of $\mathcal{U}_{4,3,1}$}
 
 \begin{eqnarray*}
\mathcal{U}_{4,3,1}&=&-\frac{1}{2\pi} \int \int\int_{0}^{\infty}\int_{0}^{\infty} \int_{0}^{\alpha}\ \sigma e^{-\gamma-\sigma} \ \Lambda^{4} \mathcal{H}f \ 
\partial_{x}D_\alpha f \ \frac{1}{\alpha^3}\sin(\frac{\gamma}{2}D_\alpha f)\cos(\frac{\gamma}{2}S_\alpha f)  s_{\eta} f_{x}     \\
&& \  \times\ {  \delta_\alpha f_x } \left(\partial_{\alpha}\tau_{\alpha}f_x\right) \sin(\sigma \tau_{\alpha} f_x )\cos(\arctan(\tau_{\alpha} f_{x})) \ d\eta  \ d\gamma \ d\sigma \ dx \ d\alpha\\
&+& \frac{1}{4\pi} \int \int\int_{0}^{\infty}\int_{0}^{\infty} \int_{0}^{\alpha}\ \sigma \gamma e^{-\gamma-\sigma} \ \Lambda^{4} \mathcal{H}f \ 
\partial_{x}D_\alpha f \ \frac{1}{\alpha^2}\left(\partial_{\alpha}D_\alpha f \right)\cos(\frac{\gamma}{2}D_\alpha f)\cos(\frac{\gamma}{2}S_\alpha f)  s_{\eta} f_{x}    \\
&& \  \times\ {  \delta_\alpha f_x } \left(\partial_{\alpha}\tau_{\alpha}f_x\right) \sin(\sigma \tau_{\alpha} f_x )\cos(\arctan(\tau_{\alpha} f_{x})) \ d\eta  \ d\gamma \ d\sigma \ dx \ d\alpha\\
&-& \frac{1}{4\pi} \int \int\int_{0}^{\infty}\int_{0}^{\infty} \int_{0}^{\alpha} \sigma \gamma e^{-\gamma-\sigma} \ \Lambda^{4} \mathcal{H}f \ \partial_{x}D_\alpha f \  \frac{1}{\alpha^2}\left(\partial_{\alpha}S_\alpha f \right)\sin(\frac{\gamma}{2}D_\alpha f)\sin(\frac{\gamma}{2}S_\alpha f)  s_{\eta} f_{x}     \\
&& \  \times\ {  \delta_\alpha f_x } \left(\partial_{\alpha}\tau_{\alpha}f_x\right) \sin(\sigma \tau_{\alpha} f_x )\cos(\arctan(\tau_{\alpha} f_{x})) \ d\eta  \ d\gamma \ d\sigma \ dx \ d\alpha\\
&+&\frac{1}{2\pi} \int \int\int_{0}^{\infty}\int_{0}^{\infty} \ \sigma e^{-\gamma-\sigma} \ \Lambda^{4} \mathcal{H}f \ \partial_{x}D_\alpha f \  \frac{1}{\alpha}\left(\partial_{\alpha}D_\alpha f \right)\sin(\frac{\gamma}{2}D_\alpha f)\cos(\frac{\gamma}{2}S_\alpha f)\\
&& \  \times\ {  \delta_\alpha f_x } \left(\partial_{\alpha}\tau_{\alpha}f_x\right) \sin(\sigma \tau_{\alpha} f_x )\cos(\arctan(\tau_{\alpha} f_{x}))   \ d\gamma \ d\sigma \ dx \ d\alpha\\
&=&\mathcal{U}_{4,3,1,1}+\mathcal{U}_{4,3,1,2}+\mathcal{U}_{4,3,1,3}+\mathcal{U}_{4,3,1,4}
\end{eqnarray*}

To estimate $\mathcal{U}_{4,3,1,1}$, we observe that
\begin{eqnarray*}
\mathcal{U}_{4,3,1,1} &\lesssim& \Vert f \Vert_{\dot H^4} \Vert f \Vert_{\dot H^{\frac{9}{4}}} \int \frac{\Vert \delta_{\alpha} f_x\Vert^2_{L^{\infty}}}{\vert\alpha\vert^4} \int_0^{\alpha}  \Vert s_{\eta} f_{x} \Vert_{L^4} \ d\eta \ d\alpha \\
&\lesssim& \Vert f \Vert_{\dot H^4} \Vert f \Vert_{\dot H^{\frac{9}{4}}} \int \frac{\Vert \delta_{\alpha} f_x\Vert^2_{L^{\infty}}}{\vert\alpha\vert^{\frac{5}{2}}} \left(\int \frac{\Vert s_{\eta} f_{x} \Vert^2_{L^4}}{\eta^2} \ d\eta\right)^{\frac{1}{2}} \ d\alpha \\
&\lesssim& \Vert f \Vert_{\dot H^4} \Vert f \Vert^3_{\dot H^{\frac{9}{4}}} \Vert f \Vert_{\dot H^{\frac{7}{4}}} \\
&\lesssim& \Vert f \Vert_{\dot H^4} \Vert f \Vert^{\frac{9}{10}}_{\dot H^{4}} \Vert f \Vert^{\frac{21}{10}}_{\dot H^{\frac{3}{2}}} \Vert f \Vert^{\frac{1}{10}}_{\dot H^{4}} \Vert f \Vert^{\frac{9}{10}}_{\dot H^{\frac{3}{2}}} \\
&\lesssim& \Vert f \Vert^2_{\dot H^4} \Vert f \Vert^{3}_{\dot H^{\frac{3}{2}}}
\end{eqnarray*}

As $\partial_{\alpha}D_\alpha f= - \frac{s_\alpha f_x}{\alpha}-\frac{1}{\alpha^2} \ {\int_0^\alpha s_\kappa f_x \ d\kappa}$, we find

\begin{eqnarray*}
\mathcal{U}_{4,3,1,2} &\lesssim& \Vert f \Vert_{\dot H^4}\Vert f \Vert_{\dot H^{2}} \sup_{\alpha \in \mathbb R} \frac{\Vert s_{\alpha} f_x \Vert_{L^\infty}}{\vert \alpha \vert} \int \frac{\Vert \delta_{\alpha} f_x\Vert^2_{L^{\infty}}}{\vert\alpha\vert^3}  \int_0^\alpha \Vert s_{\eta} f_x \Vert_{L^\infty} \ d\eta  \ d\alpha \\
&+& \Vert f \Vert_{\dot H^4} \Vert f \Vert_{\dot H^2}  \int \frac{\Vert \delta_{\alpha} f_x\Vert^2_{L^{\infty}}}{\vert\alpha\vert^4}  \int_0^\alpha \Vert s_{\eta} f_x \Vert_{L^\infty} \ d\eta \ \int_0^\alpha \Vert s_{\kappa} f_x \Vert_{L^\infty} \ d\kappa \ d\alpha \\
 &\lesssim& \Vert f \Vert_{\dot H^4}\Vert f \Vert_{\dot H^{2}} \sup_{\alpha \in \mathbb R} \frac{\Vert s_{\alpha} f_x \Vert_{L^\infty}}{\vert \alpha \vert} \int \frac{\Vert \delta_{\alpha} f_x\Vert^2_{L^{\infty}}}{\vert\alpha\vert^{\frac{3}{2}}}  d\alpha \left(\int \frac{\Vert s_{\eta} f_x \Vert^2_{L^\infty}}{\vert \eta\vert^2} \ d\eta \ \right)^{\frac{1}{2}} \\
 &+& \Vert f \Vert_{\dot H^4} \Vert f \Vert_{\dot H^2}  \int \frac{\Vert \delta_{\alpha} f_x\Vert^2_{L^{\infty}}}{\vert\alpha\vert^2} \ d\alpha \  \left(\int \frac{\Vert s_{\eta} f_x \Vert^2_{L^\infty}}{\vert \eta\vert^2} \ d\eta \ \right)^{\frac{1}{2}} \left(\int \frac{\Vert s_{\kappa} f_x \Vert^2_{L^\infty}}{\vert \kappa\vert^2} \ d\kappa \ \right)^{\frac{1}{2}}  \\
 &\lesssim& \Vert f \Vert_{\dot H^4}\Vert f \Vert_{\dot H^{2}} \left(\Vert f \Vert_{\dot B^{2}_{\infty,\infty}} \Vert f \Vert^2_{\dot B^{\frac{5}{4}}_{\infty,2}} \Vert f \Vert_{\dot B^{\frac{3}{2}}_{\infty,2}} + \Vert f \Vert^2_{\dot B^{\frac{3}{2}}_{\infty,2}}\Vert f \Vert^2_{\dot H^{2}} \right) \\
 &\lesssim& \Vert f \Vert_{\dot H^4} \left(\Vert f \Vert^2_{\dot H^{2}} \Vert f \Vert^{\frac{3}{5}}_{\dot B^{1}_{\infty,\infty}} \Vert f \Vert^{\frac{2}{5}}_{\dot B^{\frac{7}{2}}_{\infty,\infty}} \Vert f \Vert^2_{\dot H^{\frac{7}{4}}}+\Vert f \Vert^5_{\dot H^{2}}\right)\\
 &\lesssim& \Vert f \Vert_{\dot H^4} \left( \Vert f \Vert^{\frac{8}{5}}_{\dot H^{\frac{3}{2}}}  \Vert f \Vert^{\frac{2}{5}}_{\dot H^4} \Vert f \Vert^{\frac{3}{5}}_{\dot H^{\frac{3}{2}}}  \Vert f \Vert^{\frac{2}{5}}_{\dot H^4}  \Vert f \Vert^{\frac{9}{5}}_{\dot H^{\frac{3}{2}}}  \Vert f \Vert^{\frac{1}{5}}_{\dot H^4} + \Vert f \Vert_{\dot H^4} \Vert f \Vert^{4}_{\dot H^{\frac{3}{2}}}\right)
 \end{eqnarray*}
 Therefore,
 \begin{eqnarray*}
\mathcal{U}_{4,3,1,2} &\lesssim& \Vert f \Vert_{\dot H^4} \Vert f \Vert^{4}_{\dot H^{\frac{3}{2}}} 
\end{eqnarray*}

Since,
\begin{eqnarray*}
\partial_{\alpha}S_\alpha f &=&\frac{f_x(x-\alpha)-f_x(x+\alpha)}{\alpha}+\left(\frac{f(x-\alpha)+f(x+\alpha)-2f(x)}{\alpha^2}\right) \\
&=& \bar\Delta_\alpha f_x-\Delta_\alpha f_x-\frac{s_\alpha f}{\alpha^2}
\end{eqnarray*}
Then, 
\begin{eqnarray*}
\mathcal{U}_{4,3,1,3}&=&- \frac{1}{4\pi} \int \int\int_{0}^{\infty}\int_{0}^{\infty} \int_{0}^{\alpha} \sigma \gamma e^{-\gamma-\sigma} \ \Lambda^{4} \mathcal{H}f \ \partial_{x}D_\alpha f \  \frac{1}{\alpha^2}\left(\bar\Delta_\alpha f_x-\Delta_\alpha f_x-\frac{s_\alpha f}{\alpha^2} \right) \\
&& \sin(\frac{\gamma}{2}D_\alpha f)\sin(\frac{\gamma}{2}S_\alpha f)  s_{\eta} f_{x}     {  \delta_\alpha f_x } \left(\partial_{\alpha}\tau_{\alpha}f_x\right)\sin(\sigma \tau_{\alpha} f_x ) \cos(\arctan(\tau_{\alpha} f_{x}))\\
&&   \ d\eta  \ d\gamma \ d\sigma \ dx \ d\alpha\\
\end{eqnarray*}
By omitting to estimate the term involving $\bar\Delta_\alpha f_x$ as it is similar to $\Delta_\alpha f_x$, one finds
\begin{eqnarray*}
\mathcal{U}_{4,3,1,3}&\lesssim& \Vert f \Vert_{\dot H^4} \Vert f \Vert_{\dot H^2} \int \frac{\Vert \delta_{\alpha} f_x \Vert^2_{L^\infty}}{\vert \alpha \vert^{\frac{5}{2}}} \left( \int_0^\alpha \frac{\Vert s_\eta f_x \Vert^2_{L^\infty}}{\vert \eta \vert^2} \ d\eta \right)^{\frac{1}{2}} \ d\alpha \\
&+&  \Vert f \Vert_{\dot H^4} \Vert f \Vert_{\dot H^2}\int \frac{\Vert \delta_{\alpha} f_x \Vert_{L^\infty} \Vert s_{\alpha} f \Vert_{L^\infty}}{\vert \alpha \vert^{\frac{7}{2}}} \left( \int_0^\alpha \frac{\Vert s_\eta f_x \Vert^2_{L^\infty}}{\vert \eta \vert^2} \ d\eta \right)^{\frac{1}{2}} \ d\alpha \\
&\lesssim& \Vert f \Vert_{\dot H^4} \Vert f \Vert^2_{\dot H^2} \left( \Vert f \Vert^2_{\dot B^{\frac{7}{4}}_{\infty,2}} + \left( \int \frac{\Vert \delta_\alpha f_x \Vert^2_{L^\infty}}{\vert \alpha \vert^{\frac{5}{2}}} \ d\alpha \right)^{\frac{1}{2}}  \left( \int \frac{\Vert s_\alpha f \Vert^2_{L^\infty}}{\vert \alpha \vert^{\frac{9}{2}}} \ d\alpha \right)^{\frac{1}{2}} \right) \\
&\lesssim& \Vert f \Vert_{\dot H^4} \Vert f \Vert^2_{\dot H^2} \left( \Vert f \Vert^2_{\dot H^{\frac{9}{4}}} + \Vert f \Vert^2_{\dot B^{\frac{7}{4}}_{\infty,2}}   \right)\\
&\lesssim& \Vert f \Vert_{\dot H^4}  \Vert f \Vert^{\frac{8}{5}}_{\dot H^{\frac{3}{2}}}  \Vert f \Vert^{\frac{2}{5}}_{\dot H^4} \Vert f \Vert^{\frac{7}{5}}_{\dot H^{\frac{3}{2}}}  \Vert f \Vert^{\frac{3}{5}}_{\dot H^4} \\
&\lesssim& \Vert f \Vert^2_{\dot H^4}  \Vert f \Vert^{3}_{\dot H^{\frac{3}{2}}}
\end{eqnarray*}

Finally, we have obtained

\begin{eqnarray*}
\mathcal{U}_{4,3,1,3}&\lesssim&  \Vert f \Vert^2_{\dot H^4} \left(  \Vert f \Vert^{3}_{\dot H^{\frac{3}{2}}} + \Vert f \Vert^{4}_{\dot H^{\frac{3}{2}}} \right)
\end{eqnarray*}

$\bullet$ \ {Estimate of $\mathcal{U}_{4,3,2}$} \\

Recall that 
\begin{eqnarray*}
\mathcal{U}_{4,3,2}&=& \frac{1}{4\pi} \int \int\int_{0}^{\infty}\int_{0}^{\infty} \gamma \sigma e^{-\gamma-\sigma} \ \Lambda^{4}\mathcal{H}f \ \partial_{\alpha}D_\alpha f_x \ \left[\frac{1}{\alpha}\cos(\frac{\gamma}{2}D_\alpha f)\cos(\frac{\gamma}{2}S_\alpha f)\frac{1}{\alpha} \int_{0}^{\alpha} s_{\eta} f_{x}    \ d\eta \right]\\
&& \  \times\ {  \delta_\alpha f_x } \left(\partial_{\alpha}\tau_{\alpha}f_x\right) \sin(\sigma \tau_{\alpha} f_x )\cos(\arctan(\tau_{\alpha} f_{x})) \ d\eta  \ d\gamma \ d\sigma \ dx \ d\alpha \\
\end{eqnarray*}
Since we have $\partial_{\alpha}D_\alpha f_x= - \frac{s_\alpha f_{xx}}{\alpha}-\frac{1}{\alpha^2} \ {\int_0^\alpha s_\kappa f_{xx} \ d\kappa}$,
we find that 
\begin{eqnarray*}
\mathcal{U}_{4,3,2}&\lesssim& \Vert f \Vert_{\dot H^4}\Vert f \Vert_{\dot H^2} \int \frac{\Vert \delta_{\alpha} f_x \Vert_{L^\infty} \Vert s_{\alpha} f_{xx} \Vert_{L^\infty}}{\vert \alpha\vert^3} \int_0^\alpha \Vert s_{\eta}f_x \Vert_{L^\infty} \ d\eta  \\
 &+&\Vert f \Vert_{\dot H^4}\Vert f \Vert_{\dot H^2}\int \frac{\Vert \delta_{\alpha} f_x \Vert_{L^\infty} }{\vert \alpha\vert^4} \int_0^\alpha \Vert s_{\eta}f_x \Vert_{L^\infty} \ d\eta \ \int_0^\alpha \Vert s_\kappa f_{xx} \Vert_{L^\infty} \ d\kappa \\
 &\lesssim& \Vert f \Vert_{\dot H^4}\Vert f \Vert_{\dot H^2} \int \frac{\Vert \delta_{\alpha} f_x \Vert_{L^\infty} \Vert s_{\alpha} f_{xx} \Vert_{L^\infty}}{\vert \alpha\vert^{\frac{7}{4}}} \left(\int \frac{\Vert s_{\eta}f_x \Vert^2_{L^\infty}}{\vert \eta\vert^{\frac{3}{2}}} \ d\eta \right)^{\frac{1}{2}} \ d\alpha  \\
 &+&\Vert f \Vert_{\dot H^4}\Vert f \Vert_{\dot H^2}\int \frac{\Vert \delta_{\alpha} f_x \Vert_{L^\infty} }{\vert \alpha\vert^{\frac{3}{2}}} \left(\int \frac{\Vert s_{\eta}f_x \Vert^2_{L^\infty}}{\vert \eta\vert^{\frac{3}{2}}} \ d\eta \right)^{\frac{1}{2}} \ \left(\int \frac{\Vert s_\kappa f_{xx} \Vert^2_{L^\infty}}{\vert \eta\vert^{\frac{3}{2}}} \ d\kappa \right)^{\frac{1}{2}} \ d\alpha \\
 &\lesssim&\Vert f \Vert_{\dot H^4}\Vert f \Vert_{\dot H^2} \Vert f \Vert_{\dot H^\frac{7}{4}} \left( \int \frac{\Vert \delta_{\alpha} f_x \Vert^2_{L^\infty} }{\vert \alpha\vert^{2}} \ d\alpha \right)^{\frac{1}{2}}  \left( \int \frac{\Vert s_{\alpha} f_{xx} \Vert^2_{L^\infty}}{\vert \alpha\vert^{\frac{5}{2}}} \ d\alpha \right)^{\frac{1}{2}} \\
 &+&\Vert f \Vert_{\dot H^4}\Vert f \Vert_{\dot H^2} \Vert f_x \Vert_{\dot B^{\frac{1}{2}}_{\infty,1}} \Vert f \Vert_{\dot B^{\frac{5}{4}}_{\infty,2}} \Vert f \Vert_{\dot B^{\frac{9}{4}}_{\infty,2}} \\
 &\lesssim&\Vert f \Vert_{\dot H^4}\left(\Vert f \Vert_{\dot H^2}\Vert f \Vert_{\dot H^{\frac{7}{4}}}\Vert f \Vert_{\dot B^{\frac{3}{2}}_{\infty,2}} \Vert f \Vert_{\dot B^{\frac{9}{4}}_{\infty,2}} +\Vert f \Vert_{\dot H^2} \Vert f \Vert_{\dot B^{\frac{3}{2}}_{\infty,1}} \Vert f \Vert_{\dot H^{\frac{7}{4}}} \Vert f \Vert_{\dot H^{\frac{11}{4}}} \right) \\
 &\lesssim&\Vert f \Vert_{\dot H^4}\Vert f \Vert_{\dot H^{\frac{11}{4}}}\Vert f \Vert_{\dot H^{\frac{7}{4}}}\left(\Vert f \Vert^2_{\dot H^2}  +\Vert f \Vert_{\dot H^2} \Vert f \Vert^{\frac{4}{5}}_{\dot B^{1}_{\infty,\infty}} \Vert f \Vert^{\frac{1}{5}}_{\dot B^{\frac{7}{2}}_{\infty,\infty}}   \right)\\
 &\lesssim&\Vert f \Vert_{\dot H^4}\Vert f \Vert^{\frac{1}{2}}_{\dot H^4}\Vert f \Vert^{\frac{1}{2}}_{\dot H^\frac{3}{2}} \Vert f \Vert^{\frac{1}{10}}_{\dot H^4}\Vert f \Vert^{\frac{9}{10}}_{\dot H^\frac{3}{2}}\left( \Vert f \Vert^{\frac{8}{5}}_{\dot H^{\frac{3}{2}}}  \Vert f \Vert^{\frac{2}{5}}_{\dot H^4}
 +\Vert f \Vert^{\frac{4}{5}}_{\dot H^{\frac{3}{2}}}  \Vert f \Vert^{\frac{1}{5}}_{\dot H^4}
 \Vert f \Vert^{\frac{4}{5}}_{\dot H^{\frac{3}{2}}} \Vert f \Vert^{\frac{1}{5}}_{\dot H^{4}}  \right)
\end{eqnarray*}
Finally,
\begin{eqnarray*}
\mathcal{U}_{4,3,2}&\lesssim& \Vert f \Vert^2_{\dot H^4}\Vert f \Vert^{3}_{\dot H^{\frac{3}{2}}}
\end{eqnarray*}

$\bullet$ \ {Estimate of $\mathcal{U}_{4,3,3}$} \\

We have,
\begin{eqnarray*}
\mathcal{U}_{4,3,3}&=&-  \frac{1}{4\pi} \int \int\int_{0}^{\infty}\int_{0}^{\infty} \gamma \sigma e^{-\gamma-\sigma} \ \Lambda^{4}\mathcal{H}f \ \partial_{x}S_\alpha f \\
&& \partial_{\alpha}\left[\frac{1}{\alpha}\sin(\frac{\gamma}{2}D_\alpha f)\sin(\frac{\gamma}{2}S_\alpha f)\frac{1}{\alpha} \int_{0}^{\alpha} s_{\eta} f_{x}    \ d\eta \right]  \ d\gamma \ d\sigma \ dx \ d\alpha
\end{eqnarray*}

Note that this term is similar to $\mathcal{U}_{4,3,1}$ up to interchanging the term $\partial_{x}S_\alpha f$  and $\partial_{x}D_\alpha f$. As we only used the bounded $\Vert \partial_{x}D_\alpha f \Vert_{L^{\infty}} \lesssim {\Vert \delta_\alpha f_x \Vert_{L^\infty}}{\vert \alpha \vert^{-1}}$ and as this bounded is also obviously verified by $\partial_{x}S_\alpha f$, following the same estimates as $\mathcal{U}_{4,3,1}$,  we infer that
\begin{eqnarray*}
\mathcal{U}_{4,3,3}&\lesssim&  \Vert f \Vert^2_{\dot H^4} \left(  \Vert f \Vert^{3}_{\dot H^{\frac{3}{2}}} + \Vert f \Vert^{4}_{\dot H^{\frac{3}{2}}} \right)
\end{eqnarray*}

$\bullet$ {Estimate of $\mathcal{U}_{4,3,4}$}\\

We have

\begin{eqnarray*}
\mathcal{U}_{4,3,4}&=&- \frac{1}{4\pi} \int \int\int_{0}^{\infty}\int_{0}^{\infty} \gamma \sigma e^{-\gamma-\sigma} \ \Lambda^{4}\mathcal{H}f \ \partial_{\alpha}S_\alpha f_x \ \left[\frac{1}{\alpha}\sin(\frac{\gamma}{2}D_\alpha f)\sin(\frac{\gamma}{2}S_\alpha f)\frac{1}{\alpha} \int_{0}^{\alpha} s_{\eta} f_{x}    \ d\eta \right]\\
&& \  \times\ {  \delta_\alpha f_x } \left(\partial_{\alpha}\tau_{\alpha}f_x\right) \sin(\sigma \tau_{\alpha} f_x )\cos(\arctan(\tau_{\alpha} f_{x})) \ d\eta  \ d\gamma \ d\sigma \ dx \ d\alpha \\
\end{eqnarray*}
We use the identity 
$\partial_{\alpha}S_\alpha f_x = \bar\Delta_\alpha f_{xx}-\Delta_\alpha f_{xx}-\frac{s_\alpha f_x}{\alpha^2}$
and we may omit to treat the term involving $\bar\Delta_\alpha f_{xx}$ as it is similar to $\Delta_\alpha f_{xx}$, hence we need to estimate
\begin{eqnarray*}
\mathcal{U}_{4,3,4}&=& \frac{1}{4\pi} \int \int\int_{0}^{\infty}\int_{0}^{\infty} \gamma \sigma e^{-\gamma-\sigma} \ \Lambda^{4}\mathcal{H}f \ \left(\Delta_\alpha f_{xx}+\frac{s_\alpha f_x}{\alpha^2}\right) \\
&& \left[\frac{1}{\alpha}\sin(\frac{\gamma}{2}D_\alpha f)\sin(\frac{\gamma}{2}S_\alpha f)\frac{1}{\alpha} \int_{0}^{\alpha} s_{\eta} f_{x}    \ d\eta \right]\\
&& \  \times\ {  \delta_\alpha f_x } \left(\partial_{\alpha}\tau_{\alpha}f_x\right) \sin(\sigma \tau_{\alpha} f_x )\cos(\arctan(\tau_{\alpha} f_{x})) \ d\eta  \ d\gamma \ d\sigma \ dx \ d\alpha \\
\end{eqnarray*}
We find,
\begin{eqnarray*}
\mathcal{U}_{4,3,4}&\lesssim& \Vert f \Vert_{\dot H^4} \Vert f \Vert_{\dot H^2} \int \frac{\Vert \delta_{\alpha} f_{xx} \Vert_{L^\infty} \Vert \delta_{\alpha} f_x \Vert_{L^\infty}}{\vert \alpha\vert^3}  \int_{0}^{\alpha} \Vert s_{\eta} f_{x} \Vert_{L^\infty}    \ d\eta \ d\alpha \\
&+& \Vert f \Vert_{\dot H^4} \Vert f \Vert_{\dot H^2} \int \frac{\Vert s_{\alpha} f_{x} \Vert_{L^\infty} \Vert \delta_{\alpha} f_x \Vert_{L^\infty}}{\vert \alpha\vert^4}  \int_{0}^{\alpha} \Vert s_{\eta} f_{x} \Vert_{L^\infty}    \ d\eta \ d\alpha \\
&\lesssim& \Vert f \Vert_{\dot H^4} \Vert f \Vert_{\dot H^2} \int \frac{\Vert \delta_{\alpha} f_{xx} \Vert_{L^\infty} \Vert \delta_{\alpha} f_x \Vert_{L^\infty}}{\vert \alpha\vert^{\frac{3}{2}}}  \left(\int \frac{\Vert s_{\eta} f_{x} \Vert^2_{L^\infty}}{\alpha^2}    \ d\eta \right)^{\frac{1}{2}} \ d\alpha \\
&+& \Vert f \Vert_{\dot H^4} \Vert f \Vert_{\dot H^2} \int \frac{\Vert s_{\alpha} f_{x} \Vert_{L^\infty} \Vert \delta_{\alpha} f_x \Vert_{L^\infty}}{\vert \alpha\vert^\frac{5}{2}}  \left(\int \frac{\Vert s_{\eta} f_{x} \Vert^2_{L^\infty}}{\alpha^2}    \ d\eta \right)^{\frac{1}{2}} \ d\alpha \\
&\lesssim& \Vert f \Vert_{\dot H^4} \Vert f \Vert^2_{\dot H^2} \left(  \Vert f_{xx} \Vert_{\dot B^{\frac{1}{4}}_{\infty,2}} \Vert f_{x} \Vert_{\dot B^{\frac{1}{4}}_{\infty,2}} + \Vert f_{x} \Vert^2_{\dot B^{\frac{3}{4}}_{\infty,2}} \right)  
\end{eqnarray*}

Therefore,
\begin{eqnarray*}
\mathcal{U}_{4,3,4}&\lesssim& \Vert f \Vert_{\dot H^4} \Vert f \Vert^2_{\dot H^2} \left(  \Vert f\Vert_{\dot H^{\frac{11}{4}}} \Vert f\Vert_{\dot H^{\frac{7}{4}}} + \Vert f\Vert^2_{\dot H^{\frac{9}{4}}} \right)  \\
&\lesssim& \Vert f \Vert_{\dot H^4}  \Vert f \Vert^{\frac{8}{5}}_{\dot H^{\frac{3}{2}}}  \Vert f \Vert^{\frac{2}{5}}_{\dot H^4} \left(   \Vert f \Vert^{\frac{1}{2}}_{\dot H^{\frac{3}{2}}}  \Vert f \Vert^{\frac{1}{2}}_{\dot H^4} \Vert f \Vert^{\frac{1}{10}}_{\dot H^4}\Vert f \Vert^{\frac{9}{10}}_{\dot H^\frac{3}{2}} + 
\Vert f \Vert^{\frac{3}{5}}_{\dot H^4}\Vert f \Vert^{\frac{7}{5}}_{\dot H^\frac{3}{2}} \right)
\end{eqnarray*}
So that,
\begin{eqnarray*}
\mathcal{U}_{4,3,4}&\lesssim& \Vert f \Vert^2_{\dot H^4}\Vert f \Vert^{3}_{\dot H^\frac{3}{2}}
\end{eqnarray*}

$\bullet$ {Estimate of $\mathcal{U}_{4,3,5}$} \\

We have that
\begin{eqnarray*}
\mathcal{U}_{4,3,5}&=&\frac{1}{2\pi} \int \int\int_{0}^{\infty}\int_{0}^{\infty} \ \sigma e^{-\gamma-\sigma} \ \Lambda^{4}\mathcal{H}f \ \partial_{\alpha}\left[\frac{1}{\alpha}\sin(\frac{\gamma}{2}D_\alpha f)\cos(\frac{\gamma}{2}S_\alpha f)\frac{1}{\alpha} \int_{0}^{\alpha} s_{\eta} f_{xx}    \ d\eta \right]\\
&& \  \times\ {  \delta_\alpha f_x } \left(\partial_{\alpha}\tau_{\alpha}f_x\right) \sin(\sigma \tau_{\alpha} f_x )\cos(\arctan(\tau_{\alpha} f_{x})) \ d\eta  \ d\gamma \ d\sigma \ dx \ d\alpha \\
\end{eqnarray*}
When the derivatives in $\alpha$ hits  the integral in $\eta$ it is more convenient to consider its local version and we immediately  find $\partial_\alpha D_\alpha f_x$. Therefore, we find
\begin{eqnarray*}
\mathcal{U}_{4,3,5}&=&\frac{1}{2\pi} \int \int\int_{0}^{\infty}\int_{0}^{\infty} \ \sigma e^{-\gamma-\sigma} \ \Lambda^{4}\mathcal{H}f \ \partial_{\alpha}\left[\frac{1}{\alpha}\sin(\frac{\gamma}{2}D_\alpha f)\cos(\frac{\gamma}{2}S_\alpha f)\frac{1}{\alpha} \int_{0}^{\alpha} s_{\eta} f_{xx}    \ d\eta \right]\\
&& \  \times\ {  \delta_\alpha f_x } \left(\partial_{\alpha}\tau_{\alpha}f_x\right) \sin(\sigma \tau_{\alpha} f_x )\cos(\arctan(\tau_{\alpha} f_{x})) \ d\eta  \ d\gamma \ d\sigma \ dx \ d\alpha \\
&=&-\frac{1}{2\pi} \int \int\int_{0}^{\infty}\int_{0}^{\infty} \ \sigma e^{-\gamma-\sigma} \ \Lambda^{4}\mathcal{H}f \ \left[\frac{1}{\alpha^2}\sin(\frac{\gamma}{2}D_\alpha f)\cos(\frac{\gamma}{2}S_\alpha f)\frac{1}{\alpha} \int_{0}^{\alpha} s_{\eta} f_{xx}    \ d\eta \right]\\
&& \  \times\ {  \delta_\alpha f_x } \left(\partial_{\alpha}\tau_{\alpha}f_x\right) \sin(\sigma \tau_{\alpha} f_x )\cos(\arctan(\tau_{\alpha} f_{x})) \ d\eta  \ d\gamma \ d\sigma \ dx \ d\alpha \\
&+&\frac{1}{4\pi} \int \int\int_{0}^{\infty}\int_{0}^{\infty} \ \gamma\sigma e^{-\gamma-\sigma} \ \Lambda^{4}\mathcal{H}f \ \left[\frac{1}{\alpha}\left(\partial_\alpha D_\alpha f\right) \ \cos(\frac{\gamma}{2}D_\alpha f)\cos(\frac{\gamma}{2}S_\alpha f)\frac{1}{\alpha} \int_{0}^{\alpha} s_{\eta} f_{xx}    \ d\eta \right]\\
&& \  \times\ {  \delta_\alpha f_x } \left(\partial_{\alpha}\tau_{\alpha}f_x\right) \sin(\sigma \tau_{\alpha} f_x )\cos(\arctan(\tau_{\alpha} f_{x})) \ d\eta  \ d\gamma \ d\sigma \ dx \ d\alpha \\
&-&\frac{1}{4\pi} \int \int\int_{0}^{\infty}\int_{0}^{\infty} \ \gamma\sigma e^{-\gamma-\sigma} \ \Lambda^{4}\mathcal{H}f \ \left[\frac{1}{\alpha}\left(\partial_\alpha S_\alpha f\right) \sin(\frac{\gamma}{2}D_\alpha f)\sin(\frac{\gamma}{2}S_\alpha f)\frac{1}{\alpha} \int_{0}^{\alpha} s_{\eta} f_{xx}    \ d\eta \right]\\
&& \  \times\ {  \delta_\alpha f_x } \left(\partial_{\alpha}\tau_{\alpha}f_x\right) \sin(\sigma \tau_{\alpha} f_x )\cos(\arctan(\tau_{\alpha} f_{x})) \ d\eta  \ d\gamma \ d\sigma \ dx \ d\alpha \\
&+&\frac{1}{2\pi} \int \int\int_{0}^{\infty}\int_{0}^{\infty} \ \sigma e^{-\gamma-\sigma} \ \Lambda^{4}\mathcal{H}f \ 
\left[\frac{1}{\alpha}\sin(\frac{\gamma}{2}D_\alpha f)\cos(\frac{\gamma}{2}S_\alpha f) \ \partial_\alpha D_\alpha f_x \right]\\
&& \  \times\ {  \delta_\alpha f_x } \left(\partial_{\alpha}\tau_{\alpha}f_x\right) \sin(\sigma \tau_{\alpha} f_x )\cos(\arctan(\tau_{\alpha} f_{x})) \ d\eta  \ d\gamma \ d\sigma \ dx \ d\alpha \\
&=&\sum_{i=1}^4 \mathcal{U}_{4,3,5,i}
\end{eqnarray*}
In order to estimate $\mathcal{U}_{4,3,5,1}$, we write that
\begin{eqnarray*}
\mathcal{U}_{4,3,5,1}&\lesssim&\Vert f \Vert_{\dot H^4} \Vert f \Vert_{\dot H^2} \int \frac{\Vert \delta_\alpha f_x \Vert_{L^\infty}}{\vert \alpha\vert^3} \ \int_{0}^{\alpha} \Vert s_{\eta} f_{xx} \Vert_{L^\infty}   \ d\eta \ d\alpha \\
&\lesssim& \Vert f \Vert_{\dot H^4} \Vert f \Vert_{\dot H^2} \int \frac{\Vert \delta_\alpha f_x \Vert_{L^\infty}}{\vert \alpha\vert^{\frac{3}{2}}} \ \left(\int \frac{\Vert s_{\eta} f_{xx} \Vert^2_{L^\infty} }{\eta^2}  \ d\eta \right)^{\frac{1}{2}} \ d\alpha \\
&\lesssim& \Vert f \Vert_{\dot H^4} \Vert f \Vert_{\dot H^3}  \Vert f \Vert_{\dot H^2} \Vert f \Vert_{\dot B^{\frac{3}{2}}_{\infty,1}}  \\
&\lesssim& \Vert f \Vert_{\dot H^4} \Vert f \Vert^{\frac{3}{5}}_{\dot H^4} \Vert f \Vert^{\frac{2}{5}}_{\dot H^\frac{3}{2}} \Vert f \Vert^{\frac{1}{5}}_{\dot H^4} \Vert f \Vert^{\frac{4}{5}}_{\dot H^\frac{3}{2}}   \Vert f \Vert^{\frac{4}{5}}_{\dot B^{1}_{\infty,\infty}} \Vert f \Vert^{\frac{1}{5}}_{\dot B^{\frac{7}{2}}_{\infty,\infty}} 
\end{eqnarray*}
Therefore,
\begin{eqnarray*}
\mathcal{U}_{4,3,5,1}&\lesssim&\Vert f \Vert^2_{\dot H^4} \Vert f \Vert^{2}_{\dot H^\frac{3}{2}}
\end{eqnarray*}
To estimate $\mathcal{U}_{4,3,5,2}$ and $\mathcal{U}_{4,3,5,3}$ , we use the fact that $\partial_\alpha D_\alpha f\approx   \frac{s_\alpha f_x}{\alpha}+\frac{1}{\alpha^2} \ {\int_0^\alpha s_\kappa f_x \ d\kappa}$ and that $\partial_\alpha S_\alpha f\approx \Delta_\alpha f_x+\frac{s_\alpha f}{\alpha^2}$, therefore we find the following control 

\begin{eqnarray*}
\mathcal{U}_{4,3,5,2} + \mathcal{U}_{4,3,5,3} &\lesssim& \Vert f \Vert_{\dot H^4} \Vert f \Vert_{\dot H^2} \int \frac{\Vert \delta_\alpha f_x \Vert_{L^\infty} \Vert s_\alpha f_x \Vert_{L^\infty}}{\vert \alpha\vert^3} \int_{0}^{\alpha} \Vert s_{\eta} f_{xx} \Vert_{L^\infty}   \ d\eta \ d\alpha \\
&+& \Vert f \Vert_{\dot H^4} \Vert f \Vert_{\dot H^2} \int \frac{\Vert \delta_\alpha f_x \Vert_{L^\infty} }{\vert \alpha\vert^4} \int_{0}^{\alpha} \Vert s_{\eta} f_{xx} \Vert_{L^\infty}   \ d\eta \ \int_{0}^{\alpha} \Vert s_{\kappa} f_{x} \Vert_{L^\infty}   \ d\kappa \ d\alpha  \\
&+& \Vert f \Vert_{\dot H^4} \Vert f \Vert_{\dot H^2} \int \frac{\Vert \delta_\alpha f_x \Vert^2_{L^\infty} }{\vert \alpha\vert^3} \int_{0}^{\alpha} \Vert s_{\eta} f_{xx} \Vert_{L^\infty}   \ d\eta \ d\alpha \\
&+& \Vert f \Vert_{\dot H^4} \Vert f \Vert_{\dot H^2} \int \frac{\Vert \delta_\alpha f_x \Vert_{L^\infty} \Vert s_\alpha f \Vert_{L^\infty} }{\vert \alpha\vert^4} \int_{0}^{\alpha} \Vert s_{\eta} f_{xx} \Vert_{L^\infty}   \ d\eta \ d\alpha \\
&\lesssim&\Vert f \Vert_{\dot H^4} \Vert f \Vert_{\dot H^2} \int \frac{\Vert \delta_\alpha f_x \Vert_{L^\infty} \Vert s_\alpha f_x \Vert_{L^\infty}}{\vert \alpha\vert^{\frac{3}{2}}} \left(\int \frac{\Vert s_{\eta} f_{xx} \Vert^2_{L^\infty}}{\alpha^2}   \ d\eta \right)^{\frac{1}{2}} \ d\alpha \\
&+& \Vert f \Vert_{\dot H^4} \Vert f \Vert_{\dot H^2} \int \frac{\Vert \delta_\alpha f_x \Vert_{L^\infty} }{\vert \alpha\vert^\frac{3}{2}} \left(\int \frac{\Vert s_{\eta} f_{xx} \Vert^2_{L^\infty}}{\vert \eta \vert^{\frac{3}{2}}}   \ d\eta \right)^{\frac{1}{2}} \ \left(\int \frac{\Vert s_{\kappa} f_{x} \Vert^2_{L^\infty}}{\vert \kappa \vert^{\frac{3}{2}}}   \ d\kappa \right)^{\frac{1}{2}} \ d\alpha  \\
&+&\Vert f \Vert_{\dot H^4} \Vert f \Vert_{\dot H^2} \int \frac{\Vert \delta_\alpha f_x \Vert^2_{L^\infty} }{\vert \alpha\vert^{\frac{3}{2}}} \left(\int \frac{\Vert s_{\eta} f_{xx} \Vert^2_{L^\infty}}{\alpha^2}   \ d\eta \right)^{\frac{1}{2}} \ d\alpha \\
&+&\Vert f \Vert_{\dot H^4} \Vert f \Vert_{\dot H^2} \int \frac{\Vert \delta_\alpha f_x \Vert_{L^\infty} \Vert s_\alpha f \Vert_{L^\infty} }{\vert \alpha\vert^{\frac{5}{2}}} \left(\int \frac{\Vert s_{\eta} f_{xx} \Vert^2_{L^\infty}}{\alpha^2}   \ d\eta \right)^{\frac{1}{2}} \ d\alpha \\
&\lesssim&\Vert f \Vert_{\dot H^4} \Vert f \Vert_{\dot H^2} 
\Vert f \Vert_{\dot B^{\frac{5}{2}}_{\infty,2}}  \Vert f \Vert^2_{\dot B^{\frac{5}{4}}_{\infty,2}}  \\
&+& \Vert f \Vert_{\dot H^4} \Vert f \Vert_{\dot H^2} \Vert f \Vert_{\dot B^{\frac{3}{2}}_{\infty,1}} \Vert f \Vert_{\dot B^{\frac{9}{4}}_{\infty,2}}  \Vert f \Vert_{\dot B^{\frac{5}{4}}_{\infty,2}}  \\
&+&\Vert f \Vert_{\dot H^4} \Vert f \Vert_{\dot H^2}  \Vert f \Vert^2_{\dot B^{\frac{5}{4}}_{\infty,2}} \Vert_{\dot B^{\frac{5}{2}}_{\infty,2}}\\
&+&\Vert f \Vert_{\dot H^4} \Vert f \Vert_{\dot H^2} \Vert f \Vert_{\dot B^{\frac{5}{2}}_{\infty,2}} \Vert f \Vert_{\dot B^{\frac{3}{2}}_{\infty,2}} \Vert f \Vert_{\dot B^{1}_{\infty,2}} \\
&\lesssim&\Vert f \Vert_{\dot H^4} \Vert f \Vert_{\dot H^2} 
\Vert f \Vert_{\dot H^{3}}  \Vert f \Vert^2_{\dot H^{\frac{7}{4}}}  \\
&+& \Vert f \Vert_{\dot H^4} \Vert f \Vert_{\dot H^2} \Vert f \Vert_{\dot B^{\frac{3}{2}}_{\infty,1}} \Vert f \Vert_{\dot H^{\frac{11}{4}}}  \Vert f \Vert_{\dot H^{\frac{7}{4}}}  \\
&+&\Vert f \Vert_{\dot H^4} \Vert f \Vert_{\dot H^2}  \Vert f \Vert^2_{\dot H^{\frac{7}{4}}} \Vert f\Vert_{\dot H^{3}}\\
&+&\Vert f \Vert_{\dot H^4} \Vert f \Vert^2_{\dot H^2} \Vert f \Vert_{\dot H^{3}}  \Vert f \Vert_{\dot H^{\frac{3}{2}}} \\
&\lesssim&\Vert f \Vert_{\dot H^4} \Vert f \Vert^{\frac{4}{5}}_{\dot H^\frac{3}{2}} \Vert f \Vert^{\frac{1}{5}}_{\dot H^4} 
\Vert f \Vert^{\frac{2}{5}}_{\dot H^\frac{3}{2}} \Vert f \Vert^{\frac{3}{5}}_{\dot H^4}  \Vert f \Vert^{\frac{9}{5}}_{\dot H^\frac{3}{2}} \Vert f \Vert^{\frac{1}{5}}_{\dot H^4}  \\
&+& \Vert f \Vert_{\dot H^4} \Vert f \Vert^{\frac{4}{5}}_{\dot H^\frac{3}{2}} \Vert f \Vert^{\frac{1}{5}}_{\dot H^4}  \Vert f \Vert^{\frac{4}{5}}_{\dot B^{1}_{\infty,\infty}} \Vert f \Vert^{\frac{1}{5}}_{\dot B^{\frac{7}{2}}_{\infty,\infty}}  \Vert f \Vert^{\frac{1}{2}}_{\dot H^{\frac{3}{2}}} \Vert f \Vert^{\frac{1}{2}}_{\dot H^{4}}  \Vert f \Vert^{\frac{9}{10}}_{\dot H^{\frac{3}{2}}} \Vert f \Vert^{\frac{1}{10}}_{\dot H^{4}}  \\
&+&\Vert f \Vert_{\dot H^4} \Vert f \Vert^{\frac{4}{5}}_{\dot H^\frac{3}{2}} \Vert f \Vert^{\frac{1}{5}}_{\dot H^4} \Vert f \Vert^{\frac{9}{5}}_{\dot H^\frac{3}{2}} \Vert f \Vert^{\frac{1}{5}}_{\dot H^4}   \Vert f \Vert^{\frac{2}{5}}_{\dot H^\frac{3}{2}} \Vert f \Vert^{\frac{3}{5}}_{\dot H^4}\\
&+&\Vert f \Vert_{\dot H^4} \Vert f \Vert^{\frac{8}{5}}_{\dot H^\frac{3}{2}} \Vert f \Vert^{\frac{2}{5}}_{\dot H^4} \Vert f \Vert^{\frac{2}{5}}_{\dot H^\frac{3}{2}} \Vert f \Vert^{\frac{3}{5}}_{\dot H^4}   \Vert f \Vert_{\dot H^{\frac{3}{2}}}
\end{eqnarray*}
Finally,
\begin{eqnarray*}
\mathcal{U}_{4,3,5,2} + \mathcal{U}_{4,3,5,3} &\lesssim& \Vert f \Vert^2_{\dot H^4} \left(\Vert f \Vert^2_{\dot H^{\frac{3}{2}}}+\Vert f \Vert^3_{\dot H^{\frac{3}{2}}}\right)
\end{eqnarray*}
It remains to estimate $\mathcal{U}_{4,3,5,4}$, since $\partial_\alpha D_\alpha f_x \approx   \frac{s_\alpha f_{xx}}{\alpha}+\frac{1}{\alpha^2} \ {\int_0^\alpha s_\kappa f_{xx} \ d\kappa},$ we easily find that
\begin{eqnarray*}
\mathcal{U}_{4,3,5,4} &\lesssim& \Vert f \Vert_{\dot H^4}\Vert f \Vert_{\dot H^2} \int  \frac{\Vert \delta_{\alpha} f_x \Vert_{L^\infty} \Vert s_{\alpha} f_{xx} \Vert_{L^\infty}}{\vert \alpha\vert^2}  \ d\alpha \\
&+& \Vert f \Vert_{\dot H^4}\Vert f \Vert_{\dot H^2} \int  \frac{\Vert \delta_{\alpha} f_x \Vert_{L^\infty}}{\vert \alpha\vert^{\frac{3}{2}}}  \left( \int \frac{\Vert s_{\kappa} f_{xx} \Vert^2_{L^\infty}}{\alpha^2} \ d\kappa \right)^{\frac{1}{2}} \ d\alpha \\
&\lesssim& \Vert f \Vert_{\dot H^4}\left(\Vert f \Vert^2_{\dot H^2}\Vert f \Vert_{\dot H^3} + \Vert f \Vert_{\dot H^2} \Vert f \Vert_{\dot H^3} \Vert f \Vert_{\dot B^{\frac{3}{2}}_{\infty,1}} \right) \\
&\lesssim& \Vert f \Vert_{\dot H^4}\Vert f \Vert^{\frac{4}{5}}_{\dot H^\frac{3}{2}} \Vert f \Vert^{\frac{1}{5}}_{\dot H^4} 
\Vert f \Vert^{\frac{2}{5}}_{\dot H^\frac{3}{2}} \Vert f \Vert^{\frac{3}{5}}_{\dot H^4}\left(\Vert f \Vert^{\frac{4}{5}}_{\dot H^\frac{3}{2}} \Vert f \Vert^{\frac{1}{5}}_{\dot H^4}  + \Vert f \Vert^{\frac{4}{5}}_{\dot B^{1}_{\infty,\infty}} \Vert f \Vert^{\frac{1}{5}}_{\dot B^{\frac{7}{2}}_{\infty,\infty}}  \right)
\end{eqnarray*}
Finally,
\begin{eqnarray*}
\mathcal{U}_{4,3,5,4} &\lesssim& \Vert f \Vert^2_{\dot H^4}\Vert f \Vert^{2}_{\dot H^\frac{3}{2}}
 \end{eqnarray*}
 Gathering all the estimates of the $\mathcal{U}_{4,3,5,i}$  for $i=1,...4$
 \begin{eqnarray*}
\mathcal{U}_{4,3,5} \lesssim  \Vert f \Vert^2_{\dot H^4} \left(\Vert f \Vert^2_{\dot H^{\frac{3}{2}}}+\Vert f \Vert^3_{\dot H^{\frac{3}{2}}}\right)
\end{eqnarray*}

$\bullet$ {Estimate of $\mathcal{U}_{4,4}$} \\

Recall that 

\begin{eqnarray*}
\mathcal{U}_{4,4}&=&\frac{1}{2\pi} \int \int\int_{0}^{\infty}\int_{0}^{\infty} \ \sigma e^{-\gamma-\sigma} \ \Lambda^{4}\mathcal{H}f \ \partial_{\alpha}\left[\frac{1}{\alpha}\sin(\frac{\gamma}{2}D_\alpha f)\cos(\frac{\gamma}{2}S_\alpha f)\frac{1}{\alpha} \int_{0}^{\alpha} s_{\eta} f_{x}    \ d\eta \right]\\
&& \  \times\  \partial_{x}\left({  \delta_\alpha f_x } \left(\partial_{\alpha}\tau_{\alpha}f_x\right) \sin(\sigma \tau_{\alpha} f_x )\cos(\arctan(\tau_{\alpha} f_{x})) \right) \ d\eta  \ d\gamma \ d\sigma \ dx \ d\alpha\\
\end{eqnarray*}
As the term $$\partial_{\alpha}\left[\frac{1}{\alpha}\sin(\frac{\gamma}{2}D_\alpha f)\cos(\frac{\gamma}{2}S_\alpha f)\frac{1}{\alpha} \int_{0}^{\alpha} s_{\eta} f_{xx}    \ d\eta \right]$$ is one derivative more singular than $\partial_{\alpha}\left[\frac{1}{\alpha}\sin(\frac{\gamma}{2}D_\alpha f)\cos(\frac{\gamma}{2}S_\alpha f)\frac{1}{\alpha} \int_{0}^{\alpha} s_{\eta} f_{x}    \ d\eta \right]$ and the latter has been already estimated in $L^\infty$ in the control of the term $\mathcal{U}_{4,3,5}$. Hence, we may directly follow the same steps. More precisely, one may replace all the $s_{\eta} f_{xx}$ by $s_{\eta} f_{x}$ and of course the remainder term will be one derivative more singular. In order to compare the estimate of $\mathcal{U}_{4,3,5}$ and $\mathcal{U}_{4,4}$ we should first extract $\delta_\alpha f_x$ from the remainder, more precisely we write
\begin{eqnarray*}
\partial_{x}\left({  \delta_\alpha f_x } \left(\partial_{\alpha}\tau_{\alpha}f_x\right) \sin(\sigma \tau_{\alpha} f_x )\cos(\arctan(\tau_{\alpha} f_{x})) \right)&=&\delta_\alpha f_{xx} \ R_{1,\gamma} f + \delta_\alpha f_{x} \ R_{2,\gamma} f
\end{eqnarray*}
where $\delta_\alpha f_{xx}$ and $\delta_\alpha f_{xx}$ will be estimated in $L^\infty$  with the terms into brackets above, namely
 $$
 \left\Vert g \ \partial_{\alpha}\left[\frac{1}{\alpha}\sin(\frac{\gamma}{2}D_\alpha f)\cos(\frac{\gamma}{2}S_\alpha f)\frac{1}{\alpha} \int_{0}^{\alpha} s_{\eta} f_{x}    \ d\eta \right] \right \Vert_{L^\infty}
 $$
where $g$ is either $\delta_\alpha f_{xx}$ or $\delta_\alpha f_{x}$. If $g$ is $\delta_\alpha f_{xx}$, then we need to control the $L^2$ norm of $R_{1,\gamma} f$. It is not difficult to see that
$$
 \left\Vert R_{1,\gamma} f \right\Vert_{L^2} \lesssim \Vert f \Vert_{\dot H^2}
 $$
 and in this case, by adapting accordingly the estimates obtained in the control of $\mathcal{U}_{4,3,5}$. In particular, one needs to adapt the real interpolation parts as $\dot B^{\frac{3}{2}}_{\infty,1}$ becomes 
 $\dot B^{\frac{5}{2}}_{\infty,1}$ while the term in $s_\eta$ becomes one derivative less singular. Therefore, one finds  
 \begin{eqnarray*}
\mathcal{U}_{4,4,1}\lesssim  \Vert f \Vert^2_{\dot H^4} \left(\Vert f \Vert^2_{\dot H^{\frac{3}{2}}}+\Vert f \Vert^3_{\dot H^{\frac{3}{2}}}\right)
\end{eqnarray*}
where we set $\mathcal{U}_{4,4,1}$ the term corresponding to the case $g=\delta_\alpha f_{xx}$.
 
 Otherwise, $g=\delta_\alpha f_{x}$ and in this case we need to control the $L^2$ norm of $R_{2,\gamma} f$ and one may easily check that
 $$
 \left\Vert R_{2,\gamma} f \right\Vert_{L^2} \lesssim \Vert f \Vert_{\dot H^3} + \Vert f \Vert^2_{\dot W^{2,4}}
 $$
 which, by Sobolev embedding, is equivalent to 
 $$
 \left\Vert R_{2,\gamma} f \right\Vert_{L^2} \lesssim \Vert f \Vert_{\dot H^3} + \Vert f \Vert^2_{\dot H^{\frac{9}{4}}}
 $$
As $\Vert f \Vert_{\dot H^3} \lesssim \Vert f \Vert^{\frac{3}{5}}_{\dot H^4} \Vert f \Vert^{\frac{2}{5}}_{\dot H^\frac{3}{2}}$
and 
$\Vert f \Vert_{\dot H^{\frac{9}{4}}} \lesssim \Vert f \Vert^{\frac{3}{5}}_{\dot H^4} \Vert f \Vert^{\frac{2}{5}}_{\dot H^\frac{3}{2}}$, it suffices to consider one of them. Assume it is control by $\dot H^3$ then replace all the $H^2$ by $\dot H^3$ in the control already obtained for the term $\mathcal{U}_{4,3,5}$ and of course take advantage of the fact that $s_\eta$ is one derivative less singular, hence by interpolation we find again the same control as $\mathcal{U}_{4,3,5}$ that is 
 \begin{eqnarray*}
\mathcal{U}_{4,4,2}\lesssim  \Vert f \Vert^2_{\dot H^4} \left(\Vert f \Vert^2_{\dot H^{\frac{3}{2}}}+\Vert f \Vert^3_{\dot H^{\frac{3}{2}}}\right)
\end{eqnarray*}
where we denote $\mathcal{U}_{4,4,2}$  the term corresponding to the case $g=\delta_\alpha f_{x}$. Gathering the estimates of $\mathcal{U}_{4,4,1}$ and $\mathcal{U}_{4,4,2}$, we have obtained

\begin{eqnarray*}
\mathcal{U}_{4,4}\lesssim  \Vert f \Vert^2_{\dot H^4} \left(\Vert f \Vert^2_{\dot H^{\frac{3}{2}}}+\Vert f \Vert^3_{\dot H^{\frac{3}{2}}}\right)
\end{eqnarray*}

$\bullet$ {Estimate of $\mathcal{U}_{4,5}.$}\\

This term is the same as $\mathcal{U}_{4,3}$ as $\frac{1}{1+(f_x)^2}\leq 1$\\

$\bullet$ {Estimate of $\mathcal{U}_{4,6}$}\\

This term is the same as $\mathcal{U}_{4,4}$ as $\left\Vert\partial_{x}\frac{1}{1+(f_x)^2} \right\Vert_{L^p} \lesssim \Vert f_{xx} \Vert_{L^p}$ and therefore this extra term may be treated as $\partial_{x} \sin(\sigma \tau_{\alpha} f_x)$ as we always used the estimate $\partial_{x} \sin(\sigma \tau_{\alpha} f_x)\lesssim \sigma \Vert f_{xx} \Vert_{L^p}$.

We may conclude the estimate of the term $\mathcal{U}_{4}$ by collecting all the obtained estimates, namely

\begin{eqnarray*}
\mathcal{U}_{4}\lesssim  \Vert f \Vert^2_{\dot H^4} \left(\Vert f \Vert_{\dot H^{\frac{3}{2}}}+\Vert f \Vert^2_{\dot H^{\frac{3}{2}}}+\Vert f \Vert^3_{\dot H^{\frac{3}{2}}}+\Vert f \Vert^4_{\dot H^{\frac{3}{2}}}\right)
\end{eqnarray*}

\subsection{Estimate of $\mathcal{U}_{5}$ }

 Similar to $\mathcal{U}_{4}$ up to interchanging the role of sin and cos as they were bounded by 1 this gives the same estimate as $\mathcal{U}_{5}$. Note that, unlike the terms $\mathcal{U}_{3}$ and $\mathcal{U}_{6}$, we have never used the fact that the operator $S_\alpha$ is more regular than the operator $D_\alpha$, however the symmetry is crucial to estimate these terms. Therefore,
 
 \begin{eqnarray*}
\mathcal{U}_{5}\lesssim  \Vert f \Vert^2_{\dot H^4} \left(\Vert f \Vert_{\dot H^{\frac{3}{2}}}+\Vert f \Vert^2_{\dot H^{\frac{3}{2}}}+\Vert f \Vert^3_{\dot H^{\frac{3}{2}}}+\Vert f \Vert^4_{\dot H^{\frac{3}{2}}}\right)
\end{eqnarray*}
 
 Finally, it remains to estimate $\mathcal{U}_{6}$.
 
 \subsection{Estimate of $\mathcal{U}_{6}$ }
 We have
 \begin{eqnarray*}
 \mathcal{U}_{6}&=&- \frac{1}{4\pi}\int \Lambda^{4}\mathcal{H}f \ \partial_x \int\int_{0}^{\infty}\int_{0}^{\infty} \ e^{-\gamma-\sigma} \ S_{\alpha} f \left(\sin(\gamma\Delta_\alpha f)-\sin(\gamma\bar\Delta_\alpha f)\right) \\
&& \ \times \  \partial_x \left(\frac{\partial_x^2 \tau_{\alpha} f}{\alpha} \cos(\sigma \tau_{\alpha} f_x )\cos(\arctan(\tau_{\alpha}f_{x}))\right) \ d\alpha \ d\gamma \ d\sigma \ dx  \\
&=&- \frac{1}{4\pi}\int \Lambda^{4}\mathcal{H}f \  \int\int_{0}^{\infty}\int_{0}^{\infty} \ e^{-\gamma-\sigma} \ S_{\alpha} f_x \left(\sin(\gamma\Delta_\alpha f)-\sin(\gamma\bar\Delta_\alpha f)\right) \\
&& \ \times \  \partial_x \left(\frac{\partial_x^2 \tau_{\alpha} f}{\alpha} \cos(\sigma \tau_{\alpha} f_x )\cos(\arctan(\tau_{\alpha}f_{x}))\right) \ d\alpha \ d\gamma \ d\sigma \ dx \\
&-& \frac{1}{4\pi}\int \Lambda^{4}\mathcal{H}f \  \int\int_{0}^{\infty}\int_{0}^{\infty} \ e^{-\gamma-\sigma} \ S_{\alpha} f \left(\partial_x\left(\sin(\gamma\Delta_\alpha f)-\sin(\gamma\bar\Delta_\alpha f)\right)\right) \\
&& \ \times \  \partial_x \left(\frac{\partial_x^2 \tau_{\alpha} f}{\alpha} \cos(\sigma \tau_{\alpha} f_x )\cos(\arctan(\tau_{\alpha}f_{x}))\right) \ d\alpha \ d\gamma \ d\sigma \ dx \\
&-& \frac{1}{4\pi}\int \Lambda^{4}\mathcal{H}f \  \int\int_{0}^{\infty}\int_{0}^{\infty} \ e^{-\gamma-\sigma} \ S_{\alpha} f \left(\sin(\gamma\Delta_\alpha f)-\sin(\gamma\bar\Delta_\alpha f)\right) \\
&& \ \times \  \partial_{xx} \left(\frac{\partial_x^2 \tau_{\alpha} f}{\alpha} \cos(\sigma \tau_{\alpha} f_x )\cos(\arctan(\tau_{\alpha}f_{x}))\right) \ d\alpha \ d\gamma \ d\sigma \ dx \\
&=& \mathcal{U}_{6,1} +  \mathcal{U}_{6,2} +  \mathcal{U}_{6,3}
\end{eqnarray*}

$\bullet$ {Estimate of $\mathcal{U}_{6,1}$} \\

Using Sobolev embedding $\dot H^{9/4} \hookrightarrow \dot W^{2,4}$ and real interpolation we find the following estimate
\begin{eqnarray*}
\mathcal{U}_{6,1} &\lesssim& \Vert f \Vert_{\dot H^4} \left( \Vert f \Vert_{\dot H^3} + \Vert f \Vert_{\dot H^{\frac{9}{4}}}^2\right) \Vert f \Vert_{\dot B^2_{\infty,1}}\\
&\lesssim& \Vert f \Vert_{\dot H^4} \Vert f\Vert^{\frac{3}{5}}_{\dot H^4}  \Vert f \Vert^{\frac{2}{5}}_{\dot H^\frac{3}{2}} \Vert f \Vert^{\frac{2}{5}}_{\dot H^4}  \Vert f \Vert^{\frac{3}{5}}_{\dot H^\frac{3}{2}} \\
&\lesssim& \Vert f \Vert^2_{\dot H^4}   \Vert f \Vert^2_{\dot H^\frac{3}{2}}
\end{eqnarray*}

$\bullet$ {Estimate of $\mathcal{U}_{6,2}$}\\

\begin{eqnarray*}
\mathcal{U}_{6,2} &\lesssim& \Vert f \Vert_{\dot H^4} (\Vert f \Vert_{\dot H^3} + \Vert f \Vert^2_{\dot H^{\frac{9}{4}}}) \int \frac{\Vert s_{\alpha} f \Vert_{L^\infty} \Vert \delta_{\alpha} f_x \Vert_{L^\infty}}{\vert \alpha\vert^3} \ d\alpha \\
&\lesssim& \Vert f \Vert_{\dot H^4} (\Vert f \Vert_{\dot H^3} + \Vert f \Vert^2_{\dot H^{\frac{9}{4}}}) \Vert f \Vert^2_{\dot H^2} \\
&\lesssim& \Vert f \Vert_{\dot H^4} \Vert f\Vert^{\frac{3}{5}}_{\dot H^4}  \Vert f \Vert^{\frac{2}{5}}_{\dot H^\frac{3}{2}}\Vert f\Vert^{\frac{2}{5}}_{\dot H^4}  \Vert f \Vert^{\frac{3}{5}}_{\dot H^\frac{3}{2}} \\
&\lesssim& \Vert f \Vert^2_{\dot H^4} \Vert f \Vert^2_{\dot H^\frac{3}{2}}
\end{eqnarray*}

$\bullet$ {Estimate of $\mathcal{U}_{6,3}$}\\

We have that

\begin{eqnarray*}
\mathcal{U}_{6,3}&=&- \frac{1}{4\pi}\int \Lambda^{4}\mathcal{H}f \  \int\int_{0}^{\infty}\int_{0}^{\infty} \ e^{-\gamma-\sigma} \ S_{\alpha} f \left(\sin(\gamma\Delta_\alpha f)-\sin(\gamma\bar\Delta_\alpha f)\right) \\
&& \ \times \ \frac{\partial_x^4 \tau_{\alpha} f}{\alpha} \cos(\sigma \tau_{\alpha} f_x )\cos(\arctan(\tau_{\alpha}f_{x})) \ d\alpha \ d\gamma \ d\sigma \ dx \\
&-& \frac{1}{4\pi}\int \Lambda^{4}\mathcal{H}f \  \int\int_{0}^{\infty}\int_{0}^{\infty} \ e^{-\gamma-\sigma} \ S_{\alpha} f \left(\sin(\gamma\Delta_\alpha f)-\sin(\gamma\bar\Delta_\alpha f)\right) \\
&& \ \times \ \frac{\partial_x^2 \tau_{\alpha} f}{\alpha} \partial_{xx}\left(\cos(\sigma \tau_{\alpha} f_x )\cos(\arctan(\tau_{\alpha}f_{x}))\right) \ d\alpha \ d\gamma \ d\sigma \ dx \\
&-& \frac{1}{2\pi}\int \Lambda^{4}\mathcal{H}f \  \int\int_{0}^{\infty}\int_{0}^{\infty} \ e^{-\gamma-\sigma} \ S_{\alpha} f \left(\sin(\gamma\Delta_\alpha f)-\sin(\gamma\bar\Delta_\alpha f)\right) \\
&& \ \times \ \frac{\partial_x^3 \tau_{\alpha} f}{\alpha} \partial_{x}\left(\cos(\sigma \tau_{\alpha} f_x )\cos(\arctan(\tau_{\alpha}f_{x}))\right) \ d\alpha \ d\gamma \ d\sigma \ dx \\
\end{eqnarray*}
The term $\mathcal{U}_{6,3}$  is the most singular term. We have the following control
\begin{lemma} \label{lu63}
The following estimate holds
\begin{equation} \label{U63}
\mathcal{U}_{6,3} \lesssim \Vert f \Vert^2_{\dot H^4} \left( \Vert f \Vert_{\dot H^{3/2}}+\Vert f \Vert^2_{\dot H^{3/2}}  +\Vert f \Vert_{\dot H^{3/2}}  \Vert  f \Vert_{\dot B^1_{\infty,1}}\right)
\end{equation}
\end{lemma}
\noindent {\bf{Proof of Lemma \ref{lu63}}.}
As the term  $\partial_x^4 \tau_{\alpha} f$ contains too many derivatives we need to balance the regularity by writting $\partial_x^4 \tau_{\alpha} f= -\partial_\alpha \delta_{\alpha} \partial_x^3 f$ and then integrating by parts in $\alpha$ 
\begin{eqnarray*}
\mathcal{U}_{6,3} &=&-\frac{1}{2\pi} \int \int\int_{0}^{\infty}\int_{0}^{\infty} \ e^{-\gamma-\sigma} \ \Lambda^{3} f \ S_{\alpha} f  \sin(\frac{\gamma}{2} D_\alpha f)  \\
&& \  \times\  \frac{ \partial_\alpha \delta_{\alpha} \partial_x^3 f}{\alpha} \cos(\sigma \tau_{\alpha} f_x )\cos(\arctan(\tau_{\alpha} f_{x}))  \ d\alpha \ d\gamma \ d\sigma \ dx \\
&+&\frac{1}{\pi} \int \int\int_{0}^{\infty}\int_{0}^{\infty} \ e^{-\gamma-\sigma} \ \Lambda^{3} f \ S_{\alpha} f \sin^2(\frac{\gamma}{4} S_\alpha f) \sin(\frac{\gamma}{2} D_\alpha f)  \\
&& \  \times\ \frac{ \partial_\alpha \delta_{\alpha} \partial_x^3 f}{\alpha} \cos(\sigma \tau_{\alpha} f_x )\cos(\arctan(\tau_{\alpha} f_{x}))  \ d\alpha \ d\gamma \ d\sigma \ dx \\
&=& \mathcal{U}_{6,3,1}+\mathcal{U}_{6,3,2}
\end{eqnarray*}
In order to estimate $\mathcal{U}_{6,3,1}$ we first integrate by parts in $\alpha$,

\begin{eqnarray*}
\mathcal{U}_{6,3,1}&=&-\frac{1}{2\pi} \int \int\int_{0}^{\infty}\int_{0}^{\infty} \ e^{-\gamma-\sigma} \ \Lambda^{4} \mathcal{H} f \ S_{\alpha} f  \sin(\frac{\gamma}{2} D_\alpha f)  \\
&& \times\ \frac{  \delta_{\alpha} \partial_x^3 f}{\alpha^2} \cos(\sigma \tau_{\alpha} f_x )\cos(\arctan(\tau_{\alpha} f_{x}))  \ d\alpha \ d\gamma \ d\sigma \ dx \\
 &+&\frac{1}{2\pi} \int \int\int_{0}^{\infty}\int_{0}^{\infty} \ e^{-\gamma-\sigma} \ \Lambda^{4} \mathcal{H} f \ \partial_\alpha(S_{\alpha} f)  \sin(\frac{\gamma}{2} D_\alpha f)\\
&& \times\ \frac{  \delta_{\alpha} \partial_x^3 f}{\alpha} \cos(\sigma \tau_{\alpha} f_x )\cos(\arctan(\tau_{\alpha} f_{x}))  \ d\alpha \ d\gamma \ d\sigma \ dx \\
&+&\frac{1}{4\pi} \int \int\int_{0}^{\infty}\int_{0}^{\infty} \ \gamma e^{-\gamma-\sigma} \ \Lambda^{4} \mathcal{H} f \ \partial_\alpha(D_{\alpha} f)  \cos(\frac{\gamma}{2} D_\alpha f)\\
&& \times\ \frac{\delta_{\alpha} \partial_x^3 f}{\alpha} \cos(\sigma \tau_{\alpha} f_x )\cos(\arctan(\tau_{\alpha} f_{x}))  \ d\alpha \ d\gamma \ d\sigma \ dx \\
&-&\frac{1}{2\pi} \int \int\int_{0}^{\infty}\int_{0}^{\infty} \ \gamma e^{-\gamma-\sigma} \ \Lambda^{4} \mathcal{H} f \ S_{\alpha} f  \sin(\frac{\gamma}{2} D_\alpha f)\\
&& \times\ \frac{\delta_{\alpha} \partial_x^3 f}{\alpha} \partial_\alpha \left( \cos(\sigma \tau_{\alpha} f_x )\cos(\arctan(\tau_{\alpha} f_{x})) \right)  \ d\alpha \ d\gamma \ d\sigma \ dx \\
&:=&\sum_{i=1}^5 \mathcal{U}_{6,3,1,i}
\end{eqnarray*}

We shall estimates  the $\mathcal{U}_{6,3,1,i}$, $i=1,...,5$. \\

$\bullet$ {{Estimate of $\mathcal{U}_{6,3,1,1}$}} \\

It is easy to see that
\begin{eqnarray*}
\mathcal{U}_{6,1,1,1}&\lesssim& \Vert f \Vert_{\dot H^4} \int \frac{\Vert s_{\alpha}f \Vert_{L^{2}}}{\vert\alpha \vert^2} \frac{\Vert \delta_{\alpha}f_{xx} \Vert_{L^{\infty}}}{\vert\alpha \vert}\ d\alpha  \\
&\lesssim& \Vert f \Vert_{\dot H^4} \left(\int \frac{\Vert s_{\alpha}f \Vert^2_{L^{2}}}{\vert\alpha \vert^4} \ d\alpha \right)^{1/2} \left(\frac{\Vert \delta_{\alpha}\partial_x^3 f \Vert^2_{L^{\infty}}}{\vert\alpha \vert^2}\ d\alpha \right)^{1/2} \\
&\lesssim& \Vert f \Vert_{\dot H^4} \Vert f \Vert_{\dot H^{3/2}} \Vert \partial_x^3 f \Vert_{\dot B^{1/2}_{\infty,2}} \\
&\lesssim& \Vert f \Vert^2_{\dot H^4} \Vert f \Vert_{\dot H^{3/2}}    
\end{eqnarray*}

$\bullet$ {{Estimate of $\mathcal{U}_{6,3,1,2}$}} \\

By using the fact that $\partial_\alpha(S_{\alpha} f)\approx \Delta_ \alpha f_{x}+\bar\Delta_ \alpha f_{x} + \frac{s_{\alpha}f}{\alpha^2} $ and that we may omit to treat the term $\bar\Delta_ \alpha f_{x}$ term as it is similar to $\Delta_ \alpha f_{x}$, we find

\begin{eqnarray*}
\mathcal{U}_{6,3,1,2}&\lesssim& \Vert f \Vert_{\dot H^4}  \int \frac{\Vert \bar \delta_{\alpha} f_x \Vert_{L^{\infty}}}{\vert\alpha\vert} \frac{\Vert \delta_{\alpha} \partial_{x}^3 f \Vert_{L^{2}}}{\vert\alpha\vert}\ d\alpha \\
&+&\Vert f \Vert_{\dot H^3} \int \frac{\Vert s_{\alpha}f \Vert_{L^{\infty}}}{\vert\alpha \vert^2} \frac{\Vert \delta_{\alpha}\partial_{x}^3 f\Vert_{L^{\infty}}}{\vert\alpha \vert}\ d\alpha  \\
&\lesssim& \Vert f \Vert_{\dot H^4}  \Vert f_x \Vert_{\dot B^{1/2}_{\infty,2}} \Vert \partial_{x}^3 f \Vert_{\dot H^{1/2}} +\Vert f \Vert^2_{\dot H^4} \Vert f \Vert_{\dot B^{1}_{\infty,2}} \\
&\lesssim& \Vert f \Vert_{\dot H^4}  \Vert f \Vert_{\dot H^{2}} \Vert f \Vert_{\dot H^{7/2}} + \Vert f \Vert^2_{\dot H^4}   \Vert f \Vert_{\dot H^{3/2}} \\
&\lesssim& \Vert f \Vert^2_{\dot H^4}   \Vert f \Vert_{\dot H^{3/2}}
\end{eqnarray*}

$\bullet$ {{Estimate of $\mathcal{U}_{6,3,1,3}$}}\\

Since $\partial_\alpha D \approx \frac{s_\alpha f_x}{\alpha}+\frac{1}{\alpha^2} \ {\displaystyle \int_0^\alpha s_\kappa f_x \ d\kappa}$, we find that

\begin{eqnarray*}
\mathcal{U}_{6,3,1,3}&\lesssim& \Vert f \Vert_{\dot H^4} \left(\int \frac{\Vert  \delta_{\alpha} \partial_x^3 f \Vert_{L^{\infty}}}{\vert \alpha \vert} \frac{\Vert s_{\alpha} f_{x} \Vert_{L^{2}}}{\vert \alpha \vert} \ d\alpha +
\int \frac{\Vert  \delta_{\alpha} \partial_x^3 f \Vert_{L^{\infty}}}{\vert\alpha\vert^3} \int_0^\alpha s_\kappa f_x \ d\kappa \ d\alpha \right)\\
&\lesssim&\Vert f \Vert_{\dot H^4}  \left(\Vert \partial_x^3 f \Vert_{\dot B^{1/2}_{\infty,2}} \Vert f \Vert_{\dot H^{3/2}} + \int \frac{\Vert  \delta_{\alpha} \partial_x^3 f \Vert_{L^{\infty}}}{\vert\alpha\vert^{\frac{5}{4}}} \left(\int \frac{\Vert s_\kappa f_x \Vert^2_{L^2}}{\vert\kappa\vert^\frac{5}{2}} \ d\kappa \right)^{\frac{1}{2}} \ d\alpha \right)\\
&\lesssim&\Vert f \Vert_{\dot H^4}  \left(\Vert f \Vert_{\dot H^{4}} \Vert f \Vert_{\dot H^{3/2}} + \Vert  f \Vert_{\dot B^{\frac{13}{4}}_{\infty,1}} \Vert  f \Vert_{\dot H^{\frac{7}{4}}} \right) \\
\end{eqnarray*}
Since $\Vert  f \Vert_{\dot B^{\frac{13}{4}}_{\infty,1}}\lesssim  \Vert  f \Vert^{\frac{1}{10}}_{\dot B^{1}_{\infty,\infty}} \Vert  f \Vert^{\frac{9}{10}}_{\dot B^{\frac{7}{2}}_{\infty,\infty}}$ and $\Vert  f \Vert_{\dot H^{\frac{7}{4}}}\lesssim  \Vert  f \Vert^{\frac{9}{10}}_{\dot H^{\frac{3}{2}}} \Vert  f \Vert^{\frac{1}{10}}_{\dot H^{4}}$ , then 

\begin{eqnarray*}
\mathcal{U}_{6,3,1,3}&\lesssim&\Vert f \Vert^2_{\dot H^4}  \Vert f \Vert_{\dot H^{3/2}} 
\end{eqnarray*}

$\bullet$ {{Estimate of $\mathcal{U}_{6,3,1,4}$}} \\

We have 
\begin{eqnarray*}
\mathcal{U}_{6,3,1,4}&=&-\frac{1}{2\pi} \int \int\int_{0}^{\infty}\int_{0}^{\infty} \ \gamma e^{-\gamma-\sigma} \ \Lambda^{4} \mathcal{H} f \ S_{\alpha} f  \sin(\frac{\gamma}{2} D_\alpha f)\\
&& \times\ \frac{\delta_{\alpha} \partial_x^3 f}{\alpha} \partial_\alpha \left( \cos(\sigma \tau_{\alpha} f_x )\cos(\arctan(\tau_{\alpha} f_{x})) \right)  \ d\alpha \ d\gamma \ d\sigma \ dx  \\
\end{eqnarray*}
This term may be controlled as follows, 
\begin{eqnarray*}
\mathcal{U}_{6,3,1,4} &\lesssim&  \Vert f \Vert_{\dot H^4} \ \int \frac{\Vert s_{\alpha} f \Vert_{L^{\infty}}}{\vert \alpha \vert^2}
 \Vert \delta_{\alpha} \partial_{x}^3 f \Vert_{L^4}  \Vert \tau_{\alpha}f_{xx} \Vert_{L^4}   \ d\alpha \\
 &\lesssim&  \Vert f \Vert_{\dot H^4} \ \Vert \partial_{x}^3 f \Vert^2_{L^4}  \int \frac{\Vert s_{\alpha} f \Vert_{L^{\infty}}}{\vert \alpha \vert^2}    \ d\alpha \\
 &\lesssim&  \Vert f \Vert_{\dot H^4}  \ \Vert f \Vert^2_{\dot H^{\frac{13}{4}}}  \int \frac{\Vert s_{\alpha} f \Vert_{L^{\infty}}}{\vert \alpha \vert^2}    \ d\alpha \\
 &\lesssim& \Vert f \Vert^2_{\dot H^4} \ \Vert f \Vert_{\dot H^{\frac{3}{2}}}  \Vert  f \Vert_{\dot B^1_{\infty,1}}   \\
\end{eqnarray*}

 $\mathcal{U}_{6,3,1, 4}$ and  $\mathcal{U}_{6,3,1, 4}$ are similar, up to some bounded terms which are less or equal to 1. This ends the proof of Lemma.
 
 \qed
 
 \section{End of the proof of the main theorem}
 
Gathering all the estimates of the $\mathcal{U}_i, i=2,..,6$, we have obtained the following control
\begin{equation*} 
\frac{1}{2} \partial_t \Vert f \Vert^2_{\dot H^{5/2}} + \int   \ \frac{\vert\Lambda^4 f \vert^2}{(1+ \vert f_x \vert^2)^{\frac{3}{2}}} \ dx \lesssim \Vert f \Vert^2_{\dot H^4} \left(\mathcal{P}(\Vert f \Vert_{\dot H^{3/2}})+\Vert f \Vert_{\dot H^{3/2}}  \Vert  f \Vert_{\dot B^1_{\infty,1}}\right)
\end{equation*}
where $\mathcal{P}(X)=X+X^2+X^3+X^4 \approx X + X^4.$ Integrating in time gives
 \begin{eqnarray*}
\Vert f(T) \Vert^2_{\dot H^{5/2}}+ \int_0^T\int   \ \frac{\vert\Lambda^4 f \vert^2}{(1+ \vert f_x \vert^2)^{\frac{3}{2}}} \ dx \ ds &\lesssim& \Vert f_0 \Vert^2_{\dot H^{5/2}} \\
&+& \int_0^T\Vert f \Vert^2_{\dot H^4} \left( \mathcal{P}(\Vert f \Vert_{\dot H^{3/2}})+\Vert f \Vert_{\dot H^{3/2}}  \Vert  f \Vert_{\dot B^1_{\infty,1}}\right) \ ds,
\end{eqnarray*}
 Then, since $0\leq a \mapsto \frac{1}{(1+ a^2)^{\frac{3}{2}}}$ is decreasing  we infer that
$$
\frac{1}{(1+ \vert f_x\vert^2)^{\frac{3}{2}}} \geq \frac{1}{(1+  \Vert f_x \Vert_{L^\infty}^2)^{\frac{3}{2}}}
$$
Hence, we may write
\begin{equation*} 
\Vert f(T) \Vert^2_{\dot H^{5/2}}+    \ \displaystyle\int_0^T \frac{\int\vert\Lambda^4 f \vert^2  \ dx}{(1+  \displaystyle  \Vert f_x \Vert_{L^\infty_{x}}^2)^{\frac{3}{2}}} \ ds \lesssim \Vert f_0 \Vert^2_{\dot H^{5/2}} +{\displaystyle\int_0^T\Vert f \Vert^2_{\dot H^4} \left( \mathcal{P}(\Vert f \Vert_{\dot H^{3/2}})+\Vert f \Vert_{\dot H^{3/2}}  \Vert  f \Vert_{\dot B^1_{\infty,1}}\right) \ ds}
\end{equation*}
Set $L:= \Vert f_x \Vert_{L^\infty_{x}}$. Then, combining the above estimate with the $L^2$ control obtained in Lemma \ref{L2-est} we get the following {\it{a priori}} estimates for the full Muskat problem with surface tension \eqref{CP0},

\begin{eqnarray*} 
&&\Vert f(T) \Vert^2_{\dot H^{5/2}} +   \ \underbrace{\int_0^T \frac{\displaystyle \Vert f \Vert_{\dot H^4}^2}{(1+  L^2)^{\frac{3}{2}}} \ ds}_{{\rm{contribution \ when}} \ g=0}  + \underbrace{\int_0^T \frac{\displaystyle \Vert f \Vert_{\dot H^3}^2}{1+  L^2} \ ds}_{\rm contribution \ when \ \sigma=0}  \\
&& \lesssim \Vert f_0 \Vert^2_{\dot H^{5/2}} +\underbrace{\displaystyle\int_0^T\Vert f \Vert^2_{\dot H^4} \left( \mathcal{P}(\Vert f \Vert_{\dot H^{3/2}})+\Vert f \Vert_{\dot H^{3/2}}  \Vert  f \Vert_{\dot B^1_{\infty,1}}\right) \ ds}_{{\rm{contribution \ when}} \ g=0}+\underbrace{\displaystyle\int_0^T  \Vert f \Vert^2_{\dot H^{3}}\mathcal{Q}(\Vert f \Vert_{\dot H^{3/2}}) \  \ ds}_{\rm contribution \ when \ \sigma=0},
\end{eqnarray*}
where $\mathcal{Q}(X)=X+X^2$. This gives 

\begin{eqnarray} \label{H52}
&&\Vert f(T) \Vert^2_{\dot H^{5/2}} +    \ \int_0^T {\displaystyle \Vert f \Vert_{\dot H^4}^2}\frac{\left(1 - (1+  L^2)^{\frac{3}{2}}\mathcal{P}(\Vert f \Vert_{\dot H^{3/2}})-\Vert f \Vert_{\dot H^{3/2}}  \Vert  f \Vert_{\dot B^1_{\infty,1}}(1+  L^2)^{\frac{3}{2}} \right) }{(1+  L^2)^{\frac{3}{2}}} \ ds \nonumber \\
&&  \ + \int_0^T {\displaystyle \Vert f \Vert_{\dot H^3}^2}\left(\frac{1 - \mathcal{Q}(\Vert f \Vert_{\dot H^{3/2}})(1+  L^2)  }{1+  L^2} \right) \ ds \lesssim \Vert f_0 \Vert^2_{\dot H^{5/2}} 
\end{eqnarray}
Hence, if 
$$ 
(1+  L^2)^{\frac{3}{2}}\mathcal{P}(\Vert f \Vert_{\dot H^{3/2}})+\Vert f \Vert_{\dot H^{3/2}}  \Vert  f \Vert_{\dot B^1_{\infty,1}}(1+  L^2)^{\frac{3}{2}}<1,
$$
we get the wanted global control of the $H^{5/2}$ norm. This ends the proof of the main result.

\qed

\end{document}